\documentstyle[11pt]{report}
\title{Light-Cone Expansion of the Dirac Sea with Light Cone Integrals}
\author{Felix Finster}
\date{July 1997}

\makeatletter
\long\def\@makefntext#1{\parindent 1em\noindent
        \hbox to 1.8em{\hss$^{\@thefnmark}$}#1}
\makeatother

\newtheorem{Def}{Def.}[section]
\newtheorem{Thm}[Def]{Theorem}

\newtheorem{Korollar}[Def]{Korollar}
\newtheorem{Lemma}[Def]{Lemma}

\newtheorem{Bem}[Def]{Bemerkung}
\newtheorem{Satz}[Def]{Satz}
\newcommand{\Beweis}{\em{Beweis:}}
\newcommand{\QED}{\ \hfill $\Box$ \\[1em]}
\newcommand{\EndBem}{\ \hfill $\Delta$ \\[1em]}
\newcommand{\Aslsh}{\mbox{ $\!\!A$ \hspace{-1.2 em} $/$}}

\newcommand{\Equ}[1]{\begin{equation} \label{eq:#1}}
\newcommand{\EndEqu}{\end{equation}}
\newcommand{\Ref}[1]{(\ref{eq:#1})}
\newcommand{\bra}{\mbox{$< \!\!$ \nolinebreak}}
\newcommand{\ket}{\mbox{\nolinebreak $>$}}

\newcommand{\R}{\mbox{\rm I \hspace{-.8 em} R}}
\newcommand{\1}{\mbox{\rm 1 \hspace{-1.05 em} 1}}

\newcommand{\sR}{\mbox{\rm \scriptsize I \hspace{-.8 em} R}}
\newcommand{\N}{\mbox{\rm I \hspace{-.8 em} N}}

\newcommand{\Pdd}{\mbox{$\partial$ \hspace{-1.2 em} $/$}}
\newcommand{\slsh}{\mbox{ \hspace{-1.1 em} $/$}}

\newcommand{\Tr}{\mbox{Tr\/}}

\newcommand{\spc}{\;\;\;\;\;\;\;\;\;\;}
\newcommand{\lint}{\: \Diamond \hspace{-0.98 em} \int}
\newcommand{\wint}{\: \bigtriangleup \hspace{-1.09 em} \int}
\newcommand{\vint}{\: \bigtriangledown \hspace{-0.98 em} \int}
\newcommand{\olint}{\: \odot \hspace{-1.00 em} \int}
\newcommand{\xint}{\: \times \hspace{-1.00 em} \int}
\newcommand{\dint}{\: \odot \hspace{-1.00 em} \int}
\newcommand{\xdint}{\: \otimes \hspace{-1.00 em} \int}
\newcommand{\hint}{\: - \hspace{-1.00 em} \int}
\newcommand{\hinti}{\: - \hspace{-1.00 em} \int^\infty_{-\infty}}
\newcommand{\inti}{\int^\infty_{-\infty}}
\newcommand{\dotint}{\! \mbox{ \begin{picture}(10,10)
        \put(5,4){\circle*{6}}
        \end{picture}   }
        \hspace{-1.38 em} \int}

\newcommand{\veeint}{\: \vee \hspace{-1.00 em} \int}
\newcommand{\slint}{{\scriptstyle \Diamond \hspace{-0.72 em}} \int}
\newcommand{\swint}{{\scriptstyle \bigtriangleup \hspace{-0.72 em}} \int}
\newcommand{\svint}{{\scriptstyle \bigtriangledown \hspace{-0.72 em}} \int}
\newcommand{\solint}{{\scriptstyle \odot \hspace{-0.72 em}} \int}
\newcommand{\sxint}{{\scriptstyle \times \hspace{-0.72 em}} \int}
\newcommand{\sdint}{{\scriptstyle \odot \hspace{-0.72 em}} \int}
\newcommand{\sxdint}{{\scriptstyle \otimes \hspace{-0.72 em}} \int}
\newcommand{\sveeint}{{\scriptstyle \vee \hspace{-0.72 em}} \int}
\newcommand{\flint}{{\scriptstyle \Diamond \hspace{-1.07 em}} \int}
\newcommand{\sdotint}{\!\!\mbox{ \begin{picture}(10,10)
        \put(5,3){\circle*{4}}
        \end{picture}   }
        \hspace{-1.2 em} \int}
\newcommand{\Ra}{\Re}
\newcommand{\I}{\Im}
\newcommand{\Li}{{\cal L}}
\newcommand{\pn}{p^{(0)}}
\newcommand{\kn}{k^{(0)}}
\newcommand{\sn}{s^{(0)}}
\newcommand{\pe}{p^{(1)}}
\newcommand{\ke}{k^{(1)}}
\newcommand{\se}{s^{(1)}}
\newcommand{\sve}{s^\vee_{(1)}}
\newcommand{\svn}{s^\vee_0}
\newcommand{\T}{{\mbox{T}}}
\newcommand{\Texp}{{\mbox{Texp}}}

\begin{document}
\maketitle
\begin{abstract}
The Dirac sea is calculated in an expansion around the light cone. The 
method is to analyze the perturbation expansion for the Dirac sea in 
position space. This leads to integrals over expressions containing
distributions which are singular on the light cone.
We derive asymptotic formulas for these ``light cone integrals'' in 
terms of line integrals over the external potential and its partial 
derivatives.

The calculations are based on the perturbation expansion for the Dirac
sea in the preprint gr-qc/9606040 and yield the formulas listed in 
the appendix of this preprint. The results can be obtained easier 
with a combination of calculations in position and momentum space (see 
the corresponding preprints on the hep-th server).
Therefore the calculations are preliminary and will remain 
unpublished; they are intended as reference for people who encounter 
similar mathematical problems.
\end{abstract}

\tableofcontents
\newpage

\section{Einige vorbereitende Formeln}
In dieser Arbeit werden die Bezeichungen von \cite{F3} verwendet.
Dort werden die hier aufgelisteten
Formeln auch genauer abgeleitet und erkl\"art.

Die Distributionen
\begin{eqnarray}
P_{m^2}(x) &=& \int \frac{d^4k}{(2 \pi)^4} \: \delta(k^2-m^2) 
\:e^{-ikx} \nonumber \\
\label{eq:2_20}
&=& \left\{ \displaystyle \begin{array}{ll}
	\displaystyle \frac{m^2}{8 \pi^2} \frac{
	Y_1(\sqrt{m^2 x^2})}{\sqrt{m^2 x^2}} & {\mbox{f\"ur $x^2>0$}} \\[.3cm]
	\displaystyle \frac{m^2}{4 \pi^3} \frac{K_1(\sqrt{-m^2 x^2})}
	{\sqrt{-m^2 x^2}} & {\mbox{f\"ur $x^2<0$}} \end{array} \right. \\
\label{eq:2_1a}
K_{m^2}(x) &=& \int \frac{d^4k}{(2 \pi)^4} \: \delta(k^2-m^2) \:
\epsilon(k^0) \:e^{-ikx} \\
\label{eq:2_21}
&=& -\frac{i}{4 \pi^2} \: \delta(x^2) \: \epsilon(x^0)
	\;+\; \frac{i m^2}{8 \pi^2} \frac{J_1(\sqrt{m^2 x^2})}{\sqrt{m^2 x^2}}
	\: \Theta(x^2) \: \epsilon(x^0)
\end{eqnarray}
sind L\"osungen der Klein-Gordon Gleichung. Wir verwenden die Greensfunktion
\begin{eqnarray}
S_{m^2}(x) &=& \frac{1}{2} \: \sum_{\pm} \int \frac{d^4k}{(2 \pi)^4} \:
\frac{1}{k^2-m^2 \pm i \varepsilon k^0} \:e^{-ikx} \nonumber \\
\label{eq:2_22}
&=& -\frac{1}{4 \pi} \: \delta(x^2) \;+\; \frac{m^2}{8 \pi}
	\frac{J_1(\sqrt{m^2 x^2})}{\sqrt{m^2 x^2}} \: \Theta(x^2) \spc .
\end{eqnarray}
F\"ur den Diracoperator $(i \Pdd-m)$ erh\"alt man die entsprechenden 
Gr\"o{\ss}en durch Differentiation,
\begin{eqnarray}
p_m(x) &=& \frac{|m|}{m} \: (i \Pdd + m) \: P_{m^2}(x) \\
\label{eq:2_3a}
k_m(x) &=& \frac{|m|}{m} \: (i \Pdd + m) \: K_{m^2}(x) \\
\label{eq:2_4}
s_m(x) &=& (i \Pdd + m) \: S_{m^2}(x) \spc .
\end{eqnarray}

Die St\"orungsrechnung wird f\"ur den Diracoperator
\Equ{2_0}
i \Pdd + {\cal{B}}
\EndEqu
durchgef\"uhrt. In erster Ordnung erh\"alt man f\"ur die Eigenr\"aume
$\tilde{p}_m, \tilde{k}_m$ des Diracoperators
\begin{eqnarray}
\label{eq:2_28}
\tilde{p}_m &=& p_m \:-\: p_m \:{\cal{B}}\: s_m \:-\: s_m 
\:{\cal{B}}\: p_m \:+\: {\cal{O}}({\cal{B}}^2) \\
\label{eq:a1_100}
\tilde{k}_m &=& k_m \:-\: k_m \:{\cal{B}}\: s_m \:-\: s_m 
\:{\cal{B}}\: k_m \:+\: {\cal{O}}({\cal{B}}^2) \spc .
\end{eqnarray}
F\"ur die St\"orungsrechnung h\"oherer Ordnung braucht man die 
avancierte und retardierte Greensfunktion
\begin{eqnarray}
\label{eq:2_f7}
s_m^\vee &=& s_m \:+\: i \pi \: k_m \\
\label{eq:2_f8}
s_m^\wedge &=& s_m \:-\: i \pi \: k_m \spc .
\end{eqnarray}
Die formale St\"orungsrechnung kann mit der unit\"aren 
Transformation
\Equ{2_69}
U \;=\; \int_{\sR \cup i \sR} dm \; \sum_{l=0}^\infty \:
	(-s_m \: {\cal{B}})^l \: p_m
\EndEqu
durchgef\"uhrt werden. Allerdings treten bei dieser St\"orungsrechnung 
nichtlokale Linienintegrale auf; um sie zu vermeiden, muss man mit der
nichtunit\"aren St\"ortransformation
\begin{eqnarray}
V &=& \int_{\sR \cup i \sR} dm \;
	\sum_{l=0}^\infty \: (-1)^l \!\!\! \sum_{ \scriptsize
	\begin{array}{cc} \scriptsize Q \in {\cal{P}}(l) , \\
		\scriptsize \# Q \; {\mbox{gerade}}
	\end{array} }  \frac{(\#Q-1)!!}{(\#Q/2)! \cdot 2^{\#Q/2} }
	\; (i \pi)^{\#Q} \nonumber \\
\label{eq:2_62}
&& \hspace*{1cm} \times \;
	C_m(Q,1) \: {\cal{B}} \: C_m(Q,2) \cdots C_m(Q,l-1) \: {\cal{B}}
	\: C_m(Q,l) \: {\cal{B}} \: p_m \spc
\end{eqnarray}
arbeiten, dabei ist
\[ C_m(Q, n) \;=\; \left\{ \begin{array}{ll}
	k_m & {\mbox{falls $n \in Q$}} \\
	s_m & {\mbox{falls $n \not \in Q$}} \end{array} \right.
	\;\;\;, \spc Q \subset \N \spc . \]
(Zur Vollst\"andigkeit sei erw\"ahnt, dass in der Arbeit
hep-th/9705006 eine etwas andere St\"ortransformation verwendet
wird. F\"ur die Rechnungen in dieser Arbeit spielt dieser Unterschied 
aber keine Rolle, weil die St\"orungsrechnung direkt mit (\ref{eq:2_69})
durchgef\"uhrt wird.)
Die gest\"orten Eigenr\"aume sind dann durch
\begin{eqnarray}
	\tilde{p}_m & = & V \: p_m \: V^*
	\label{eq:2_62a}  \\
	\tilde{k}_m & = & V \: k_m \: V^*
	\label{eq:2_123a}
\end{eqnarray}
definiert. F\"ur die gest\"orten Greensfunktionen hat man
\Equ{2_t2}
\tilde{s}_m^\vee \;=\; \sum_{k=0}^\infty \left(- s_m^\vee \: {\cal{B}}
	\right)^k \: s_m^\vee \;\;\;,\spc
\tilde{s}_m^\wedge \;=\; \sum_{k=0}^\infty \left(- s_m^\wedge \: {\cal{B}}
	\right)^k \: s_m^\wedge \spc .
\EndEqu
Der Eigenraum $\tilde{k}_m$ l\"a{\ss}t sich auch direkt mit den 
gest\"orten Greensfunktionen ausdr\"ucken,
\begin{eqnarray}
\label{eq:2_tm}
\tilde{k}_m &=& \frac{1}{2 \pi i} \: \left(\tilde{s}_m^\vee -
	\tilde{s}_m^\wedge \right) \spc .
\end{eqnarray}
Schlie{\ss}lich hat $\tilde{k}$ auch die explizitere Form
\begin{eqnarray}
V \: k_m \: V^* &=& \sum_{l_1, l_2 =0}^\infty \; (-1)^{l_1+l_2}
\sum_{ \scriptsize
	\begin{array}{cc} \scriptsize Q_1 \in {\cal{P}}(l_1) , \\
		\scriptsize \# Q_1 \; {\mbox{gerade}}
	\end{array} }
\sum_{ \scriptsize
	\begin{array}{cc} \scriptsize Q_2 \in {\cal{P}}(l_2) , \\
		\scriptsize \# Q_2 \; {\mbox{gerade}}
	\end{array} } \!\!\!
	(i \pi)^{\#Q_1 + \#Q_2} \; c\left( \frac{\#Q_1}{2} \right)
	\: c \left( \frac{\#Q_2}{2} \right) \nonumber \\
\label{eq:2_65}
&& \hspace*{-2cm} \times \; C_m(Q_1, 1) \: {\cal{B}} \cdots {\cal{B}} \:
	C_m(Q_1, l_1) \:
	{\cal{B}} \: k_m \: {\cal{B}} \: C_m(Q_2, 1) \: {\cal{B}} \cdots
	{\cal{B}} \: C_m(Q_2, l_2) \spc .
\end{eqnarray}

\appendix

\chapter{St\"orungsrechnung f\"ur $k_0$ im Ortsraum}
\label{anh1}
In diesem Kapitel werden wir die Distribution $\tilde{k}_0$ f\"ur
verschiedene St\"orungen des Di\-rac\-ope\-ra\-tors perturbativ in erster
Ordnung berechnen.

Die Distribution $\tilde{k}_m$ ist in erster Ordnung durch \Ref{a1_100}
gegeben, dabei ist ${\cal{B}}$ die St\"orung des Diracoperators \Ref{2_0}.
Im Impulsraum k\"onnen wir diese Gleichung mit Hilfe von \Ref{2_3a}, \Ref{2_4}
direkt auswerten
\begin{eqnarray}
\label{eq:a1_1}
   \tilde{k}_m (p,q) &=& \frac{|m|}{m} (p \slsh + m) \: \delta(p^2-m^2)
      \epsilon(p^0) \: \delta^4(p-q) \nonumber \\
   && - \lambda \frac{|m|}{m} \left(
      (p \slsh + m) \: \delta(p^2-m^2) \: \epsilon(p^0) \; {\cal{B}}_{p,q} \;
	 \frac{1}{q \slsh - m}   \right. \nonumber\\
   && \spc + \left.
      \frac{1}{p \slsh - m} \; {\cal{B}}_{p,q} \; (q \slsh + m) \: \delta(q^2-m^2) \:
	 \epsilon(q^0)  \right) \spc .
\end{eqnarray}
Das ist f\"ur unsere Zwecke aber nicht ausreichend, wir m\"ussen eine explizite
Formel f\"ur $\tilde{k}_m$ im Ortsraum ableiten.
Dazu k\"onnten wir die Fouriertransformation von \Ref{a1_1} berechnen, was aber recht
aufwendig und schwierig ist. Etwas bequemer ist es, die St\"orungsrechnung
direkt im Ortsraum durchzuf\"uhren. Darum wollen wir diesen Rechenweg w\"ahlen.

Nach \Ref{2_21}, \Ref{2_22} und \Ref{2_3a}, \Ref{2_4}
haben wir im Spezialfall $m=0$
\begin{eqnarray}
   K_0(x) &=& \int \frac{d^4 k}{(2 \pi)^4} \: \delta(k^2) \: \epsilon(k^0)
      \; e^{-i k x} \nonumber \\
\label{eq:a1_75}
      &=& - \frac{i}{4 \pi^2} \: \delta(x^2)
	 \: \epsilon(x^0)    \\
   S_0(x) &=& \int \frac{d^4 k}{(2 \pi)^4} \:\frac{1}{k^2}
	    \; e^{-i k x} \nonumber \\
\label{eq:a1_76}
      &=& - \frac{1}{4 \pi} \: \delta(x^2)    \\
\label{eq:a1_101}
   (k_0 {\cal{B}} s_0 + s_0 {\cal{B}} k_0) (x,y) &=&
	(i \Pdd_x) \; (K_0 {\cal{B}} S_0 + S_0 {\cal{B}} K_0)(x,y)
      \; (i \Pdd_y)  \spc .
\end{eqnarray}
Damit treten Ausdr\"ucke der Form
\Equ{a1_2}
   (K_0 {\cal{B}} S_0 + S_0 {\cal{B}} K_0)(x,y)
\EndEqu
sowie Ableitungen von~\Ref{a1_2} nach $x$ und $y$ auf.
Ist ${\cal{B}}$ eine Funktion, so mu{\ss} in~\Ref{a1_2} ein Integral
\"uber den Schnitt der Lichtkegel um die Punkte $x$ und $y$ ausgewertet
werden:
\begin{eqnarray*}
   (K_0 f S_0 + S_0 f K_0)(x,y) &=& \frac{i}{16 \pi^3} \int d^4 z \; 
      \delta((x-z)^2) \:\delta((z-y)^2) \\
      && \spc \;\;\;\; f(z)  \left( \epsilon(x^0-z^0) + \epsilon(z^0-y^0)
      \right)
\end{eqnarray*}
Zur technischen Vorbereitung werden wir zun\"achst solche sogenannten
Lichtkegelintegrale untersuchen:

\section{Lichtkegelintegrale}
\label{lichtkegelint}
In diesem Abschnitt gehen wir recht ``elementar'' vor, berechnen also vieles
explizit in speziellen Koordinatensytemen.
Das ist aus mathematischer Sicht sicher nicht der eleganteste Zugang, wir
kommen so aber direkter ans Ziel.
F\"ur eine systematische Behandlung von Randwertproblemen der Wellengleichung
siehe \cite{Fl}.

Wir bezeichnen im folgenden den Minkowski-Raum mit $M$ und den Rand bzw. das
Inneren und \"Au{\ss}ere des Lichtkegels mit
\begin{eqnarray*}
   \Li &=& \left\{ y \in M \:|\: \bra y,y \ket=0 \right\} \\
   \I &=& \left\{ y \in M \:|\: \bra y,y \ket>0 \right\} \\
   \Ra &=& \left\{ y \in M \:|\: \bra y,y \ket<0 \right\} \spc .
\end{eqnarray*}
F\"ur den oberen und unteren Lichtkegel verwenden wir die Schreibweise
\[ \begin{array}{rclcrcl}
   	\Li^\vee &=& \left\{ y \in \Li \;|\; y^0 > 0 \right\}   &\;\;\;,\spc&
   	\Li^\wedge &=& \left\{ y \in \Li \;|\; y^0 < 0 \right\}  \\
   	\I^\vee &=& \left\{ y \in \I \;|\; y^0 > 0 \right\}	&\;\;\;,\spc&
   	\I^\wedge &=& \left\{ y \in \I \;|\; y^0 < 0 \right\} \spc .
\end{array} \]

\begin{Def}
   Definiere die gem\"a{{\ss}}igten Distributionen $l^\vee$ und $l^\wedge$
   durch
\begin{eqnarray*}
   l^\vee(y) &=& \delta(y^2) \: \Theta(y^0)  \\
   l^\wedge(y) &=& \delta(y^2) \: \Theta(-y^0)  \spc .
\end{eqnarray*}
\end{Def}

Wir wollen das Lichtkegelintegral einer Funktion $f$ als das Integral
von $f$ \"uber den Schnitt des nach oben ge\"offneten Lichtkegels um
den Ursprung mit dem nach unten ge\"offneten Lichtkegel um $y$
eingef\"uhren,
also in symbolischer Schreibweise
\[   (\lint f)(y) \;=\; \int d^4z \; l^\wedge(z-y) \: l^\vee(z) \; f(z)
      \spc . \]
Da zun\"achst nicht klar ist,
ob dieses Integral punktweise existiert, soll
$\slint f$ als Distribution definiert werden.

F\"ur $f \in C^\infty_c(M)$ ist $l^\vee f$ eine Distribution mit kompaktem
Tr\"ager, und wir k\"onnen
\Equ{a1_26}
   \lint f \;=\; l^\vee * (f l^\vee)
\EndEqu
setzen\footnote{Die Faltung
   einer Distribution $a$ mit
   einer Distribution $b$, die kompakten Tr\"ager besitzt, ist definiert
   durch
\[    (a * b)(g) \;=\; a (\hat{b} * g) \;\;\;, \spc g \in C^\infty_c(M)
	 \;\;\;, \]
   wobei $\hat{b}(x) = b(-x)$. Diese Definition ist sinnvoll, weil
   $\hat{b} * g \in C^\infty_c(M)$ ist.}.

Falls $f$ keinen kompakten Tr\"ager besitzt, m\"ussen wir eine etwas andere
Definition verwenden:
W\"ahle eine Funktion $\eta \in C^\infty_c(\R)$ mit
\Equ{a1_130}
      \eta_{|[-0.5, 0.5]}=1 \spc 0 \leq \eta \leq 1
	 \spc \mbox{supp } \eta \subset [-1,1] \spc ,
\EndEqu
und setze (f\"ur beliebiges $\varepsilon > 0$)
\Equ{a1_24}
    \eta_\varepsilon(x) = \eta(\frac{t-r}{\varepsilon}) \spc .
\EndEqu
F\"ur $g \in C^\infty_c(M)$ ist die Funktion
\begin{eqnarray*}
   h(z) &=& (\hat{f} \hat{l}^\vee) * g (z)   \\
	&=& \int d^4 y \; f(y-z) \; l^\vee(y-z) \: g(y)  \spc
	 \in C^\infty(M)
\end{eqnarray*}
im Halbraum $t \geq 0$ lediglich auf einer beschr\"ankten Menge von Null
verschieden. Daher besitzt $h \: \eta_\varepsilon$ kompakten Tr\"ager, und
man kann $\slint f$ durch
\[   (\lint f)(g) \;=\; l^\vee (h \: \eta_\varepsilon)     \]
sinnvoll definieren.
Beachte, da{\ss} diese Definition f\"ur $f \in C^\infty_c(M)$ mit~\Ref{a1_26}
\"ubereinstimmt. In diesem Fall ist n\"amlich der Tr\"ager von $h$ kompakt, aus
${\mbox{supp }} l^\vee = \Li^\vee$ und $\eta_{\varepsilon | \Li^\vee} \equiv 1$
folgt
\[  l^\vee(\eta_\varepsilon h) \;=\; l^\vee(h) \;=\; \left( l^\vee *
	(f l^\vee) \right) (g) \spc . \]
Als abk\"urzende Schreibweise verwenden wir
\Ref{a1_26} auch f\"ur allgemeines $f \in C^\infty(M)$.

\begin{Def}
F\"ur $f \in C^\infty(M)$ hei{\ss}t die Distribution
\begin{eqnarray}
\label{eq:a1_lint}
      (\lint f) (y) &=& l^\vee * (f l^\vee) (y)
\end{eqnarray}
   das {\bf Lichtkegelintegral} von $f$.
\end{Def}

\begin{Satz}
\label{a1_lemma1}
   $\slint f$ ist eine Funktion, die au{\ss}erhalb von $\I^\vee$ verschwindet.
   Auf $\I^\vee$ ist $\slint f$ glatt und harmonisch,
\[ \Box_y({\scriptstyle \Diamond \hspace{-0.9 em}} \int
 f)(y)=0 \spc {\mbox{f\"ur $y \in \I^\vee$.}} \]
   F\"ur $z \in \Li^\vee$ gilt
\Equ{a1_8}
      \lim_{\I^\vee \ni y \rightarrow z} (\lint f)(y) \;=\;
	 \frac{\pi}{2} \int_0^1 f(\alpha z) \: d \alpha \spc .
\EndEqu
   Das Verhalten von $\slint f$ kann also in der N\"ahe des Lichtkegels
   durch ein Linienintegral beschrieben werden.
\end{Satz}
{\Beweis}
F\"ur $g \in C^\infty_c(M)$ hat man nach Definition von $\slint f$:
\begin{eqnarray}
\label{eq:a1_39}
  ( \lint f) (g) &=& \int d^4z \; l^\vee(z) \int d^4y \; l^\vee(y-z)
	\: f(y-z) \; g(y) \\
\label{eq:a1_40}
   &=& \int d^4z \: l^\vee(z) \int d^4y \; l^\vee(y) \: f(y) \: g(y+z)
\end{eqnarray}
In~\Ref{a1_40} geht nur der Funktionswert $g(y+z)$ mit $y, z \in \Li^\vee$
ein. Aus $y, z \in \Li^\vee$ folgt $(y+z) \in \I^\vee \cup \Li^\vee$, daher
verschwindet $\slint f$ au{\ss}erhalb von $\I^\vee \cup \Li^\vee$.

Sei nun $\mbox{supp } g \subset \I^\vee$. In~\Ref{a1_39} kann das Integral
\"uber $z$ f\"ur festes $y \in \mbox{supp } g$ ausgef\"uhrt werden, es f\"uhrt
auf das Integral von $f$ \"uber ein zweidimensionales
Ellipsoid\footnote{Um die
  Normierungsfaktoren zu bestimmen, kann man den Fall $y^0=t_0$,
  $y^\alpha=0$ betrachten:
\begin{eqnarray*}
   (\lint f)(y) &:=& \int d^4z \; l^\vee(y-z) \: l^\vee(z) \: f(z) \\
    &=& \int d^4z \; \delta((y-z)^2) \: \delta(z^2)
	\; \Theta(y^0-z^0) \: \Theta(z^0) \; f(z)  \\
    &=& \int_{-\infty}^{\infty} dt \int_0^\infty r^2 dr \int_{S^2} d \omega
	\; \delta((t_0-t)^2-r^2) \: \delta(t^2-r^2) \: \Theta(t^0-t)
	\: \Theta(t) \; f(t,r,\omega) \spc ,
\end{eqnarray*} 
  wobei $\omega$ die \"ublichen
  Winkelkoordinaten $(\vartheta, \varphi)$ bezeichnet.
\begin{eqnarray}
   &=& \int_{-\infty}^\infty dt \int_{S^2} d \omega \; \frac{|t|}{2}
	\; \delta(-2 t_0 t + t_0^2) \; f(t, |t|, \omega)  \nonumber \\
\label{eq:a1_6}
   &=& \frac{1}{8} \int_{S^2} d \omega \; f(\frac{t_0}{2}, \frac{|t_0|}{2},
	\omega)
\end{eqnarray}
Die Normierung wurde also so gew\"ahlt, da{\ss} $(\flint 1)(y) = \pi/2$ f\"ur
$y \in \I^\vee$.}.
Da dieses Integral gleichm\"a{\ss}ig in $y$ beschr\"ankt ist, kann man die
$y$- und $z$-Integration vertauschen:
\[  (\lint f)(g) \;=\; \int d^4z \; (\lint f)(y) \; g(y)  \]
mit
\begin{eqnarray}
\label{eq:a1_41}
    (\lint f)(y) &=& \int d^4z \; l^\wedge(z-y) \: \l^\vee(z) \: f(z) \\
	&=& \int d^4z \; \delta((z-y)^2) \: \delta(z^2) \: f(z) \spc .
	\nonumber 
\end{eqnarray}
Also ist die Distribution $\slint f$ auf $\I^\vee$ eine Funktion.

Da das Ellipsoid, \"uber das in~\Ref{a1_41} integriert werden mu{\ss},
differenzierbar von $y$ abh\"angt, ist $\slint f$ auf $\I^\vee$ glatt.

Es bleibt noch zu zeigen:
\begin{enumerate}
\item $\slint f$ ist auf $\I^\vee$ eine harmonische Funktion:

  Beachte zun\"achst, da{\ss}
\[  \Box l^\vee(y) \;=\; \Box l^\wedge(y) \;=\; 2 \pi \delta^4(y) \spc
	\mbox{im Distributionssinne} \spc .     \]
  W\"ahle $g \in C^\infty_c (\I^\vee)$. Es gilt
\begin{eqnarray}
  ( \lint f) (\Box g) &=& (l^\vee * f l^\vee)(\Box g) \;=\;
	(f l^\vee * l^\vee)(\Box g) \nonumber \\
  &=& \int d^4 z \; l^\vee(z) \: f(z) \int d^4y \;
	l^\wedge(z-y) \: (\Box g(y)) \nonumber \\
  &=& 2 \pi \int d^4z \; l^\vee(z) \: f(z) \int d^4y \; \delta^4(z-y) \: g(y)
   \nonumber \\
\label{eq:a1_78}
  &=& 2 \pi \int d^4 y \; l^\vee(y) \: f(y) \: g(y) \;=\; 0 \spc ,
\end{eqnarray}
da $g$ auf dem Tr\"ager von $l^\vee$ verschwindet.

  Weil $(\slint f)$ auf $\I^\vee$ glatt ist, folgt $\Box ( \slint f) (y)=0$
  punktweise.
\item Es gilt
\[  \lim_{\I^\vee \ni y \rightarrow v} (\lint f)(y) \;=\;
     \frac{\pi}{2}  \int_0^1 f(\alpha v) d \alpha \spc :    \]
  Diese Gleichung kann man anschaulich leicht einsehen: N\"ahert sich $y$
  dem Lichtkegel an, so wird das Ellipsoid, \"uber das in~\Ref{a1_41}
  integriert werden mu{\ss}, immer schmaler, so da{\ss}~\Ref{a1_41} immer besser
  durch ein Linienintegral approximiert werden kann.
  F\"ur den Beweis m\"ussen wir die Differenz zwischen dem Lichtkegel- und dem
  Linienintegral geeignet absch\"atzen.

  Sei $v \in \Li^\vee$, erweitere $v$ zu einer Basis $(v, w, a, b)$ von $M$
  mit
\[  \bra v,w \ket =1 \spc \bra a,a \ket = \bra b,b \ket = -1 \spc ,     \]
  alle anderen Skalarprodukte zwischen den Basisvektoren sollen verschwinden.

  W\"ahle f\"ur $z \in \I^\vee$ die Basisdarstellung
\Equ{a1_10}
  z \;=\; \mu v + \alpha_1 w + \alpha_2 a + \alpha_3 b \spc .
\EndEqu
  Im folgenden bezeichnen $\mu$, $\alpha_j$ stets die Komponenten von $z$ in
  dieser Basis.
  Definiere die Funktion $\tilde{f}$ durch
\[      \tilde{f}(z) \;=\; f( \mu v) \spc .     \]
  Gehe nun in zwei Schritten vor:
  \begin{description}
  \item[(a)]
$\displaystyle    \lim_{\I^\vee \ni y \rightarrow v} (\lint \tilde{f})(y) \;=\;
	 \frac{\pi}{2} \int_0^1 f(\alpha v) d \alpha \spc :  \\$
    F\"ur $y \in \I^\vee$ gilt:
\begin{eqnarray*}
  (\lint \tilde{f})(y) &=& \int d^4z \; \delta( (y-z)^2) \:
	\delta(z^2) \: \tilde{f}(z) \\
    &=& \int_{-\infty}^\infty d\mu \; f(\mu v)
	\int_{-\infty}^\infty d\alpha_1  \int_{-\infty}^\infty d\alpha_2 
	\int_{-\infty}^\infty d\alpha_3 \; \delta((y-z)^2) \; \delta(z^2) \\
    &=& \int_{-\infty}^\infty d\lambda \; f(\lambda v) \int d^4z \;
	\delta( \mu - \lambda) \; \delta((y-z)^2) \; \delta(z^2)  \\
    &=& \int_{-\infty}^\infty d\lambda \; f(\lambda v) \int d^4z \;
	\delta( \bra z,w \ket - \lambda) \; \delta((y-z)^2) \; \delta(z^2)  \\
    &=& \int_{-\infty}^\infty d \lambda \;
	f( \bra y,w \ket \lambda v) \int d^4z \;
	\delta \left( \frac{\bra z,w \ket}{\bra y,w \ket} - \lambda \right) \;
	\delta((y-z)^2) \; \delta(y^2)
\end{eqnarray*}
    Wir wenden Lemma~\ref{lemma_schnitt} an, da{\ss} im Anschlu{\ss} an diesen
    Satz bewiesen wird.
\[  \;=\; \frac{\pi}{2} \int_0^1 f( \bra y,w \ket \alpha v) \: d \alpha \]
    Im Grenzfall $y \rightarrow v$ hat man $\bra y,w \ket \rightarrow 1$ und
    somit
\[   (\lint f)(y) \;\rightarrow\;   \frac{\pi}{2}
      \int_0^1 f( \alpha v) \: d\alpha    \spc .                     \]
  \item[(b)]
$\displaystyle  \lim_{\I^\vee \ni y \rightarrow v} \left( \lint |f - \tilde{f}| \right)
	\;=\; 0  \spc : \\$
    Sei $\varepsilon > 0$ vorgegeben, $y \in \I^\vee$. Es gibt es ein
    $\delta>0$, so da{\ss}
\Equ{a1_17}
    | f(\mu v + \alpha_1 w + \alpha_2 a + \alpha_3 b) \;-\; f(\mu v)|
	\;<\; \frac{2 \varepsilon}{\pi} \spc ,
\EndEqu
    falls $|\mu|<2$, $|\alpha_1|, \alpha_2^2, \alpha_3^2 < \delta$.
    Wir wollen nun in der Gleichung 
\Equ{a1_9}
    \left( \lint |f-\tilde{f}| \right)(y) \;=\; \int d^4z \; \delta((y-z)^2)
	\: \delta(z^2) \; |f(z) - \tilde{f}(z)|
\EndEqu
    die Differenz $|f(z)-\tilde{f}(z)|$ absch\"atzen.

    In~\Ref{a1_9} geht nur der Funktionswert von $f-\tilde{f}$ an den
    Punkten $z$ mit
\Equ{a1_11}
	z^2=0 \spc (y-z)^2 =0
\EndEqu
    ein. In einer Basisdarstellung von $z$, $y$ in der Form~\Ref{a1_10}
    bzw.
\Equ{a1_15}
	y \;=\; \tau v + \beta_1 w + \beta_2 a + \beta_3 b
\EndEqu
    haben die Bedingungen~\Ref{a1_11} die Form
\begin{eqnarray}
\label{eq:a1_12}
    2 \mu \alpha_1 - \alpha_2^2 - \alpha_3^2 &=& 0 \\
\label{eq:a1_13}
    2 (\tau - \mu) (\beta_1 - \alpha_1) - (\beta_2 - \alpha_2)^2
	- (\beta_3 - \alpha_3)^2 &=& 0 \spc .
\end{eqnarray}
    Da wir den Grenzfall $y \rightarrow v$ betrachten wollen,
    kann man
\Equ{a1_19}
   \frac{1}{2} < \tau < 2  \spc
\EndEqu
   annehmen.

   Aus~\Ref{a1_12}, \Ref{a1_13} und~\Ref{a1_19} folgt
   $0 \leq \mu \leq \tau$, $\alpha_1 \geq 0$. Dies ist auch direkt
   einsichtig, weil sich die Lichtkegel nur in dem Bereich zwischen
   den Null-Hyperebenen $\mu=0$, $\mu=\tau$ schneiden.
   
    Aufl\"osen von~\Ref{a1_12} nach $\alpha_1$ und Einsetzen in~\Ref{a1_13}
    liefert
\[  2 (\tau - \mu) \beta_1 - \frac{\tau}{\mu} ( \alpha_2^2 + \alpha_3^2 )
    - \beta_2^2 + 2 \beta_2 \alpha_2 - \beta_3^2 + 2 \beta_3 \alpha_3
	\;=\; 0  \]
    Mit Hilfe der Ungleichung
\[      2 \beta_j \alpha_j \; \leq \; \frac{2 \mu}{\tau} \beta_j^2 +
	\frac{\tau}{2 \mu} \alpha_j^2   \spc , j=1,2    \]
    erh\"alt man
\begin{eqnarray*}
    0 &\leq& 2 (\tau - \mu) \beta_1 - \frac{\tau}{2 \mu} (\alpha_2^2
    + \alpha_3^2) + \left( \frac{2 \mu}{\tau} - 1 \right) (\beta_2^2 +
	\beta_3^2)      \\
      &\leq& 4 (\beta_1 + \beta_2^2 + \beta_3^2) - \frac{1}{4 \mu}
	(\alpha_2^2 + \alpha_3^2)    \\
      &\leq& 4 (\beta_1 + \beta_2^2 + \beta_3^2) - \frac{1}{8}
	(\alpha_2^2 + \alpha_3^2)
\end{eqnarray*}
    und nach Einsetzen in~\Ref{a1_12}
\[      \alpha_1 \;\leq\; 8 (\beta_1 + \beta_2^2 + \beta_3^2)   \spc . \]
    W\"ahlt man $y$ so dicht an $v$, da{\ss} in der Basisdarstellung~\Ref{a1_15}
\Equ{a1_16}
	\beta_1, \beta_2^2, \beta_3^2 < \frac{\delta}{32}       \spc ,
\EndEqu
    so folgt also
\[      |\alpha_1|, \alpha_2^2, \alpha_3^2 \;<\; \delta \spc .  \]

    Nun kann man~\Ref{a1_9} unter
    Verwendung von~\Ref{a1_17} weiter absch\"atzen:
\begin{eqnarray*}
    \left( \lint |f - \tilde{f}| \right) (y) &\leq&
	\frac{2 \varepsilon}{\pi} ( \lint 1)(y)  \;=\; \varepsilon
\end{eqnarray*}
  \end{description}
\item $\slint f$ ist eine Funktion: \\
   Definiere die Funktion $h$ durch
\[    h(y) \;=\;  \left\{ \begin{array}{ll}
		  (\slint f)(y) & \mbox{f\"ur } y \in \I^\vee   \\
		  0            & \mbox{sonst}
		  \end{array} \right.           \]
   Offensichtlich stimmen $\slint f$ und $h$ auf $M \setminus
   \Li^\vee$ \"uberein.
   Wir wollen zeigen, da{\ss} sogar $\slint f = h$ im Distributionssinne.
   
   Sei dazu $g \in C^\infty_c(M)$ vorgegeben. Da $h$ auf $\I^\vee$ glatt
   und auf $\Li^\vee$ lokal beschr\"ankt ist, ist das Integral
\[    \int h(z) \: g(z) \: d^4z  \]
   wohldefiniert und endlich.
   W\"ahle $\eta_\varepsilon$ wie in~\Ref{a1_24} mit $\varepsilon>0$.
   In einer Umgebung von $\Li^\vee$ ist $\eta_\varepsilon \equiv 1$.
   Daher gilt:
\[    (\lint f)(g) - \int h g \;=\; (\lint f)(g \eta_\varepsilon) -
	 \int h (g \eta_\varepsilon)         \]
   Wegen
\[    \lim_{\varepsilon \rightarrow 0} \int h (g \eta_\epsilon) \;=\; 0 \]
   gen\"ugt es zu zeigen, da{\ss}
\[    (\lint f)(\eta_\varepsilon \: g) \;\rightarrow\;0  \spc, \]
   denn dann folgt $(\slint f)(g) = \int h g$ f\"ur beliebiges
   $g \in C^\infty_c(M)$ und somit die Behauptung.

   Forme dazu~\Ref{a1_40} weiter um:
\[ (\lint f) (\eta_\varepsilon g) \;=\: \int \frac{d^3 \vec{z}}{2 |\vec{z}|}
     \int \frac{d^3 \vec{y}}{2 |\vec{y}|} \; f(y) \: 
     (\eta_\varepsilon g)(y+z)          \]
   mit $y=(|\vec{y}|,\vec{y})$, $z=(|\vec{z}|, \vec{z})$.

   Da $g$ kompakten Tr\"ager besitzt, gibt es $R>0$ (unabh\"angig von
   $\varepsilon$), so da{\ss}
\[  \left| (\lint f)(\eta_\varepsilon g) \right|
    \;\leq\;  \int_{B_R(0)} d^3 \vec{z} \int_{B_R(0)} d^3 \vec{y} \;
    \frac{1}{4 \: |\vec{z}| \: |\vec{y}|} \; |f(y) \: (\eta_\varepsilon g)
      (y+z)| \spc .     \]
   Da die einzelnen Integrale \"uber $y$, $z$ 
   fast \"uberall existieren und bez\"uglich der jeweils
   anderen Variablen integrierbar sind, erh\"alt man nach dem Satz
   von Fubini
\[  \;=\; \int_{B_R(0) \times B_R(0)} d^3 \vec{y} \: d^3 \vec{z} \;
    \frac{1}{4 \: |\vec{z}| \: |\vec{y}|} \; |f(y) \:
      (\eta_\varepsilon g)(y+z)|
       \spc .   \]
   Im Grenzfall $\varepsilon \rightarrow 0$ hat man $(\eta_\varepsilon g)
   \rightarrow 0$ auf $\R^6 \setminus \{ (\vec{y}, \vec{z}) \;|\; y = \lambda z,
	\lambda \in \R \}$, also $(\eta_\varepsilon g) \rightarrow 0$ fast
   \"uberall.

   Nach Lebesgues monotonem Konvergenzsatz folgt
\[   \lim_{\varepsilon \rightarrow 0} \; \int_{B_R(0) \times B_R(0)}
    d^3 \vec{y} \: d^3 \vec{z} \;
    \frac{1}{4 \: |\vec{z}| \: |\vec{y}|} \; |f(y) \:
      (\eta_\varepsilon g)(y+z)|
       \;=\; 0 \spc ,   \]
   was den Beweis abschlie{\ss}t.   \QED
\end{enumerate}

Es bleibt noch das folgende kleine Lemma nachzutragen:

\begin{Lemma}
\label{lemma_schnitt}
  Sei $y \in \I^\vee$ und $w \in \Li^\vee$. Dann ist f\"ur die
  Distribution
\[    f(z) \;=\; 
    \delta \left( \frac{\bra z,w \ket}{\bra y,w \ket} - \lambda \right) \]
  der Ausdruck $\slint f$ sinnvoll definiert, und es gilt
\Equ{a1_5}
  ( \lint f )(y)
       \;=\; \left\{ \begin{array}{ll}
			\frac{\pi}{2} & \mbox{f\"ur $0<\lambda<1$} \\
			0   & \mbox{f\"ur $\lambda<0 \; \vee \; \lambda > 1$}
		      \end{array} \right. \spc .
\EndEqu
\end{Lemma}
  Gleichung~\Ref{a1_5} ist das Integral \"uber den Schnitt der Lichtkegel
  um die Punkte $0$, $y$ mit der Null-Hyperebene $\bra z,w \ket = \lambda 
  \: \bra y,w \ket$. 

{\Beweis}
  W\"ahle ein Bezugssystem\footnote{Ein Bezugssystem ist ein Koordinatensystem
  mit $g_{ij} = \eta_{ij}$.}, so da{\ss} $y^0=t$, $y^\alpha=0$, $w^0=w^1=w$,
  $w^2=w^3=0$.
  Nach~\Ref{a1_6} gilt:
\begin{eqnarray*}
  \lint \delta \left( \frac{\bra z,w \ket}{\bra y,w \ket} - \lambda \right)
   _{|y} &=& \frac{1}{8} \int_{S^2} d \omega \;
   \delta \left( \frac{\bra z,w \ket}{\bra y,w \ket} - \lambda \right) \\
   &=& \frac{1}{8} \int_{-1}^1 d \cos \vartheta \int_{0}^{2 \pi} d \varphi
      \; \delta \left( \frac{1}{2} (1 - \cos \vartheta) - \lambda \right)  \\
   &=& \frac{\pi}{2} \int_{-1}^{1} d \cos \vartheta \; \delta(  1
      - 2 \lambda - \cos \vartheta)
\end{eqnarray*}
   Hieraus folgt unmittelbar die Behauptung.
\QED
Die Funktion $\slint f$ kann nach ihrer Definitionsgleichung
\[   (\lint f)(g) \;=\; \int d^4y \; (\lint f)(y) \: g(y)
	\;\;\;\;, \; g \in C^\infty_c(M)        \]
auf einer Nullmenge beliebig abge\"andert werden.
Auf $\I \cup \Ra \cup \Li^\wedge$ kann $\slint f$ stetig gew\"ahlt werden,
wodurch die Definition eindeutig wird. Auf $\Li^\vee$ verwenden wir die
Konvention
\Equ{a1_56}
   (\lint f)(y) \;=\; \frac{\pi}{2} \int_0^1 f(\alpha y) \: d \alpha
	\;\;\;\;, \; y \in \Li^\vee \spc .
\EndEqu
Dadurch ist $\slint f$ auf $\I^\vee \cup \Li^\vee$ stetig.
\\[1em]
Wir wollen nun f\"ur die partiellen Ableitungen von $\slint f$ explizite
Formeln abgeleiten. Ein erster Schritt ist das folgende Lemma.
\begin{Lemma}
\label{a1_lemma5}
F\"ur $y \in \I^\vee$ gilt
\Equ{a1_abl1}
   \partial_j (\lint f)(y) \;=\; \left( \lint h_j^{(y)} \right) (y)
\EndEqu
mit
\Equ{a1_hjy}
   h_j^{(y)}(z) \;=\; \partial_j f_{|z} + 2 \: \frac{z_j (z-y)^l}{y^2}
	\: \partial_l f_{|z}
   - 2 \: \frac{(y - 2 z)_j}{y^2} \: f_{|z} \spc .
\EndEqu
\end{Lemma}
{\Beweis}
   Nehme zun\"achst an, da{\ss} $y$ in der $tx$-Ebene liegt, also
   $y = (t_0, x_0,0,0)$, weiterhin sei $t_0 > 0$.
   
   Nach einer Lorentztransformation in der $t$- und $x$-Koordinaten
   gem\"a{\ss}
\begin{eqnarray*}
   \tilde{t} &=& \frac{1}{\gamma} (t - \beta x)   \\
   \tilde{x} &=& \frac{1}{\gamma} (x - \beta t)
\end{eqnarray*}
   mit $\beta = x_0/t_0$, $\gamma = \sqrt{1-\beta^2} $  hat man
   $y=(2 \tau, 0,0,0)$
   mit $\tau = \gamma t_0 / 2$.
   Wir bezeichnen die Funktion $f$ in den gestrichenen Koordinaten mit
   $\tilde{f}$.
   Eine Rechnung in Polarkoordinaten gem\"a{\ss}~\Ref{a1_6} liefert
\[ (\lint f)(y) \;=\; \frac{1}{8} \int_{S^2} d \omega \;
      \tilde{f} \left( \tau , \tau \cos \vartheta,
      \tau \sin \vartheta \cos \varphi,
      \tau \sin \vartheta \sin \varphi \right)   \]
   oder in den ungestrichenen Koordinaten
\[    \;=\;  \frac{1}{8} \int_{S^2} d \omega \;
      f \left( \frac{1}{\gamma} ( \tau + \beta \tau \cos \vartheta ),
	 \frac{1}{\gamma} ( \tau \cos \vartheta + \beta \tau ),
	 \tau \sin \vartheta \cos \varphi,
	 \tau \sin \vartheta \sin \varphi \right)   \spc .   \]
    Nun k\"onnen wir die partiellen Ableitungen nach $t_0$, $x_0$ an der
    Stelle $x_0=0$ direkt ausrechnen. In
    Polarkoordinaten $(t,r, \vartheta, \varphi)$ erh\"alt man
\begin{eqnarray}
\label{eq:a1_30}
   \frac{\partial}{\partial t_0} ( \lint f )(y)_{|x_0=0} &=&
      \frac{1}{16} \int_{S^2} d \omega \;
      ( \partial_t + \partial_r ) f(\tau, \tau, \omega)   \\
\label{eq:a1_31}
   \frac{\partial}{\partial x_0} ( \lint f )(y)_{|x_0=0} &=&
      \frac{1}{16} \int_{S^2} d \omega \;
      ( \cos \vartheta \; \partial_t + \partial_x ) f(\tau, \tau, \omega)
      \spc .
\end{eqnarray}
Wir wollen Gleichung~\Ref{a1_31} weiter umformen. Dabei nutzen wir aus, da{\ss}
die zu $S^2$ tangentiale Komponente der Ableitung partiell integriert werden kann:
Wir haben
\begin{eqnarray*}
   \int_{S^2} d \omega \; \partial_x f
   &=& \int_0^{2 \pi} d \varphi \int_{-1}^{1}
    d \alpha \left( \frac{x}{r} \frac{\partial}{\partial r} +
      \frac{1-\alpha^2}{r} \frac{\partial}{\partial \alpha} \right) f \\
   &=& \int_{S^2} d \omega \left( \frac{x}{r}
      \frac{\partial}{\partial r}
      + \frac{2 x}{r^2} \right) f
\end{eqnarray*}
und somit
\begin{eqnarray}
   \frac{\partial}{\partial x_0} ( \lint f )(y) &=&
      \frac{1}{16} \int_{S^2} d \omega \;
      \left( 2 \partial_x + \cos \vartheta \; \partial_t
	 -  \frac{x}{r} \partial_r
	 -  \frac{2 x}{r^2} \right) f(\tau, \tau, \omega) \nonumber \\
\label{eq:a1_32}
&=& \frac{1}{16} \int_{S^2} d \omega \;
      \left( 2 \partial_x + \frac{x}{r} (\partial_t - \partial_r )
	 -  \frac{2 x}{r^2} \right) f(\tau, \tau, \omega)   \spc .
\end{eqnarray}
Schreibt man die Gleichungen~\Ref{a1_30} und~\Ref{a1_32} in
koordinateninvarianter Form, so erh\"alt man gerade~\Ref{a1_abl1}.

Damit ist auch das Lemma bewiesen, denn durch geeignete Wahl des
Bezugssystems kann man immer erreichen, da{\ss}
der Vektor $y$ in $t$-Richtung zeigt sowie
die Ableitungsrichtung in der $tx$-Ebene liegt.
\QED
Wir k\"onnen auch die Gleichungen~\Ref{a1_30} und~\Ref{a1_31} kovariant in der
Form~\Ref{a1_abl1} schreiben mit
\Equ{a1_43}
   h_j^{(y)} \;=\; \frac{1}{2} \partial_j f
   + \frac{ z_j (z-y)^l - (z-y)_j z^l}{y^2}
	\; \partial_l f \spc .
\EndEqu
Nat\"urlich ist~\Ref{a1_43} \"aquivalent zu~\Ref{a1_hjy}, f\"ur das weitere
Vorgehen ist jedoch~\Ref{a1_hjy} g\"unstiger.
\\[1em]
Es ist unsch\"on an Gleichung~\Ref{a1_abl1}, da{\ss} die Funktion
$h_j^{(y)}$ von $y$ abh\"angt.
W\"un\-schens\-wert w\"are eine Gleichung der Form\label{a1_ank}
\Equ{a1_36}
   \partial_j (\lint f) \;=\; (\lint h_j)
\EndEqu
mit geeigneten Funktionen $h_j$. Eine solche Relation scheint auf den
ersten Blick sinnvoll zu sein, weil in~\Ref{a1_36} sowohl die linke Seite
(als partielle Ableitung einer harmonischen Funktion) als auch die
rechte Seite harmonisch sind.

Der wesentliche Vorteil von~\Ref{a1_36}
gegen\"uber~\Ref{a1_abl1} besteht darin, da{\ss}~\Ref{a1_36} leicht iteriert werden kann,
wodurch man unmittelbar auch Formeln f\"ur die h\"oheren Ableitungen von
$\slint f$ erh\"alt.

Die folgenden drei Lemmata \ref{a1_lemma6}, \ref{a1_lemma8}, \ref{a1_lemma9}
dienen als Vorbereitung f\"ur die Konstruktion der Funktionen $h_j$ in
Satz \ref{a1_parabl}.

\begin{Lemma}
\label{a1_lemma6}
   F\"ur $y \in \I^\vee$ und die Funktion
\[       h(z) \;=\; \int_a^b \alpha^n \; f(\alpha z) \; d \alpha     \]
   (mit $a,b \in \R^+_0$, $n \in \N_0$, $f \in C^\infty(M)$) gilt
\[    (\lint h)(y) \;=\; \int_a^b \alpha^n \; (\lint f)_{|\alpha y} \;
	d \alpha  \spc . \]
\end{Lemma}  
{\Beweis}
\begin{eqnarray*}
   (\lint f)(y) &=& \int d^4 z \; l^\wedge(z-y) \: l^\vee(z)
      \int_a^b d \alpha \; \alpha^n f(\alpha z)    \\
   &=& \int_a^b d \alpha \; \alpha^n \int d^4 z \; l^\wedge(z-y) \: l^\vee(z)
      \: f(\alpha z) \\
   &=& \int_a^b d \alpha \; \alpha^n \int d^4 z \; l^\wedge(z-\alpha y)
      \: l^\vee(z) \: f(z) \\
   &=& \int_a^b \alpha^n \: (\lint f)_{|\alpha y} \; d \alpha \spc ,
\end{eqnarray*}
   wobei die Beziehung $\delta(\lambda x)=\delta(x)/\lambda$
   verwendet wurde.
\QED
\begin{Def}
\label{def_a7}
   Bei den bisherigen Definitionen war der Lichtkegel um den Ursprung
   ausgezeichnet.
   Wir bezeichnen die durch Parallelverschiebung um $x \in M$ entstehenden
   Mengen und Distributionen
   durch einen zus\"atzlichen Index $x$, also z.B.
   $\Li_x$, $\I_x$, $\Ra_x$ f\"ur den Lichtkegel um den Punkt $x$,
   entsprechend
\[    l_x^\vee(y)\;:=\;l^\vee(y-x) \spc l_x^\wedge(y)\;:=\;l^\wedge(y-x)
      \spc . \]
   Schreibe f\"ur das verschobene Lichtkegelintegral
\begin{eqnarray}
    \lint_x^y f &=& l^\vee_x * (f l_x^\vee) (y) \nonumber \\
\label{eq:a1_57}
    &=& \int d^4z \; l_y^\wedge(z) \: l_x^\vee(z) \: f(z)   \spc .
\end{eqnarray}
   Wenn die Integrationsvariable angegeben werden soll, verwenden wir f\"ur
   \Ref{a1_57} auch die Schreibweise
\[  \lint_x^y f(z) \: dz  \spc . \]
   Die bisherigen Ergebnisse \"ubertragen sich unmittelbar
   auf diesen etwas allgemeineren Fall.
\end{Def}

\begin{Lemma}
\label{a1_lemma8}
   F\"ur die Distributionsableitung von $\slint f$ gilt
\Equ{a1_45}
      \lint_x^y \partial_j f \;=\; \left( \frac{\partial}{\partial y^j} +
   \frac{\partial}{\partial x^j} \right) \lint_x^y f \spc .
\EndEqu
   F\"ur $(y-x) \in \I^\vee$ gilt~\Ref{a1_45} sogar punktweise.
\end{Lemma}
{\Beweis}
Die zweite Aussage folgt unmittelbar aus der ersten, da die linke
Seite von~\Ref{a1_45} f\"ur $(y-x) \in \I^\vee$ eine glatte Funktion ist.

F\"ur $g \in C^\infty_c(M)$ hat man nach Definition der Distributionsableitung
\begin{eqnarray*}
   (\lint_x \partial_j f)(g) &=& \int d^4z \:l^\vee(z-x) \: \partial_j f(z)
   \int d^4y \: l^\vee(y-z) \: g(y) \\
   &=& - \int d^4z \: \left( \frac{\partial}{\partial z^j} l^\vee(z-x)
      \right) f(z) \int d^4y \: l^\vee(y-z) \: g(y) \\
   && - \int d^4 z \: l^\vee(z-x) \: f(z) \int d^4y \: \left(
      \frac{\partial}{\partial z_j} l^\vee(y-z) \right) g(y) \\
   &=& \int d^4z \: \left( \frac{\partial}{\partial x^j} l^\vee(z-x)
      \right) f(z) \int d^4y \: l^\vee(y-z) \: g(y) \\
   && + \int d^4 z \: l^\vee(z-x) \: f(z) \int d^4y \: \left(
      \frac{\partial}{\partial y_j} l^\vee(y-z) \right) g(y) \\
   &=& \frac{\partial}{\partial x^j} (\lint_x f)(g) - (\lint_x f)
      (\partial_j g) \\
   &=& \int d^4y \; g(y) \: \left( \frac{\partial}{\partial x^j} +
      \frac{\partial}{\partial y^j} \right) \left( \lint_x^y f \right)
	\spc .
\end{eqnarray*}
\QED
Zwei Spezialf\"alle von Ableitungen der Lichtkegelintegrale k\"onnen wir durch
Einsetzen in~\Ref{a1_abl1} und~\Ref{a1_45} direkt ausrechnen:
\begin{eqnarray}
\label{eq:a1_46}
  (y-x)^j  \frac{\partial}{\partial y^j} \; \lint_x^y f &=&
	\lint_x^y (z-x)^j \; \partial_j f(z) \: dz \\
\label{eq:a1_47}
  (y-x)^j  \frac{\partial}{\partial x^j} \; \lint_x^y f &=&
	\lint_x^y (y-z)^j \; \partial_j f(z) \: dz
\end{eqnarray}

\begin{Lemma}
\label{a1_lemma9}
   F\"ur $(y-x) \in \I^\vee$ gilt
\Equ{a1_29}
   \lint_x^y (\Box f) \;=\; \frac{4}{|y-x|^3} \; (y-x)^j (y-x)^k \;
   \frac{\partial}{\partial y^j} \frac{\partial}{\partial x^k}
   \left( |y-x| \lint_x^y f \right)
\EndEqu
\end{Lemma}
{\Beweis}
Beachte, da{\ss} in~\Ref{a1_29} nur Ableitungen in Richtung des Vektors
$(y-x)$ auftreten.

Man kann annehmen, da{\ss} $(y-x)$ in $t$-Richtung zeigt.
Durch eine anschlie{\ss}ende Verschiebung des Koordinatensystems kann man
zus\"atzlich erreichen, da{\ss} $x=(x^0,0,0,0)$, $y=(y^0,0,0,0)$ mit $y^0 > x^0 > 0$.

Rechne nun in Polarkoordinaten und wende~\Ref{a1_6} an:
\begin{eqnarray*}
  \lint_x^y \Box f &=& \frac{1}{8} \int_{S^2} d \omega \; \Box f \left(
	\frac{y^0+x^0}{2}, \frac{y^0+x^0}{2}, \omega \right) \\
  &=& \frac{1}{8} \int_{S^2} d \omega \; \left( \partial_t^2 - \frac{1}{r^2}
	\partial_r ( r^2 \partial_r ) - \Delta_s
	\right) f \left( \frac{y^0+x^0}{2},
	\frac{y^0+x^0}{2}, \omega \right)
\end{eqnarray*}
Das Integral \"uber $\Delta_s f$ verschwindet nach dem Satz von Gau{\ss}:
\[ \;=\; \frac{1}{8} \int_{S^2} d \omega \; \left( \partial_t^2 - 
	\partial_r^2 - \frac{2}{r} \partial_r
	\right) f \left( \frac{y^0+x^0}{2},
	\frac{y^0+x^0}{2}, \omega \right)       \]
Nach Lemma~\ref{a1_lemma5} und Lemma~\ref{a1_lemma8} gilt weiterhin
\begin{eqnarray*}
(y-x)^j \frac{\partial}{\partial y^j} \; \lint_x^y &=&
  \frac{1}{8} (y^0-x^0) \int_{S^2} d \omega \; \frac{1}{2} \left( \partial_t
	+ \partial_r \right) f \\
(y-x)^j (y-x)^k \frac{\partial}{\partial y^j} \frac{\partial}{\partial x^k}
  \; \lint_x^y &=&
  \frac{1}{8} (y^0-x^0)^2 \int_{S^2} d \omega \; \frac{1}{4}
  \left( \partial_t^2 - \partial_r^2 \right) f \\
(y-x)^j \left( \frac{\partial}{\partial y^j} - \frac{\partial}{\partial x^j}
  \right) \lint_x^y f &=&
  \frac{1}{8} (y^0-x^0) \int_{S^2} d \omega \; \partial_r f \\
  &=& \frac{1}{8} (y^0-x^0)^2 \int_{S^2} d \omega \; \frac{1}{2 r} \partial_r
  f \spc ,
\end{eqnarray*}
also insgesamt
\begin{eqnarray*}
  \lint_x^y \Box f &=& \frac{4}{(y-x)^2} \left( (y-x)^j (y-x)^k
  \frac{\partial}{\partial y^j} \frac{\partial}{\partial x^k}
    - (y-x)^j (\frac{\partial}{\partial y^j} - \frac{\partial}{\partial x^j})
  \right) \lint_x^y f   \\
  &=& \frac{4}{|y-x|^3} \; (y-x)^j (y-x)^k \frac{\partial}{\partial y^j}
     \frac{\partial}{\partial x^k} \left( |y-x| \lint_x^y f \right) \spc .
\end{eqnarray*}
\QED

\begin{Satz}
\label{a1_parabl}
  F\"ur die partiellen Ableitungen von $\slint f$ auf $\I^\vee$ gilt
\Equ{a1_50}
  \partial_j \lint f \;=\; \lint h_j
\EndEqu
mit
\Equ{a1_55}
   h_j[f](y) \;=\; \partial_j f(y) - \frac{1}{2} \int_0^1
      \alpha \;\: \Box_z
      (z_j f(z))_{|z = \alpha y} \; d \alpha  \spc .
\EndEqu
\end{Satz}
{\Beweis}
  Schreibe zur Abk\"urzung $b_j(z)=z_j f(z)$.
  Nach Lemma~\ref{a1_lemma6} gilt f\"ur $(y-x) \in \I^\vee$
\[  \lint_x^y dz \; \int_0^1 d \alpha \; \alpha \: \Box b_j(\alpha z +
     (1-\alpha) x) \;=\; \int_0^1 d \alpha \; \alpha \; \lint_x^{\alpha y +
     (1-\alpha) x} \: \Box b_j  \spc . \]
  Setze $\tilde{y}=\alpha y + (1-\alpha) x$, $\tilde{x}=x$ und fasse die
  Variablen $x$, $y$ (f\"ur festes $\alpha$) als von $\tilde{x}$, $\tilde{y}$
  abh\"angige Variablen auf.
  Lemma~\ref{a1_lemma9} liefert
\[  \;=\; \int_0^1 d \alpha \; \alpha \; \frac{4}{|\tilde{y}-\tilde{x}|^3} \;
   (\tilde{y}-\tilde{x})^k (\tilde{y}-\tilde{x})^l
   \frac{\partial}{\partial \tilde{x}^k} \frac{\partial}{\partial \tilde{y}^l}
   \left( |\tilde{y}-\tilde{x}| \lint_{\tilde{x}}^{\tilde{y}} b_j \right)
	\spc . \]
  Um diese Relation in den Variablen $x$ und $y$ auszudr\"ucken, setzen wir die
  Beziehungen
\begin{eqnarray*}
  |\tilde{y}-\tilde{x}| &=& \alpha \: |y-x| \\
  \frac{\partial}{\partial \tilde{x}^k} &=& \frac{\partial}{\partial x^k}
     - \frac{1-\alpha}{\alpha} \frac{\partial}{\partial y^k} \\
  \frac{\partial}{\partial \tilde{y}^l} &=& \frac{1}{\alpha}
	\frac{\partial}{\partial y^l}
\end{eqnarray*}
  ein und erhalten
\begin{eqnarray}
  &=& \frac{4}{(y-x)^2} \int_0^1 d \alpha \; (y-x)^k \; \frac{(y-x)^l}{|y-x|}
    \left( \frac{\partial}{\partial x^k} - \frac{1-\alpha}{\alpha}
	\; \frac{\partial}{\partial y^k} \right) \nonumber \\
    && \spc \frac{\partial}{\partial y^l} \left( |y-x|
	\lint_x^{\alpha y + (1-\alpha) x} b_j \right)   \nonumber \\
  &=& \frac{4}{(y-x)^2} \int_0^1 d \alpha \; (y-x)^k
    \left( \frac{\partial}{\partial x^k} - \frac{1-\alpha}{\alpha}
	\frac{\partial}{\partial y^k} \right) \frac{(y-x)^l}{|y-x|}
	\nonumber \\
    && \spc \frac{\partial}{\partial y^l} \left( |y-x|
	\lint_x^{\alpha y + (1-\alpha) x} b_j \right) \nonumber \\
\label{eq:a1_42aa}
  &=& \frac{4}{(y-x)^2} \int_0^1 d \alpha \; (y-x)^k
    \left( \frac{\partial}{\partial x^k} - \frac{1-\alpha}{\alpha}
	\frac{\partial}{\partial y^k} \right)  \frac{d}{d \alpha} \;
	\left( \alpha \;
	\lint_x^{\alpha y + (1-\alpha) x} b_j \right) \;\; . \spc
\end{eqnarray}
Nun kann in $\alpha$ partiell integriert werden. F\"ur die Randwerte bei $\alpha=1$
erh\"alt man
\[   \frac{4}{(y-x)^2} \; (y-x)^k \frac{\partial}{\partial x^k} \lint_x^{y}
	b_j \spc . \]
Bei den Randwerten f\"ur $\alpha=0$ mu{\ss} man etwas aufpassen. Beachte dazu,
da{\ss}
\begin{eqnarray}
\lefteqn{ \lim_{\alpha \rightarrow 0} \alpha \: (y-x)^k
	\frac{\partial}{\partial x^k}
	\lint_x^{\alpha y + (1-\alpha) x} b_j } \nonumber \\
&=& \lim_{\alpha \rightarrow 0} \alpha
	\lint_x^{\alpha y + (1-\alpha) x} dz \; (y - z)^k \:
	\partial_k b_j(z) \;=\; 0 \\
\lefteqn{ \lim_{\alpha \rightarrow 0} (1-\alpha) \: (y-x)^k \:
	\frac{\partial}
	{\partial y^k} \lint_x^{\alpha y + (1-\alpha) x} b_j } \nonumber \\
&=& \lim_{\alpha \rightarrow 0} (1-\alpha) \lint_x^{\alpha y +
	(1-\alpha) x} dz \; (z-x)^k \: \partial_k b_j(z) \nonumber \\
\label{eq:a1_44b}
&=& \lim_{\alpha \rightarrow 0} (1-\alpha) \alpha \lint_x^{\alpha y +
	(1-\alpha) x} dz \; (z-x)^k \: \partial_k b_j(\alpha z + (1-\alpha) x)
	\;=\; 0 \;\;\; , \spc
\end{eqnarray}
wobei \Ref{a1_46}, \Ref{a1_47} angewendet und ausgenutzt wurde,
da{\ss} das Lichtkegelintegral eine beschr\"ankte Funktion ist.
Man erh\"alt also
\begin{eqnarray*}
\lefteqn{ \lint_x^y dz \; \int_0^1 d \alpha \; \alpha \: \Box b_j(\alpha z +
     (1-\alpha) x) } \\
&=& \frac{4}{(y-x)^2} \: (y-x)^k \frac{\partial}{\partial x^k} \:
	\lint_x^y b_j \;-\; \frac{4}{(y-x)^2} \: (y-x)^k \int_0^1 d \alpha \;
	\;\frac{1}{\alpha} \; \frac{\partial}{\partial y^k}
	\lint_x^{\alpha y - (1-\alpha) x} b_j \;\;\; .
\end{eqnarray*}
  Setze nun die spezielle Form von $b_j$ ein, wende wiederum
  Lemma~\ref{a1_lemma6} an
\begin{eqnarray*}
  &=& \frac{4}{(y-x)^2} \: (y-x)^k \frac{\partial}{\partial x^k} \:
	\lint_x^y dz \: z_j f(z) \\
  && \spc - \frac{4}{(y-x)^2} \: (y-x)^k \frac{\partial}{\partial y^k} \:
	\lint_x^y dz \: z_j \int_0^1 f(\alpha z)
	d \alpha
\end{eqnarray*}
  und verwende die Beziehungen~\Ref{a1_46} und~\Ref{a1_47}
\begin{eqnarray}
  &=& \frac{4}{(y-x)^2} \lint_x^y dz \: (y-z)^k \partial_k \left( z^j
		f(z) \right)  \nonumber \\
\label{eq:a1_52}
  && \spc - \frac{4}{(y-x)^2} \lint_x^y dz \: (z-x)^k \partial_k
	\left( z_j \int_0^1 f(\alpha z) d \alpha \right) \spc .
\end{eqnarray}
  Betrachte jetzt den Spezialfall $x=0$. Im zweiten Summanden von~\Ref{a1_52}
  kann die Umformung
\begin{eqnarray*}
  \lint_0^y dz \; z^k \partial_k \left( z_j \int_0^1 f(\alpha z) d \alpha
	\right) &=&
     \lint_0^y dz \; z_j \int_0^1  \left(1 + \alpha \frac{d}{d \alpha} \right)
	f(\alpha z) d \alpha \\
     &=& \lint_0^y dz \; z_j f(z)
\end{eqnarray*}
  angewendet werden, es ergibt sich insgesamt
\[  \lint_0^y dz \; \int_0^1 d \alpha \; \alpha \: \Box b_j(\alpha z)
     \;=\; \frac{4}{y^2} \lint_0^y dz \; \left( (y-2z)_j f(z)
	+ z_j (y-z)^k \partial_k f(z) \right) \;\;\; .  \]
  Nun kann man~\Ref{a1_50} unter Verwendung von Lemma~\ref{a1_lemma5}
  verifizieren:
\begin{eqnarray*}
  (\lint h_j)(y) &=& \lint_0^y dz \: \left( \partial_j f(z) - \frac{1}{2}
	\int_0^1 d \alpha \; \alpha \: \Box b_j(\alpha z) \right) \\
  &=& \lint_0^y dz \: \left( \partial_j f + 2 \: \frac{z_j (z-y)^l}{y^2} \:
	\partial_l f
	- 2 \: \frac{(y-2z)_j}{y^2} \: f \right) \\
  &=& \partial_j \left( \lint f \right)(y)
\end{eqnarray*}
\QED

\begin{Korollar}
\label{a1_korollar}
Mit Konvention \Ref{a1_56} ist $\slint f \in C^\infty(\I^\vee \cup \Li^\vee)$.
Die partiellen Ableitungen beliebiger Ordnung von $\slint f$ in $\I^\vee$
sind also stetig auf $\I^\vee \cup \Li^\vee$ fortsetzbar.
\end{Korollar}
{\Beweis}
Nach Satz~\ref{a1_parabl} und Satz~\ref{a1_lemma1} sind die ersten
Ableitungen auf $\I^\vee \cup \Li^\vee$ stetig fortsetzbar.

F\"ur die h\"oheren Ableitungen folgt die Behauptung durch Induktion; beachte
dazu, da{\ss} $h_j$ nach~\Ref{a1_55} eine glatte Funktion ist.
\QED

Wir wollen nun durch Iteration von \Ref{a1_50} einen expliziten Ausdruck f\"ur
die zweiten Ableitungen der Lichtkegelintegrale ableiten. Als Vorbereitung
ben\"otigen wir folgendes Lemma, mit dem sich geschachtelte Linienintegrale
vereinfachen lassen:
\begin{Lemma}
\label{integr_umf}
Sei $f \in C^\infty([0,1])$. Dann gilt f\"ur $n, m \geq 0$
\[  \int_0^1 d\alpha \; \alpha^m \int_0^1 d\beta \; \beta^n \: f(\alpha \:
	\beta) \;=\; - \frac{1}{m-n} \: \int_0^1 d\alpha \;
	(\alpha^m-\alpha^n) \: f(\alpha)        \spc .  \]
\end{Lemma}
{\Beweis}
Durch partielle Integration und Anwendung der Beziehung
$\alpha \: \partial_\alpha f(\alpha \beta) = \beta \: \partial_\beta f(\alpha
\beta)$ erh\"alt man
\begin{eqnarray*}
\lefteqn{ (n+1) \int_0^1 d\alpha \; \alpha^m \int_0^1 d\beta \; \beta^n \;
	f(\alpha \beta) \;=\; \int_0^1 d\alpha \; \alpha^m \int_0^1
	d\beta \; f(\alpha \beta) \: \frac{d}{d\beta} \beta^{n+1} } \\
&=& \int_0^1 d\alpha \; \alpha^m \; f(\alpha) \;-\; \int_0^1 d\alpha \;
	\alpha^m \int_0^1 d\beta \; \beta^{n+1} \:
	\frac{\partial}{\partial \beta} f(\alpha \beta) \\
&=& \int_0^1 d\alpha \; \alpha^m \; f(\alpha) \;-\; \int_0^1 d\beta \;
	\beta^n \int_0^1 d\alpha \; \alpha^{m+1} \:
	\frac{\partial}{\partial \alpha} f(\alpha \beta) \\
&=& \int_0^1 d\alpha \; (\alpha^m-\alpha^n) \; f(\alpha) \;+\;
	(m+1) \int_0^1 d\beta \;
	\beta^n \int_0^1 d\alpha \; \alpha^{m} \:
	f(\alpha \beta) \spc . \\
\end{eqnarray*}
\QED
Als kleine Anwendung berechnen wir die Randwerte von $\partial_j \slint f$ auf
dem Lichtkegel:
\begin{Lemma}
\Equ{a1_41a}
\lim_{\I^\vee \ni y \rightarrow z \in \Li^\vee} \partial_j \left(
	\lint f \right) (y) \;=\;
	\frac{\pi}{2} \int_0^1 \alpha \: \partial_j f(\alpha z) \; d\alpha
	+ \frac{\pi}{4} \int_0^1 (\alpha^2-\alpha) \: y_j \;
	(\Box f)_{|\alpha y} \; d \alpha
\EndEqu
\end{Lemma}
{\Beweis}
\begin{eqnarray}
\lefteqn{\lim_{\I^\vee \ni y \rightarrow z \in \Li^\vee} \partial_j \left(
	\lint f \right) (y)
  \;=\; \lim_{\I^\vee \ni y \rightarrow z \in \Li^\vee}
	\left(\lint h_j \right) (y)}
  \;=\; \frac{\pi}{2} \int_0^1 h_j(\beta z) \; d \beta \nonumber \\
  &=& \frac{\pi}{2} \int_0^1 \partial_j f(\beta z) \; d\beta - \frac{\pi}{2}
	\int_0^1 d\beta \int_0^1 d\alpha \; \alpha \: \partial_j f_{| \alpha
	\beta y}
  - \frac{\pi}{4} \int_0^1 d\beta \int_0^1 d\alpha \;
	\alpha^2 \: y_j \; \Box f_{|\alpha \beta y} \nonumber \\
  &=& \frac{\pi}{2} \int_0^1 \alpha \: \partial_j f(\alpha z) \; d\alpha
	+ \frac{\pi}{4} \int_0^1 (\alpha^2-\alpha) \: y_j \;
	(\Box f)_{|\alpha y} \; d \alpha \nonumber
\end{eqnarray}
\QED

\begin{Satz}
  Auf $\I^\vee$ gilt
\[ \partial_{jk} \; \lint f \;=\; \lint h_{jk}  \]
  mit
\begin{eqnarray}
  h_{jk}[f](y) &=& \partial_{jk} f(y) - \int_0^1 \left(
    2 \alpha^2 \: \partial_{jk} f_{|\alpha y}
    + \frac{1}{2} (2 \alpha^2 - \alpha) \: g_{jk} \:
   (\Box f)_{|\alpha y} \right) d \alpha \nonumber \\
  && \;\;\; - \int_0^1 \frac{1}{2} (2\alpha^3 - \alpha^2)
    \left( y_k \; (\Box \partial_j f)_{|\alpha y}
   + y_j \; (\Box \partial_k f)_{|\alpha y} \right) d \alpha \nonumber \\
\label{eq:a1_60}
   && \;\;\; - \int_0^1 \frac{1}{4} (\alpha^4 - \alpha^3) \; y_j \: y_k
   \; (\Box^2 f)_{|\alpha y}
   \; d \alpha \spc .
\end{eqnarray}
\end{Satz}
{\Beweis}
Nach Satz~\ref{a1_parabl} gilt auf $\I^\vee$
\[ \partial_{jk} \lint f \;=\; \partial_j \lint h_k \;=\; \lint h_{jk} \]
mit
\begin{eqnarray}
  h_k(y) &=& \partial_k f(y) - \int_0^1 \left( \alpha \; \partial_k
	f_{|\alpha y} + \frac{1}{2} \alpha^2 \; y_k \: (\Box f)_{|\alpha y}
	\right) d \alpha \\
\label{eq:a1_61}
  h_{jk}(y) &=& \partial_j h_k(y) - \int_0^1 \left( \beta \; \partial_j
	h_{k|\beta y} + \frac{1}{2} \beta^2 \; y_k \: (\Box h_k)_{|\beta y}
	\right) d \beta \spc .
\end{eqnarray}
Hieraus kann $h_{jk}$ direkt bestimmt werden:
\begin{eqnarray*}
  \partial_j h_k(y) &=& \partial_{jk} f_{|y} - \int_0^1 \left( \alpha^2 \:
     \partial_{jk} f_{|\alpha y} + \frac{1}{2} \alpha^2 \: g_{jk} \: (\Box
     f)_{|\alpha y} + \frac{1}{2} \alpha^3 \: y_k \: (\Box
     \partial_j f)_{|\alpha y}
     \right) d \alpha \\
  \partial_j h_k(\beta y) &=& \partial_{jk} f_{|\beta y}
	- \int_0^1 \left( \alpha^2 \:
     \partial_{jk} f_{|\alpha \beta y} + \frac{1}{2} \alpha^2 \: g_{jk} \:
     (\Box f)_{|\alpha \beta y} \right) d \alpha \\
     && \spc - \frac{1}{2} \int_0^1
     \alpha^3 \beta \: y_k \: (\Box \partial_j f)_{|\alpha \beta y}
     \: d \alpha \\
  \Box h_k(\beta y) &=& \Box \partial_k f_{|\beta y} - \int_0^1 \left(
     2 \alpha^3 \Box \partial_k f_{|\alpha \beta y} + \frac{1}{2} \alpha^4
	\beta \: y_k \: \Box^2 f_{|\alpha \beta y} \right) d \alpha
\end{eqnarray*}
Wir wenden die Integralumformungen von Lemma \ref{integr_umf} an
\begin{eqnarray*}
  \int_0^1 \beta \: (\partial_j h_k)_{|\beta y} \; d \beta
   &=& \int_0^1 \left( \alpha^2 \: \partial_{jk} f_{|\alpha y}
   + \frac{1}{2} (\alpha^2 - \alpha) \: g_{jk} \: (\Box f)_{|\alpha y}
   \right) d \alpha \\
   && \spc + \frac{1}{2} \int_0^1 (\alpha^3 - \alpha^2) \: y_k \: (\Box
	\partial_j f)_{|\alpha y} \: d \alpha \\
  \int_0^1 \beta^2 \: (\Box h_k)_{|\beta y} \; d \beta
   &=& \int_0^1 \left( (2 \alpha^3 - \alpha^2) \; (\Box
   \partial_k f)_{|\alpha y}
	+ \frac{1}{2} (\alpha^4 - \alpha^3) \: y_k \: (\Box^2 f)_{|\alpha y}
	\right) d \alpha \;\;\;.
\end{eqnarray*}
und setzen in \Ref{a1_61} ein.
\QED

Abschlie{\ss}end leiten wir Formeln f\"ur die Distributionsableitungen von
$\slint f$ ab:
\begin{Satz}
\label{a1_dis_abl}
  F\"ur die ersten Ableitungen von $\slint f$ gilt (im Distributionssinne)
\Equ{a1_58}
  \partial_j \: (\lint f)(y) \;=\; (\lint h_j)(y) \:+\:
   \pi \: y_j \: l^\vee(y)  \: \int_0^1 f(\alpha y)
	\: d \alpha 
\EndEqu
   mit
\[   h_j(y) \;=\; (\partial_j f)(y) - \frac{1}{2} \int_0^1 \alpha \;\:
	\Box_z(z_j f(z))_{|z=\alpha y} \; d \alpha \spc . \]
\end{Satz}
Gleichung~\Ref{a1_58} l\"a{\ss}t sich anhand einer formalen Rechnung leicht
einsehen: Dazu nimmt man an, da{\ss} die glatte Funktion $(\slint f)_{|\I^\vee}$
auf ganz $M$ glatt fortgesetzt werden kann, es also eine Funktion
$\gamma \in C^\infty(M)$ gibt mit
\[  \gamma_{| \I^\vee} \;=\; (\lint f)_{|\I^\vee}  \spc .       \]
Nach Satz~\ref{a1_lemma1} und Satz~\ref{a1_parabl} hat man
\begin{eqnarray}
\label{eq:a1_83}
  \gamma(y) &=& \frac{\pi}{2} \int_0^1 f(\alpha y) \; d \alpha \;\;\;, \spc
	y \in \Li^\vee  \\
\label{eq:a1_84a}
   \partial_j \gamma(y) &=& (\lint h_j)(y) \;\;\;, \spc \spc y \in \I^\vee \\
\label{eq:a1_85a}
   (\lint f)(y) &=& \gamma(y) \; \Theta(y^2) \: \Theta(y^0) \spc .
\end{eqnarray}
Formales Ableiten von~\Ref{a1_85a} liefert
\begin{eqnarray*}
  \partial_j (\lint f)(y) &=& \partial_j \gamma(y) \; \Theta(y^2) \;
    \Theta(y^0) \;+\; 2 \: y_j \: \gamma(y) \; \delta(y^2) \; \Theta(y^0) \\
    && +\; \gamma(y) \; \Theta(y^2) \; g_{j0} \: \delta(y^0) \spc .
\end{eqnarray*}
Im ersten Summanden kann man~\Ref{a1_83}, im zweiten~\Ref{a1_84a} einsetzen,
der dritte Summand verschwindet. Auf diese Weise erh\"alt man~\Ref{a1_58}
und kann leicht die Vorfaktoren und Vorzeichen verifizieren.

Die technischen Probleme dieser Herleitung kann man folgenderma{\ss}en umgehen:\\
{\Beweis}
  Da $\slint f$ au{\ss}erhalb von $\I^\vee \cup \Li^\vee$ verschwindet, hat man
  f\"ur $g \in C^\infty_c(M)$
\[  (\lint f)(g) \;=\; \int_{\I^\vee \cup \Li^\vee} (\lint f)(y) \: g(y) \:
	d^4 y   \spc . \]
  Bei der Berechnung der schwachen Ableitung von $\slint f$ kann man den
  Satz von Gau{\ss} anwenden\footnote{Beachte, da{\ss} $\flint f$ aufgrund von
  Konvention~\Ref{a1_56} bis an den Rand des Integrationsgebietes stetig
  ist.}
\begin{eqnarray*}
    (\lint f)(\partial_j g) &=&
       \int_{\I^\vee \cup \Li^\vee} (\lint f)_{|y} \: \partial_j g_{|y} \:
	d^4 y   \\
  &=& \int_{\Li^\vee} (-y_j) \: (\lint f)_{|y} \: g_{|y} \: d \mu \;-\;
	\int_{\I^\vee \cup \Li^\vee} (\partial_j \lint f)_{|y} \: g_{|y} \: d^4y
	\spc ,
\end{eqnarray*}
  wobei $d \mu$ das kanonische Ma{\ss} $d^3 \vec{x}/|\vec{x}|$
  auf $L^\vee$ bezeichnet.

  Einsetzen von~\Ref{a1_50} und~\Ref{a1_56} liefert
\begin{eqnarray*}
  &=& -\int_{\Li^\vee} y_j \left( \frac{\pi}{2} \int_0^1 f(\alpha y) \:d \alpha
	\right) \: g(y) \; d\mu \;-\;
	\int_{\I^\vee \cup \Li^\vee} (\lint h_j)_{|y} \: g_{|y} \: d^4y \\
  &=& - \int_M \left( (\lint h_j)(y) + \pi \: y_j \: l^\vee(y)
	\int_0^1 f(\alpha y) \: d \alpha \right) g(y) \: d^4y \spc .
\end{eqnarray*}
\QED
Die zweiten Ableitungen erh\"alt man wiederum durch Iteration:
\begin{Satz}
\label{a1_dis_abl2}
F\"ur die zweiten Ableitungen von $\slint f$ gilt (im Distributionssinne)
\begin{eqnarray}
  \partial_{jk} \; (\lint f)(y) &=& 2 \pi \left( \int_0^1 f(\alpha y) \:
    d \alpha \right) \: y_j \: y_k \: m^\vee(y)  \nonumber \\
  && + \pi \left( \int_0^1 
   f(\alpha y) \: d \alpha \right) \: g_{jk} \:
   l^\vee(y) \nonumber \\
  && + \pi \left( \int_0^1 \alpha \left( y_j \: \partial_k
   f_{|\alpha y} + y_k \: \partial_j f_{|\alpha y} \right) d \alpha \right)
   \: l^\vee(y) \nonumber \\
  && + \frac{\pi}{2} \left( \int_0^1 (\alpha^2 - \alpha)
   (\Box f)_{|\alpha y} \: d \alpha \right) \: y_j \: y_k \: l^\vee(y)
   \nonumber \\
\label{eq:a1_70}
  && + (\lint h_{jk})(y) \spc ,
\end{eqnarray}
  wobei $h_{jk}$ die Funktion~\Ref{a1_60} und $m^\vee$ die
  Distribution
\Equ{51a}
    m^\vee(y) \;=\; \delta^\prime(y^2) \: \Theta(y^0)
\EndEqu
  bezeichnet.
\end{Satz}
{\Beweis}
Nach Satz~\ref{a1_dis_abl} gilt
\[  \partial_k \: (\lint h_j)(y) \;=\; (\lint h_{jk})(y) + \pi
   \: y_k \: l^\vee(y) \int_0^1 h_j(\alpha y) \: d \alpha \spc .  \]
   Einsetzen der Umformung
\begin{eqnarray*}
  \int_0^1 h_j(\alpha y) \: d \alpha &=& \int_0^1 \partial_j f_{|\alpha y}
   \: d \alpha - \int_0^1 d \alpha \int_0^1 \beta \: d \beta \;
   \partial_j f_{|\alpha \beta y} \\
  && - \frac{1}{2} \: y_j \int_0^1 \alpha \: d\alpha \int_0^1 \beta^2 \:
  d\beta \; (\Box f)_{|\alpha \beta y} \\
  &=& \int_0^1 \left( \alpha \: \partial_j f_{|\alpha y} + \frac{1}{2}
  (\alpha^2 - \alpha) \: y_j \: (\Box f)_{|\alpha y} \right) d \alpha
\end{eqnarray*}
  f\"uhrt auf
\begin{eqnarray}
   \partial_k \: (\lint h_j) (y) &=& (\lint h_{jk})(y) + \frac{\pi}{2}
   \left( \int_0^1 (\alpha^2 - \alpha) \: (\Box f)_{|\alpha y} \: d \alpha
   \right) \: y_j \: y_k \: l^\vee(y) \nonumber \\
\label{eq:a1_71}
   && + \pi \left( \int_0^1 \alpha \: y_k \: \partial_j
   f_{|\alpha y} \: d \alpha \right) \: l^\vee(y) \spc .
\end{eqnarray}
Weiterhin hat man wegen $\partial_k l^\vee(y) = 2 y_k \: m^\vee(y)$
\begin{eqnarray}
   \partial_k \left( y_j \: l^\vee(y) \: \int_0^1 f(\alpha y) \: d \alpha
   \right) &=& g_{jk} \: l^\vee(y) \: \int_0^1 f(\alpha y) \: d \alpha
   + l^\vee(y) \int_0^1 \alpha \: y_j \: \partial_k f_{|\alpha y} \: d \alpha
   \nonumber \\
\label{eq:a1_72}
   && + 2 \: y_j \: y_k \: \left( \int_0^1 f(\alpha y) \: d \alpha \right)
   m^\vee(y) \spc .
\end{eqnarray}
Partielles Ableiten von~\Ref{a1_58} nach $y^k$ sowie
Einsetzen von~\Ref{a1_71} und~\Ref{a1_72} liefert~\Ref{a1_70}.
\QED
\begin{Bem}\
\em Aus Rechnung~\Ref{a1_78} folgt die Beziehung
\Equ{a1_77}
   \Box \: \lint f \;=\; 2 \pi f \: l^\vee \spc .
\EndEqu
   Zur besseren Kontrolle der abgeleiteten Formeln kann
   man~\Ref{a1_77} auch aus~\Ref{a1_70} durch \"Uberschieben mit der Metrik
   ableiten:
   
   Aus~\Ref{a1_60} erh\"alt man unmittelbar $\slint h_{jk} g^{jk}=0$. Beachte
   dazu, da{\ss} Terme, die $y^2$ enthalten, auf $\Li^\vee$
   verschwinden. Au{\ss}erdem gilt $y^j \partial_j (\Box f)(\alpha y) =
   (d/d \alpha) (\Box f)(\alpha y)$, so da{\ss} man partiell integrieren
   kann.
   
   Daher liefert Gleichung~\Ref{a1_70} unter Verwendung von
   $y^2 \: l^\vee(y)=0$
\begin{eqnarray*}
   \Box \: (\lint f)(y) &=& \pi \left( \int_0^1 f(\alpha y) \: d\alpha
   \right) y^j \partial_j \: l^\vee(y) \\
   && + 4 \pi \left( \int_0^1 f(\alpha y) \: d\alpha \right) \: l^\vee(y) \\
   && + 2 \pi \left( \int_0^1 \alpha \: \frac{d}{d \alpha} f(\alpha y) \:
   d \alpha \right) \: l^\vee(y) \spc .
\end{eqnarray*}
   Anwendung der Identit\"at\footnote{Das sieht man folgenderma{\ss}en:
\begin{eqnarray*}
   \int g(y) \; y^j \partial_j l^\vee(y) \; d^4y &=&
   - \int l^\vee(y) \: (4 + y^j \partial_j) g(y) \: d^4 y \\
   &=& - \int \frac{d^3\vec{y}}{2 |\vec{y}|} (4 + y^j \partial_j)g_{|y=
      (|\vec{y}|, \vec{y})} \; d^4y
\end{eqnarray*}
   Gehe nun in Polarkoordinaten

\[ \;=\; - \int_{S^2} d \omega \int_0^\infty \frac{r}{2} \: dr
   \:\left( 4 + r \frac{\partial}{\partial r} \right) g(r,r,\omega)  \]
   und integriere partiell
\begin{eqnarray*}
   &=& - \int_{S^2} d \omega \int_0^\infty \frac{r}{2} \: dr \:
      2 g(r,r,\omega) \;=\; - 2 \int l^\vee(y) \: g(y) \: d^4y \spc .
\end{eqnarray*}}
\Equ{a1_54b}
   y^j \partial_j l^\vee(y) = - 2 l^\vee
\EndEqu
und partielle Integration im dritten Summanden f\"uhrt auf~\Ref{a1_77}.
\em \EndBem \end{Bem}
Mit Hilfe der Formeln~\Ref{a1_58} und~\Ref{a1_70} f\"ur die Ableitungen von
Lichtkegelintegralen kann man die St\"orungsrechnung f\"ur $k_0$ relativ
direkt durchf\"uhren. Es ist g\"unstig, vorher die Notation etwas zu
vereinfachen:

\section{Einige Bezeichungen}
F\"ur die St\"orungsrechnung im Ortsraum werden oft Integrale l\"angs der
Verbindungsstrecke zweier Punkte in $M$ ben\"otigt. Wir wollen eine
Kurzschreibweise vereinbaren:
\begin{Def}
   Schreibe f\"ur das Linienintegral zwischen den Punkten $x, y \in M$
\[    \int_x^y f \; := \; \int_0^1 f(\alpha y + (1-\alpha)x) \: d \alpha 
   \spc . \]
   Als Integrationsvariable soll dabei stets $\alpha$ verwendet
   werden, also
\[    \int_x^y \alpha^n \: f \;:=\; \int_0^1 \alpha^n \: 
	 f(\alpha y + (1-\alpha)x) \: d \alpha  \spc .   \]
\end{Def}

Wenn wir Distributionen
$p_m$, $s_m$, $k_m$, $\tilde{p}_m$, $\tilde{k}_m$,
\dots~ nach $m$ entwickeln, lassen wir zur besseren \"Ubersicht oft
den Index $m$ weg.
Den Beitrag $\sim m^k$ dieser Ausdr\"ucke bezeichnen wir durch einen
zus\"atzlichen Index $^{(k)}$. Insbesondere haben wir
\begin{eqnarray}
\label{eq:a5_a}
\pn(x,y)  &=& p_0(x,y) \;=\; -\frac{i}{2 \pi^3} \frac{\xi \slsh}{\xi^4} \\
\label{eq:a5_b}
\pe(x,y)  &=& - \frac{1}{4 \pi^3} \frac{1}{\xi^2} \\
\label{eq:a5_cc}
p^{(2)}(x,y) &=& -\frac{i}{8 \pi^3} \: \frac{\xi \slsh}{\xi^2} \\
p^{(3)}(x,y) &=& \frac{1}{16 \pi^3} \: (\ln(|\xi^2|) + C_e) \\
p^{(4)}(x,y) &=& \frac{i}{64 \pi^3} \:\xi\slsh\; (\ln(|\xi^2|) + C_e) \\
\kn(x,y) &=& k_0(x,y) \;=\; \frac{1}{2 \pi^2} \; \xi\slsh \left( m^\vee(\xi) -
m^\wedge(\xi) \right) \\
\label{eq:a5_dd}
\ke(x,y) &=& \frac{i}{4 \pi^2} \left( l^\vee(\xi) - l^\wedge(\xi)
	\right) \\
k^{(2)}(x,y) &=& -\frac{1}{8 \pi^2} \; \xi\slsh \:
	(l^\vee(\xi)-l^\wedge(\xi)) \\
k^{(3)}(x,y) &=& -\frac{i}{16 \pi^2} \; \Theta(\xi^2) \:
	\epsilon(\xi^0) \\
k^{(4)}(x,y) &=& \frac{1}{64 \pi^2} \:\xi\slsh \; \Theta(\xi^2) \: \epsilon(\xi^0)
	\spc ,
\end{eqnarray}
dabei ist $C_e$ die Eulersche Konstante.

Bei zusammengesetzten Termen mu{\ss} der gesamte Ausdruck in Potenzen
von $m$ entwickelt werden, also beispielsweise
\[  (s \Aslsh p)^{(1)} \;=\; \sn \: \Aslsh \: \pe + \se \: \Aslsh \: \pn
	\spc . \]

Oftmals gen\"ugt es, das Verhalten einzelner Beitr\"age der St\"orungsrechnung
nur in der N\"ahe des Lichtkegels zu untersuchen.
In diesem Fall f\"uhren wir eine Entwicklung um den Lichtkegel durch:
Wir nennen eine Funktion $f(y)$, $y \in M$, von der Ordnung $y^{2n}$, $n \in
\N^0$, falls
\[  \lim_{\I \ni y \rightarrow z} \left| y^{-2n} \;
  f(y) \right| \;<\; \infty \;\;\;\;\;\;
	 {\mbox{f\"ur alle $z \in \Li$.}} \]
Entsprechend hei{\ss}t eine Funktion $f(x,y)$ von der Ordnung $(y-x)^{2n}$,
falls $f_x(\xi) := f(x,x+\xi)$ f\"ur jedes $x \in M$ von der Ordnung
$\xi^{2n}$ ist.
M\"ochte man die Beitr\"age der Ordnung $(y-x)^{2n}$ berechnen, so mu{\ss}
man die Randwerte der Ableitungen bis zur $n$-ten Ordnung auf dem
Lichtkegel bestimmen.
Das ist unmittelbar einleuchtend, wir wollen es aber trotzdem
exemplarisch f\"ur eine Entwicklung bis zur Ordnung $y^4$ beweisen:
\begin{Lemma}
\label{a2_lemma2}
Sei $f \in C^\infty(\Li^\vee \cup \I^\vee)$ und gelte
\[ \lim_{\I^\vee \ni y \rightarrow z} f(y) \;=\; \lim_{\I^\vee \ni y
\rightarrow z} \partial_j f(y) \;=\; 0 \;\;\;\;\; {\mbox{f\"ur alle $z \in
	\Li^\vee$.}}            \]
Dann ist $f(y)$ von der Ordnung $y^4$.
\end{Lemma}
{\Beweis}
Nach Voraussetzung sind $f$ und dessen partielle Ableitungen
auf $\I^\vee \cup \Li^\vee$ stetig fortsetzbar. Bezeichne f\"ur den
Beweis der Einfachheit halber den Grenzwert auf dem Lichtkegel durch den
entsprechenden Funktionswert, also\footnote{Diese Bezeichnung wird
ansonsten bewu{\ss}t vermieden, weil sie zu Verwechslungen zwischen der
Distributionsableitung $(\partial_j f)(z)$, $z \in \Li^\vee$ und dem 
Grenzwert der partiellen Ableitung $(\partial_j f)(y)$ bei
$\I^\vee \ni y \rightarrow z \in \Li^\vee$ f\"uhren kann.}
\[ f(z) := \lim_{\I^\vee \ni y \rightarrow z} f(y), \spc
   \partial_j f(z) := \lim_{\I^\vee \ni y \rightarrow z} \partial_j f(y), \;\;
	{\mbox{u.s.w..}}  \]
Betrachte eine Folge $\I^\vee \ni (y_n) \rightarrow z \in \Li^\vee$.
W\"ahle einen zeitartigen Vektor mit $\bra v,v \ket = 1$, $\: v^0 > 0$.
Die Folgenglieder $y_n$ k\"onnen eindeutig in der Form
\[  y_n = z_n + \tau_n v \;\;,\;\;\;\;\; z_n \in \Li^\vee, \;\tau_n > 0 \]
dargestellt werden.
Im Grenzfall $n \rightarrow \infty$ hat man
$z_n \rightarrow z$, $\tau_n \rightarrow 0$.

Es m\"ussen zwei F\"alle getrennt betrachtet werden:
\begin{enumerate}
\item $z \neq 0$: \\
  Man kann eine Umgebung $U \subset \Li^\vee$ von $z$ und
  Konstanten $c_1, c_2 > 0$ w\"ahlen, so da{\ss}
\begin{eqnarray*}
  \left| v^i v^j \: \partial_{ij} f(x) \right| &<& c_1  \\
  \bra x,v \ket &>& c_2
\end{eqnarray*}
f\"ur alle $x \in U$. F\"ur gen\"ugend gro{\ss}es $n$ ist $z_n \in U$ und
somit
\begin{eqnarray*}
|f(y_n)|
  &=& \left| f(z_n) + \tau_n \: v^j \partial_j f(z_n) + \frac{1}{2} \:
	\tau_n^2 \: v^i v^j
	\partial_{ij} f(z_n) \right| + {\cal{O}}(\tau_n^3) \\
  &\leq& \frac{1}{2} c_1 \tau_n^2 + {\cal{O}}(\tau_n^3)  \\
y_n^2 &=& 2 \tau_n \: \bra z_n, v \ket + \tau_n^2 \: \bra v,v \ket \\
y_n^4 &\geq& 4 c_2^2 \: \tau_n^2 + {\cal{O}}(\tau_n^3) \spc .
\end{eqnarray*}
Es folgt
\[  \lim_{n \rightarrow \infty} \left| y_n^{-4} \: f(y_n) \right| \;<\; \infty
	\spc . \]
\item $z=0$:\\
Wegen $\partial_k f_{|\Li^\vee}=0$ folgt f\"ur alle $v \in \Li^\vee$
\begin{eqnarray}
\label{eq:a2_100}
v^j \partial_{jk} f(0) &=& 0 \\
\label{eq:a2_101}
v^i v^j \partial_{ijk} f(0) &=& 0 \spc .
\end{eqnarray}
Aus~\Ref{a2_100} erh\"alt man $\partial_{jk} f(0)=0$. Aus~\Ref{a2_101} folgt
zun\"achst
\[  \partial_{ijk} f(0) \;=\; g_{ij} \: w_k \]
f\"ur einen geeignetes $w \in M$ und wegen der Symmetrie in den Indizes $i$,
$j$ und $k$ sogar
\[  \partial_{ijk} f(0) \;=\; 0 \spc . \]
Die partiellen Ableitungen von $f$ verschwinden also am Ursprung bis zur
dritten Ordnung.
Man kann eine Umgebung $U \subset \Li^\vee$ von $z$ und eine Konstante
$c_3>0$ w\"ahlen, so da{\ss} f\"ur alle $x \in U$ die Ungleichungen
\begin{eqnarray*}
\left| \partial_{jk} f(x) \right| &=& \left| \frac{1}{2} \: x^l x^m
   \partial_{jklm} f(0) \right| + {\cal{O}}(x^3)  \\
  &\leq& c_3 \: \bra x,v \ket^2 + {\cal{O}}( \bra x,v \ket^3 ) \\
\left| \partial_{ijk} f(x) \right| &=& \left| x^l \partial_{ijkl} f(0)
   \right| + {\cal{O}}(x^2) \\
  &\leq& c_3 \: | \bra x,v \ket | + {\cal{O}}( \bra x,v \ket^2 )
\end{eqnarray*}
gelten.
Setze $\kappa_n = \bra z_n,v \ket \geq 0$. Man hat $y_n^2 = 2 \tau_n
\kappa_n + \tau_n^2$. F\"ur gen\"ugend gro{\ss}es $n$ ist $z \in U$ und
somit
\begin{eqnarray*}
|f(y_n)| &=& \left| \frac{1}{2} \: \tau_n^2 \: v^i v^j \partial_{ij} f(z_n)
+ \frac{1}{3!} \: \tau_n^3 \: v^i v^j v^k \partial_{ijk} f(z_n) \right| + 
{\cal{O}}(\tau_n^4) \\
& \leq & c_4 \left( \tau_n^2 \kappa_n^2 + \tau_n^3 \kappa_n \right)
	+ {\cal{O}}(\tau_n^4 + \tau_n^2 \kappa_n^3 + \tau_n^3 \kappa_n^2) \\
y_n^4 &=& 4 \left( \tau_n^2 \kappa_n^2 + \tau_n^3 \kappa_n \right) + \tau_n^4
\spc .
\end{eqnarray*}
Es folgt
\[ \lim_{n \rightarrow \infty} \left| y_n^{-4} \:f_n(y_n) \right| \;<\; \infty
	\spc . \]
\end{enumerate}
\QED

\section{St\"orungsrechnung f\"ur das elektromagnetische Feld}
\label{elek_k0}
Wir betrachten die St\"orung des freien Diracoperators durch ein \"au{\ss}eres
elektromagnetisches Feld
\Equ{a1_54a}
    G \;=\; i \Pdd + e \Aslsh \spc .
\EndEqu
In erster Ordnung St\"orungstheorie hat man nach~\Ref{a1_100}
und~\Ref{a1_101}
\[    \tilde{k}_0 \;=\; k_0 + \Delta k_0  \]
mit
\begin{eqnarray}
\label{eq:a1_205}
  \Delta k_0 (x,y) &=& - e \: (k_0 \Aslsh s_0 + s_0 \Aslsh k_0)(x,y)  \\
   &=& - e \: (i \Pdd)_x \left( K_0 \Aslsh S_0 + S_0 \Aslsh K_0 \right)(x,y)
   \: (i \Pdd)_y  \spc . \nonumber
\end{eqnarray}
Die Gleichungen~\Ref{a1_75} und~\Ref{a1_76} kann man unter Verwendung
der Distributionen $l^\vee$ und $l^\wedge$ schreiben als
\begin{eqnarray*}
   K_0(x,y) &=& - \frac{i}{4 \pi^2} \left( l^\vee(x-y) - l^\wedge(x-y)
	\right) \\
   S_0(x,y) &=& - \frac{1}{4 \pi} \left( l^\vee(x-y) + l^\wedge(x-y) \right)
   \spc .
\end{eqnarray*}
Folglich gilt
\begin{eqnarray}
\lefteqn{(K_0 \Aslsh S_0 + S_0 \Aslsh K_0)(x,y)} \nonumber \\
   &=&
   \frac{i}{16 \pi^3} \int d^4z \: \left( l^\vee(x-z) - l^\wedge(x-z) \right)
      \Aslsh(z)\: 
	  \left( l^\vee(z-y) + l^\wedge(z-y) \right) \nonumber \\
   && + \frac{i}{16 \pi^3} \int d^4z \: \left( l^\vee(x-z) + l^\wedge(x-z)
      \right) \Aslsh(z)\:
	  \left( l^\vee(z-y) - l^\wedge(z-y) \right) \nonumber \\
   &=& \frac{i}{8 \pi^3} \int d^4z \: \left( l^\vee(x-z) \: \Aslsh (z) \:
      l^\vee(z-y)
      \;-\; l^\wedge(x-z) \:\Aslsh(z)\: l^\wedge(z-y) \right) \nonumber \\
   &=& \frac{i}{8 \pi^3} \int d^4z \: \left(l^\wedge_x(z) \: l^\vee_y(z)
      - l^\wedge_y(z) \: l^\vee_x(z) \right) \: \Aslsh (z) \nonumber \\
\label{eq:a1_80}
   &=& - \frac{i}{8 \pi^3} \left( \lint_x^y \Aslsh \;-\; \lint_y^x \Aslsh
      \right)  \spc .
\end{eqnarray}
Den zu berechnenden Operator $\Delta k_0$ kann man also mit
Lichtkegelintegralen aus\-dr\"ucken\footnote{Beachte wegen der Vorzeichen,
da{\ss} f\"ur den Operator $B(x,y) = \flint_x^y \Aslsh - \flint_y^x \Aslsh$ und
eine Wellenfunktion $\Psi$ gilt:
\begin{eqnarray*}
   \left( (i \Pdd) \: B \: (i \Pdd) \; \Psi \right)(x) &=&
   \int d^4y \; (i \Pdd)_x  B(x,y) (i \Pdd)_y \;\Psi(y) \\
   &=& - \int d^4y \; \left[ 
	(i \Pdd)_x \left(i \frac{\partial}{\partial y^j} \right)
      B(x,y) \: \gamma^j \right] \; \Psi(y)
\end{eqnarray*}   }:
\begin{eqnarray}
\label{eq:a1_56a}
   \Delta k_0(x,y) &=& \frac{i e}{8 \pi^3} (i \Pdd)_x \left(
      \lint_x^y \Aslsh - \lint_y^x \Aslsh \right) (i \Pdd)_y \\
\label{eq:a1_81}
   &=& \frac{i e}{8 \pi^3} \: \gamma^i \gamma^j \gamma^k \:
      \frac{\partial}{\partial x^i} \frac{\partial}{\partial y^k}
      \left( \lint_x^y A_j - \lint_y^x A_j \right)
\end{eqnarray}
Die Ableitungen k\"onnen nun mit Hilfe von Satz~\ref{a1_dis_abl} und
Satz~\ref{a1_dis_abl2} ausgewertet werden.

\begin{Thm}
\label{a1_theorem1}
   In erster Ordnung St\"orungstheorie gilt
\begin{eqnarray}
\label{eq:a1_111}
 \Delta k_0(x,y) &=&  - i e \left( \int_x^y A_j \right) \: \xi^j \; k_0(x,y) \\
\label{eq:a1_112}
   && + \frac{i e}{8 \pi^2} \left( \int_x^y (\alpha^2-\alpha) \; \xi \slsh
	\: \xi^k \: j_k \right) \;\; \left( l^\vee(\xi) - l^\wedge(\xi)
	\right) \\
\label{eq:a1_113}
   && - \frac{i e}{8 \pi^2} \left( \int_x^y (2 \alpha - 1) \; \xi^j \:
	\gamma^k \: F_{kj} \right) \;\; \left( l^\vee(\xi) - \l^\wedge(\xi)
	\right) \\
\label{eq:a1_114}
   && + \frac{e}{16 \pi^2} \left( \int_x^y \varepsilon^{ijkl} \; F_{ij} \:
	\xi_k \; \rho \gamma_l \right) \;\; \left( l^\vee(\xi) -l^\wedge(\xi)
	\right) \\
\label{eq:a1_115}
   && - \frac{i e}{16 \pi^3} \left( \lint_x^y - \lint_y^x \right) dz
	\int_x^z (\alpha^4 - \alpha^3) \; \zeta \slsh \;
	\zeta_k \; \Box j^k \\
\label{eq:a1_116}
   && + \frac{i e}{16 \pi^3} \left( \lint_x^y - \lint_y^x \right) dz
	\int_x^z (4 \alpha^3 - 3 \alpha^2) \; \zeta^j \: \gamma^k
	\; \Box F_{kj} \\
\label{eq:a1_117}
   && - \frac{e}{32 \pi^3} \left( \lint_x^y - \lint_y^x \right) dz
	\int_x^z \alpha^2 \; \varepsilon^{ijkl} \; (\Box F_{ij}) \:
	\zeta_k \; \rho \gamma_l \\
\label{eq:a1_118}
   && + \frac{i e}{4 \pi^3} \left( \lint_x^y - \lint_y^x \right) dz
	\int_x^z (2 \alpha^2 - \alpha) \; \gamma^k \: j_k \spc ,
\end{eqnarray}
   wobei $F_{jk}=\partial_j A_k - \partial_k A_j$ den elektromagnetischen
   Feldst\"arketensor und $j^k = \partial_l F^{kl}$ den Maxwell-Strom
   bezeichnet. Zur Abk\"urzung wurde $\xi = y-x$ und $\zeta = z-x$ gesetzt.
\end{Thm}
{\Beweis}
   Betrachte zun\"achst den Fall $\xi^0>0$:

   Da $\slint_y^x \Aslsh$ au{\ss}erhalb von $\xi \in \I^\wedge \cup \Li^\wedge$
   verschwindet, vereinfacht sich~\Ref{a1_81} unter Verwendung
   von Lemma~\ref{a1_lemma8} zu
\begin{eqnarray}
   \Delta k_0 (x,y) &=& \frac{i e}{8 \pi^3} \: \gamma^i \gamma^j \gamma^k
   \: \frac{\partial}{\partial x^i} \frac{\partial}{\partial y^k}
   \lint_x^y A_j  \nonumber \\
\label{eq:a1_95}
   &=& - \frac{i e}{8 \pi^3} \: \gamma^i \gamma^j \gamma^k \: \left(
      \frac{\partial^2}{\partial y^i \partial y^k} \lint_x^y A_j
      - \frac{\partial}{\partial y^k} \lint_x^y \partial_i A_j \right)
   \spc .
\end{eqnarray}
   
   Berechne nun nacheinander die einzelnen Terme:
\begin{eqnarray*}
   \gamma^i \gamma^j \gamma^k \frac{\partial^2}{\partial y^i \partial y^k}
   \lint_x^y A_j &=& \left( \gamma^i \gamma^j + \gamma^j \gamma^i \right)
      \gamma^k \frac{\partial^2}{\partial y^i \partial y^k}
      \lint_x^y A_j
      \;-\; \gamma^j \gamma^i \gamma^k \frac{\partial^2}{\partial y^i
      \partial y^k} \lint_x^y A_j   \\
    &=& 2 \gamma^k \frac{\partial^2}{\partial y^j \partial y^k}
      \lint_x^y A^j \;-\; \Box_y \lint_x^y \Aslsh   \spc ,
\end{eqnarray*}
   wobei die Symmetrie der zweiten partiellen Ableitungen ausgenutzt
   wurde.
   
   Wende nun Satz~\ref{a1_dis_abl2} und~\Ref{a1_77} an:
\begin{eqnarray}
   &=& 2 \: \lint h_{jk}[A^j] \: \gamma^k
      \;+\; 4 \pi \: (\int_x^y A^j) \: \xi_j \: \xi \slsh \: m^\vee(\xi) \\
\label{eq:a1_84}
   && - 2 \pi l^\vee(\xi) \left( \Aslsh(y) - \int_x^y \Aslsh \right)
   \;+\; 2 \pi l^\vee(\xi) \left( \int_x^y \alpha \: \left( \xi_j \:
      \Pdd A^j + \xi \slsh \: \partial_j A^j \right) \right)   \\
\label{eq:a1_85}
   && + \pi l^\vee(\xi) \left( \int_x^y (\alpha^2 - \alpha) \: \xi \slsh
      \xi_j \: \Box A^j \right)
\end{eqnarray}
   F\"uhre in~\Ref{a1_85} die Ersetzung
\[    \xi_j \: \Box A^j \;=\; \xi^j \partial_j \: \partial_k A^k
      - \xi^k \: j_k \]
   durch. Die Ableitung in Richtung des Vektors $\xi$ kann wieder
   in eine Ableitung nach der Variablen $\alpha$ umgewandelt und partiell
   integriert werden.
   
   Auf analoge Weise kann man im zweiten Summanden
   von~\Ref{a1_84} den Ausdruck
\[    \xi_j \: \Pdd A^j \;=\; \xi^j \partial_j \Aslsh + \xi^j \gamma^k \:
	 F_{kj}   \]
   behandeln. Die nun bei der partiellen Integration auftretenden Terme
   kompensieren gerade den ersten Summanden in~\Ref{a1_84}.
   
   Man erh\"alt auf diese Weise den Ausdruck
\begin{eqnarray*}
   \gamma^i \gamma^j \gamma^k \frac{\partial^2}{\partial y^i \partial y^k}
   \lint_x^y A_j
   &=& 2 \: \lint h_{jk}[A^j] \: \gamma^k
      \;+\; 4 \pi \: (\int_x^y A_j) \: \xi^j \: \xi \slsh \: m^\vee(\xi) \\
   && + 2 \pi l^\vee(\xi) \int_x^y \alpha \: \xi^j \gamma^k \: F_{kj}
      \;-\; \pi l^\vee(\xi) \left( \int_x^y (\alpha^2 - \alpha) \: \xi \slsh
      \xi^k \: j_k \right)  \\
   && + \pi l^\vee(\xi) \int_x^y \xi \slsh \: \partial_j A^j \spc .
\end{eqnarray*}
   Weiterhin gilt
\begin{eqnarray*}
   \gamma^i \gamma^j \gamma^k \: \frac{\partial}{\partial y^k}
   \lint_x^y \partial_i A_j
   &=& \gamma^i \gamma^j \gamma^k \left( \lint_x^y h_k[\partial_i A_j]
   \right) \;+\; \pi \: \gamma^i \gamma^j \: \xi \slsh \: l^\vee(\xi)
   \left( \int_x^y \partial_i A_j \right)   \\
   &=& \gamma^i \gamma^j \gamma^k \left( \lint_x^y h_k[\partial_i A_j] \right)
      \;+\; \frac{\pi}{2} l^\vee(\xi)
      \left( \int_x^y F_{ij} \gamma^i \gamma^j \xi \slsh \right)  \\
   && + \pi l^\vee(\xi) \; \int_x^y \xi \slsh \: \partial_j A^j
	  \spc .
\end{eqnarray*}
   Man kann die Identit\"at
\Equ{a1_70a}
    F_{ij} \: \gamma^i \gamma^j \xi \slsh \;=\; 2 \gamma^i F_{ij} \: \xi^j
   - i \: \varepsilon^{ijkl} \: F_{ij} \: \xi_k \; \rho \gamma_l
\EndEqu
   anwenden und erh\"alt nach Einsetzen in~\Ref{a1_95} insgesamt:
\begin{eqnarray}
\label{eq:a1_87}
   \Delta k_0(x,y) &=& -\frac{i e}{2 \pi^2} \: (\int_x^y A_j) \: \xi^j
      \xi \slsh \: m^\vee(\xi) \;-\; \frac{i e}{8 \pi^2} \: l^\vee(\xi)
      \left( \int_x^y (2 \alpha-1) \xi^j \gamma^k F_{kj} \right) \\
   && + \frac{i e}{8 \pi^2} \: l^\vee(\xi) \left( \int_x^y (\alpha^2-\alpha)
      \: \xi \slsh \xi^k \: j_k \right) \\
\label{eq:a1_88}
   && + \frac{e}{16 \pi^2} \:
      l^\vee(\xi) \left( \int_x^y \varepsilon^{ijkl} F_{ij} \xi_k
      \; \rho \gamma_l \right) \\
\label{eq:a1_89}
   && - \frac{i e}{4 \pi^3} \lint_x^y  \left( \gamma^j \: h_{jk}[A^k]
      \;-\; \frac{1}{2} \gamma^i \gamma^j \gamma^k \:
      h_k[\partial_i A_j] \right)
\end{eqnarray}
   Beachtet man, da{\ss} f\"ur $\xi^0>0$
\[    k_0(x,y) \;=\; (i \Pdd)_x \: K_0(x,y) \;=\; \frac{1}{2 \pi^2} \:
   \xi \slsh \: m^\vee(\xi)   \]
   gilt, so sieht man, da{\ss} die Summanden in den Zeilen~\Ref{a1_87}
   bis~\Ref{a1_88} mit den Aus\-dr\"ucken~\Ref{a1_111} bis~\Ref{a1_114}
   \"ubereinstimmen.
   
   Es soll nun gezeigt werden, da{\ss} auch die verbleibenden Terme von
   $\tilde{k}_0$ korrekt sind, also insbesondere, da{\ss}~\Ref{a1_89} mit
   den Summanden~\Ref{a1_114} bis~\Ref{a1_118} \"ubereinstimmt:

   Aus~\Ref{a1_60} erh\"alt man
\begin{eqnarray}
\label{eq:a1_102}
   \gamma^j \: h_{jk}[A^k]_{|z} &=& \Pdd \partial_k A^k
	- \int_x^z \left( 2 \alpha^2
   \: \Pdd \partial_k A^k + \frac{1}{2} (2 \alpha^2 - \alpha) \:
   \Box \Aslsh \right) \\
\label{eq:a1_103}
   && - \frac{1}{2} \int_x^z (2 \alpha^3 - \alpha^2) \left(
   \zeta_k \: \Box \Pdd A^k + \zeta \slsh \: \Box \partial_k A^k \right) \\
\label{eq:a1_104}
   && - \frac{1}{4} \int_x^z (\alpha^4 - \alpha^3) \zeta \slsh \:
	\zeta_k \; \Box^2 A^k  \spc .
\end{eqnarray}
   F\"uhre in~\Ref{a1_103} bzw.~\Ref{a1_104} die Ersetzungen
\begin{eqnarray*}
   \zeta_k \: \Box \Pdd A^k &=& \zeta^k \partial_k \Box \Aslsh +
      \zeta^k \gamma^j \: \Box F_{jk} \\
   \zeta_k \: \Box^2 A^k &=& \zeta^k \partial_k \Box \partial_j A^j -
      \zeta^k \: \Box j_k
\end{eqnarray*}
   durch, schreibe die Ableitung in Richtung des Vektors $\zeta$
   als Ableitung nach $\alpha$ und integriere partiell.
   Man erh\"alt dann
\begin{eqnarray*}
   \gamma^j \: h_{jk}[A^j]_{|z} &=& \frac{1}{4} \int_x^z (\alpha^4 - \alpha^3)
      \zeta \slsh \: \zeta^k \: \Box j_k   \\
    && - \frac{1}{2} \int_x^z (2 \alpha^3 - \alpha^2) \zeta^k \: \gamma^j
      \: \Box F_{jk} \\
    && - \int_x^z 2 \alpha^2 \: \Pdd \partial_k A^k + \int_x^z
      (2 \alpha^2 - \frac{\alpha}{2}) \Box \Aslsh \\
    && + \Pdd \partial_k A^k - \frac{1}{2} \Box \Aslsh 
       - \frac{1}{4} \int_x^z \alpha^2 \zeta \slsh \: \Box \partial_k
       A^k \spc . \\
\end{eqnarray*}
   Nach~\Ref{a1_55} gilt weiterhin:
\begin{eqnarray*}
   \gamma^i \gamma^j \gamma^k \: h_k[\partial_i A_j]_{|z} &=&
   \gamma^i \gamma^j \gamma^k \left( \partial_{ik} A_j(z) - \int_x^z
      \alpha \: \partial_{ik} A_j + \frac{1}{2} \alpha^2 \zeta_k \: \Box
      \partial_i A_j \right) \\
    &=& 2 \: \Pdd \partial_k A^k - \Box \Aslsh - \int_x^z
      \left( 2 \alpha \: \Pdd \partial_k A^k -
      \alpha \: \Box \Aslsh \right) \\
    && - \frac{1}{2}
      \int_x^z \left( \alpha^2 \: (\Box \partial_i A_j) \: \gamma^i \gamma^j
      \: \zeta \slsh \right) \\
    &=& - \int_x^z \left( 2 \alpha \: \Pdd \partial_k A^k 
      - \alpha \: \Box \Aslsh \right) \\
    && - \frac{1}{2} \int_x^z \left( \alpha^2 \: \zeta^j \gamma^k \:
      \Box F_{kj} \right) + \frac{i}{4} \int_x^z \left( \alpha^2
      \: \varepsilon^{ijkl} \: (\Box F_{ij}) \: \zeta_k \; \rho \gamma_l
      \right)  \\
    && + 2 \: \Pdd \partial_k A^k - \Box \Aslsh - \frac{1}{2} \int_x^z \alpha^2
	\zeta \slsh \: \Box \partial_k A^k
\end{eqnarray*}
   Insgesamt ergibt sich
\begin{eqnarray}
\lefteqn{\gamma^j \: h_{jk}[A^k]_{|z} - \frac{1}{2} \gamma^i \gamma^j \gamma^k
      \: h_k[\partial_i A_j]_{|z}} \nonumber \\
   &=& - \int_x^z (2 \alpha^2 - \alpha) \;
      \gamma^k \: j_k 
      \;-\; \int_x^z (\alpha^3 - \frac{3}{4} \alpha^2) \: \zeta^j \: \gamma^k
      \; \Box F_{kj} \nonumber \\
\label{eq:a3_80a}
   && + \frac{1}{4} \int_x^z (\alpha^4 - \alpha^3) \: \zeta\slsh \: \zeta_k
      \; \Box j^k
      \;-\; \frac{i}{8} \int_x^z \alpha^2 \: \varepsilon^{ijkl} \:
      (\Box F_{ij}) \: \zeta_k \; \rho \gamma_l \spc .
\end{eqnarray}
   Durch Einsetzen in~\Ref{a1_89} erh\"alt man gerade die Terme \Ref{a1_115}
   bis \Ref{a1_118}. Damit ist der Fall $\xi^0>0$ bewiesen.
\\[1em]
   Den Fall $\xi^0<0$ kann man durch Punktspiegelung des Minkowski-Raumes
   am Ursprung auf den Fall $\xi^0>0$ zur\"uckf\"uhren:
   Beachte dazu zun\"achst, da{\ss}\footnote{Dabei bezeichnet $\hat{f}$ die
   Funktion $\hat{f}(x) = f(-x)$.}
\Equ{a1_128}
     \lint_x^y f \;=\; \lint_{-y}^{-x} \hat{f}
\EndEqu
   und, nach Einsetzen in~\Ref{a1_81},
\Equ{a1_129}
     \Delta k_0(x,y) \;=\; - \Delta k_0[\hat{A}] (-x,-y)
\EndEqu
   gilt. Dabei bedeutet die eckige Klammer $[\hat{A}]$, da{\ss} das
   elektromagnetische Potential $A$ in der Definitionsgleichung von
   $\Delta k_0$ durch $\hat{A}$ zu ersetzen ist.
   Auf der rechten Seite von~\Ref{a1_129} kann man die Gleichung
   \Ref{a1_111} bis~\Ref{a1_118} anwenden und erh\"alt
   die Aussage des Theorems f\"ur $\xi^0<0$.
   Verwende dabei die Identit\"aten~\Ref{a1_128} und
\begin{eqnarray*}
	F_{ij}[\hat{A}] &=& - \hat{F}_{ij}  \\
	j_k[\hat{A}] &=& \hat{\jmath}_k \\
	\int_{-x}^{-y} \alpha^n \hat{f} &=& \int_x^y \alpha^n f \spc .
\end{eqnarray*}
  Um den Fall $\xi^0=0$ zu behandeln, gen\"ugt die Feststellung, da{\ss}
  im bisherigen Beweis die Distributionen
\begin{eqnarray*}
  -\frac{i e}{8 \pi^3} \gamma^i \gamma^j \gamma^k
    \frac{\partial}{\partial x^i} \frac{\partial}{\partial y^k}
    \: \lint_x^y A_j  \\
  \frac{i e}{8 \pi^3} \gamma^i \gamma^j \gamma^k
    \frac{\partial}{\partial x^i} \frac{\partial}{\partial y^k}
    \: \lint_y^x A_j    
\end{eqnarray*}
  vollst\"andig (also ohne die Einschr\"ankungen $\xi^0>0$ bzw. $\xi^0<0$)
  berechnet wurden und da{\ss} $\Delta k_0$ nach~\Ref{a1_81} die Summe dieser
  Operatoren ist. Wir haben also alle Beitr\"age zu $\tilde{k}_0$ gefunden,
  es gibt keine zus\"atzlichen Terme, die f\"ur $\xi^0=0$ beitragen.
\QED
Wir wollen
die gewonnene Gleichung f\"ur $\tilde{k}_0$ kurz diskutieren:

Zun\"achst einmal f\"allt auf, da{\ss} $\tilde{k}_0(x,y)$ f\"ur raumartiges
$\xi$ verschwindet, was die Kausalit\"at widerspiegelt.

Die Summanden~\Ref{a1_111} bis~\Ref{a1_114} tragen nur auf dem Lichtkegel
(also f\"ur $\xi \in \Li$) bei und werden daher die {\em auf dem Lichtkegel
lokalisierten Terme} genannt. Die Summanden~\Ref{a1_115} bis~\Ref{a1_118}
f\"uhren dagegen zu einem Beitrag im Innern des Lichtkegels (also f\"ur
$\xi \in \I$) und hei{\ss}en {\em delokalisierte Terme}.

Wir wollen nun den Spezialfall betrachten, da{\ss} $\tilde{k}_0$
durch eine Eichtransformation aus $k_0$ hervorgeht, also da{\ss}
\[     e \Aslsh \;=\; \Pdd \Lambda     \]
f\"ur eine geeignete reelle Funktion $\Lambda$.
Da diese Eichtransformation der Phasentransformation
\[      \Psi(x) \; \rightarrow \; e^{i \Lambda(x)} \: \Psi(x) \]
der Wellenfunktionen entspricht, kann man $\tilde{k}_0$ sehr einfach
exakt angeben:
\begin{eqnarray}
  \tilde{k}_0(x,y) &=& e^{i \Lambda(x)} \; k_0(x,y) \; e^{-i \Lambda(y)}
	    \nonumber \\
   &=& e^{i \Lambda(x)} \; (i \Pdd)_x K_0(x,y) \; e^{-i \Lambda(y)}
	    \nonumber \\
\label{eq:a1_120}
   &=& e^{i \Lambda(x)- i \Lambda(y)} \: k_0(x,y)
\end{eqnarray}
Auf der anderen Seite kann man $\tilde{k}_0(x,y)$ mit Hilfe von
Theorem~\ref{a1_theorem1} berechnen. Da $F_{ij}=j_k=0$ ist, verschwinden
die Summanden~\Ref{a1_112} bis~\Ref{a1_118}. Man erh\"alt
\begin{eqnarray*}
   \tilde{k}_0(x,y) &=& k_0(x,y) - i \left( \int_x^y \partial_j \Lambda 
	\right) \: \xi^j \; k_0(x,y)  \\
   &=& k_0(x,y) - i \left( \int_0^1 d \alpha \; \frac{d}{d \alpha}
	\Lambda(\alpha y + (1-\alpha) x) \right) \; k_0(x,y) \\
   &=& k_0(x,y) + i \left( \Lambda(x)-\Lambda(y) \right) \; k_0(x,y) \spc ,
\end{eqnarray*}
was in erster Ordnung in $\Lambda$ mit~\Ref{a1_120} \"ubereinstimmt.
Man sieht auf diese Weise, da{\ss} der Summand~\Ref{a1_111} f\"ur das 
richtige Eichtransformationsverhalten von $\tilde{k}_0$ verantwortlich ist.
Er wird daher {\em Eichterm} genannt.

Aus Satz~\ref{a1_lemma1} folgt, da{\ss} $\tilde{k}_0(x,y)$
f\"ur $y-x \in \I$ eine glatte Funktion ist. Dar\"uber hinaus kann
man den Grenzfall untersuchen, da{\ss} sich $y-x$ dem Lichtkegel ann\"ahert:
\begin{Satz}
\label{a1_randwert}
   F\"ur $y-x \in \Li$ gilt
\begin{eqnarray}
\label{eq:a1_121}
  \lim_{\I_x \ni u \rightarrow y} \tilde{k}_0(x,u) &=&
  + \frac{i e}{64 \pi^2} \: \epsilon(\xi^0) \int_x^y (\alpha^4 - 2 \alpha^3
	+ \alpha^2) \: \xi \slsh \: \xi_k \; \Box j^k \\
  && - \frac{i e}{64 \pi^2} \: \epsilon(\xi^0) \int_x^y (4 \alpha^3
	- 6 \alpha^2 + 2 \alpha) \: \xi^j \: \gamma^k \: (\Box F_{kj}) \\
  && + \frac{e}{64 \pi^2} \: \epsilon(\xi^0) \int_x^y (\alpha^2 - \alpha) \;
	\varepsilon^{ijkl} \: (\Box F_{ij}) \: \xi_k \; \rho \gamma_l \\
\label{eq:a1_122}
  && - \frac{i e}{8 \pi^2} \: \epsilon(\xi^0) \int_x^y (\alpha^2 - \alpha)
	\: \gamma^k \: j_k \spc ,
\end{eqnarray}
wobei wieder $\xi = y-x$ gesetzt wurde.
\end{Satz}
{\Beweis}
  Die Behauptung folgt mit Hilfe von~\Ref{a1_8} durch eine direkte Rechnung.
  Dabei wendet man die folgenden Integralumformungen von Lemma \ref{integr_umf}
  an:
\begin{eqnarray*}
  \int_x^y dz \int_x^z (\alpha^4 - \alpha^3) \; \zeta \slsh \: \zeta_k \;
	\Box  j^k
  &=& \int_0^1 d \beta \int_0^1 d \alpha \; (\alpha^4 - \alpha^3) \;
	\beta^2 \: \xi \slsh \: \xi_k \;
	(\Box j^k)(\alpha \beta y + (1\!-\!\alpha
	\beta) x)  \\
  &=& - \frac{1}{2} \int_x^y (\alpha^4 - 2 \alpha^3 + \alpha^2) \:
	\xi \slsh \: \xi_k \; \Box j^k  \\
  \int_x^y dz \int_x^z (4 \alpha^3 - 3 \alpha^2) \; \zeta^j \gamma^k \:
	(\Box F_{kj}) &=& - \frac{1}{2} \int_x^y (4 \alpha^3 - 6 \alpha^2
	+ 2 \alpha) \; \xi^j \gamma^k \: (\Box F_{kj}) \\
  \int_x^y dz \int_x^z \alpha^2 \: \varepsilon^{ijkl} \: (\Box F_{ij}) \:
	\zeta_k \; \rho \gamma_l &=&
	\int_x^y (\alpha - \alpha^2) \: \varepsilon^{ijkl} \:
	(\Box F_{ij}) \: \xi_k \; \rho \gamma_l \\
  \int_x^y dz \int_x^z (2 \alpha^2 - \alpha) \: \gamma^k \: j_k &=&
	\int_x^y (\alpha - \alpha^2) \: \gamma^k \: j_k
\end{eqnarray*}
(Die Integrationsvariable $\alpha$ bezieht sich immer auf das innere
Linienintegral.)
\QED
Nun kann das Verhalten der einzelnen Summanden~\Ref{a1_111} bis
\Ref{a1_118} in der N\"ahe des Lichtkegels untersucht werden.
Bei den delokalisierten Termen betrachten wir zur Einfachheit nur die
N\"aherungsausdr\"ucke \Ref{a1_121} bis \Ref{a1_122}.
Die Summanden~\Ref{a1_112} und~\Ref{a1_115} sind (als $4 \times 4$-Matrizen
bei festem $x, y$) proportional zu $\xi \slsh$ und hei{\ss}en
{\em radiale Beitr\"age}.
Aufgrund der Antisymmetrie des elektromagnetischen Feldst\"arketensors
sind die Terme~\Ref{a1_113}, \Ref{a1_114}, \Ref{a1_116} und~\Ref{a1_117}
senkrecht zu $\xi \slsh$ und werden daher {\em orthogonale Beitr\"age}
genannt.
Die Summanden \Ref{a1_112}, \Ref{a1_118} enthalten den Maxwellstrom.
Sie spielen f\"ur die Feldgleichungen eine entscheidende Rolle und werden
{\em Stromterme} genannt.

\section{St\"orungsrechnung f\"ur das Gravitationsfeld}
\label{grav_k0}

Der Diracoperator im Gravitationsfeld wird in~\cite{Physdip}
ausf\"uhrlich behandelt. Es wird dort insbesondere gezeigt, da{\ss} die
Diracmatrizen $\gamma^j$ bei Einf\"uhrung des Gra\-vi\-ta\-tions\-fel\-des
durch dynamische Matrixfelder $G^j(x)$ zu ersetzen sind, die mit der
Lorentzmetrik \"uber die Antikommutatorrelation
\Equ{a1_199}
      g^{jk} \;=\; \frac{1}{2} \left\{ G^j, G^k \right\}
\EndEqu
verkn\"upft sind. Der Diracoperator nimmt in allgemeinen Eichungen die
Form
\Equ{a1_200}
   G \;=\; i G^j \partial_j + G^j E_j
\EndEqu
an, wobei die Matrizen $E_j$, die sogenannten eichfixierenden Terme,
von $G^j$ abh\"angige Matrixfelder sind:
\Equ{eich_fix}
    E_j = \frac{i}{2} \rho \rho,_j - \frac{i}{16} \Tr
    \left( G^m G^n;_j \right) G_m G_n + \frac{i}{8} \Tr \left( \rho
    G_j G^m;_m \right) \rho
\EndEqu
Die St\"orungsrechnung f\"ur das elektromagnetische Feld wurde im vorigen
Abschnitt ohne Wahl einer speziellen Eichung (z.B. Lorentzeichnung)
durchgef\"uhrt. Dadurch bot das Eichverhalten der verschiedenen Terme
bei den Rechnungen eine wichtige Kontrollm\"oglichkeit.
Um den Rechenaufwand so gering wie m\"oglich zu halten, soll in diesem
Abschnitt hingegen eine g\"unstige Eichung gew\"ahlt werden:
Zun\"achst einmal kann man die Matrixfelder $G^j$
durch eine geeignete Eichtransformation in die
Form\footnote{Das sieht man folgenderma{\ss}en: Nach Definition des
Diracoperators im Gravitationsfeld (siehe~\cite{Physdip}) gibt es
zu jedem Raum-Zeit-Punkt $x$ eine spezielle Eichung und ein spezielles
Koordinatensystem, so da{\ss} $G^j(x)=\gamma^j$. In dieser Eichung hat
man in beliebigen Koordinaten
\Equ{a1_201}
   G^j(x) = \sum_{a=0}^3 k^{ja} \gamma_a
\EndEqu
mit geeigneten Koeffizienten $k^{ja}$.
Da $G^j(x)$ nur von der Eichung im Punkte $x$ abh\"angt,
kann man durch eine geeignete Eichtransformation
erreichen, da{\ss}~\Ref{a1_201} sogar auf
dem ganzen Definitionsbereich der Karte mit glatten Funktionen
$k^{ja}(x)$ gilt.}
\[ G^j(x) \;=\; \gamma^j + \sum_{a=0}^3 h^{ja}(x) \: \gamma_a  \]
bringen. Die Funktionen $h^{ja}$ beschreiben dabei
die St\"orung durch das Gra\-vi\-ta\-tions\-feld\footnote{Beachte,
da{\ss} $h^{ja}$ kein Tensorfeld ist, sondern sich unter
Koordinatentransformationen gem\"a{\ss}
\[   \tilde{h}^{ja} \;=\; \frac{\partial \tilde{x}^j}{\partial x^k}
	 \: h^{ka}      \]
transformiert.}.
Da die St\"orungsrechnung f\"ur $k_0$ in
erster Ordnung durchgef\"uhrt werden soll,
gen\"ugt es im folgenden, die linearen Terme in $h^{ja}$ mitzunehmen.

Bei einer Eichtransfomation
\[  \tilde{\Psi}(x) = (1 + i \: a_{kl}(x) \: \sigma^{kl} + {\cal O }(a^2)) \:
      \Psi(x)     \]
verhalten sich die Matrixfelder $G^j$ wie
\begin{eqnarray*}
   \tilde{G}^j &=& \gamma^j + \sum_{a=0}^3 h^{ja} \gamma_a +
      i \: a_{kl} \left[ \sigma^{kl}, \gamma^j \right] + {\cal O}(a^2,h a) \\
   &=& \gamma^j + \sum_{a=0}^3 (h^{ja} + 4 a^{ja}) \: \gamma_a
      + {\cal O}(a^2, h a) \spc .
\end{eqnarray*}
Man kann also durch Wahl der Eichung erreichen, da{\ss}
$h^{ja}$ symmetrisch ist, also $h^{ja}=h^{aj}$. Eine solche Eichung
wird {\bf symmetrische Eichung} genannt. Wir werden im folgenden
stets in symmetrischer Eichung arbeiten und auch das Ergebnis der
St\"orungsrechnung nur in symmetrischer Eichung angeben.

Die Funktionen $h^{ja}$ sollen von nun an formal wie ein Tensorfeld
behandelt werden. Wie in der linearisierten Gravitationstheorie
\"ublich\footnote{Siehe z.B. Landau/Lifshitz \cite{landau},
Band 2, \S 105.}, soll das Heben und Senken von Tensorindizes mit der
Minkowski-Metrik $\eta_{jk}$ durchgef\"uhrt werden.

F\"ur die Metrik erh\"alt man nach~\Ref{a1_199}
\begin{eqnarray*}
   g^{jk} &=& \eta^{jk} + 2 h^{jk}  \\
   g_{jk} &=& \eta_{jk} - 2 h_{jk}  \spc .
\end{eqnarray*}
Insbesondere sieht man, da{\ss} die $h^{jk}$ durch die Lorentzmetrik eindeutig
bestimmt sind.
Man hat nun noch die Freiheit, kleine Koordinatentransformationen
$\tilde{x}^j = x^j + \kappa^j(x)$ durchzuf\"uhren. Diese Willk\"ur kann man
ausnutzen, um die Zusatzbedingung
\Equ{a1_210}
      h_{jk},^k \;=\; \frac{1}{2} h,_j \spc\spc {\mbox{mit}}
	 \spc h=h^k_k
\EndEqu
zu erf\"ullen\footnote{Siehe Landau/Lifshitz \cite{landau},
Band 2, \S 107.}.
Man hat dann
\begin{eqnarray*}
   \Gamma^i_{jk} &=& -h^i_j,_k - h^i_k,_j + h_{jk},^i    \\
   \sqrt{|\det g|} &=& 1 - h \\
   R_{jk} &=& - \Box h_{jk}   \\
   E_j &=& - \frac{i}{16} \Tr ( G^m G^n;_j ) \: G_m G_n  \\
   &=& \frac{i}{4} \left( h_{jn},_m - h_{jm},_n \right) \:
      \gamma^m \gamma^n    \\
   G^j E_j &=& \frac{i}{4} \left( h_{jn},_m \: \gamma^j \gamma^m \gamma^n
      - h_{jm},_n \: \gamma^j \gamma^m \gamma^n \right) \\
   &=& \frac{i}{2} \: \gamma^j \: \left( h_{jk},^k - h,_j \right) \\
   &=& - \frac{i}{4} \: \gamma^j \: h_{,j} \spc .
\end{eqnarray*}
Durch Einsetzen in~\Ref{a1_200} erh\"alt man\footnote{Durch partielle
Integration kann man direkt \"uberpr\"ufen, da{\ss} $G$ hermitesch ist, also
(in erster Ordnung in $h_{jk}$) die Gleichung
\[  \int (G \Psi) \overline{\Phi} \sqrt{|g|} d^4 x \;=\;
	\int \Psi \: \overline{G \Phi} \: \sqrt{|g|} d^4 x      \]
erf\"ullt.}
\begin{eqnarray}
\label{eq:a1_208}
   G &=& i \Pdd + i \gamma_j h^{jk} \frac{\partial}{\partial x^k}
      - \frac{i}{4} (\Pdd h) \spc .
\end{eqnarray}
Bei der Durchf\"uhrung der St\"orungsrechnung tritt eine kleine Schwierigkeit
auf: Das Gravitationsfeld f\"uhrt nicht nur zu einer St\"orung des
Diracoperators in der Form~\Ref{a1_208}, sondern auch zu einer \"Anderung
des Integrationsma{\ss}es $\sqrt{-g} \: d^4x$. Gleichung~\Ref{a1_100} kann
jedoch nur angewendet werden, wenn das Integrationsma{\ss} unver\"andert
bleibt.

Der einfachste Ausweg besteht darin, auf ein System mit dem Ma{\ss}
$d^4x$ zu transformieren, die St\"orungsrechnung in diesem System
durchzuf\"uhren und anschlie{\ss}end auf das Ma{\ss} $\sqrt{-g} \: d^4x$
zur\"uckzutransformieren.
Beachte dazu, da{\ss} sich eine Wellenfunktion $\Psi$, der Diracoperator $G$
und ein Integraloperator $A(x,y)$ beim \"Ubergang von dem Ma{\ss} $a(x) \: d^4x$
zu $b(x) \: d^4x$ mit $0 < a,b \in C^\infty(M)$ gem\"a{\ss}
\begin{eqnarray}
   \Psi^{(a)}(x) &=& \sqrt{\frac{b(x)}{a(x)}} \: \Psi^{(b)}(x) \nonumber \\
   G^{(a)} &=& \sqrt{\frac{b}{a}} \:G^{(b)}\: \sqrt{\frac{a}{b}} \nonumber \\
\label{eq:a1_209}
   A^{(a)}(x,y) &=& \sqrt{\frac{b(x)}{a(x)}} \sqrt{\frac{b(y)}{a(y)}} \:
	A^{(b)}(x,y)
\end{eqnarray}
transformieren, also insbesondere
\begin{eqnarray*}
   G^{(d^4x)} &=& \left( 1 - \frac{h}{2} \right) \: G \: \left(1+\frac{h}{2}
	\right) \\
   &=& i \Pdd + i \gamma_j h^{jk} \frac{\partial}{\partial x^k} +
	\frac{i}{4} \: (\Pdd h) \spc .
\end{eqnarray*}
Aus~\Ref{a1_100} folgt nun
\begin{eqnarray*}
  \Delta k^{(d^4x)}_0 &=& - i \left( k_0 \:\gamma^j \: (h_j^k \partial_k
	+ \frac{1}{4} h,_j) \: s_0
	+ s_0 \: \gamma^j \: (h_j^k \partial_k + \frac{1}{4} h,_j) \: k_0
	\right)
\end{eqnarray*}
und nach~\Ref{a1_209}
\begin{eqnarray*}
   \tilde{k}_0(x,y) &=& \left( 1+\frac{h(x)}{2} \right)
	\left( 1+\frac{h(y)}{2} \right) \:
	\tilde{k}^{(d^4x)}_0(x,y) \\
   &=& k_0(x,y) + i \frac{\partial}{\partial y^k}
	\left( k_0 \: \gamma^j h^k_j \: s_0
	+ s_0 \: \gamma^j h^k_j \: k_0 \right)(x,y)  \\&& - \frac{1}{4}
	\left( k_0 \:
	(i \Pdd h) \: s_0 + s_0 \: (i \Pdd h) \: k_0 \right)(x,y)
	+ \frac{1}{2} (h(x)+h(y)) \: k_0(x,y) \spc .
\end{eqnarray*}
Einsetzen der Beziehung\footnote{Dies folgt aus
\begin{eqnarray*}
   k_0 \: (i \Pdd h) \: s_0 &=& k_0 \: [i \Pdd, h] \: s_0 \;=\; (i \Pdd k_0)
	\: h \: s_0 - k_0 \: h \: (i \Pdd s_0) \\
   &=& - k_0 \: h \spc ,
\end{eqnarray*}
wobei $h$ als Multiplikationsoperator aufgefa{\ss}t wird.}
\[ \left( k_0 \: (i \Pdd h) \: s_0 + s_0 \: (i \Pdd h) \: k_0 \right)(x,y)
	\;=\; \left( h(x)-h(y) \right) \: k_0(x,y)      \]
liefert
\begin{eqnarray*}
   \tilde{k}_0(x,y) &=& k_0(x,y) + \left( \frac{1}{4} h(x) + \frac{3}{4}
	h(y) \right) \: k_0(x,y) \\
   && + i \frac{\partial}{\partial y^k} \left( k_0 \: \gamma^j h^k_j \: s_0
	+ s_0 \: \gamma^j h^k_j \: k_0 \right)(x,y) \spc .
\end{eqnarray*}
Ein Vergleich mit~\Ref{a1_205} zeigt, da{\ss} man die Rechnungen zum Teil auf
die Ergebnisse des vorigen Abschnittes zur\"uckf\"uhren kann:
\begin{eqnarray}
\label{eq:a1_223}
   \;=\; k_0(x,y) + \left( \frac{1}{4} h(x) + \frac{3}{4}h(y) \right)
	\: k_0(x,y)
     - \frac{i}{e} \frac{\partial}{\partial y^k} \:
	\Delta k_0[\gamma^j h_j^k](x,y)
\end{eqnarray}
Dabei bezeichnet $\Delta k_0[\gamma^j h_j^k]$ den Beitrag der
St\"orung durch das elektromagnetische Potential $A_j = h_j^k$.

Mit Hilfe von Theorem~\ref{a1_theorem1} kann man den letzten Summanden
in~\Ref{a1_223} direkt auswerten und erh\"alt auf diese Weise:
\begin{Thm}
\label{a1_theorem2}
In erster Ordnung St\"orungstheorie gilt in symmetrischer Eichung
\begin{eqnarray}
\label{eq:a1_212}
\lefteqn{\Delta k_0(x,y) \;=\; - \left( \int_x^y h^k_j \right) \xi^j
	\frac{\partial}{\partial y^k} \;\; k_0(x,y) } \\
\label{eq:a1_213}
  && - \frac{1}{4 \pi^2} \left( \int_x^y (2 \alpha -1) \; \gamma^i \: \xi^j \:\xi^k
	\; (h_{jk},_i - h_{ik},_j)  \right) \;\; (m^\vee(\xi) - m^\wedge(\xi)
	) \\
\label{eq:a1_214}
  && + \frac{i}{8 \pi^2} \left( \int_x^y \varepsilon^{ijlm} \; (h_{jk},_i
	- h_{ik},_j) \: \xi^k \; \xi_l \: \rho \gamma_m \right) \;\;
	(m^\vee(\xi) - m^\wedge(\xi)) \\
\label{eq:a1_215}
  && + \frac{1}{2} \left( \int_x^y (\alpha^2 - \alpha) \; \xi^j \; \xi^k \;
	R_{jk} \right) \;\; k_0(x,y) \\
\label{eq:a1_216}
  && + \frac{1}{32 \pi^2} \left( \int_x^y (\alpha^4 - 2 \alpha^3 + \alpha^2)
	\; \xi \slsh \; \xi^j \; \xi^k \; \Box R_{jk} \right) \;\;
	(l^\vee(\xi)-l^\wedge(\xi)) \\
\label{eq:a1_217}
  && - \frac{1}{32 \pi^2} \left( \int_x^y (6 \alpha^2 - 6 \alpha + 1) \;
	\xi \slsh \; R \right) \;\; (l^\vee(\xi)-l^\wedge(\xi)) \\
\label{eq:a1_218}
  && + \frac{1}{32 \pi^2} \left( \int_x^y (4 \alpha^3 - 6 \alpha^2 + 2 \alpha)
	\; \xi^j \: \xi^k \: \gamma^l \; R_{j[k},_{l]} \right) \;\;
	(l^\vee(\xi)-l^\wedge(\xi)) \\
\label{eq:a1_219}
  && - \frac{i}{16 \pi^2} \left( \int_x^y (\alpha^2 - \alpha) \;
	\varepsilon^{ijlm} \: R_{ki},_j \: \xi^k \: \xi_l \: \rho \gamma_m
	\right) \;\; (l^\vee(\xi)-l^\wedge(\xi)) \\
\label{eq:a1_220}
  && - \frac{1}{8 \pi^2} \left( \int_x^y (\alpha^2 - \alpha) \; \xi^j \:
	\gamma^k \: G_{jk} \right) \;\; (l^\vee(\xi)-l^\wedge(\xi)) \\
\label{eq:a1_221}
  && + {\cal{O}}(\xi^0) \spc , \nonumber
\end{eqnarray}
wobei $R_{jk}$ den Ricci- und $G_{jk}=R_{jk} - \frac{1}{2} R \: g_{jk}$
den Einstein-Tensor bezeichnet.\\
($m^\vee$, $m^\wedge$ sind die Distributionen $m^\vee(y)=
\delta^\prime(y^2) \:\Theta(y^0)$, $m^\wedge(y)=\delta^\prime(y^2) \:
\Theta(-y^0)$,
ferner wurde $\xi=y-x$ gesetzt.)
\end{Thm}
{\Beweis}
Setzt man in Theorem~\ref{a1_theorem1} das elektromagnetische Potential
$A_j=h_j^k$ ein, so ergibt sich f\"ur den elektromagnetischen
Feldst\"arketensor und den Maxwellstrom
\begin{eqnarray*}
  F_{ij} &=& h^k_j,_i - h^k_i,_j \\
  j_i &=& h^{kj},_{ij} - \Box h^k_i \;=\;
  \frac{1}{2} h,^k_i - \Box h^k_i \spc ,
\end{eqnarray*}
wobei \Ref{a1_210} angewendet wurde.
Bei der Berechnung von 
\[  -\frac{i}{e} \frac{\partial}{\partial y^k} \: \Delta
	k_0[\gamma^j h_j^k](x,y)        \]
f\"uhren die Summanden~\Ref{a1_111} bis~\Ref{a1_118} zu folgenden Beitr\"agen:
\begin{eqnarray*}
\Ref{a1_111} &:& \Ref{a1_212} - \frac{1}{2} h(y) \: k_0(x,y) -
	\frac{1}{2} \left(\int_x^y h \right) k_0(x,y) \\
\Ref{a1_112} &:& \Ref{a1_215} - \frac{1}{4} (h(y)+h(x)) \: k_0(x,y)
	+ \frac{1}{2} \left(\int_x^y h \right) k_0(x,y) \\
  && + \frac{1}{8 \pi^2} \left( \int_x^y (\alpha^2-\alpha) \;
	\xi^j \: \gamma^k \:
	(R_{jk} + \frac{1}{2} \:R\: g_{jk} ) \right) \;\; (l^\vee(\xi)-
	l^\wedge(\xi)) \\
  && - \frac{1}{16 \pi^2} \left( \int_x^y (2 \alpha-1) \; (\Pdd h) \right)\;\;
	(l^\vee(\xi)-l^\wedge(\xi)) \\
\Ref{a1_113} &:& \Ref{a1_213} - \frac{1}{16 \pi^2} \left(\int_x^y (2 \alpha-1)
	\; (\Pdd h) \right) \;\; (l^\vee(\xi)-l^\wedge(\xi)) \\
\Ref{a1_114} &:& \Ref{a1_214}
\end{eqnarray*}
Bei der Ableitung der Summanden~\Ref{a1_115} bis~\Ref{a1_118} kann man
Satz~\ref{a1_dis_abl} anwenden. Da wir die delokalisierten Terme
(also diejenigen Terme, die f\"ur $y \in \I_x$ beitragen und in einer Umgebung
von $\Li_x$ beschr\"ankt sind) nicht berechnen wollen, braucht man nur den
zweiten Summanden in~\Ref{a1_58} zu betrachten. 
Man kann sich die Rechenarbeit erleichtern, wenn man beachtet, da{\ss} die
dabei auftretenden Linienintegrale bereits im Beweis von
Satz~\ref{a1_randwert} bestimmt wurden. Man erh\"alt als Beitr\"age:
\begin{eqnarray*}
\Ref{a1_115} &:& \Ref{a1_216} + \Ref{a1_217} \\
\Ref{a1_116} &:& \Ref{a1_218} \\
\Ref{a1_117} &:& \Ref{a1_219} \\
\Ref{a1_118} &:& - \frac{1}{4 \pi^2} \left(\int_x^y (\alpha^2-\alpha) \;
	\xi^j \: \gamma^k \; R_{jk} \right) \;\; (l^\vee(\xi)-l^\wedge(\xi))\\
  && +\frac{1}{8 \pi^2} \left( \int_x^y (2 \alpha-1) \; (\Pdd h)
	\right) \;\; (l^\vee(\xi)-l^\wedge(\xi))
\end{eqnarray*}
Aufsummieren aller Terme und Einsetzen in~\Ref{a1_223} liefert die
Behauptung.
\QED
Die abgeleitete Formel f\"ur $\tilde{k}_0$ soll noch kurz diskutiert
werden:

Wir wollen dazu zun\"achst den Spezialfall untersuchen, da{\ss} die St\"orung
durch eine Koordinatentransformation $\tilde{x}^i=x^i + \kappa^i(x)$
hervorgerufen wird, wobei $\kappa$ als kleines Vektorfeld angenommen wird.
Aus dem Transformationsverhalten von Integraloperatoren erh\"alt man
unmittelbar in erster Ordnung in $\kappa$:
\begin{eqnarray}
  \tilde{k}_0(x,y) &=& k_0(x\!-\!\kappa(x), y\!-\!\kappa(y)) \nonumber \\
\label{eq:a1_224}
  &=& k_0(x,y) - \kappa^i(x) \frac{\partial}{\partial x^i} k_0(x,y)
	- \kappa^i(y) \frac{\partial}{\partial y^i} k_0(x,y)
\end{eqnarray}
Man kann \"uberpr\"ufen, da{\ss} Theorem~\ref{a1_theorem2} auf das gleiche
Resultat f\"uhrt: F\"ur die St\"orung der Metrik erh\"alt man
\[  h_{ij} \;=\; \frac{1}{2} \left( \kappa_i,_j + \kappa_j,_i \right) \spc . \]
Die Kr\"ummungsgr\"o{\ss}en $R_{jk}$, $G_{jk}$, $R$
verschwinden, so da{\ss} nur die Terme~\Ref{a1_212}
bis~\Ref{a1_214} ausgewertet m\"ussen, also
\begin{eqnarray}
\lefteqn{ \tilde{k}_0(x,y) \;=\; k_0(x,y) - \frac{1}{2} \int_x^y (\kappa^k_{\;,j} +
	\kappa_j,^{\:k}) \: \xi^j \: \frac{\partial}{\partial y^k}
	k_0(x,y) } \nonumber \\
&&- \frac{1}{8 \pi^2} \int_x^y (2\alpha-1) \: \gamma^i \: \xi^j \: \xi^k
	\: (\kappa_{j,ki} - \kappa_{i,kj}) \: (m^\vee(\xi)-m^\wedge(\xi))
	\nonumber \\
&&+ \frac{i}{16 \pi^2} \int_x^y \varepsilon^{ijlm} \: (\kappa_{j,ki} -
	\kappa_{i,kj}) \: \xi^k \: \xi_l \; \rho \gamma_m \;
	(m^\vee(\xi)-m^\wedge(\xi)) \nonumber \\
&=& k_0(x,y) - \kappa^i(x) \frac{\partial}{\partial x^i} k_0(x,y) -
	\kappa^i(y) \frac{\partial}{\partial y^i} k_0(x,y) \nonumber \\
&&- \frac{1}{2} \: \left( \int_x^y \kappa^{[j,k]} \right) \: \xi_j \:
	\frac{\partial}{\partial y^k} k_0(x,y) \nonumber \\
&&- \frac{1}{8 \pi^2} \; \xi^i \: \gamma^k \: (\kappa_{[j,k]}(y) +
	\kappa_{[j,k]}(x)) \:
	(m^\vee(\xi) - m^\wedge(\xi)) \nonumber \\
&&+ \frac{1}{4 \pi^2} \: \left( \int_x^y \kappa_{[j,k]} \: \xi^j \: \gamma^k
	\right) \: (m^\vee(\xi) - m^\wedge(\xi)) \nonumber \\
&&- \frac{i}{16 \pi^2} \; \varepsilon^{jklm} \: (\kappa_{[j,k]}(y) -
	\kappa_{[j,k]}(x)) \: \xi_l \: \rho \gamma_m
	\; (m^\vee(\xi) - m^\wedge(\xi)) \nonumber \\
&=& k_0(x,y) - \kappa^i(x) \frac{\partial}{\partial x^i}
	k_0(x,y) - \kappa^i(y) \frac{\partial}{\partial y^i} k_0(x,y)
	\nonumber \\
&& - \frac{1}{8 \pi^2} \left( \xi^j \: \gamma^k \: (\kappa_{[j},_{k]}(x)
	+ \kappa_{[j},_{k]}(y)) \right) \;\;
	(m^\vee(\xi)-m^\wedge(\xi)) \nonumber \\
\label{eq:a1_113a}
&& + \frac{i}{16 \pi^2} \left(\varepsilon^{jklm} \: (\kappa_{[j},_{k]}(x)
	- \kappa_{[j},_{k]}(y)) \: \xi_l \: \rho \gamma_m \right) \;\;
	(m^\vee(\xi) - m^\wedge(\xi)) \;\;\; ,
\end{eqnarray}
denn es gilt
\begin{eqnarray*}
\left(\int_x^y \kappa^{[j,k]}\right) \xi_j \: \frac{\partial}{\partial y^k}
	k_0(x,y) &=& \frac{1}{2 \pi^2} \: \left(\int_x^y \kappa^{[j,k]}\right)
	\xi_j \: \frac{\partial}{\partial y^k} \; \xi\slsh \:
	(m^\vee(\xi)-m^\wedge(\xi)) \\
&=& \frac{1}{2 \pi^2} \: \left(\int_x^y \kappa_{[j,k]} \right) \:
	\xi^j \: \gamma^k \: (m^\vee(\xi)-m^\wedge(\xi)) \spc .
\end{eqnarray*}
F\"uhrt man nun die Eichtransformation $\tilde{\Psi}(x)= U(x) \: \Psi(x)$ mit
\[  U(x) \;=\; 1 + \frac{1}{8} \: \kappa_{[j},_{k]}(x) \: \gamma^j \: \gamma^k
	+ {\cal{O}}(\kappa^2)   \]
durch, so erh\"alt man aus~\Ref{a1_113a} gerade~\Ref{a1_224}.

Man sieht also, da{\ss} die Summanden~\Ref{a1_212} bis~\Ref{a1_214} f\"ur das
richtige Verhalten bei Koordinatentransformationen verantwortlich sind.
Sie haben \"Ahnlichkeit mit dem Eichterm~\Ref{a1_111} bei einer
elektromagnetischen St\"orung, der auf das korrekte Verhalten bei
$U(1)$-Eichtransformationen f\"uhrt.

Die Summanden~\Ref{a1_215} bis~\Ref{a1_217} sind {\em{tangentiale Terme}},
\Ref{a1_218} und~\Ref{a1_219} {\em{orthogonale Terme}}.
Die entscheidende Rolle spielen die Summanden \Ref{a1_215} und \Ref{a1_220}.
Sie werden {\em{Kr\"ummungsterme}} genannt. Falls $y-x$ klein wird, sind sie
proportional zum Ricci- bzw. Einstein-Tensor.

\section{Skalare St\"orung}
Wir betrachten die St\"orung durch ein skalares Potential
\Equ{a1_s1}
	G \;=\; i \Pdd + \Xi \spc .
\EndEqu
F\"ur $\Delta k_0$ hat man in erster Ordnung in $\Xi$
\[      \Delta k_0 \;=\; - \left( s_0 \: \Xi \: k_0 + k_0 \: \Xi \: s_0
	\right) \spc . \]
\begin{Thm}
\label{theorem_sk0}
In erster Ordnung St\"orungstheorie gilt
\begin{eqnarray}
\label{a1_s2}
\Delta k_0(x,y) &=& -\frac{1}{2} \: (\Xi(y)+\Xi(x)) \; k^{(1)}(x,y) \\
\label{a1_s3}
&&+ \frac{1}{8 \pi^2} \: (l^\vee(\xi)-l^\wedge(\xi)) \;
	\int_x^y (\partial_j \Xi) \; \xi_k \; \sigma^{jk} \\
\label{a1_s4}
&&- \frac{1}{16 \pi^3} \left(\lint_x^y - \lint_y^x \right) dz \int_x^z
	\alpha^2 \: (\partial_j \Box \Xi) \; \zeta_k \; \sigma^{jk} \\
\label{a1_s5}
&&+ \frac{i}{16 \pi^3} \left(\lint_x^y - \lint_y^x \right) \Box \Xi \spc .
\end{eqnarray}
\end{Thm}
{\Beweis}
Es gen\"ugt wieder, den Fall $\xi^0 > 0$ zu betrachten.
\begin{eqnarray}
\Delta k_0(x,y) &=& \frac{i}{8 \pi^3} \: \Pdd_x \: \Pdd_y \lint_x^y \Xi
\;=\; \frac{i}{8 \pi^3} \left( \gamma^i \: \Pdd_y \lint_x^y \partial_i \Xi
	\;-\; \Box_y \lint_x^y \Xi \right) \nonumber \\
&=& - \frac{i}{4 \pi^2} \: \Xi(y) \: l^\vee(\xi) \;+\; \frac{i}{8 \pi^2} \:
	l^\vee(\xi) \int_x^y (\Pdd \Xi) \; \xi\slsh \nonumber \\
\label{eq:a1_127a}
&& + \frac{i}{8 \pi^3} \lint_x^y dz \left( \Box \Xi - \int_x^z \alpha \;
	\Box \Xi - \frac{1}{2} \int_x^z \alpha^2
	\: (\Pdd \Box \Xi) \: \zeta \slsh \right) \\
&=& -\frac{i}{8 \pi^2} \: l^\vee(\xi) \; (\Xi(y) + \Xi(x))
	+ \frac{1}{8 \pi^2} \: l^\vee(\xi) \int_x^y (\partial_j \Xi) \;
	\xi_k \; \sigma^{jk} \nonumber \\
&&- \frac{1}{16 \pi^3} \lint_x^y dz \int_x^z \alpha^2 \: (\partial_j \Box \Xi)
	\: \zeta_k \; \sigma^{jk}
	+ \frac{i}{16 \pi^3} \lint_x^y \Box \Xi \nonumber
\end{eqnarray}
\QED
F\"ur die Randwerte auf dem Lichtkegel hat man
\begin{Satz}
F\"ur $y-x \in \Li$ gilt
\begin{eqnarray*}
\lim_{\I_x \ni u \rightarrow y} \Delta k_0(x,u) &=&
\frac{1}{32 \pi^2} \: \epsilon(\xi^0) \int_x^y (\alpha^2-\alpha) \;
	(\partial_j \Box \Xi) \; \xi_k \; \sigma^{jk} \\
&&+ \frac{i}{32 \pi^2} \:
	\epsilon(\xi^0) \int_x^y \Box \Xi \;+\; {\cal{O}}(\xi^2) \spc .
\end{eqnarray*} 
\end{Satz}
{\Beweis}
Folgt direkt mit Hilfe von \Ref{a1_8} und den Integralumformungen
\ref{integr_umf}.
\QED

\section{Bilineare St\"orung}
Wir betrachten die St\"orung durch ein bilineares Potential, also
\Equ{a1_b00}
	G \;=\; i \Pdd + B_{jk} \: \sigma^{jk}
\EndEqu
mit einem antisymmetrischen Tensorfeld $B_{jk} = - B_{kj}$.
F\"ur $\Delta k_0$ hat man in erster Ordnung in $B$
\Equ{a1_b0}
  \Delta k_0 \;=\; - \left( s_0 \: B_{jk} \: \sigma^{jk} \; k_0 +
	k_0 \: B_{jk} \: \sigma^{jk} \: s_0 \right) \spc .
\EndEqu
\begin{Thm}
\label{theorem_b0}
In erster Ordnung St\"orungstheorie gilt
\begin{eqnarray}
\label{eq:a1_b1}
\lefteqn{\Delta k_0(x,y) \;=\; -\frac{i}{\pi^2} \: (m^\vee(\xi)-m^\wedge(\xi))
	\int_x^y B_{ij} \; \xi^i \: \xi_k \; \sigma^{jk} } \\
\label{eq:a1_b3}
&&- \frac{i}{8 \pi^2} \: (l^\vee(\xi)-l^\wedge(\xi)) \; (B_{jk}(y) +
	B_{jk}(x)) \; \sigma^{jk} \\
\label{eq:a1_b4}
&&+ \frac{i}{2 \pi^2} \: (l^\vee(\xi)-l^\wedge(\xi)) \int_x^y B_{jk} \;
	\sigma^{jk} \\
\label{eq:a1_b2}
&&+ \frac{1}{4 \pi^2} \: (l^\vee(\xi)-l^\wedge(\xi)) \int_x^y \xi_j \;
	B^{jk}_{\;\;\:,k} \\
\label{eq:a1_b5}
&&- \frac{i}{4 \pi^2} \: (l^\vee(\xi)-l^\wedge(\xi)) \int_x^y (2\alpha-1) \;
	(\xi^k \: B_{jk,i} + \xi_i \: B_{jk,}^{\;\;\;\;k}) \;
	\sigma^{ij} \\
\label{eq:a1_b6}
&&- \frac{i}{4 \pi^2} \: (l^\vee(\xi)-l^\wedge(\xi)) \int_x^y
	(\alpha^2-\alpha) \; (\Box B_{ij}) \: \xi^i \: \xi_k \; \sigma^{jk} \\
\label{eq:a1_b7}
&&+ \frac{i}{8 \pi^2} \: (l^\vee(\xi)-l^\wedge(\xi)) \int_x^y
	\varepsilon^{ijkl} \; B_{ij,k} \; \xi_l \; \rho \\
&&+ {\cal{O}}(\xi^0) \spc . \nonumber
\end{eqnarray}
\end{Thm}
{\Beweis}
Es gen\"ugt wieder, den Fall $\xi^0>0$ zu betrachten.
\begin{eqnarray*}
\Delta k_0(x,y) &=& \frac{i}{8 \pi^3} \: \Pdd_x \: \sigma^{jk} \: \Pdd_y
	\lint_x^y B_{jk}
\;=\; -\frac{1}{8 \pi^3} \: \frac{\partial}{\partial x^a}
	\frac{\partial}{\partial y^b} \lint_x^y \gamma^a \; \gamma^j
	\gamma^k \: \gamma^b \; B_{jk}
\end{eqnarray*}
Wir setzen die Relation
\Equ{a1_b8}
\gamma^a \: \gamma^j \gamma^k \: \gamma^b \; B_{jk} \;=\; 2 B^{ab} - i
	B_{jk} \: \sigma^{jk} \; g^{ab} - 2i \: (B^a_{\;\:j} \: \sigma^{jb}
	+ B^b_{\;\:j} \: \sigma^{ja} ) - i \rho \: \varepsilon^{ijab}
	\: B_{ij}
\EndEqu
ein und erhalten
\begin{eqnarray*}
\Delta k_0(x,y) &=& - \frac{1}{4 \pi^3} \: \frac{\partial}{\partial
x^j} \frac{\partial}{\partial y^k} \lint_x^y B^{jk} \\
&&+ \frac{i}{8 \pi^3} \: \frac{\partial}{\partial x^a}
	\frac{\partial}{\partial y_a} \lint_x^y B_{jk} \; \sigma^{jk} \\
&&+ \frac{i}{4 \pi^3} \: \frac{\partial}{\partial x^a}
	\frac{\partial}{\partial y^b} \lint_x^y (B^a_{\;\:j} \: \sigma^{jb}
	+ B^b_{\;\:j} \: \sigma^{ja}) \\
&&+ \frac{i}{8 \pi^3} \: \rho \; \varepsilon^{ijab} \;
	\frac{\partial}{\partial x^a} \frac{\partial}{\partial y^b}
	\lint_x^y B_{ij} \spc .
\end{eqnarray*}
Mit Hilfe von Satz \ref{a1_dis_abl} und Satz \ref{a1_dis_abl2} berechnen wir
nun die auftretenden Terme nacheinander, wobei wir alle Beitr\"age der
Ordnung $\xi^0$ weglassen.
\begin{eqnarray*}
\lefteqn{ \frac{\partial}{\partial x^j} \frac{\partial}{\partial y^k}
	\lint_x^y B^{jk} \;=\; \frac{\partial}{\partial y^k} \lint_x^y
	\partial_j B^{jk} - \frac{\partial^2}{\partial y^j \: \partial y^k}
	\lint_x^y B^{jk} } \\
&=& - \frac{\partial}{\partial y^j} \lint_x^y B^{jk}_{\;\;,k} \;=\;
	- \pi \: l^\vee(\xi) \int_x^y \xi_j \: B^{jk}_{\;\;,k} +
	{\cal{O}}(\xi^0) \\
\lefteqn{ \frac{\partial}{\partial x^a} \frac{\partial}{\partial y_a}
	\lint_x^y B_{jk} \;=\; \frac{\partial}{\partial y^a} \lint_x^y
	\partial^a B_{jk} - \Box_y \lint_x^y B_{jk} } \\
&=& \pi \: l^\vee(\xi) \int_x^y \xi^a B_{jk,a} - 2 \pi \:
	l^\vee(\xi) \; B_{jk}(y) + {\cal{O}}(\xi^0) \\
&=& - \pi \: l^\vee(\xi) \; (B_{jk}(y) + B_{jk}(x)) + {\cal{O}}(\xi^0) \\
\lefteqn{ \frac{\partial}{\partial x^a} \frac{\partial}{\partial y^b}
	\lint_x^y \left( B^a_{\;\:j} \: \sigma^{jb} + B^b_{\;\:j} \:
	\sigma^{ja} \right)} \\
&=& \frac{\partial}{\partial y^b} \lint_x^y \left(\partial_a B^a_{\;\:j}
	\: \sigma^{jb} + \partial_a B^b_{\;\:j} \: \sigma^{ja} \right)
	- 2 \: \frac{\partial^2}{\partial y^a \: \partial
	y^b} \lint_x^y B^a_{\;\:j} \: \sigma^{jb} \\
&=& - 4 \pi \: m^\vee(\xi) \int_x^y B_{aj} \: \xi^a \: \xi_b \; \sigma^{jb}
	\;+\; 2 \pi \: l^\vee(\xi) \int_x^y B_{jk} \; \sigma^{jk} \\
&&- \pi \: l^\vee(\xi) \int_x^y (\alpha^2-\alpha) \; (\Box B_{aj}) \;
	\xi^a \: \xi_b \; \sigma^{jb} \\
&&- \pi \: l^\vee(\xi) \int_x^y (2\alpha-1) \; \left( \xi^a \:
	B_{aj,b} + \xi_b \:  B_{aj,}^{\;\;\;a} \right) \; \sigma^{jb}
	+ {\cal{O}}(\xi^0) \\
\lefteqn{ \varepsilon^{ijab} \: \frac{\partial}{\partial x^a}
	\frac{\partial}{\partial y^b} \lint_x^y B_{ij} \;=\;
	\varepsilon^{ijab}
	\left( \frac{\partial}{\partial y^b} \lint_x^y \partial_a B_{ij} -
	\frac{\partial^2}{\partial y^a \: \partial y^b} \lint_x^y B_{ij}
	\right) } \\
&=& \varepsilon^{ijkl} \: \frac{\partial}{\partial y^l} \lint_x^y \partial_k
	B_{ij}
	\;=\; \pi \: l^\vee(\xi) \int_x^y \varepsilon^{ijkl} \; B_{ij,k} \;
	\xi_l + {\cal{O}}(\xi^0)
\end{eqnarray*}
\QED

\section{Differentialst\"orung durch Vektorpotential}
Wir betrachten jetzt die St\"orung des Diracoperators
\Equ{a1_va}
G \;=\; i \Pdd + i \: L^j \partial_j + \frac{i}{2} \: L^j_{\;,j}
\EndEqu
mit einem reellen Vektorfeld $L$. F\"ur $\Delta k_0$ hat man in 
erster Ordnung in $L$
\begin{eqnarray}
\Delta k_0(x,y) &=& -i \left( s_0 \: L^j \partial_j \: k_0 + k_0 \: L^j
	\partial_j \: s_0 \right)(x,y)
	- \frac{i}{2} \left( s_0 \: L^j_{\;,j} \: k_0 + k_0 \: L^j_{\;,j} \:
	s_0 \right)(x,y) \nonumber \\
&=& i \frac{\partial}{\partial y^j} \left( s_0 \: L^j \: k_0 + k_0 \: L^j \:
	s_0 \right)(x,y)
	- \frac{i}{2} \left( s_0 \: L^j_{\;,j} \: k_0 + k_0 \: L^j_{\;,j} \:
	s_0 \right)(x,y) \nonumber \\
\label{eq:a1_v0}
&=& -i \frac{\partial}{\partial y^j} \: \Delta k_0[L^j](x,y) + \frac{i}{2}
	\: \Delta k_0[L^j_{\;,j}](x,y) \spc .
\end{eqnarray}
Damit k\"onnen wir die Rechnung zum Teil auf diejenige f\"ur skalare
St\"orungen, Theorem~\ref{theorem_sk0}, zur\"uckf\"uhren.

\begin{Thm}
\label{theorem_vk0}
In erster Ordnung St\"orungstheorie gilt
\begin{eqnarray}
\label{eq:a1_v1}
\Delta k_0(x,y) &=& - \frac{1}{4 \pi^2} \; (m^\vee(\xi)-m^\wedge(\xi)) \;
	(L^j(y) + L^j(x)) \: \xi_j \\
\label{eq:a1_v2}
&&- \frac{i}{4 \pi^2} \; (m^\vee(\xi)-m^\wedge(\xi)) \int_x^y L_{m,j} \;
	\xi^m \: \xi_k \; \sigma^{jk} \\
\label{eq:a1_v3}
&&- \frac{1}{16 \pi^2} \; (l^\vee(\xi)-l^\wedge(\xi)) \; (L^j_{\;,j}(y)
	- L^j_{\;,j}(x)) \\
\label{eq:a1_v4}
&&- \frac{i}{16 \pi^2} \; (l^\vee(\xi)-l^\wedge(\xi)) \int_x^y
	(2\alpha-1) \: L^m_{\;,mj} \: \xi_k \; \sigma^{jk} \\
\label{eq:a1_v5}
&&- \frac{i}{8 \pi^2} \; (l^\vee(\xi)-l^\wedge(\xi)) \int_x^y L_{k,j} \;
	\sigma^{jk} \\
\label{eq:a1_v6}
&&- \frac{i}{16 \pi^2} \; (l^\vee(\xi)-l^\wedge(\xi)) \int_x^y 
	(\alpha^2-\alpha) \: (\Box L_{m,j}) \; \xi^m \: \xi_k \; \sigma^{jk}\\
\label{eq:a1_v7}
&&+ \frac{1}{16 \pi^2} \; (l^\vee(\xi)-l^\wedge(\xi)) \int_x^y (\Box L_j) \:
	\xi^j \\
&&+ {\cal{O}}(\xi^0) \spc . \nonumber
\end{eqnarray}
\end{Thm}
{\Beweis}
Nach Theorem~\ref{theorem_sk0} hat man
\begin{eqnarray*}
-i \frac{\partial}{\partial y^j} \: \Delta k_0[L^j](x,y) &=&
	-\frac{1}{4 \pi^2} \: (m^\vee(\xi)-m^\wedge(\xi)) \;
	(L^j(y)+L^j(x)) \: \xi_j \\
&&- \frac{1}{8 \pi^2} \: (l^\vee(\xi)-l^\wedge(\xi)) \; L^j_{\;,j}(y) \\
&&- \frac{i}{4 \pi^2} \: (m^\vee(\xi)-m^\wedge(\xi)) \int_x^y (\partial_j
	L_m) \: \xi^m \: \xi_k \; \sigma^{jk} \\
&&- \frac{i}{8 \pi^2} \: (l^\vee(\xi)-l^\wedge(\xi)) \int_x^y \left(
	\alpha \: (\partial_{jm} L^m) \: \xi_k + (\partial_j L_k) \right)
	\: \sigma^{jk} \\
&&- \frac{i}{16 \pi^2} \: (l^\vee(\xi)-l^\wedge(\xi)) \int_x^y (\alpha^2
	-\alpha) \; (\Box L_{m,j}) \; \xi^m \: \xi_k \; \sigma^{jk} \\
&&+ \frac{1}{16 \pi^2} \: (l^\vee(\xi)-l^\wedge(\xi)) \int_x^y (\Box L_j)
	\: \xi^j \;+\; {\cal{O}}(\xi^0) \\
\frac{i}{2} \: \Delta k_0[L^j_{\;,j}](x,y) &=& \frac{1}{16 \pi^2}
	\: (l^\vee(\xi)-l^\wedge(\xi)) \; \left( L^j_{\;,j}(y)
	+ L^j_{\;,j}(x) \right) \\
&&+ \frac{i}{16 \pi^2} \int_x^y (\partial_{jm} L^m) \: \xi_k \; \sigma^{jk}
	\;+\; {\cal{O}}(\xi^0) \spc ,
\end{eqnarray*}
setze nun in~\Ref{a1_v0} ein.
\QED

\section{Bilineare Differentialst\"orung durch Vektorpotential}
Wir betrachten jetzt die St\"orung des Diracoperators
\Equ{a1_bva}
G \;=\; i \Pdd + i L_j \:\sigma^{jk} \: \partial_k + \frac{i}{2} \: L_{j,k}
	\sigma^{jk}
\EndEqu
mit einem reellen Vektorfeld $L$. F\"ur $\Delta k_0$ hat man in 
erster Ordnung in $L$
\begin{eqnarray*}
\Delta k_0(x,y) &=& -i \left( s_0 \:  L_j \:\sigma^{jk} \: \partial_k \: k_0
	+ k_0 \: L_j \: \sigma^{jk} \:
	\partial_k \: s_0 \right)(x,y) \\
&&- \frac{i}{2} \left( s_0 \: L_{j,k} \sigma^{jk} \: k_0 + k_0 \: L_{jk}
	\sigma^{jk} \: s_0 \right)(x,y)
\end{eqnarray*}
Wir setzen die Relation
\Equ{a1_345a}
L_j \: \sigma^{jk} \: \partial_k \;=\; -i L^j \partial_j \:+\: L \slsh \:
	(i \Pdd)
\EndEqu
ein und erhalten
\begin{eqnarray}
\Delta k_0(x,y) &=& -i \Delta k_0 \left[ i L^j \partial_j + \frac{i}{2} \: L^j_{\;,j}
	\right](x,y) \:-\: i k_0(x,y) \: L \slsh(y) \nonumber \\
\label{eq:a1_u}
&&+ \frac{i}{2} \: \Delta k_0[L_{j,k} \sigma^{jk}](x,y) \:-\: \frac{1}{2} \:
	\Delta k_0[L^j_{\;,j}](x,y) \spc .
\end{eqnarray}
Damit k\"onnen wir die Rechnung zum Teil auf diejenige f\"ur Differentialst\"orungen,
bilineare St\"orungen und skalare St\"orungen zur\"uckf\"uhren.

\begin{Thm}
\label{theorem_bvk0}
In erster Ordnung St\"orungstheorie gilt
\begin{eqnarray}
\label{eq:a1_bv1}
\Delta k_0(x,y) &=& \frac{1}{4 \pi^2} \; (m^\vee(\xi)-m^\wedge(\xi)) \;
	(L_j(y) + L_j(x)) \: \xi_k \: \sigma^{jk} \\
\label{eq:a1_bv2}
&&- \frac{i}{4 \pi^2} \; (m^\vee(\xi)-m^\wedge(\xi)) \; (L_j(y) - L_j(x)) \: \xi^j \\
\label{eq:a1_bv3}
&& \frac{i}{16 \pi^2} (l^\vee(\xi)-l^\wedge(\xi)) \; (L^k_{\;,k}(y)+L^k_{\;,k}(x)) \\
\label{eq:a1_bv4}
&&- \frac{1}{16 \pi^2} (l^\vee(\xi)-l^\wedge(\xi)) \;
	\int_x^y L^m_{\;,mj} \; \xi_k \; \sigma^{jk} \;+\; {\cal{O}}(\xi^0) \;\;\;,
	\spc
\end{eqnarray}
\end{Thm}
{\Beweis}
Aus Theorem~\ref{theorem_b0} erh\"alt man durch partielle Integration
\begin{eqnarray*}
\lefteqn{ \frac{i}{2} \: \Delta k_0[L_{j,k} \sigma^{jk}](x,y) \;=\;
	\frac{1}{4 \pi^2} \: (m^\vee(\xi)-m^\wedge(\xi)) \;
	\int_x^y L_{i,j} \: \xi^i \xi_k \: \sigma^{jk} } \\
&&-\frac{1}{4 \pi^2} \: (m^\vee(\xi)-m^\wedge(\xi)) \;
	(L_j(y) - L_j(x)) \: \xi_k \: \sigma^{jk} \\
&&+\frac{1}{16 \pi^2} \: (l^\vee(\xi)-l^\wedge(\xi)) \;
	(L_{i,j}(y) + L_{i,j}(x)) \: \sigma^{ij} \\
&&-\frac{1}{4 \pi^2} \: (l^\vee(\xi)-l^\wedge(\xi)) \; \int_x^y L_{i,j} \:
	\sigma^{ij} \\
&&+\frac{i}{16 \pi^2} \: (l^\vee(\xi)-l^\wedge(\xi)) \; \int_x^y (\Box L_j) \: \xi^j \\
&&-\frac{i}{16 \pi^2} \: (l^\vee(\xi)-l^\wedge(\xi)) \;
	(L^k_{\;,k}(y)-L^k_{\;,k}(x)) \\
&&+\frac{1}{16 \pi^2} \: (l^\vee(\xi)-l^\wedge(\xi)) \; (L_{j,i}(y) + L_{j,i}(x))
	\: \sigma^{ij} \\
&&-\frac{1}{8 \pi^2} \: (l^\vee(\xi)-l^\wedge(\xi)) \; \int_x^y L_{j,i} \: \sigma^{ij}
	\\
&&+\frac{1}{16 \pi^2} \: (l^\vee(\xi)-l^\wedge(\xi)) \; \int_x^y (2\alpha-1) \:
	\xi_i \: (\Box L_j) \: \sigma^{ij} \\
&&-\frac{1}{16 \pi^2} \: (l^\vee(\xi)-l^\wedge(\xi)) \; \int_x^y (2\alpha-1) \:
	\xi_i \: L^k_{\;,jk} \: \sigma^{ij} \\
&&+\frac{1}{16 \pi^2} \: (l^\vee(\xi)-l^\wedge(\xi)) \; \int_x^y (\alpha^2-\alpha)
	\: \Box (L_{i,j}) \: \xi^i \: \xi_k \: \sigma^{jk} \\
&&+\frac{1}{16 \pi^2} \: (l^\vee(\xi)-l^\wedge(\xi)) \; \int_x^y (2\alpha-1) \:
	(\Box L_j) \: \xi_k \: \sigma^{jk} \;+\; {\cal{O}}(\xi^0) \spc .
\end{eqnarray*}
Durch Vergleich mit Theorem~\ref{theorem_vk0} erh\"alt man
\begin{eqnarray*}
\lefteqn{ -i \Delta k_0 \left[ i L^j \partial_j + \frac{i}{2} \: L^j_{\;,j}
	\right](x,y) \:-\: i k_0(x,y) \: L \slsh(y)
	\:+\: \frac{i}{2} \: \Delta k_0[L_{j,k} \sigma^{jk}](x,y) } \\
&=& \frac{1}{4 \pi^2} \; (m^\vee(\xi)-m^\wedge(\xi)) \;
	(L_j(y) + L_j(x)) \: \xi_k \: \sigma^{jk} \\
&&- \frac{i}{4 \pi^2} \; (m^\vee(\xi)-m^\wedge(\xi)) \; (L_j(y) + L_j(x)) \: \xi^j
	\;+\; {\cal{O}}(\xi^0) \spc .
\end{eqnarray*}
Die Behauptung folgt mit \Ref{a1_u} und Theorem \ref{theorem_sk0}.
\QED

\chapter{St\"orungsrechnung f\"ur $k_m$ im Ortsraum}
\label{anh2}

In diesem Kapitel werden wir die St\"orungsrechnung f\"ur $k_m$ durchf\"uhren.
Dabei werden wir \"ahnlich wie im Fall $m=0$ in Anhang~\ref{anh1} vorgehen
und Gleichung~\Ref{a1_100} im Ortsraum auswerten.

Der Einfachheit halber nehmen wir $m>0$ an, den allgemeinen Fall
erh\"alt man durch eine PCT-Transformation\footnote{F\"ur eine St\"orung
durch ein elektromagnetisches Feld oder ein Gravitationsfeld hat man
insbesondere
\[      \rho B \rho \;=\; - B   \]
und wegen $k_{-m}=\rho k_m \rho$, $s_{-m}=-\rho s_m \rho$
\begin{eqnarray*}
  \tilde{k}_{-m} &=& k_{-m} - \lambda k_{-m} B s_{-m} -
	\lambda s_{-m} B k_{-m} \\
  &=& \rho k_m \rho + \lambda (\rho k_m \rho) B (\rho s_m \rho) +
	\lambda (\rho s_m \rho) B (\rho k_m \rho)  \\
  &=& \rho \tilde{k}_m \rho \spc .
\end{eqnarray*}}.
Man hat dann
\begin{eqnarray}
\label{eq:a2_1}
  k_m(x) &=& \left( i \Pdd + m \right) K_{m^2}(x) \\
\label{eq:a2_2}
  s_m(x) &=& \left( i \Pdd + m \right) S_{m^2}(x)
\end{eqnarray}
mit
\begin{eqnarray}
  K_{m^2}(x) &=& \int \frac{d^4k}{(2 \pi)^4} \: \delta(k^2-m^2)
    \:\epsilon(k^0) \; e^{- i k x} \nonumber \\
\label{eq:a2_3}
   &=& - \frac{i}{4 \pi^2} \left( \delta(x^2) - m^2 \: \frac{J_1(m \tau)}
	{2 m \tau}
	\Theta(x^2) \right) \epsilon(x^0) \\
  S_{m^2}(x) &=& \int \frac{d^4k}{(2 \pi)^4} \: \frac{1}{k^2-m^2} \; e^{-i k x}
	\nonumber \\
\label{eq:a2_4}
   &=& - \frac{1}{4 \pi} \left( \delta(x^2) - m^2 \: \frac{J_1(m \tau)}
	{2 m \tau}      
	\Theta(x^2) \right) \spc .
\end{eqnarray}
Dabei wurde $\tau = \sqrt{x^2}$ gesetzt, $J_1$ ist eine Besselfunktion
erster Art.

Im Unterschied zu~\Ref{a1_75}, \Ref{a1_76} tritt im Innern des Lichtkegels
ein zus\"atzlicher Beitrag auf, der f\"ur gro{\ss}es $\tau$ asymptotisch wie
\[ K_{m^2}(x), S_{m^2}(x) \;\sim \; \tau^{-\frac{3}{2}}
	\cos( \frac{3}{4} \pi - m \tau) \]
abf\"allt und in der N\"ahe des Lichtkegels mit
\begin{eqnarray*}
  \lim_{\I \ni x \rightarrow y \in \Li} K_{m^2}(x) &=&
	- \frac{i}{4 \pi^2} \left( - \frac{m^2}{4} \right) \epsilon(y^0) \\
  \lim_{\I \ni x \rightarrow y \in \Li} S_{m^2}(x) &=&
	- \frac{1}{4 \pi} \left( - \frac{m^2}{4} \right) \\
\end{eqnarray*}
beschr\"ankt ist.

Durch Einsetzen von~\Ref{a2_3}, \Ref{a2_4} in \Ref{a2_1}, \Ref{a2_2} erh\"alt
man unter Verwendung von
\[  \frac{d}{dx} \left( \frac{J_1(x)}{x} \right) \;=\; - \frac{J_2(x)}{x} \]
f\"ur die Distributionen $k_m$, $s_m$\footnote{Um die Gleichungen~\Ref{a2_5}
und~\Ref{a2_6} auf Vorzeichen und Vorfaktoren zu verifizieren, kann man
direkt nachrechnen, da{\ss} f\"ur $y \neq 0$ tats\"achlich $(i \Pdd - m) k_m(y) =
(i \Pdd - m) s_m(y) = 0$ gilt.
Dazu verwendet man die Relation
\[  \Pdd_y \left( y \slsh l^\vee(y) \right) \;=\; 2 l^\vee(y) \spc , \]
die man durch
\begin{eqnarray*}
  \Pdd_y \left( y \slsh \delta(y^2) \right) &=& 4 \delta(y^2) + 2 y \slsh
	y \slsh \delta^\prime(y^2) \\
   &=& 4 \delta(y^2) + 2 y^2 \delta^\prime(y^2) \;=\; 2 \delta(y^2)
\end{eqnarray*}
herleiten kann, sowie die Beziehungen zwischen Besselfunktionen
\begin{eqnarray*}
  \frac{d}{dx} \left( \frac{J_2(x)}{x^2} \right) &=& - \frac{J_3(x)}{x^2} \\
  x J_1(x) + x J_3(x) &=& 4 J_2(x) \spc .
\end{eqnarray*}}
\begin{eqnarray}
\label{eq:a2_5}
  k_m(x) &=& k_0(x) + \left( m - \frac{i}{2} \: m^2 \; x \slsh \right) K_0(x)
	- \frac{i}{4 \pi^2} \Upsilon(x) \; \epsilon(x^0)  \\
\label{eq:a2_6}
  s_m(x) &=& s_0(x) + \left( m - \frac{i}{2} \: m^2 \; x \slsh \right) S_0(x)
	- \frac{1}{4 \pi} \Upsilon(x)
\end{eqnarray}
mit
\Equ{a2_7}
  \Upsilon(x) \;=\; - \frac{1}{2} m^3 \frac{J_1(m \tau)}{m \tau} \Theta(x^2)
    + \frac{i}{2} m^4 \: x \slsh \:\frac{J_2(m \tau)}{m^2 \tau^2} \Theta(x^2)
   \spc .
\EndEqu
F\"ur die Randwerte von $\Upsilon$ auf dem Lichtkegel hat man
\[ \lim_{\I \ni x \rightarrow y \in \Li} \Upsilon(x) \;=\; -\frac{1}{4}
   m^3 + \frac{i}{16} m^4 \; y \slsh  \spc . \]
F\"ur das weitere Vorgehen wollen wir~\Ref{a2_5}, \Ref{a2_6} in~\Ref{a1_100}
einsetzen und explizit ausrechnen. Zus\"atzlich zu den Lichtkegelintegralen von
Abschnitt~\ref{lichtkegelint} tritt dabei ein weiterer Typ von Integralen
auf, beispielsweise das Integral
\[  \int d^4z \; l^\vee_x(z) \: \Theta((y-z)^2) \; f(z) \]
der Funktion $f$ \"uber den Schnitt des Lichtkegels um den Punkt $x$ mit
dem Inneren des Lichtkegels um den Punkt $y$.

Es ist hilfreich, solche sogenannten inneren Lichtkegelintegrale zun\"achst
allgemeiner zu untersuchen:

\section{Innere Lichtkegelintegrale}

\begin{Def}
Setze
\begin{eqnarray*}
  \Theta^\vee(y) &=& \Theta(y^2) \: \Theta(y^0) \\
  \Theta^\wedge(y) &=& \Theta(y^2) \: \Theta(-y^0)
\end{eqnarray*}
sowie $\Theta^\vee_x(y) = \Theta^\vee(y-x)$, $\Theta^\wedge_x(y)
= \Theta^\wedge(y-x)$.

Definiere f\"ur $f \in C^\infty(M)$ die {\bf inneren Lichtkegelintegrale}
durch
\begin{eqnarray}
\label{eq:a2_11}
  \left( \vint f \right) (y) &=& \int d^4z \; l^\vee(z) \: \Theta^\wedge_y(z)
	\; f(z) \\
\label{eq:a2_12}
  \left( \wint f \right) (y) &=& \int d^4z \; \Theta^\vee(z) \: l^\wedge_y(z)
	\; f(z) \\
\label{eq:a2_13}
  \left( \olint f \right) (y) &=& \int d^4z \; \Theta^\vee(z)
	\: \Theta^\wedge_y(z) \; f(z) \spc .
\end{eqnarray}
\end{Def}
Offensichtlich sind die Integrale~\Ref{a2_11} bis~\Ref{a2_13} wohldefiniert
und endlich. Auf $y \in \I^\vee$ sind sie glatte Funktionen.
F\"ur $y \not\in \Li^\vee \cup \I^\vee$ verschwinden die inneren
Lichtkegelintegrale, weil sich die entsprechenden Lichtkegel
um den Ursprung und den Punkt $y$ dann nicht schneiden.

Die Bezeichungen von Definition~\ref{def_a7} werden in analoger Weise auch
f\"ur die inneren Lichtkegelintegrale verwendet.

Die inneren Lichtkegelintegrale k\"onnen auch durch gew\"ohnliche
Lichtkegelintegrale ausgedr\"uckt werden:
\begin{Lemma}
\label{lemma_b1}
Es gilt:
\begin{eqnarray}
\label{eq:a2_14a}
  \left( \vint f \right) (y) &=& y^2 \int_0^1 d \alpha \; \alpha
    \lint_0^{\alpha y} f \\
\label{eq:a2_14}
  \left( \wint f \right) (y) &=& y^2 \int_0^1 d \alpha \; (1 - \alpha)
    \lint_{\alpha y}^y f \\
\label{eq:a2_15}
  \left( \olint f \right) (y) &=& y^4 \int_0^1 d \alpha \int_\alpha^1 d \beta \;
    (\beta - \alpha)^2 \; \lint_{\alpha y}^{\beta y} f
\end{eqnarray}
\end{Lemma}
{\Beweis}
Man hat
\[ y^2 \int_0^1 d\alpha \; \alpha \lint_0^{\alpha y} f \;=\; y^2 \int_0^1
  d\alpha \; \alpha \int d^4z \; l^\vee(z) \: l^\wedge
	(z-\alpha y) \; f(z) \spc . \]
Nach Satz~\ref{a1_lemma1} existiert das innere Integral und ist in $\alpha$
stetig, man kann also die beiden Integrale vertauschen:
\begin{eqnarray*}
  &=& y^2 \int d^4z f(z) l^\vee(z) \int_0^1 d\alpha \; \alpha \;
    \delta((z-\alpha y)^2) \: \Theta(\alpha y^0 - z^0) \\
  &=& y^2 \int d^4z f(z) l^\vee(z) \int_0^1 d\alpha \; \alpha \;
    \delta(-2 \alpha \bra y,z \ket + \alpha^2 y^2) \;
	\Theta(\alpha y^0 - z^0) \spc .
\end{eqnarray*}
Die Integration \"uber $\alpha$ kann ausgef\"uhrt werden. Man erh\"alt einen
Beitrag bei $\alpha = 2 \: y^{-2} \bra y,z \ket$ falls $z \in \I^\wedge_y$,
genauer
\[  \;=\; \int d^4z \; f(z) \; l^\vee(z) \: \Theta^\wedge_y(z) \;=\; 
	\left( \vint f \right) (y) \spc . \]
Gleichung~\Ref{a2_14} erh\"alt man ganz analog. Weiterhin gilt
\begin{eqnarray*}
\lefteqn{y^4 \int_0^1 d\alpha \int_0^1 d\beta \; (\beta - \alpha)^2
	\lint_{\alpha y}^{\beta y} f} \\
  &=& y^4 \int_0^1 d\alpha \int_0^1 d\beta \; (\beta - \alpha)^2 \int d^4z
    \; l^\vee(z-\alpha y) \: l^\wedge(z-\beta y) \; f(z) \spc.
\end{eqnarray*}
Die Integration \"uber $\alpha$, $\beta$ kann ausgef\"uhrt werden und liefert
Beitr\"age bei
\begin{eqnarray*}
  \alpha &=& \frac{1}{y^2} \left( \bra z,y \ket - \sqrt{ \bra z,y \ket^2 -
   z^2 y^2} \right) \\
  \beta &=& \frac{1}{y^2} \left( \bra z,y \ket + \sqrt{ \bra z,y \ket^2 -
   z^2 y^2} \right)
\end{eqnarray*}
falls $z \in \I^\vee \cap \I^\wedge_y$.
Man erh\"alt so gerade~\Ref{a2_15}.
\QED

\begin{Korollar}
Es gilt $\svint f, \swint f, \solint f \in C^\infty(\I^\vee \cup \Li^\vee)$.

Die partiellen Ableitungen beliebiger Ordnung der inneren Lichtkegelintegrale
sind also stetig auf $\Li^\vee \cup \I^\vee$ fortsetzbar.
\end{Korollar}
{\Beweis}
Nach Lemma~\ref{a1_lemma6} gilt
\[  (\vint f)(y) \;=\; y^2 \lint_0^y dz \int_0^1 \alpha f(\alpha z) \: d\alpha
	\spc .   \]
Nach Korollar~\ref{a1_korollar} ist das Lichtkegelintegral, und somit
auch der ganze Ausdruck, in $C^\infty(\I^\vee \cup \Li^\vee)$.

Die Randwerte der partiellen Ableitungen kann man direkt ausrechen: Nach
Satz~\ref{a1_parabl} folgt f\"ur $y \in \I^\vee$
\[ \partial_j \left(\vint f\right)(y) \;=\; 2 \: y_j \lint_0^y dz \int_0^1
 \beta \: f(\beta z) \: d\beta + y^2 \lint_0^y
	h_j \left[ \int_0^1 \beta f(\beta z) d\beta \right]
	\]
und nach Satz~\ref{a1_lemma1} und \ref{integr_umf}
\begin{eqnarray}
\lim_{\I^\vee \ni y \rightarrow z \in \Li^\vee} \partial_j (\vint f)(y) &=&
  \pi \: z_j \int_0^1 d\alpha \int_0^1 d\beta \; \beta \: f(\alpha \beta z)
	\nonumber \\
\label{eq:a2_ka}
&=& \pi \: z_j \int_0^1 (1-\alpha) \: f(\alpha z) \: d\alpha \spc .
\end{eqnarray}
F\"ur die inneren Lichtkegelintegrale $\swint f$ und $\solint f$ kann
man ganz analog argumentieren, man erh\"alt explizit
\begin{eqnarray}
\lim_{\I^\vee \ni y \rightarrow z \in \Li^\vee} \partial_j (\wint f)(y)
\label{eq:a2_kb}
&=& \pi \: z_j \int_0^1 \alpha f(\alpha z) \: d\alpha \\
\lim_{\I^\vee \ni y \rightarrow z \in \Li^\vee} \partial_j (\olint f)(y)
\label{eq:a2_kc}
&=& 0 \spc .
\end{eqnarray}   \QED
Die Gleichungen \Ref{a2_ka}-\Ref{a2_kc} werden sp\"ater noch ben\"otigt.

Mit Hilfe von Lemma~\ref{lemma_b1} kann man in Analogie zu
Satz~\ref{a1_lemma1}
das Verhalten der inneren Lichtkegelintegrale in der N\"ahe des Lichtkegels
$\Li^\vee$ untersuchen.
\begin{Satz}
\label{satz_b2}
F\"ur $y \in \Li^\vee$ gilt
\begin{eqnarray}
\label{eq:a2_16}
\lim_{\I^\vee \ni x \rightarrow y} \frac{1}{x^2} \left( \vint f \right)(x) &=&
  \frac{\pi}{2} \int_0^1 (1 - \alpha) \: f(\alpha y) \; d\alpha \\
\label{eq:a2_17}
\lim_{\I^\vee \ni x \rightarrow y} \frac{1}{x^2} \left( \wint f \right)(x) &=&
  \frac{\pi}{2} \int_0^1 \alpha \: f(\alpha y) \; d\alpha \\
\label{eq:a2_18}
\lim_{\I^\vee \ni x \rightarrow y} \frac{1}{x^4} \left( \olint f \right)(x)&=&
  \frac{\pi}{4} \int_0^1 (\alpha - \alpha^2) \: f(\alpha y) \; d\alpha
\end{eqnarray}
\end{Satz}
{\Beweis}
  Nach Lemma~\ref{lemma_b1} und Lemma~\ref{a1_lemma6} gilt
\begin{eqnarray*}
  \frac{1}{x^2} \left( \vint f \right) (x) &=& \int_0^1 d \alpha \; \alpha
   \lint_0^{\alpha x} f \;=\; \lint_0^x dz \int_0^1 d \alpha \; \alpha
   f(\alpha z) \spc ,
\end{eqnarray*}
  nach Satz~\ref{a1_lemma1} folgt hieraus
\begin{eqnarray*}
  \lim_{\I^\vee \ni x \rightarrow y} \frac{1}{x^2} \left( \vint f \right)(x)
   &=& \frac{\pi}{2} \int_0^1 d\beta \int_0^1 d\alpha \; \alpha
       f(\alpha \beta y) \\
  &=& \frac{\pi}{2} \int_0^1 d\alpha \; (1-\alpha) f(\alpha y) \spc ,
\end{eqnarray*}
wobei die Integralumformungen \ref{integr_umf} angewendet wurden.
Gleichung~\Ref{a2_17} erh\"alt man ganz analog.

Weiterhin hat man nach Lemma~\ref{a1_lemma6}
\begin{eqnarray*}
  \frac{1}{x^4} \left( \olint f \right) (x) &=& \int_0^1 d \alpha\int_\alpha^1
   d\beta \; (\beta - \alpha)^2 \lint_{\alpha x}^{\beta x} f \\
  &=& \lint_0^x dz \int_0^1 d\alpha \int_\alpha^1 d\beta \: (\beta-\alpha)^2
   \; f(\alpha x + (\beta - \alpha) z)
\end{eqnarray*}
und nach Satz~\ref{a1_lemma1}
\begin{eqnarray*}
\lim_{\I^\vee \ni x \rightarrow y} \frac{1}{x^4} \left( \olint f \right)(x)
&=& \frac{\pi}{2} \int_0^1 d\lambda \int_0^1 d\alpha \int_\alpha^1 d\beta
\; (\beta - \alpha)^2 \: f(\alpha y + (\beta - \alpha) \lambda y) \spc .
\end{eqnarray*}
F\"uhre die Variablentransformation $\tilde{\beta} = (\beta-\alpha)/(1-\alpha)$
durch
\[ \;=\; \frac{\pi}{2} \int_0^1 d\alpha \; (1-\alpha)^3 \int_0^1 d\lambda
  \int_0^1 d\tilde{\beta} \; \tilde{\beta}^2 \;
  f(\alpha y + \lambda \tilde{\beta} (1-\alpha) y)      \]
und wende die Integralumformung~\ref{integr_umf} in $\lambda$, $\tilde{\beta}$
an
\[ \;=\; \frac{\pi}{4} \int_0^1 d\alpha \; (1-\alpha)^3 \int_0^1 d\beta \;
  (1-\beta^2) \; f(\alpha y + \beta (1-\alpha) y) \spc . \]
In den neuen Koordinaten $\tilde{\alpha}=1-\alpha$, $\tilde{\beta}=1-\beta$
\[ \;=\; \frac{\pi}{4} \int_0^1 d\tilde{\alpha} \; \tilde{\alpha}^3
  \int_0^1 d\tilde{\beta} \; (2 \tilde{\beta}-\tilde{\beta}^2) \; f(y -
	\tilde{\alpha} \tilde{\beta} y) \]
kann man wieder~\ref{integr_umf} anwenden und erh\"alt
\[ \;=\; \frac{\pi}{4} \int_0^1 d\alpha \; (\alpha - \alpha^2) \; f((1-\alpha) y)
   \;=\; \frac{\pi}{4} \int_0^1 d\alpha \; (\alpha - \alpha^2) \; f(\alpha y)
	\spc . \]
\QED
Es werden noch einfache Formeln f\"ur die Ableitungen der inneren
Lichtkegelintegrale ben\"otigt.
\begin{Lemma}
\label{lemma_b3}
Es gilt im Distributionssinne:
\begin{eqnarray}
\label{eq:a2_21}
  \frac{\partial}{\partial y^j} \: \vint_x^y f &=& 2 \lint_x^y
	dz \; (y-z)_j \; f(z) \\
\label{eq:a2_22}
  \frac{\partial}{\partial x^j} \: \wint_x^y f &=& -2 \lint_x^y
	dz \; f(z) \;(z-x)_j \\
\label{eq:a2_23}
  \frac{\partial}{\partial x^j} \: \vint_x^y f &=& \vint_x^y \partial_j f
	\:-\: 2 \lint_x^y dz \; (y-z)_j \; f(z) \\
\label{eq:a2_24}
  \frac{\partial}{\partial y^j} \: \wint_x^y f &=& \wint_x^y \partial_j f
	\:+\:2  \lint_x^y  dz \; f(z) \; (z-x)_j
\end{eqnarray}
F\"ur $(y-x) \not \in \Li^\vee$ gelten diese Gleichungen sogar punktweise.
\end{Lemma}
{\Beweis}
Berechne zun\"achst die Distributionsableitungen:
\[ - \int d^4y \; (\partial_j g(y)) \: \left(\vint_x^y f\right) \;=\;
  - \int d^4y \; (\partial_j g(y)) \: \int d^4z \; l^\vee_x(z) \:
    \Theta^\wedge_y(z) \; f(z) \]
Man kann die beiden Integrationen vertauschen, weil die einzelnen Integrale
existieren und bez\"uglich der anderen Variablen stetig sind:
\[ \;=\; - \int d^4z \; l^\vee_x(z) \: f(z) \: \int_{\I^\vee_z} d^4y \;
    (\partial_j g)(y)  \]
Jetzt kann man partiell integrieren
\[ \;=\; \int d^4z \; l^\vee_x(z) \: f(z) \: \int d^4y \; g(y) \:
    2 (y-z)_j \; l^\vee_z(y) \]
und nach Satz~\ref{a1_lemma1} wiederum die Integrale vertauschen
\[ \;=\; \int d^4y \: g(y) \: \lint_x^y dz \; 2(y-x)_j f(z)  \spc . \]
Gleichung~\Ref{a2_22} zeigt man auf die gleiche Weise.
\Ref{a2_23} und \Ref{a2_24}
folgen daraus analog zum Beweis von Lemma~\ref{a1_lemma8}.

Da die inneren Lichtkegelintegrale f\"ur $(y-x) \not \in \Li^\vee$ glatte
Funktionen sind, gelten die abgeleiteten Gleichungen dort sogar punktweise.
\QED

\section{St\"orungsrechnung f\"ur das elektromagnetische Feld}
Betrachte nun speziell die St\"orung~\Ref{a1_54a} durch ein
elektromagnetisches Feld. In erster Ordnung St\"orungstheorie hat man
\[  \tilde{k}_m \;=\; k_m + \Delta k_m \]
mit
\Equ{a2_a10}
  \Delta k_m \;=\; - e \left( k_m \:\Aslsh\: s_m + s_m \:\Aslsh\: k_m \right)
   \spc .
\EndEqu
Schreibt man~\Ref{a2_5} und~\Ref{a2_6} in der Form
\begin{eqnarray*}
  k_m(x) &=& - \frac{i}{4 \pi^2} \left\{ (i \Pdd + m - \frac{i}{2} \: m^2\; x
    \slsh) (l^\vee(x) - l^\wedge(x)) \;+\; \Upsilon(x)
    (\Theta^\vee(x)-\Theta^\wedge(x)) \right\}   \\
  s_m(x) &=& - \frac{1}{4 \pi} \left\{ (i \Pdd + m - \frac{i}{2} \: m^2 \; x
    \slsh) (l^\vee(x) + l^\wedge(x)) \;+\; \Upsilon(x)
    (\Theta^\vee(x)+\Theta^\wedge(x)) \right\}
\end{eqnarray*}
und setzt in~\Ref{a2_a10} ein, ergibt sich
\begin{eqnarray}
\Delta k_m(x,y) &=& - \frac{ie}{8 \pi^3} \int d^4z \; \left( i \Pdd_x + m
- \frac{i}{2} m^2 \: (x-z)^j \gamma_j \right) \;\times 
\nonumber\\&& \spc\spc\;\;\;\; \left( l^\vee(x-z) l^\vee(z-y)
- l^\wedge(x-z) l^\wedge(z-y) \right) \Aslsh(z) \;\times
\nonumber\\&& \spc\spc\spc
	\left(i \Pdd_y + m - \frac{i}{2} m^2 \: (z-y)^j \gamma_j \right) \nonumber\\
&& \;\;\;+\; \left( i \Pdd_x + m - \frac{i}{2} m^2 \: (x-z)^j \gamma_j \right)
\; \times
\nonumber\\&& \spc\;\;
\left( l^\vee(x-z) \Theta^\vee(z-y) - l^\wedge(x-z) \Theta^\wedge(z-y) \right)
\Aslsh(z) \Upsilon(z-y) \nonumber\\
&& \;\;\;+\; \Upsilon(x-z) \left(\Theta^\vee(x-z) l^\vee(z-y)
- \Theta^\wedge(x-z) l^\wedge(z-y) \right) \Aslsh(z) \; \times
\nonumber\\&& \spc\;\;
\left(i \Pdd_y + m - \frac{i}{2} m^2 \: (z-y)^j \gamma_j \right) \nonumber\\
&& \;\;\;+\; \Upsilon(x-z) (\Theta^\vee(x-z) \Theta^\vee(z-y) \nonumber\\
&& \spc\spc\spc\;\;\;\; - \Theta^\wedge(x-z)
	\Theta^\wedge(z-y) ) \Aslsh(z) \Upsilon(z-y) \nonumber\\
&=& \frac{ie}{8 \pi^3} \;\;\;\left(i \Pdd_x + m \right) \left(\lint_x^y \Aslsh
	- \lint_y^x \Aslsh \right) \left(i \Pdd_y + m \right) \nonumber\\
&& + \frac{ie}{8 \pi^3} \left( -\frac{i}{2} m^2 \right) \;\;
	\left(i \Pdd_x + m \right) \left(  (\lint_x^y-\lint_y^x) dz \;
	\Aslsh(z) \; (z-y)^j \gamma_j \right) \nonumber\\
&& + \frac{ie}{8 \pi^3} \left( -\frac{i}{2} m^2 \right) \;\;
	\left( (\lint_x^y-\lint_y^x) dz \; (x-z)^j \gamma_j \;
	\Aslsh(z) \right) \left(i \Pdd_y + m \right)  \nonumber\\
&& + \frac{ie}{8 \pi^3} \left( - \frac{1}{4} m^4 \right) \;\;
	(\lint_x^y-\lint_y^x) dz \; (x-z)^j \gamma_j \;
	\Aslsh(z) \; (z-y)^j \gamma_j \nonumber\\
&& + \frac{ie}{8 \pi^3} \;\; \left(i \Pdd_x+m \right) \left(
	(\vint_x^y-\wint_y^x) dz \; \Aslsh(z) \; \Upsilon(z-y) \right) \nonumber\\
&& + \frac{ie}{8 \pi^3} \;\; \left( (\wint_x^y - \vint_y^x) dz \;
	\Upsilon(x-z) \: \Aslsh(z) \right) \left(i \Pdd_y + m \right) \nonumber\\
&& + \frac{ie}{8 \pi^3} \left(- \frac{i}{2} m^2 \right) \;\; \left(
	(\vint_x^y - \wint_y^x) dz \; (x-z)^j \gamma_j \; \Aslsh(z) \;
	\Upsilon(z-y)\right) \nonumber\\
&& + \frac{ie}{8 \pi^3} \left(- \frac{i}{2} m^2 \right) \;\; \left(
	(\wint_x^y - \vint_y^x) dz \; \Upsilon(x-z) \; \Aslsh(z) \; (z-y)^j
	\gamma_j \right) \nonumber\\
\label{eq:a2_a11}
&& + \frac{ie}{8 \pi^3} \;\; \left( (\olint_x^y - \olint_y^x) dz
	\; \Upsilon(x-z) \; \Aslsh(z) \; \Upsilon(z-y) \right) \spc .
\end{eqnarray}
In unseren Anwendungen sind in $\tilde{k}_m(x,y)$
Terme der Ordnung $(y-x)^4$ vernachl\"assigbar.
Daher f\"uhren wir eine Entwicklung um den Lichtkegel durch:

\begin{Thm}
\label{a2_theorem2}
In erster Ordnung St\"orungstheorie gilt
\begin{eqnarray}
\label{eq:a2_c1}
\lefteqn{ \Delta k_m(x,y) \;=\; \Delta k_0(x,y) } \\
\label{eq:a2_c2}
&& - i e \left( \int_x^y A_j \right) \xi^j \; (k_m-k_0)(x,y) \\
\label{eq:a2_c3}
&& + \frac{ie}{16 \pi^3} \: m \; \left( \lint_x^y - \lint_y^x \right) \:
   F_{ij} \: \sigma^{ij} \\
\label{eq:a2_c4}
&& + \frac{e}{8 \pi^3} \: m \; \left( \lint_x^y - \lint_y^x \right)
   dz \: \int_x^z \alpha^2 \; j_k \: \zeta^k \\
\label{eq:a2_c5}
&& + \frac{ie}{16 \pi^3} \: m^2 \; \left( \lint_x^y - \lint_y^x \right)
   dz \; F_{ij} \; \gamma^i \: (2 z -x-y)^j \\
\label{eq:a2_c6}
&& - \frac{e}{32 \pi^3} \: m^2 \; \left( \lint_x^y - \lint_y^x \right)
   \varepsilon^{ijkl} \: F_{ij} \: \xi_k \; \rho \gamma_l \\
\label{eq:a2_c7}
&& + \frac{ie}{16 \pi^3} \: m^2 \; \left(\lint_x^y - \lint_y^x\right)
   dz \; \int_x^z \alpha^2 \; j_k \: \zeta^k \; \xi \slsh \\
\label{eq:a2_c8}
&& - \frac{i e}{128 \pi^2} \: m^3 \; \left(\int_x^y F_{ij} \: \sigma^{ij}
   \right)
   \; \left(\Theta^\vee(\xi) - \Theta^\wedge(\xi)\right) \; \xi^2 \\
\label{eq:a2_c9}
&& + \frac{e}{64 \pi^2} \: m^3 \; \left(\int_x^y (\alpha^2-\alpha) \;
   j_k \: \xi^k \right) \; \left(\Theta^\vee(\xi) -
   \Theta^\wedge(\xi)\right) \; \xi^2 \\
\label{eq:a2_c10}
&& + \frac{e}{512 \pi^2} \: m^4 \; \left(\int_x^y
   \varepsilon^{ijkl} \: F_{ij} \: \xi_k \; \rho \gamma_l \right) \;
   \left(\Theta^\vee(\xi) - \Theta^\wedge(\xi)\right) \; \xi^2 \\
\label{eq:a2_c11}
&& + \frac{ie}{256 \pi^2} \: m^4 \; \left(\int_x^y (1-2\alpha) \: F_{ij}
   \: \gamma^i \: \xi^j \right) \;
   \left(\Theta^\vee(\xi) - \Theta^\wedge(\xi)\right) \; \xi^2 \\
\label{eq:a2_c12}
&& + \frac{ie}{256 \pi^2} \: m^4 \; \left(\int_x^y (\alpha^2-\alpha) \;
   j_k \: \xi^k \; \xi \slsh \right) \;
   \left(\Theta^\vee(\xi) - \Theta^\wedge(\xi)\right) \; \xi^2 \\
&& + {\cal{O}}(\xi^4) \spc , \nonumber
\end{eqnarray}
wobei $\xi=y-x$, $\zeta=z-x$ gesetzt wurde.
\end{Thm}
{\Beweis}
Es soll nur der Fall $\xi^0>0$ betrachtet werden, die Formel f\"ur
allgemeines $\xi$ erh\"alt man genau wie im Beweis von
Theorem~\ref{a1_theorem1} durch Punktspiegelung des Minkowski-Raumes.

Betrachte in~\Ref{a2_a11} die Terme verschiedener Ordnung in $m$
nacheinander. Beachte dabei, da{\ss} $\Upsilon$ nach~\Ref{a2_7} erst ab
dritter Ordnung in $m$ beitr\"agt.
\begin{enumerate}
\item Terme $\sim m^0$:\\
F\"uhrt auf die St\"orungsrechnung f\"ur $m=0$, vergleiche~\Ref{a1_56a}.
\item Terme $\sim m$:\\ Man hat den Beitrag
\Equ{a2_b4}
   - \frac{e}{8 \pi^3} \: m \left\{ \Pdd_x \left(\lint_x^y \Aslsh \right)
   \;+\; \left(\lint_x^y \Aslsh \right) \Pdd_y \right\}  \spc .
\EndEqu
Forme unter Verwendung von Lemma~\ref{a1_lemma8} und Satz~\ref{a1_dis_abl}
weiter um\footnote{Dabei ist $\Pdd \Aslsh=\gamma^i \gamma^j
(\partial_i A_j)$.}:
\begin{eqnarray*}
\lefteqn{\Pdd_x \left(\lint_x^y \Aslsh \right) \;+\; \left(\lint_x^y
  \Aslsh \right) \Pdd_y
\;=\; \lint_x^y (\Pdd \Aslsh) \;-\; \frac{\partial}{\partial y^j}
   \lint_x^y (\gamma^j \Aslsh + \Aslsh \gamma^j ) } \\
&=& \lint_x^y (\Pdd \Aslsh) \;-\; 2 \frac{\partial}{\partial y^j}
   \lint_x^y A^j \\
&=& - 2 \pi l^\vee(\xi) \; \left(\int_x^y A_j \right) \xi^j
 + \lint_x^y (\Pdd \Aslsh)
 - 2 \lint_x^y (\partial_j A^j) \\
&& \;\; + 2 \lint_x^y dz \int_x^z
	\alpha \; \partial_j A^j
 + \lint_x^y dz \int_x^z \alpha^2 \; (\Box A_j) \; \zeta^j
\end{eqnarray*}
Ersetze nun $\Box A_j = - j_j + \partial_j \partial_k A^k$, integriere
in $\alpha$ partiell
\begin{eqnarray*}
&=& - 2 \pi l^\vee(\xi) \; \left(\int_x^y A_j \right) \xi^j
	- \lint_x^y dz \int_x^z \alpha^2 \; j_k \: \zeta^k
	+ \lint_x^y (\Pdd \Aslsh - \partial_j A^j)
\end{eqnarray*}
und wende die Beziehung $F_{ij} \gamma^i \gamma^j = - i F_{ij} \:
\sigma^{ij}$ an
\begin{eqnarray}
\label{eq:a2_b1}
&=& - 2 \pi l^\vee(\xi) \; \left(\int_x^y A_j \right) \xi^j
\label{eq:a2_b2}
- \lint_x^y dz \int_x^z \alpha^2 \; j_k \: \zeta^k
\label{eq:a2_b3}
- \frac{i}{2} \lint_x^y F_{ij} \: \sigma^{ij} \;\; .
\end{eqnarray}
Bei Einsetzen in~\Ref{a2_b4} liefert der erste Summand von
\Ref{a2_b1} den Term erster Ordnung
in $m$ von~\Ref{a2_c2}; der zweite und dritte Summand f\"uhrt
auf~\Ref{a2_c4} bzw. \Ref{a2_c3}.
\item Terme $\sim m^2$:\\ Man hat die Beitr\"age
\[ \frac{ie}{8 \pi^3} \: m^2 \left\{ \lint_x^y \Aslsh + \frac{1}{2} \Pdd_x
   \left(\lint_x^y dz \; \Aslsh(z) \; (z-y)^j \gamma_j \right) \right. \]
\Equ{a2_b6}
\spc\spc\;\;\; \left. + \frac{1}{2} \left(\lint_x^y dz \; (x-z)^j \gamma_j \;
   \Aslsh(z) \right) \Pdd_y \right\} \spc .
\EndEqu
Forme zun\"achst um:
\begin{eqnarray}
\lefteqn{ \Pdd_x \left( \lint_x^y dz \; \Aslsh(z) \; (z-y)^j \gamma_j
   \right) } \nonumber \\
&=& \lint_x^y (\Pdd \Aslsh) \; (z-y)^j \gamma_j \;-\; \Pdd_y
   \lint_x^y dz \; \Aslsh(z) \; (z-y)^j \gamma_j \nonumber \\
&=& \lint_x^y (\Pdd \Aslsh) \; (z-y)^j \gamma_j \;-\; 2
   \frac{\partial}{\partial y^k}
   \lint_x^y dz \; A^k(z) \; (z-y)^j \gamma_j \nonumber \\
&& \spc + \frac{\partial}{\partial y^k} \lint_x^y dz \; \Aslsh(z) \; \gamma^k
   \; (z-y)^j \gamma_j \nonumber \\
\label{eq:a2_155a}
&=& \lint_x^y (\Pdd \Aslsh) \; (z-y)^j \gamma_j \;-\; 2
   \frac{\partial}{\partial y^k}
   \lint_x^y dz \; A^k(z) \; (z-y)^j \gamma_j
   - 2 \lint_x^y \Aslsh \;\;,\spc
\end{eqnarray}
wobei Lemma~\ref{a2_lemma10} angewendet wurde. Auf analoge Weise erh\"alt
man
\begin{eqnarray}
\lefteqn{ \left(\lint_x^y dz \; (x-z)^j \gamma_j \; \Aslsh(z) \right)
   \Pdd_y } \nonumber \\
\label{eq:a2_155b}
&=& \lint_x^y dz \; (x-z)^j \gamma_j \; (\Pdd \Aslsh) - 2
   \frac{\partial}{\partial y^k} \lint_x^y dz \; (x-z)^j \gamma_j \;
   A^k(z) - 2 \lint_x^y \Aslsh \spc
\end{eqnarray}
und somit insgesamt
\begin{eqnarray*}
\lefteqn{\frac{1}{2} \: \Pdd_x \left(\lint_x^y dz \; \Aslsh(z) \; (z-y)^j
   \gamma_j \right) + \frac{1}{2} \left(\lint_x^y dz \; (x-z)^j \gamma_j
   \; \Aslsh(z) \right) \Pdd_y + \lint_x^y \Aslsh } \\
&=& \frac{1}{2} \lint_x^y dz \; \left( (\Pdd \Aslsh) \; (z-y)^j \gamma_j +
   (x-z)^j \gamma_j \; (\Pdd \Aslsh) \right)
	+ \frac{\partial}{\partial y^k} \lint_x^y A^k \: \xi \slsh 
	- \lint_x^y \Aslsh \;\; .
\end{eqnarray*}
Satz~\ref{a1_dis_abl} liefert
\begin{eqnarray*}
&=& \pi l^\vee(\xi) \; \left(\int_x^y A_j \right) \xi^j \: \xi \slsh +
   \lint_x^y (\partial_j A^j) \; \xi \slsh \\
&& - \lint_x^y dz \int_x^z \alpha \: \partial_j A^j \; \xi \slsh -
   \frac{1}{2} \lint_x^y dz \int_x^z \alpha^2 \;
   (\Box A_j) \; \zeta^j \; \xi\slsh \\
&& + \frac{1}{2} \lint_x^y dz \; \left( (\Pdd \Aslsh) \; (z-y)^j  \gamma_j
	+ (x-z)^j \gamma_j \; (\Pdd \Aslsh) \right) \spc .
\end{eqnarray*}
Setze nun wieder $\Box A_j = - j_j + \partial_j \partial_k A^k$ und
integriere partiell, ersetze au{\ss}erdem $\Pdd \Aslsh = \partial_j A^j
+ \frac{1}{2} F_{ij} \: \gamma^i \gamma^j$
\begin{eqnarray*}
&=& \pi l^\vee(\xi) \; \left(\int_x^y A_j \right) \xi^j \: \xi \slsh
+ \frac{1}{2} \lint_x^y dz \; \int_x^z \alpha^2 \; j_k \: \zeta^k \;
   \xi \slsh \\
&& + \frac{1}{4} \lint_x^y dz \; \left(  F_{jk} \: \gamma^j \gamma^k \;
   (z-y)^l \gamma_l + (x-z)^l \gamma_l \; F_{jk} \: \gamma^j \gamma^k
   \right)
\end{eqnarray*}
und wende~\Ref{a1_70a} an:
\begin{eqnarray*}
&=& \pi l^\vee(\xi) \; \left(\int_x^y A_j \right) \xi^j \; \xi \slsh
\;+\; \frac{1}{2} \lint_x^y dz \int_x^z \alpha^2 \; j_k \: \zeta^k \;
   \xi \slsh \\
&& + \frac{1}{2} \lint_x^y dz \; \gamma^i F_{ij} \: (2z-x-y)^j
\;+\; \frac{i}{4} \lint_x^y dz \; \varepsilon^{ijkl} \; F_{ij} \: \xi_k
   \; \rho \gamma_l
\end{eqnarray*}
Bei Einsetzen in~\Ref{a2_b6} liefert der erste Summand den Term zweiter
Ordnung in $m$ von~\Ref{a2_c2}; die restlichen Summanden f\"uhren
auf~\Ref{a2_c5} bis~\Ref{a2_c7}.

\item Terme $\sim m^3$:\\
Durch Taylorentwicklung der Funktion $\Upsilon$ nach $m$ erh\"alt man den
Beitrag
\begin{eqnarray}
-\frac{e}{8 \pi^3} m^3 \left\{ - \frac{1}{2} \lint_x^y dz \: \left( (x-z)^j
  \gamma_j \: \Aslsh(z) + \Aslsh(z) \: (z-y)^j \gamma_j \right) \right.
	\nonumber \\
\label{eq:a2_120}
\spc \left. -\frac{1}{4} \Pdd_x \left( \vint_x^y \Aslsh \right) - \frac{1}{4} 
  \left(\wint_x^y \Aslsh \right) \Pdd_y \;\;\; \right\} \spc .
\end{eqnarray}
Setzt man die Umformungen nach Lemma~\ref{lemma_b3}
\begin{eqnarray}
\label{eq:a2_156z}
  \Pdd_x \left(\vint_x^y \Aslsh \right) &=& \vint_x^y \Pdd \Aslsh
	+ 2 \lint_x^y dz \; (z-y)^j \gamma_j \: \Aslsh(z) \\
\label{eq:a2_156a}
  \left(\wint_x^y \Aslsh \right) \Pdd_y &=& - \wint_x^y (\partial_j \Aslsh)
	\gamma^j + 2 \lint_x^y dz \; \Aslsh(z) \: (x-z)^j \gamma_j
\end{eqnarray}
ein, erh\"alt man
\Equ{a2_bb6}
\Ref{a2_120} \;=\;  -\frac{e}{8 \pi^3} m^3 \left\{\lint_x^y A_j \: \xi^j -
  \frac{1}{4} \vint_x^y \Pdd \Aslsh + \frac{1}{4} \wint_x^y
  (\partial_j \Aslsh) \gamma^j \right\} \spc .
\EndEqu
In nullter Ordnung in $\xi$ f\"allt der zweite und dritte Summand
von~\Ref{a2_bb6} nach
Satz~\ref{satz_b2} weg, f\"ur den ersten Summanden kann man~\Ref{a1_8}
einsetzen und erh\"alt
\Equ{a2_b7}
  =\; -\frac{e}{16 \pi^2} m^3 \int_x^y A_j \; \xi^j \;
  \Theta^\vee(\xi) \spc {\mbox{in der Ordnung $\xi^0$}.}
\EndEqu
Dies ist auch das richtige Ergebnis, denn~\Ref{a2_b7} stimmt gerade mit dem
Term $\sim m^3$ von~\Ref{a2_c3} \"uberein, w\"ahrend~\Ref{a2_c8}
und~\Ref{a2_c9} in der Ordnung $\xi^0$ verschwinden.

Nach Lemma~\ref{a2_lemma2} bleibt zu \"uberpr\"ufen, da{\ss} die Grenzwerte der
ersten Ableitungen f\"ur $\I_x^\vee \ni y \rightarrow \bar{y} \in \Li^\vee_x$
von \Ref{a2_bb6} und $\Ref{a2_b7}+\Ref{a2_c8}+\Ref{a2_c9}$ \"ubereinstimmen,
also da{\ss}
\[  \lim_{\I_x^\vee \ni y \rightarrow \bar{y} \in \Li^\vee_x}
	\frac{\partial}{\partial y^k} \left(
	\Ref{a2_bb6} - \Ref{a2_b7} - \Ref{a2_c8} - \Ref{a2_c9} \right)
	\;=\; 0 \spc .  \]
Offensichtlich hat man
\begin{eqnarray}
\lefteqn{ \lim_{\I_x^\vee \ni y \rightarrow \bar{y} \in \Li^\vee_x}
  \frac{\partial}{\partial y^k} \left( \Ref{a2_c8} + \Ref{a2_c9} \right) }
	\nonumber \\
&=& - \frac{ie}{64 \pi^2} \: m^3 \int_x^{\bar{y}} F_{ij} \: \sigma^{ij} \;
	\bar{\xi}_k \nonumber \\
\label{eq:a2_110}
&&  - \frac{e}{32 \pi^2} \: m^3 \int_x^{\bar{y}} (\alpha - \alpha^2) \:
	j_m \: \bar{\xi}^m \; \bar{\xi}_k \;\;,
\end{eqnarray}
wobei zur Abk\"urzung $\bar{\xi} = \bar{y} - x$ gesetzt wurde.
Unter Anwendung von~\Ref{a1_41a}, \Ref{a2_ka} und \Ref{a2_kb}
erh\"alt man
\begin{eqnarray*}
\lefteqn{\lim_{\I_x^\vee \ni y \rightarrow {\bar{y}} \in \Li^\vee_x}
\frac{\partial}{\partial y^k} \left(\Ref{a2_bb6} - \Ref{a2_b7} \right)} \\
&=&
- \frac{e}{16 \pi^2} \: m^3 \left\{ -\frac{1}{2} \int_x^{\bar{y}} (\alpha -
\alpha^2) \: (\Box A_j) \: \bar{\xi}^j \; \bar{\xi}_k \right.\\
&& \left. \spc\;\;\;\;\;\;\; - \frac{1}{2} \int_x^{\bar{y}}
(1-\alpha) \: (\Pdd \Aslsh)\; \bar{\xi}_k + \frac{1}{2} \int_x^{\bar{y}}
\alpha\: (\partial_j \Aslsh) \gamma^j \; \bar{\xi}_k \right\}  \\
&=& - \frac{e}{16 \pi^2} \: m^3 \left\{ - \frac{1}{2} \int_x^{\bar{y}} \:
  (\alpha-\alpha^2) \: (\Box A_j) \: \bar{\xi}^j \; \bar{\xi}_k \right.\\
&& \spc\;\;\;\;\;\;\; - \frac{1}{2} \int_x^{\bar{y}} (1-2\alpha) \:
	(\partial_j A^j) \; \bar{\xi}_k
 \left. -\frac{1}{4} \int_x^{\bar{y}} F_{ij} \: \gamma^i \gamma^j \;
	\bar{\xi}_k \right\} \spc .
\end{eqnarray*}
Ersetzt man nun wieder $F_{ij} \gamma^i \gamma^j = -i F_{ij} \sigma^{ij}$,
$\Box A_j = - j_j + \partial_{jk}A^k$ und integriert partiell, erh\"alt
man gerade~\Ref{a2_110}.
\item Terme $\sim m^4$:
Man hat den Beitrag
\begin{eqnarray}
  \frac{ie}{8 \pi^3} m^4 \left\{ -\frac{1}{4} \lint_x^y dz \; (x-z)^j \gamma_j
   \: \Aslsh(z) \: (z-y)^k \gamma_k \nonumber \right. \\
\spc - \frac{1}{16} \Pdd_x \left( \vint_x^y dz \; \Aslsh(z) \: (z-y)^j
	\gamma_j \right) \nonumber \\
\spc - \frac{1}{16} \left(\wint_x^y dz \; (x-z)^j \gamma_j \: \Aslsh(z)
	\right) \Pdd_y \nonumber \\
\label{eq:a2_10a} \left.
\spc - \frac{1}{4} \vint_x^y \Aslsh - \frac{1}{4} \wint_x^y \Aslsh
	\spc \right\} \spc .
\end{eqnarray}
F\"uhre zun\"achst die Ersetzungen
\begin{eqnarray}
\lefteqn{\Pdd_x \left(\vint_x^y dz \; \Aslsh(z) \: (z-y)^j \gamma^j \right)}
	\nonumber \\
&=& \vint_x^y dz \; (\Pdd \Aslsh)(z) \: (z-y)^j \gamma_j + \vint_x^y
	\gamma^j \Aslsh \gamma_j
  - \Pdd_y \vint_x^y dz \; \Aslsh(z) \: (z-y)^j \gamma_j \nonumber \\
&=& \vint_x^y dz \; (\Pdd \Aslsh)(z) \: (z-y)^j \gamma_j - 2 \vint_x^y
	\Aslsh \nonumber \\
\label{eq:a2_194a}
&& \spc + 2 \lint_x^y dz \; (z-y)^j \gamma_j \: \Aslsh(z) \: (z-y)^k \gamma_k
	\\
\lefteqn{ \left( \wint_x^y dz \; (x-z)^j \gamma_j \: \Aslsh(z) \right) \Pdd_y
\;=\; - \frac{\partial}{\partial y^k} \wint_x^y dz \; (x-z)^j \gamma_j \:
	\Aslsh(z) \: \gamma^k } \nonumber \\
&=& - \wint_x^y dz \; (x-z)^j \gamma_j \: (\partial_k \Aslsh)\gamma^k
	-2 \wint_x^y \Aslsh \nonumber \\
\label{eq:a2_194b}
&& \spc + 2 \lint_x^y dz \; (x-z)^j \gamma_j \:
	\Aslsh(z) \: (x-z)^k \gamma_k \\
\lefteqn{\lint_x^y dz \; (x-z)^j \gamma_j \: \Aslsh(z) (\: x-z)^k \gamma_k }
	\nonumber \\
&=& 2 \lint_x^y dz \; (x-z)^j \: A_j(z) \; (x-z)^k \gamma_k - \lint_x^y
	dz \; \Aslsh(z) \; (x-z)^2 \nonumber \\
\label{eq:a2_194c}
&=&  2 \lint_x^y dz \; (x-z)^j \: A_j(z) \; (x-z)^k \gamma_k \\
\label{eq:a2_194d}
\lefteqn{\lint_x^y dz \; (z-y)^j \gamma_j \: \Aslsh(z) \: (z-y)^k \gamma_k
\;=\; 2 \lint_x^y dz \; (z-y)^j \: A_j(z) \; (z-y)^k \gamma_k }
\end{eqnarray}
durch:
\begin{eqnarray}
\label{eq:a2_f1}
\Ref{a2_10a} &=& \; \frac{ie}{8 \pi^3} m^4 \left\{
-\frac{1}{4} \lint_x^y dz \; 
	(x-z)^j \: A_j(z) \: (x-z)^k \gamma_k \nonumber \right. \\
&& \spc\;\;\;\; -\frac{1}{4} \lint_x^y dz \;  (z-y)^j \: A_j(z) \: 
	(z-y)^k \gamma_k \nonumber \\
&& \spc\;\;\;\; -\frac{1}{4} \lint_x^y dz \;
	(x-z)^j \gamma_j \: \Aslsh(z) \: (z-y)^k
	\gamma_k  \\
&& \spc\;\;\;\; -\frac{1}{16} \vint_x^y dz \; (\Pdd \Aslsh)(z) \: (z-y)^j
	\gamma_j \nonumber \\
&& \spc\;\;\;\; + \frac{1}{16} \wint_x^y dz \; (x-z)^j \gamma_j \:
	(\partial_k \Aslsh)(z) \: \gamma^k \nonumber \\
\label{eq:a2_f2}
&&  \spc\;\;\;\; \left. -\frac{1}{8} \vint_x^y \Aslsh - \frac{1}{8}
	\wint_x^y \Aslsh \spc \right\}
\end{eqnarray}
In nullter Ordnung in $\xi$ tragen nur die Lichtkegelintegrale $\slint$
bei; nach~\Ref{a1_8} erh\"alt man
\begin{eqnarray}
&=& - \frac{ie}{64 \pi^2} m^4 \int_x^y \left(\alpha^2 + (1-\alpha)^2 + 2
	\alpha \: (1-\alpha) \right) \: A_j \: \xi^j \; \xi \slsh \;
	\Theta^\vee(\xi) \nonumber \\
\label{eq:a2_f4}
&=& -\frac{ie}{64 \pi^2} m^4 \int_x^y A_j \; \xi^j \; \xi \slsh \;
	\Theta^\vee(\xi) \spc {\mbox{in der Ordnung $\xi^0$.}}
\end{eqnarray}
Dies ist erwartungsgem\"a{\ss} der Beitrag $\sim m^4$ von~\Ref{a2_c2}; beachte,
da{\ss} \Ref{a2_c10} und \Ref{a2_c11} in nullter Ordnung in $\xi$ nicht
beitragen.

Nun kann man ganz analog wie unter 4. vorgehen: Es ist zu \"uberpr\"ufen,
da{\ss}
\[ \lim_{\I^\vee_x \ni y \rightarrow \bar{y} \in \Li_x^\vee}
	\frac{\partial}{\partial y^k} \left(
	\Ref{a2_f1}+\Ref{a2_f2}-\Ref{a2_f4}-\Ref{a2_c10}
	-\Ref{a2_c11}-\Ref{a2_c12} \right) \;=\; 0 \spc . \]
Offensichtlich hat man
\begin{eqnarray}
\lefteqn{\lim_{\I^\vee_x \ni y \rightarrow \bar{y} \in \Li_x^\vee}
  \frac{\partial}{\partial y^k} \left(\Ref{a2_c10} + \Ref{a2_c11} +
	\Ref{a2_c12} \right)}
	\nonumber\\
&=& \frac{e}{256 \pi^2} \: m^4 \; \bar{\xi}_k \left(\int_x^y
   \varepsilon^{ijkl} \: F_{ij} \: \bar{\xi}_k \; \rho \gamma_l \right)
	\nonumber \\
&& + \frac{ie}{128 \pi^2} \: m^4 \; \bar{\xi}_k \left(\int_x^y (1-2 \alpha)
	\; F_{ij} \: \gamma^i \: \bar{\xi}^j \right)
	\nonumber \\
\label{eq:a2_G}
&& - \frac{ie}{128 \pi^2} \: m^4 \; \bar{\xi}_k \left(\int_x^y
   (\alpha-\alpha^2) \;
   j_k \: \bar{\xi}^k \; \bar{\xi} \slsh \right) \spc ,
\end{eqnarray}
wobei wiederum $\bar{\xi} = \bar{y} - x$ gesetzt wurde.
Nach \Ref{a2_ka} und \Ref{a2_kb} folgt au{\ss}erdem:
\begin{eqnarray}
\lefteqn{\lim_{\I^\vee_x \ni y \rightarrow \bar{y} \in \Li_x^\vee}
  \frac{\partial}{\partial y^k} \Ref{a2_f2} } \nonumber \\
&=& \frac{ie}{8 \pi^2} m^4 \: \bar{\xi}_k \: \left\{ \frac{1}{16}
	\int_x^{\bar{y}} (1-\alpha)^2 (\Pdd \Aslsh) \: \bar{\xi}\slsh
	- \frac{1}{16} \int_x^{\bar{y}} \alpha^2 \: \bar{\xi}\slsh \:
	(\partial_k \Aslsh) \gamma^k
	- \frac{1}{8} \int_x^{\bar{y}} \Aslsh \right\} \nonumber \\
&=& \frac{ie}{8 \pi^2} m^4 \: \bar{\xi}_k \: \left\{ \frac{1}{16}
	\int_x^{\bar{y}} (1-2 \alpha) (\partial_j A^j) \: \bar{\xi}\slsh
	- \frac{1}{8} \int_x^{\bar{y}} \Aslsh \right. \nonumber \\
\label{eq:a2_E}
&& \spc\spc \left. + \frac{1}{32} \int_x^{\bar{y}} (1-\alpha)^2 \:
	F_{ij} \gamma^i \gamma^j \: \bar{\xi}\slsh
	+ \frac{1}{32} \int_x^{\bar{y}} \alpha^2 \: \bar{\xi} \slsh \:
	F_{ij} \gamma^i \gamma^j  \right\} \spc .
\end{eqnarray}
Mit Hilfe von \Ref{a1_41a} und~\ref{integr_umf} kann man berechnen:
\begin{eqnarray}
\lefteqn{\lim_{\I^\vee_x \ni y \rightarrow \bar{y} \in \Li_x^\vee}
  \frac{\partial}{\partial y^k} (\Ref{a2_f1}-\Ref{a2_f4}) } \nonumber \\
&=& \frac{ie}{8 \pi^2} m^4 \: \bar{\xi}_k \: \left\{ \frac{3}{4}
	\int_x^{\bar{y}} (\alpha-\alpha^2) \: \Aslsh \right. \nonumber \\
&& \spc\spc + \frac{1}{16} \int_x^{\bar{y}} (\alpha - \alpha^2) \:
	(\Box A_j) \: \bar{\xi}^j \; \bar{\xi} \slsh \nonumber \\
&& \spc\spc + \frac{1}{8} \int_x^{\bar{y}} (-\alpha + 3\alpha^2 - 2\alpha^3)
	\: \bar{\xi}^j \: (\Pdd A_j) \nonumber \\
&& \spc\spc + \frac{1}{16} \int_x^{\bar{y}} (\alpha - 2\alpha^2 + \alpha^3)
	\: F_{ij} \: \gamma^i \gamma^j \: \bar{\xi}\slsh \nonumber \\
\label{eq:a2_F}
&& \left. \spc\spc + \frac{1}{16} \int_x^{\bar{y}} (\alpha^2-\alpha^3) \;
	\bar{\xi}\slsh \: F_{ij} \: \gamma^i \gamma^j \spc \right\}
\end{eqnarray}
Ersetze wieder $\Pdd A_j = \gamma^i F_{ij} + \partial_j \Aslsh$, $\Box
A_j=-j_j + \partial_{jk} A^k$ und integriere jeweils partiell:
\begin{eqnarray}
\Ref{a2_F} &=& \frac{ie}{8 \pi^2} m^4 \: \bar{\xi}_k \: \left\{ -\frac{1}{16}
	\int_x^{\bar{y}} (\alpha-\alpha^2) \: j_j \: \bar{\xi}^j \;
	\bar{\xi}\slsh \right. \nonumber \\
&& \spc\spc - \frac{1}{16} \int_x^{\bar{y}} (1-2\alpha) \; (\partial_j A^j)
	\; \bar{\xi}\slsh \nonumber \\
&& \spc\spc + \frac{1}{8} \int_x^{\bar{y}} (-\alpha + 3 \alpha^2 - 2 \alpha^3)
	\; F_{ij} \: \gamma^i \: \bar{\xi}^j \nonumber \\
&& \spc\spc + \frac{1}{8} \int_x^{\bar{y}} \Aslsh \nonumber \\
&& \spc\spc + \frac{1}{16} \int_x^{\bar{y}} (\alpha - 2\alpha^2 + \alpha^3)
	\: F_{ij} \: \gamma^i \gamma^j \: \bar{\xi}\slsh \nonumber \\
\label{eq:a2_H}
&& \left. \spc\spc + \frac{1}{16} \int_x^{\bar{y}} (\alpha^2-\alpha^3) \;
	\bar{\xi}\slsh \: F_{ij} \: \gamma^i \gamma^j \spc \right\}
\end{eqnarray}
Addiert man die Beitr\"age \Ref{a2_E} und \Ref{a2_H} auf, so fallen alle
von der Eichung abh\"angigen Terme weg, und es ergibt sich
\begin{eqnarray*}
\lefteqn{\lim_{\I^\vee_x \ni y \rightarrow \bar{y} \in \Li_x^\vee}
  \frac{\partial}{\partial y^k} \left(\Ref{a2_f1} + \Ref{a2_f2} -
	\Ref{a2_f4} \right)} \\
&=& \frac{ie}{8 \pi^2} m^4 \: \bar{\xi}_k \: \left\{ -\frac{1}{16}
	\int_x^{\bar{y}}
	(\alpha-\alpha^2) \: j_j \: \bar{\xi}^j \; \bar{\xi}\slsh \right. \\
&& \spc\spc + \frac{1}{8} \int_x^{\bar{y}} (-\alpha +3\alpha^2 - 2\alpha^3)
	\: F_{ij} \: \gamma^i \: \bar{\xi}^j \\
&& \spc\spc + \frac{1}{32} \int_x^{\bar{y}} (1 - 3 \alpha^2 + 2 \alpha^3)
	\: F_{ij} \: \gamma^i \gamma^j \: \bar{\xi}\slsh \\
&& \left. \spc\spc + \frac{1}{32} \int_x^{\bar{y}} (3\alpha^2-2\alpha^3) \;
	\bar{\xi}\slsh \: F_{ij} \: \gamma^i \gamma^j \spc \right\} \spc .
\end{eqnarray*}
Bei Anwendung von~\Ref{a1_70a} ergibt sich gerade~\Ref{a2_G}.
\item Alle h\"oheren Terme in $m$ von $\tilde{k}_m$ sind von der Ordnung
$\xi^4$ und brauchen daher nicht berechnet zu werden.
\end{enumerate}
\QED
Es mu{\ss} noch folgendes Lemma bewiesen werden:
\begin{Lemma}
\label{a2_lemma10}
Es gilt im Distributionssinne
\begin{eqnarray}
\lefteqn{ \Pdd_y \; \lint_x^y dz \: f(z) \; (z-y)^j \gamma_j \;=\;
\frac{\partial}{\partial y^k} \; \lint_x^y dz \: f(z) \; (z-y)^j \gamma_j
	\; \gamma^k } \nonumber \\
\label{eq:a2_333} 
&=& \frac{\partial}{\partial y^k} \; \lint_x^y dz \: f(z) \;
	(z-y)^k \;=\; - 2 \lint_x^y f \\
\lefteqn{ \Pdd_x \; \lint_x^y dz \: (x-z)^j \gamma_j \; f(z) \;=\;
\frac{\partial}{\partial x^k} \; \lint_x^y dz \: (x-z)^j \gamma_j \;
\gamma^k \; f(z) } \nonumber \\
\label{eq:a2_334}
&=& \frac{\partial}{\partial y^k} \; \lint_x^y dz \: (x-z)^k \; f(z)
	\;=\;  2 \lint_x^y f \spc .
\end{eqnarray}
\end{Lemma}
{\Beweis}
Nach Definition des Lichtkegelintegrals hat man f\"ur $g \in C^\infty_c(M)$
\begin{eqnarray*}
\lefteqn{ \left( \Pdd_y \lint_x^y dz \: f(z) \; (z-y)^j \gamma_j \right)
	(g)}\\
&=& - \int d^4z \; f(z) \: l^\vee_x(z) \int d^4y \; l^\vee(y-z) \; \gamma^k \;
	(z-y)^j \gamma_j \; \frac{\partial}{\partial y^k} g(y) \\
&=& - \int d^4z \: f(z) \: l^\vee_x(z) \int d^4y \; \left( 4 l^\vee(y-z) +
	2 m^{\vee}(y-z) \; (y-z)^2 \right) \; g(y) \\
&=& - \int d^4z \: f(z) \: l^\vee_x(z) \int d^4y \; \left( 4 l^\vee(y-z) +
	(y-z)^j \frac{\partial}{\partial y^j} l^\vee(y-z) \right) \; g(y) \\
&=& - 2 \int d^4z \: f(z) \: l^\vee_x(z) \int d^4y \;
	l^\vee(y-z) \: g(y) \;=\; -2 (\lint_x^y f)(g) \spc ,
\end{eqnarray*}
wobei die Bezeichung~\Ref{51a} und Relation~\Ref{a1_54b} verwendet wurden.

Die anderen Gleichungen folgen analog.
\QED
Das Ergebnis der St\"orungsrechnung soll nun wieder kurz diskutiert werden:

Der Summand~\Ref{a2_c2} ist der Eichterm, er f\"uhrt zu einer
Phasenverschiebung von $k_m(x,y)$.

Allgemein sieht man, da{\ss} die Beitr\"age ungerader Ordnung in $m$ proportional
zu den skalaren und bilinearen Kovarianten $\1$, $\sigma^{ij}$ sind; die
Beitr\"age gerader Ordnung in $m$ sind dagegen $\sim \gamma^j, \rho \gamma^j$.

Die Summanden~\Ref{a2_c3} bis~\Ref{a2_c7} sind in der Umgebung des
Lichtkegels beschr\"ankt. F\"ur die Randwerte auf dem Lichtkegel hat man
\begin{Satz}
\label{a2_satzb9}
F\"ur $(y-x) \in \Li$ gilt
\begin{eqnarray*}
\lefteqn{\lim_{\I_x \ni u \rightarrow y} \left( \Delta k_m(x,u) -
\Delta k_0(x,u) \right) }\\
&=& \frac{ie}{32 \pi^2} \; m \: \epsilon(\xi^0) \int_x^y F_{ij} \sigma^{ij} \\
&& - \frac{e}{16 \pi^2} \; m \: \epsilon(\xi^0) \int_x^y (\alpha^2-\alpha) \:
	j_k \: \xi^k \\
&& + \frac{ie}{32 \pi^2} \; m^2 \: \epsilon(\xi^0) \int_x^y (2 \alpha -1) \:
	\gamma^i \: F_{ij} \: \xi^j \\
&& - \frac{e}{64 \pi^2} \; m^2 \: \epsilon(\xi^0) \int_x^y \varepsilon^{ijkl}
	\; F_{ij} \:  \xi_k \; \rho \gamma_l \\
&& - \frac{ie}{32 \pi^2} \; m^2 \: \epsilon(\xi^0) \int_x^y (\alpha^2-\alpha)
	\; j_k \: \xi^k \: \xi \slsh \;+\; {\cal{O}}(m^3) \spc .
\end{eqnarray*}
\end{Satz}
{\Beweis} Die Behauptung folgt unmittelbar mit Hilfe von~\Ref{a1_8}
und~\ref{integr_umf}.
\QED

In unseren Anwendungen spielen die Stromterme~\Ref{a2_c4}
und~\Ref{a2_c7} eine wichtige Rolle. Die Summanden \Ref{a2_c3}, \Ref{a2_c5},
\Ref{a2_c6} fallen werden wegfallen, \Ref{a2_c8}
bis \Ref{a2_c12} werden vernachl\"assigbar sein.

\section{St\"orungsrechnung f\"ur das Gravitationsfeld}
\label{grav_km}
Wie in Abschnitt~\ref{grav_k0} arbeiten wir mit der linearisierten
Gravitationstheorie in symmetrischer Eichung.
Das Koordinatensystem w\"ahlen wir wieder so, da{\ss} die Bedingung
\Ref{a1_210} erf\"ullt ist. Die St\"orung des Diracoperators
\Ref{a1_208} f\"uhrt auf ein $\Delta k_m$ der Form
\begin{eqnarray}
\Delta k_m(x,y) &=& i \frac{\partial}{\partial y^k} \left( k_m \:
	\gamma^j h_j^k \: s_m + s_m \: \gamma^j h_j^k \: k_m \right)(x,y)
	\nonumber \\
&&- \frac{1}{4} \left( k_m \: (i \Pdd h) \: s_m + s_m \: (i \Pdd h) \: k_m
	\right)(x,y) \nonumber
	+ \frac{1}{2} \: (h(x) + h(y)) \; k_m(x,y) \nonumber \\
\label{eq:a2_g10}
&=& \left(\frac{1}{4} \: h(x) + \frac{3}{4} \: h(y) \right) \: k_m(x,y)
	- \frac{i}{e} \frac{\partial}{\partial y^k} \:
	\Delta k_m[\gamma^j h_j^k](x,y) \spc .
\end{eqnarray}
Dabei haben wir die Umformungen
\begin{eqnarray*}
k_m \: (i \Pdd h) \: s_m &=& k_m \: [(i \Pdd - m), h] \: s_m \;=\;
	- k_m \: h \\
s_m \: (i \Pdd h) \: k_m &=&  h \: k_m
\end{eqnarray*}
eingesetzt.
Gleichung \Ref{a2_g10} erm\"oglicht es, die Rechnung wieder zum Teil auf
diejenige f\"ur das elektromagnetische Feld, Theorem \ref{a2_theorem2},
zur\"uckzuf\"uhren.

\begin{Thm}
\label{thm_gpm}
In erster Ordnung St\"orungstheorie gilt
\begin{eqnarray}
\lefteqn{ \Delta k_m(x,y) \;=\; \Delta k_0(x,y) } \\
\label{eq:a2_g11}
&&- \left( \int_x^y h^k_j \right) \: \xi^j \: \frac{\partial}{\partial y^k}
	\left( k_m(x,y) - k_0(x,y) \right) \\
\label{eq:a2_g12}
&&+ \frac{i}{2} \: m \left( \int_x^y h_{ki,j} \right) \: \xi^k \; \sigma^{ij}
	\; \ke(x,y) \\
\label{eq:a2_g13}
&&+ \frac{i}{8 \pi^2} \: m \; (l^\vee(\xi) - l^\wedge(\xi)) \int_x^y
	(\alpha^2-\alpha) \; R_{jk} \: \xi^j \: \xi^k \\
\label{eq:a2_g14}
&&- \frac{i}{16 \pi^3} \: m \left( \lint_x^y - \lint_y^x \right) dz
	\int_x^z (2 \alpha^2 - \alpha) \; R \\
\label{eq:a2_g15}
&&+\frac{i}{32 \pi^3} \: m \left( \lint_x^y - \lint_y^x \right) dz \;
	\zeta^j \: \zeta^k
	\int_x^z (\alpha^4-\alpha^3) \; (R_{\:,jk} - 2 \: \Box R_{jk}) \\
\label{eq:a2_g16}
&&- \frac{1}{16 \pi^3} \: m \left( \lint_x^y - \lint_y^x \right)dz
	\; \zeta^k \int_x^z \alpha^2 \; R_{ki,j} \; \sigma^{ij} \\
\label{eq:a2_g17}
&&+ \frac{1}{16 \pi^2} \: m^2 \; (l^\vee(\xi)-l^\wedge(\xi)) \int_x^y
	(2\alpha-1) \; (h_{jk,i} - h_{ik,j}) \: \gamma^i \;
	\xi^j \: \xi^k \\
\label{eq:a2_g18}
&&+ \frac{i}{16 \pi^2} \: m^2 \; (l^\vee(\xi)-l^\wedge(\xi)) \int_x^y
	\varepsilon^{ijlm} \; h_{jk,i} \; \xi^k \: \xi_l \; \rho \gamma_m \\
\label{eq:a2_g19}
&&- \frac{1}{16 \pi^2} \: m^2 \; (l^\vee(\xi)-l^\wedge(\xi)) \int_x^y
	(\alpha^2-\alpha) \; R_{jk} \; \xi^j \: \xi^k \; \xi\slsh \\
\label{eq:a2_g20}
&&- \frac{1}{16 \pi^3} \: m^2 \left( \lint_x^y - \lint_y^x \right)dz
	\; \zeta^j \int_x^z \alpha^2 \: R_{jk} \; \gamma^k \\
\label{eq:a2_g21}
&&+ \frac{1}{32 \pi^3} \: m^2 \left( \lint_x^y - \lint_y^x \right)dz
	\int_x^z (2\alpha^2-\alpha) \; R \; \xi\slsh \\
\label{eq:a2_g22}
&&- \frac{1}{64 \pi^3} \: m^2 \left( \lint_x^y - \lint_y^x \right)dz
	\; \zeta^j \: \zeta^k \int_x^z (\alpha^4-\alpha^3) \: (
	R_{\:,jk} - 2 \: \Box R_{jk}) \; \xi\slsh \\
\label{eq:a2_g23}
&&+ \frac{1}{32 \pi^3} \: m^2 \left( \lint_x^y - \lint_y^x \right)dz
	\; \zeta^j \int_x^z \alpha^2 \; (R_{jk,i} - R_{ik,j}) \:
	(2 \alpha \zeta^k - \xi^k) \; \gamma^i \\
\label{eq:a2_g24}
&&+ \frac{i}{32 \pi^3} \: m^2 \left( \lint_x^y - \lint_y^x \right)dz
	\; \zeta^k \int_x^z \alpha^2 \; \varepsilon^{ijlm} \; R_{jk,i} \;
	\xi_l \; \rho \gamma_m \\
&&+ {\cal{O}}(m^3) \nonumber \spc .
\end{eqnarray}
\end{Thm}
{\Beweis}
Unter Verwendung von Theorem \ref{a2_theorem2} berechnen wir zun\"achst
\[  -\frac{i}{e} \frac{\partial}{\partial y^k} k_m[\gamma^j h_j^k](x,y)
	\spc . \]
Der Eichterm \Ref{a2_c2} f\"uhrt dabei auf den Ausdruck
\[ - \left( \int_x^y h^k_j \right) \: \xi^j \: \frac{\partial}{\partial y^k}
	k_m(x,y) - \frac{1}{2} \: h(y) \: k_m(x,y) - \left( \int_x^y h
	\right) \: k_m(x,y) \spc . \]
F\"ur die restlichen Beitr\"age betrachten wir die Terme verschiedener
Ordnung in $m$ nacheinander:
\begin{enumerate}
\item Terme $\sim m$:
Die Summanden \Ref{a2_c2} und \Ref{a2_c3} ergeben die Beitr\"age
\begin{eqnarray*}
\Ref{a2_c3} &:& \Ref{a2_g12} + \Ref{a2_g16} \\
\Ref{a2_c4} &:& -\frac{1}{4} \: (h(y)+h(x)) \; k_m(x,y) + \frac{1}{2}
	\int_x^y h \; k_m(x,y)\\
&& + \Ref{a2_g13} + \Ref{a2_g14} + \Ref{a2_g15}
	\spc .
\end{eqnarray*}
\item Terme $\sim m^2$:
Die Summanden \Ref{a2_c5} bis \Ref{a2_c7} f\"uhren auf
\begin{eqnarray*}
\Ref{a2_c5} &:& -\frac{1}{8\pi^3} \: m^2 \left( \lint_x^y - \lint_y^x \right)
	dz \; \zeta^j \int_x^z \alpha^2 \; R_{jk} \: \gamma^k
	+ \Ref{a2_g18} + \Ref{a2_g23} \\
&&- \frac{1}{32 \pi^3} \: m^2 \left( \lint_x^y - \lint_y^x \right)
	\Pdd h \;+\; \frac{1}{16 \pi^3} \: m^2 \left( \lint_x^y -
	\lint_y^x \right) dz \int_x^z \alpha \; \Pdd h \\
\Ref{a2_c6} &:& \Ref{a2_g19} + \Ref{a2_g24} \\
\Ref{a2_c7} &:& -\frac{1}{4} \: (h(y) + h(x)) \; k_m(x,y) + \frac{1}{2}
	\int_x^y h \; k_m(x,y) \\
&&+ \frac{1}{32 \pi^3} \: m^2 \left( \lint_x^y - \lint_y^x \right)
	\Pdd h \;-\; \frac{1}{16 \pi^3} \: m^2 \left( \lint_x^y -
	\lint_y^x \right) dz \int_x^z \alpha \; \Pdd h \\
&& +\frac{1}{16 \pi^3} \: m^2 \left( \lint_x^y - \lint_y^x \right)
	dz \; \zeta^j \int_x^z \alpha^2 \; R_{jk} \: \gamma^k \\
&& + \Ref{a2_g21} + \Ref{a2_g22} + \Ref{a2_g19} \spc .
\end{eqnarray*}
\end{enumerate}
Setzen wir diese Gleichungen in \Ref{a2_g10} ein, folgt die Behauptung.
\QED

\section{Axiale St\"orung}
Wir betrachten jetzt die St\"orung durch ein axiales Potential
\Equ{a2_x0}
	G \;=\; i \Pdd + e \: \rho \Aslsh \spc .
\EndEqu
F\"ur $\Delta k_m$ hat man in St\"orungstheorie erster Ordnung
\[  \Delta k_m(x,y) \;=\; -e \left( s_m \: \rho \Aslsh \: p_m + p_m \:
	\rho \Aslsh
\: s_m \right)(x,y) \spc . \]
Diesen Beitrag erh\"alt man auch, wenn man in Gleichung \Ref{a2_a11} das
elektromagnetische Potential $\Aslsh$ durch $\rho \Aslsh$ ersetzt.

\begin{Thm}
\label{a2_thm13}
In erster Ordnung St\"orungstheorie gilt
\begin{eqnarray}
\label{eq:a2_x1}
\Delta k_m(x,y) &=& - \rho \: \Delta k_m[\Aslsh](x,y) \\
\label{eq:a2_x2}
&& + \frac{e}{4 \pi^2} \: m \: (l^\vee(\xi)-l^\wedge(\xi)) \: \int_x^y \rho
	\Aslsh \: \xi\slsh \\
\label{eq:a2_x3}
&& + \frac{e}{8 \pi^3} \: m \left( \lint_x^y - \lint_y^x \right) dz \;
	\int_x^y \alpha^2 \: j_k \: \zeta^k \; \rho \\
\label{eq:a2_x4}
&& + \frac{e}{8 \pi^3} \: m \left( \lint_x^y - \lint_y^x \right) \partial_k
	A^k \; \rho \\
\label{eq:a2_x5}
&& + \frac{ie}{4 \pi^3} \: m \left( \lint_x^y - \lint_y^x \right) h_j[A_k] \:
	\rho \sigma^{jk} \\
\label{eq:a2_x6}
&& + \frac{ie}{4 \pi^3} \: m^2 \left( \lint_x^y - \lint_y^x \right) \rho
	\Aslsh \\
\label{eq:a2_x7}
&& - \frac{e}{8 \pi^3} \: m^3 \left( \lint_x^y - \lint_y^x \right) \rho
	\Aslsh \: \xi\slsh \\
\label{eq:a2_x8}
&& - \frac{e}{16 \pi^3} \: m^3 \left( \wint_x^y - \wint_y^x \right)
	(\partial_j \Aslsh) \: \gamma^j \\
\label{eq:a2_x9}
&& - \frac{ie}{16 \pi^3} \: m^4 \left( \vint_x^y + \wint_x^y - \vint_y^x
	- \wint_y^x \right) \: \rho \Aslsh \\
&& + {\cal{O}}(m^5) \spc . \nonumber
\end{eqnarray}
\end{Thm}
{\Beweis}
Es gen\"ugt wieder, den Fall $\xi^0>0$ zu betrachten. Untersuche die Beitr\"age
verschiedener Ordnung in $m$ nacheinander:
\begin{enumerate}
\item Terme $\sim m^0$:
\begin{eqnarray*}
  \Delta k_m &=& -e \: (s_0 \: \rho \Aslsh \: p_0 + p_0 \: \rho \Aslsh
	\: s_0 ) \\
&=& \rho \: e \: (s_0 \: \Aslsh \: p_0 + p_0 \: \Aslsh \: s_0) \;=\; - \rho \:
	\Delta k_m[\Aslsh]
\end{eqnarray*}
\item Terme $\sim m$:
\begin{eqnarray}
\label{eq:a2_229z}
  \Delta k_m(x,y) &=& - \frac{e}{8 \pi^3} \: m \left\{ \Pdd_x \left( \lint_x^y  
	\rho \Aslsh \right) + \left( \lint_x^y \rho \Aslsh \right) \Pdd_y
	\right\} \\
&=& - \rho \: \Delta k_m[\Aslsh](x,y) + \frac{em}{4 \pi^3} \: \rho \:
	\frac{\partial}{\partial y^j} \lint_x^y \Aslsh \: \gamma^j \nonumber \\
\label{eq:a2_229a}
&=& - \rho \: \Delta k_m[\Aslsh](x,y) + \frac{em}{4 \pi^2} \: l^\vee(\xi)
	\int_x^y \rho \Aslsh \: \xi\slsh
	+ \frac{em}{4 \pi^3} \: \rho \lint_x^y h_j[\Aslsh] \: \gamma^j \spc
\end{eqnarray}
Setze nun noch die Umformungen
\begin{eqnarray}
\label{eq:a2_x10}
h_j[\Aslsh] \: \gamma^j &=& h_j[A^j] + i h_j[A_k] \: \sigma^{jk} \\
h_j[A^j](z) &=& \partial_j A^j(z) - \int_x^z \alpha \: \partial_j A^j
	- \frac{1}{2} \: \zeta_j \int_x^z \alpha^2 \: \Box A^j \nonumber \\
\label{eq:a2_x11}
&=& \frac{1}{2} \: \partial_k A^k(z) + \frac{1}{2} \: \zeta_k \int_x^z
	\alpha^2 \: j^k
\end{eqnarray}
ein.
\item Terme $\sim m^2$:
\begin{eqnarray*}
  \Delta k_m(x,y) &=& \frac{ie}{8 \pi^3} \: m^2 \left\{ \lint_x^y \rho \Aslsh
	+ \frac{1}{2} \Pdd_x \left(\lint_x^y dz \; \rho \Aslsh(z) \; (z-y)^j
	\gamma_j \right) \right. \\
&& \spc  \left. + \frac{1}{2} \left(\lint_x^y dz \; (x-z)^j \gamma_j \;
  \rho \Aslsh(z) \right) \Pdd_y \right\} \\
&=& - \rho \: \Delta k_m[\Aslsh](x,y) + \frac{ie}{4 \pi^3} \: m^2 \lint_x^y
	\rho \Aslsh
\end{eqnarray*}
\item Terme $\sim m^3$:
\begin{eqnarray*}
  \Delta k_m(x,y) &=&
-\frac{e}{8 \pi^3} \: m^3 \left\{ - \frac{1}{2} \lint_x^y dz \left( (x-z)^j
  \gamma_j \: \rho \Aslsh(z) + \rho \Aslsh(z) \: (z-y)^j \gamma_j \right) 
	\right. \nonumber \\
&& \spc \left. -\frac{1}{4} \Pdd_x \left( \vint_x^y \rho \Aslsh \right) -
	\frac{1}{4} \left(\wint_x^y \rho \Aslsh \right) \Pdd_y \;\;\; \right\}
	\\
&=& - \rho \: \Delta k_m[\Aslsh](x,y) \\
&& + \frac{e}{8 \pi^3} \: m^3 \: \rho
	\left( \lint_x^y dz \; \Aslsh(z) \: (z-y)^j \gamma_j
	+ \frac{1}{2} \left( \wint_x^y \Aslsh \right) \Pdd_y \right)
\end{eqnarray*}
Forme mit \Ref{a2_156a} weiter um
\begin{eqnarray*}
\lint_x^y dz \; \Aslsh(z) \: (z-y)^j \: \gamma_j + \frac{1}{2}
	\left( \wint_x^y \Aslsh \right) \Pdd_y
&=& - \lint_x^y \Aslsh \: \xi\slsh - \frac{1}{2} \wint_x^y (\partial_j
	\Aslsh) \: \gamma^j \spc .
\end{eqnarray*}
\item Terme $\sim m^4$:
\begin{eqnarray*}
\Delta k_m(x,y) &=&
  \frac{ie}{8 \pi^3} m^4 \left\{ -\frac{1}{4} \lint_x^y dz \; (x-z)^j \gamma_j
   \: \rho \Aslsh(z) \: (z-y)^k \gamma_k \nonumber \right. \\
&& \spc - \frac{1}{16} \Pdd_x \left( \vint_x^y dz \; \rho \Aslsh(z) \: (z-y)^j
	\gamma_j \right) \nonumber \\
&& \spc - \frac{1}{16} \left(\wint_x^y dz \; (x-z)^j \gamma_j \: \rho \Aslsh(z)
	\right) \Pdd_y \nonumber \\
&& \left. \spc - \frac{1}{4} \vint_x^y \rho \Aslsh - \frac{1}{4} \wint_x^y \rho
	\Aslsh \spc \right\} \\
&=& - \rho \: \Delta k_m[\Aslsh](x,y) - \frac{ie}{16 \pi^3} \: m^4 \left(
	\vint_x^y + \wint_x^y \right) \: \rho \Aslsh
\end{eqnarray*}
\end{enumerate}
\QED
Wir kommen nun zur Interpretation der abgeleiteten Formel: An die Stelle der
Eichterme bei elektromagnetischen St\"orungen treten nun die Beitr\"age
\begin{eqnarray}
\label{eq:a2_x19}
ie \int_x^y \rho \: A_j \xi^j \; k_m(x,y) && \spc {\mbox{gerade Ordnung in
	$m$}} \\
\label{eq:a2_x20}
ie \int_x^y \rho \: \frac{1}{2}  \left[ \xi\slsh, \Aslsh \right]
	\; k_m(x,y) && \spc {\mbox{ungerade Ordnung in $m$},} 
\end{eqnarray}
die man durch Addition der Eichterme von \Ref{a2_x1} mit \Ref{a2_x2} und
\Ref{a2_x7} erh\"alt. Wir nennen diese Terme {\bf{Pseudoeichterme}}. Sie
f\"uhren nicht nur zu einer Phasenverschiebung von $k_m$, was die Tatsache
widerspiegelt, da{\ss} axiale St\"orungen nicht lokal durch Eichtransformationen
kompensiert werden k\"onnen.

Die Summanden \Ref{a2_x3} bis \Ref{a2_x5} sind Korrekturen zu den
Stromtermen.
Die Beitr\"age \Ref{a2_x6}, \Ref{a2_x7}, \Ref{a2_x9} hei{\ss}en
{\bf{Massenterme}}; sie f\"uhren auf die Ruhemasse der $W-$ und $Z-$Eichbosonen.
\\[1em]
F\"ur die Randwerte auf dem Lichtkegel hat man
\begin{Satz}
F\"ur $y-x \in \Li$ gilt
\begin{eqnarray*}
\lefteqn{  \lim_{\I_x \ni u \rightarrow y} (\Delta k_m(x,u) - \Delta k_0(x,u))
	} \\
&=& \frac{ie}{32 \pi^2} \: m \: \epsilon(\xi^0) \int_x^y (2\alpha -1) \:
	F_{jk} \; \rho \sigma^{jk} \\
&& + \frac{e}{16 \pi^2} \: m \: \epsilon(\xi^0) \int_x^y \partial_k A^k \;
	\rho \\
&& - \frac{ie}{16 \pi^2} \: m \: \epsilon(\xi^0) \int_x^y (\alpha^2-\alpha) \;
	\Box A_j \: \xi_k \; \rho \sigma^{jk} \\
&& + \frac{ie}{8 \pi^2} \: m^2 \: \epsilon(\xi^0) \int_x^y \rho \Aslsh \\
&& - \frac{ie}{32 \pi^2} \: m^2 \: \epsilon(\xi^0) \int_x^y (2\alpha-1) \; 
	F_{ij} \: \xi^j \; \rho \gamma^i \\
&& + \frac{e}{64 \pi^2} \: m^2 \: \epsilon(\xi^0) \int_x^y \varepsilon^{ijkl}
	\; F_{ij} \: \xi_k \; \gamma_l \\
&& + \frac{ie}{32 \pi^2} \: m^2 \: \epsilon(\xi^0) \int_x^y (\alpha^2-\alpha)
	\; j_k \: \xi^k \; \rho \xi\slsh \\
&& + {\cal{O}}(m^3) \spc .
\end{eqnarray*}
\end{Satz}
{\Beweis}
Folgt direkt mit Hilfe von \Ref{a1_8}, \ref{integr_umf} und
Satz \ref{a2_satzb9}.
\QED

\section{Skalare St\"orung}
Wir betrachten wieder die skalare St\"orung~\Ref{a1_s1}. F\"ur
$\Delta k_m$ hat man in erster Ordnung in $\Xi$
\[      \Delta k_m \;=\; - \left( s_m \: \Xi \: k_m + k_m \: \Xi \: s_m
	\right) \spc . \]

\begin{Thm}
\label{a2_theorem_sm}
In erster Ordnung St\"orungstheorie gilt
\begin{eqnarray}
\Delta k_m(x,y) &=& \Delta k_0(x,y) \\
\label{eq:a21_s2}
&&- 2 m \left( \int_x^y \Xi \right) \; k^{(2)}(x,y) \\
&&+ \frac{1}{8 \pi^3} \: m \left(\lint_x^y - \lint_y^x \right) dz
	\; \left( (\Pdd \Xi)(z) - 2 \int_x^z \alpha \; (\Pdd \Xi) \right) \\
&&- \frac{1}{8 \pi^3} \: m \left(\lint_x^y - \lint_y^x \right) dz \;
	\zeta \slsh \int_x^z \alpha^2 \; (\Box \Xi) \\
&&+ \frac{3i}{8 \pi^3} \: m^2 \left(\lint_x^y - \lint_y^x \right) \Xi \\
\label{eq:a2_449}
&&+ \frac{i}{16 \pi^3} \: m^2 \left(\lint_x^y - \lint_y^x \right) dz\;
	(\partial_j \Xi) \: (2 \zeta^j - \xi^j) \\
&&- \frac{1}{16 \pi^3} \: m^2 \left(\lint_x^y - \lint_y^x \right) dz \;
	(\partial_j \Xi) \: \xi_k \; \sigma^{jk} \\
&&- \frac{1}{8 \pi^3} \: m^3 \left(\lint_x^y - \lint_y^x \right) \Xi \; \xi\slsh \\
&&+\frac{1}{32 \pi^3} \: m^3 \left(\vint_x^y - \wint_x^y + \vint_y^x -
	\wint_y^x \right) \Pdd \Xi \\
&&+ {\cal{O}}(m^4) \spc . \nonumber
\end{eqnarray}
\end{Thm}
{\Beweis}
Es gen\"ugt, den Fall $\xi^0>0$ zu betrachten. Untersuche die Terme
verschiedener Ordnung in $m$ nacheinander:
\begin{enumerate}
\item Terme $\sim m$:
\begin{eqnarray*}
\Delta k_m(x,y) &=& - \frac{1}{8 \pi^3} \: m \left\{ \Pdd_x (\lint_x^y \Xi)
	+ (\lint_x^y \Xi) \Pdd_y \right\} \\
&=& -\frac{1}{8 \pi^3} \: m \left\{ \lint_x^y (\Pdd \Xi) - 2 \Pdd_y
	\lint_x^y \Xi \right\} \\
&=& \frac{m}{4 \pi^2} \: l^\vee(\xi) \int_x^y \Xi \; \xi\slsh \\
&& + \frac{m}{8 \pi^3} \lint_x^y dz \: \left( (\Pdd \Xi)(z) - 2 \int_x^z
	\alpha \: (\Pdd \Xi) - \zeta \slsh \int_x^z \alpha^2 \: (\Box \Xi)
	\right)
\end{eqnarray*}
\item Terme $\sim m^2$:
\begin{eqnarray*}
\Delta k_m(x,y) &=& \frac{i}{8 \pi^3} \: m^2 \left\{ \lint_x^y \Xi +
	\frac{1}{2} \: \Pdd_x \left( \lint_x^y dz \; \Xi(z) \; (z-y)^j
	\gamma_j \right) \right. \\
&& \left. \spc\spc + \frac{1}{2} \left( \lint_x^y dz \; (x-z)^j \gamma_j \;
	\Xi(z) \right) \Pdd_y \right\}
\end{eqnarray*}
Unter Verwendung von Lemma~\ref{a2_lemma10} kann man den zweiten und dritten
Summanden umformen
\begin{eqnarray*}
\Pdd_x \left( \lint_x^y dz \; \Xi(z) \; (z-y)^j \gamma_j \right) &=&
	\lint_x^y dz \; (\Pdd \Xi)(z) \; (z-y)^j \gamma_j \;+\; 2
	\lint_x^y \Xi \\
\left( \lint_x^y dz \; (x-z)^j \gamma_j \; \Xi(z) \right) \Pdd_y &=&
	- \lint_x^y dz \; (x-z)^j \gamma_j \; (\Pdd \Xi)(z) \;+\; 2
	\lint_x^y \Xi
\end{eqnarray*}
und erh\"alt insgesamt
\begin{eqnarray*}
\lefteqn{ \Delta k_m(x,y) \;=\; \frac{3i}{8 \pi^3} \: m^2 \lint_x^y \Xi } \\
&&+ \frac{i}{16 \pi^3} \: m^2 \lint_x^y dz \left( (\Pdd \Xi)(z) \; (z-y)^j
	\gamma_j - (x-z)^j \gamma_j \; (\Pdd \Xi)(z) \right) \spc .
\end{eqnarray*}
\item Terme $\sim m^3$:
\begin{eqnarray*}
\Delta k_m(x,y) &=& -\frac{1}{8 \pi^3} \: m^3 \left\{
	\frac{1}{2} \: \xi\slsh \lint_x^y \Xi - \frac{1}{4} \: \Pdd_x
	\left( \vint_x^y \Xi \right) - \frac{1}{4} \left( \wint_x^y \Xi
	\right) \Pdd_y \right\}
\end{eqnarray*}
Setze nun die Umformungen
\begin{eqnarray}
\label{eq:a2_400}
\Pdd_x \left( \vint_x^y \Xi \right) &=& \vint_x^y \Pdd \Xi \:-\:
	2 \lint_x^y dz \; (y-z)^j \gamma_j \; \Xi(z) \\
\label{eq:a2_401}
\left( \wint_x^y \Xi \right) \Pdd_y &=& -\wint_x^y \Pdd \Xi \:-\:
	2 \lint_x^y dz \; (z-x)^j \gamma_j \; \Xi(z)
\end{eqnarray}
ein.
\end{enumerate}
\QED
Bei der Berechnung der Randwerte auf dem Lichtkegel kann man in~\Ref{a2_449}
partiell integrieren, so da{\ss} sich die Gleichung vereinfacht:
\begin{Satz}
F\"ur $(y-x) \in \Li$ gilt
\begin{eqnarray*}
\lefteqn{\lim_{\I_x \ni u \rightarrow y} \left( \Delta k_m(x,u) -
\Delta k_0(x,u) \right) } \\
&=& \frac{1}{16 \pi^2} \: m \: \epsilon(\xi^0) \int_x^y (2\alpha-1)
	\; (\Pdd \Xi) \\
&&+ \frac{1}{16 \pi^2} \: m \: \epsilon(\xi^0) \int_x^y (\alpha^2-\alpha)
	\; (\Box \Xi) \; \xi\slsh \\
&&+ \frac{i}{32 \pi^2} \: m^2 \: \epsilon(\xi^0) \; (\Xi(y) + \Xi(x)) \\
&&+ \frac{i}{8 \pi^2} \: m^2 \: \epsilon(\xi^0) \int_x^y \Xi \\
&&- \frac{1}{32 \pi^2} \: m^2 \: \epsilon(\xi^0) \; \int_x^y (\partial_j \Xi)
	\: \xi_k \; \sigma^{jk} \\
&&-\frac{1}{16 \pi^2} \: m^3 \: \epsilon(\xi^0) \; \int_x^y \Xi \; \xi\slsh \\
&&+ {\cal{O}}(m^4) \spc .
\end{eqnarray*}
\end{Satz}
{\Beweis} Die Behauptung folgt unmittelbar mit Hilfe von~\Ref{a1_8}
und~\ref{integr_umf}.
\QED

\section{Pseudoskalare St\"orung}
Wir betrachten die St\"orung durch ein pseudoskalares Potential
\Equ{a1_ps1}
	G \;=\; i \Pdd + i \rho \Xi \spc .
\EndEqu
F\"ur $\Delta k_m$ hat man in erster Ordnung in $\Xi$
\[      \Delta k_m \;=\; - i \: \left( s_m \: \rho \Xi \: k_m + k_m \: \rho \Xi
	\: s_m
	\right) \spc . \]
\begin{Thm}
\label{a2_theorem_psm}
In erster Ordnung St\"orungstheorie gilt
\begin{eqnarray}
\Delta k_m(x,y) &=& -i \rho \: \Delta k_0[\Xi](x,y) \\
&&+ \frac{i}{8 \pi^3} \: m \: \rho \left(\lint_x^y-\lint_y^x \right) \Pdd \Xi \\
&&+ \frac{1}{8 \pi^3} \: m^2 \left(\lint_x^y - \lint_y^x \right) \Xi \; \rho \\
\label{eq:a2_p449}
&&+ \frac{1}{16 \pi^3} \: m^2 \left(\lint_x^y - \lint_y^x \right) dz\;
	(\partial_j \Xi) \: (2 \zeta^j - \xi^j) \; \rho \\
&&+ \frac{i}{16 \pi^3} \: m^2 \left(\lint_x^y - \lint_y^x \right) dz \;
	(\partial_j \Xi) \: \xi_k \; \rho \sigma^{jk} \\
&&-\frac{i}{32 \pi^3} \: m^3 \left( \vint_x^y + \wint_x^y - \vint_y^x
	- \wint_y^x \right) \rho (\Pdd \Xi) \\
&&+ {\cal{O}}(m^4) \spc . \nonumber
\end{eqnarray}
\end{Thm}
{\Beweis}
Es gen\"ugt, den Fall $\xi^0>0$ zu betrachten. Untersuche die Terme
verschiedener Ordnung in $m$ nacheinander:
\begin{enumerate}
\item Terme $\sim m$:
\begin{eqnarray*}
\Delta k_m(x,y) &=& - \frac{i}{8 \pi^3} \: m \left\{ \Pdd_x (\lint_x^y
	\rho \Xi) + (\lint_x^y \rho \Xi) \Pdd_y \right\} \\
&=& \frac{i}{8 \pi^3} \: m \: \rho \left\{ \Pdd_x (\lint_x^y
	\Xi) + \Pdd_y (\lint_x^y \Xi) \right\} \\
&=& \frac{i}{8 \pi^3} \: m \: \rho \lint_x^y \Pdd \Xi
\end{eqnarray*}
\item Terme $\sim m^2$:
\begin{eqnarray*}
\Delta k_m(x,y) &=& -\frac{1}{8 \pi^3} \: m^2 \left\{ \lint_x^y \rho \Xi +
	\frac{1}{2} \: \Pdd_x \left( \lint_x^y dz \; \rho \Xi(z) \; (z-y)^j
	\gamma_j \right) \right. \\
&& \left. \spc\spc + \frac{1}{2} \left( \lint_x^y dz \; (x-z)^j \gamma_j \;
	\rho \Xi(z) \right) \Pdd_y \right\} \\
&=& - i \rho \: \Delta k_m[\Xi](x,y) - \frac{m^2}{4 \pi^3} \: \rho
	\lint_x^y \Xi
\end{eqnarray*}
\item Terme $\sim m^3$:
\begin{eqnarray*}
\Delta k_m(x,y) &=& -\frac{i}{8 \pi^3} \: m^3 \left\{ -\frac{1}{2}
	\lint_x^y dz \left( (x-z)^j \gamma_j \: \rho \Xi(z) +
	\rho \Xi(z) \: (z-y)^j \gamma_j \right) \right. \\
&&\hspace*{3cm} \left. -\frac{1}{4} \: \Pdd_x \left( \vint_x^y \rho
	\Xi \right) - \frac{1}{4} \left( \wint_x^y \rho \Xi \right) \Pdd_y
	\right\} \\
&=&-\frac{i \rho}{8 \pi^3} \: m^3 \left\{ \frac{1}{2} \lint_x^y dz \;
	(y-2z+x)^j \gamma_j \; \Xi(z) \right. \\
&&\hspace*{3cm} \left. +\frac{1}{4} \: \Pdd_x \left( \vint_x^y \Xi \right)
	-\frac{1}{4} \left( \wint_x^y \Xi \right) \Pdd_y \right\}
\end{eqnarray*}
Setze nun die Umformungen \Ref{a2_400}, \Ref{a2_401} ein.
\end{enumerate}
\QED

\section{Bilineare St\"orung}
Wir betrachten wieder die bilineare St\"orung~\Ref{a1_b00}.
F\"ur $\Delta k_m$ hat man in erster Ordnung in $B$
\Equ{a2_b0}
  \Delta k_m \;=\; - \left( s_m \: B_{jk} \: \sigma^{jk} \; k_m +
	k_m \: B_{jk} \: \sigma^{jk} \: s_m \right) \spc .
\EndEqu
\begin{Thm}
\label{a2_thm_b}
In erster Ordnung St\"orungstheorie gilt
\begin{eqnarray}
\lefteqn{ \Delta k_m(x,y) \;=\; \Delta k_0(x,y) } \\
&&+ \frac{1}{4 \pi^2} \: m \: (l^\vee(\xi)-l^\wedge(\xi)) \int_x^y
	\varepsilon^{ijkl} \; B_{ij} \: \xi_k \; \rho \gamma_l \\
&&+ \frac{i}{4 \pi^3} \: m \left(\lint_x^y - \lint_y^x \right)
	B_{jk,}^{\;\;\;k} \; \gamma^j \\
&&+ \frac{1}{8 \pi^3} \: m \left(\lint_x^y - \lint_y^x \right) dz \;
	\varepsilon^{ijkl} \: \left( B_{ij,k}(z) - 2 \int_x^z
	\alpha \: B_{ij,k} \right) \: \rho \gamma_l \\
&&- \frac{1}{8 \pi^3} \: m \left(\lint_x^y - \lint_y^x \right) dz \;
	\varepsilon^{ijkl} \; \zeta_i \int_x^z \alpha^2 \: (\Box B_{jk})
	\; \rho \gamma_l \\
&&+ \frac{i}{4 \pi^2} \: m^2 \: (l^\vee(\xi)-l^\wedge(\xi)) \int_x^y
	B_{ij} \: \xi^i \: \xi_k \; \sigma^{jk} \\
&&+ \frac{i}{32 \pi^2} \: m^2 \; \Theta(\xi^2) \: \epsilon(\xi^0) \;
	(B_{jk}(y) + B_{jk}(x)) \; \sigma^{jk} \\
&&- \frac{1}{16 \pi^2} \: m^2 \; \Theta(\xi^2) \: \epsilon(\xi^0) \int_x^y
	\xi_j \; B^{jk}_{\;\;\:,k} \\
&&+ \frac{i}{16 \pi^2} \: m^2 \; \Theta(\xi^2) \: \epsilon(\xi^0)
	\int_x^y (2\alpha-1) \; (\xi^k \: B_{jk,i} + \xi_i \:
	B_{jk,}^{\;\;\;k}) \; \sigma^{ij} \\
&&+ \frac{i}{16 \pi^2} \: m^2 \; \Theta(\xi^2) \: \epsilon(\xi^0)
	\int_x^y (\alpha^2-\alpha) \; (\Box B_{ij}) \: \xi^i \: \xi_k \;
	\sigma^{jk} \\
&&- \frac{i}{32 \pi^2} \: m^2 \; \Theta(\xi^2) \: \epsilon(\xi^0)
	\int_x^y \varepsilon^{ijkl} \; B_{ij,k} \: \xi_l \; \rho \\
&&- \frac{1}{16 \pi^2} \: m^3 \; \Theta(\xi^2) \: \epsilon(\xi^0)
	\int_x^y \varepsilon^{ijkl} \: B_{ij} \: \xi_k \;
	\rho \gamma_l \\
&&+ {\cal{O}}(\xi^2) + {\cal{O}}(m^4) \spc . \nonumber
\end{eqnarray}
\end{Thm}
{\Beweis}
Es gen\"ugt, den Fall $\xi^0>0$ zu betrachten. Untersuche die Beitr\"age
verschiedener Ordnung in $m$ nacheinander:
\begin{enumerate}
\item Terme $\sim m$:
\[ \Delta k_m(x,y) \;=\; -\frac{1}{8 \pi^3} \: m \: \left\{ \Pdd_x
	(\lint_x^y B_{jk} \: \sigma^{jk}) + (\lint_x^y B_{jk} \: \sigma^{jk})
	\Pdd_y \right\} \]
Wir setzen die Relationen
\begin{eqnarray}
\label{eq:a2_gam1}
\gamma^i \: \sigma^{jk} &=& i \: ( g^{ij} \: \gamma^k - g^{ik} \: \gamma^j)
	+ \varepsilon^{ijkl} \; \rho \gamma_l \\
\label{eq:a2_gam2}
\sigma^{jk} \: \gamma^i  &=& - i \: ( g^{ij} \: \gamma^k - g^{ik} \: \gamma^j)
	+ \varepsilon^{ijkl} \; \rho \gamma_l
\end{eqnarray}
ein und erhalten
\begin{eqnarray*}
&=& -\frac{im}{4 \pi^3} \: \frac{\partial}{\partial x^j} \lint_x^y B^{jk}
	\: \gamma_k - \frac{m}{8 \pi^3} \: \varepsilon^{ijkl} \:
	\frac{\partial}{\partial x^i} \lint_x^y B_{jk} \; \rho \gamma_l \\
&&-\frac{im}{4 \pi^3} \: \frac{\partial}{\partial y^j} \lint_x^y B^{jk}
	\: \gamma_k + \frac{m}{8 \pi^3} \: \varepsilon^{ijkl} \:
	\frac{\partial}{\partial y^i} \lint_x^y B_{jk} \; \rho \gamma_l \\
&=& - \frac{im}{4 \pi^3} \lint_x^y B_{jk,}^{\;\;\;j} \: \gamma^k 
	- \frac{m}{8 \pi^3} \: \varepsilon^{ijkl} \lint_x^y B_{ij,k} \;
	\rho \gamma_l + \frac{m}{4 \pi^3} \: \varepsilon^{ijkl} \:
	\frac{\partial}{\partial y^i} \lint_x^y B_{jk} \; \rho \gamma_l \\
&=& \frac{im}{4 \pi^3} \lint_x^y B_{jk,}^{\;\;\;k} \: \gamma^j \;+\;
\frac{m}{4 \pi^2} \: l^\vee(\xi) \int_x^y \varepsilon^{ijkl} \:
	B_{ij} \: \xi_k \; \rho \gamma_l \\
&&+ \frac{m}{8 \pi^3} \: \varepsilon^{ijkl} \lint_x^y dz \; \left(
	B_{ij,k}(z) - 2 \int_x^z \alpha \; B_{ij,k} \right) \: \rho \gamma_l
	\\
&&- \frac{m}{8 \pi^3} \: \varepsilon^{ijkl} \lint_x^y dz \; \zeta_i
	\int_x^z \alpha^2 \: (\Box B_{jk}) \; \rho \gamma_l \spc .
\end{eqnarray*}
\item Terme $\sim m^2$:
\begin{eqnarray}
\Delta k_m(x,y) &=& \frac{i}{8 \pi^3} \: m^2 \left\{ \lint_x^y
	B_{jk} \: \sigma^{jk} + \frac{1}{2} \: \Pdd_x \left( \lint_x^y dz
	\; B_{jk}(z) \: \sigma^{jk} \; (z-y)^m \gamma_m \right)
	\right. \nonumber \\
\label{eq:a2_bi0}
&&\left. \spc\spc + \frac{1}{2} \left( \lint_x^y dz \;
	(x-z)^m \gamma_m \; B_{jk}(z) \: \sigma^{jk} \right) \Pdd_y
	\right\}
\end{eqnarray}
Wir berechnen nun den zweiten und dritten Summanden. Dazu l\"osen wir
zun\"achst das Produkt der Dirac-Matrizen mit Hilfe von~\Ref{a1_b8}
auf. An\-schlie{\ss}end berechnen wir die Ableitungen und f\"uhren eine
Entwicklung um den Lichtkegel durch, wobei wir alle Terme der Ordnung
$\xi^2$ weglassen. Am einfachsten berechnet man die Randwerte der
Ableitungen der Lichtkegelintegrale direkt mit~\Ref{a1_41a}. 
Beachte, da{\ss} sich die Terme
\[  \frac{\partial}{\partial y^m} \lint_x^y B_{jk} \: \sigma^{jk} \; (z-y)^m
	\;\;\;, \spc \frac{\partial}{\partial x^m} \lint_x^y (x-z)^m \;
	B_{jk} \: \sigma^{jk}   \]
mit Hilfe von Lemma~\ref{a2_lemma10} auswerten lassen.
Man erh\"alt auf diese Weise
\begin{eqnarray*}
\lefteqn{ \Pdd_x \left( \lint_x^y B_{jk} \: \sigma^{jk} \; (z-y)^i \gamma_i
	\right)
\;=\; 2i \: \frac{\partial}{\partial x^j} \lint_x^y dz \; B^{jk}(z) \;
	(z-y)_k } \\
&&+ \frac{\partial}{\partial x^m} \lint_x^y dz \; B_{jk}(z) \: \sigma^{jk}
	\; (z-y)^m
\;+\; 2 \: \frac{\partial}{\partial x^j} \lint_x^y dz \; B^{jk}(z) \;
	\sigma_{km} \; (z-y)^m \\
&&+ 2 \: \frac{\partial}{\partial x^j} \lint_x^y dz \; B_{mk}(z) \;
	\sigma^{kj} \; (z-y)^m
\;+\; \frac{\partial}{\partial x^k} \lint_x^y dz \; \varepsilon^{ijkl} \:
	B_{ij}(z) \; (z-y)_l \; \rho \\
&=& - 4\pi \: l^\vee(\xi) \int_x^y (\alpha-1) \; B_{jk} \: \xi^j \: \xi_m
	\; \sigma^{km} \\
&&+ \frac{\pi}{2} \int_x^y (\alpha-1) \; \frac{d}{d\alpha} B_{jk} \:
	\sigma^{jk}
\;+\; \pi \int_x^y (2\alpha-1) \; B_{jk} \: \sigma^{jk} \\
&&-i \pi \int_x^y (\alpha-1) \; \xi_j \: B^{jk}_{\;\;\:,k}
\;+\; \frac{\pi}{2} \int_x^y (\alpha-1) \: \varepsilon^{ijkl} \:
	B_{ij,k} \: \xi_l \; \rho \\
&&- \pi \int_x^y (\alpha^3 - 2\alpha^2 + \alpha) \; (\Box B_{jk}) \;
	\xi^j \; \xi_m \; \sigma^{km} \\
&&- \pi \int_x^y (2\alpha^2 - 3\alpha +1) \: (B_{jk,}^{\;\;\:j} \; \xi_m
	+ B_{jk,m} \: \xi^j) \; \sigma^{km} \\
\lefteqn{ \left( \lint_x^y (x-z)^i \gamma_i \; B_{jk} \: \sigma^{jk}
	\right) \Pdd_y
\;=\; -2i \: \frac{\partial}{\partial y^k} \lint_x^y dz \; (x-z)_j
	\;B^{jk}(z)} \\
&&- \frac{\partial}{\partial y^m} \lint_x^y dz \; (x-z)^m \; B_{jk}(z) \:
	\sigma^{jk}
\;-\; 2 \: \frac{\partial}{\partial y^j} \lint_x^y dz \; (x-z)^m \; B^{jk}(z)
	\; \sigma_{km} \\
&&- 2 \: \frac{\partial}{\partial y^j} \lint_x^y dz \; (x-z)^m \; B_{mk}(z)
	\; \sigma^{kj}
\;-\; \frac{\partial}{\partial y^l} \lint_x^y dz \; \varepsilon^{ijkl} \:
	B_{ij}(z) \; (x-z)_k \; \rho \\
&=& 4\pi \: l^\vee(\xi) \int_x^y \alpha \; B_{jk} \: \xi^j \: \xi_m
	\; \sigma^{km} \\
&&+ \frac{\pi}{2} \int_x^y \alpha \; \frac{d}{d\alpha} B_{jk} \:
	\sigma^{jk}
\;-\; \pi \int_x^y (2\alpha-1) \; B_{jk} \: \sigma^{jk} \\
&&+i \pi \int_x^y \alpha \; \xi_j \: B^{jk}_{\;\;\:,k}
\;-\; \frac{\pi}{2} \int_x^y \alpha \: \varepsilon^{ijkl} \:
	B_{ij,k} \: \xi_l \; \rho \\
&&+ \pi \int_x^y (\alpha^3 - \alpha^2) \; (\Box B_{jk}) \;
	\xi^j \; \xi_m \; \sigma^{km} \\
&&+ \pi \int_x^y (2\alpha^2 - \alpha) \: (B_{jk,}^{\;\;\:j} \; \xi_m
	+ B_{jk,m} \: \xi^j) \; \sigma^{km}
\end{eqnarray*}
und durch Einsetzen in~\Ref{a2_b0} die Behauptung.
\item Terme $\sim m^3$:
\begin{eqnarray*}
\lefteqn{ \Delta k_m(x,y) \;=\; -\frac{1}{8 \pi^3} \: m^3 \left\{
	-\frac{1}{2} \lint_x^y dz \left( (x-z)^j \gamma_j \:
	B_{kl}(z) \: \sigma^{kl} \right. \right. } \\
&& \left. \left. + B_{kl}(z) \: \sigma^{kl} \:
	(z-y)^j \gamma_j \right) -\frac{1}{4} \: \Pdd_x \left( \vint_x^y
	B_{kl} \sigma^{kl} \right) - \frac{1}{4} \left( \wint_x^y
	B_{kl} \sigma^{kl} \right) \Pdd_y \right\}
\end{eqnarray*}
Mit den Relationen
\begin{eqnarray*}
\Pdd_x \left( \vint_x^y B_{kl} \sigma^{kl} \right) &=&
	2 \lint_x^y dz \; (z-y)^j \gamma_j \: B_{kl}(z) \: \sigma^{kl}
	\;+\; {\cal{O}}(\xi^2) \\
\left( \wint_x^y B_{kl} \sigma^{kl} \right) \Pdd_y &=&
	2 \lint_x^y dz \; B_{kl}(z) \: \sigma^{kl} \: (x-z)^j \gamma_j
	\;+\; {\cal{O}}(\xi^2)
\end{eqnarray*}
vereinfacht sich der Beitrag zu
\begin{eqnarray*}
\Delta k_m(x,y) &=& -\frac{1}{16 \pi^3} \: m^3 \lint_x^y \left( \xi\slsh
	\: B_{kl} \sigma^{kl} + B_{kl} \sigma^{kl} \: \xi\slsh \right)
	\;+\; {\cal{O}}(\xi^2) \\
&=&-\frac{1}{8 \pi^3} \: m^3 \lint_x^y \varepsilon^{ijkl} \: B_{ij}
	\: \xi_k \; \rho \gamma_l \;+\; {\cal{O}}(\xi^2) \\
&=&-\frac{1}{16 \pi^2} \: m^3 \: \Theta^\vee(\xi) \int_x^y
	\varepsilon^{ijkl} \: B_{ij} \: \xi_k \; \rho \gamma_l
	\;+\; {\cal{O}}(\xi^2) \spc .
\end{eqnarray*}
\end{enumerate}
\QED

\section{Differentialst\"orung durch Vektorpotential}
Wir betrachten wieder die St\"orung des Diracoperators~\Ref{a1_va}.
F\"ur $\Delta k_m$ hat man in  erster Ordnung in $L$
\begin{eqnarray}
\Delta k_m(x,y) &=& -i \left( s_m L^j \partial_j \: k_m + k_m L^j
	\partial_j s_m \right)(x,y)
	- \frac{i}{2} \left( s_m L^j_{\;,j} k_m + k_m L^j_{\;,j}
	s_m \right)(x,y) \nonumber \\
\label{eq:a2_v0}
&=& -i \frac{\partial}{\partial y^j} \: \Delta k_m[L^j](x,y) + \frac{i}{2}
	\: \Delta k_m[L^j_{\;,j}](x,y) \spc ,
\end{eqnarray}
so da{\ss} wir die Rechnung teilweise auf das Ergebnis f\"ur skalare St\"orungen,
Theorem~\ref{a2_theorem_sm}, zur\"uckf\"uhren k\"onnen.

\begin{Thm}
\label{theorem_vkm}
In erster Ordnung St\"orungstheorie gilt
\begin{eqnarray}
\Delta k_m(x,y) &=& \Delta k_0(x,y) \\
\label{eq:a2_v1}
&&- im \left(\int_x^y L_j \right) \: \xi^j \; k_0(x,y) \\
&&- \frac{i}{8 \pi^2} \: m \; (l^\vee(\xi)-l^\wedge(\xi)) \int_x^y (2\alpha-1)
	\; (L^j_{\;,j} \: \xi\slsh + (\Pdd L_j) \: \xi^j ) \\
\label{eq:a2_v3}
&&- \frac{i}{4 \pi^2} \: m \; (l^\vee(\xi)-l^\wedge(\xi)) \int_x^y L \slsh \\
&&- \frac{i}{8 \pi^2} \: m \; (l^\vee(\xi)-l^\wedge(\xi)) \int_x^y
	(\alpha^2-\alpha) \; (\Box L_j) \; \xi^j \: \xi\slsh \\
\label{eq:a2_v5}
&&+ \frac{1}{16 \pi^2} \: m^2 \; (l^\vee(\xi)-l^\wedge(\xi)) \; (L_j(y) +
	L_j(x)) \; \xi^j \\
\label{eq:a2_v6}
&&+ \frac{1}{4 \pi^2} \: m^2 \; (l^\vee(\xi)-l^\wedge(\xi))
	\int_x^y L_j \: \xi^j \\
&&+ \frac{i}{16 \pi^2} \: m^2 \; (l^\vee(\xi)-l^\wedge(\xi)) \int_x^y
	(\partial_j L_m) \; \xi^m \: \xi_k \; \sigma^{jk} \\
&&+ {\cal{O}}(\xi^0) + {\cal{O}}(m^3) \spc . \nonumber
\end{eqnarray}
\end{Thm}
{\Beweis}
Nach Theorem~\ref{a2_theorem_sm} hat man
\begin{eqnarray*}
-i \frac{\partial}{\partial y^j} \: \Delta k_m[L^j](x,y) &=&
	-i \frac{\partial}{\partial y^j} \: \Delta k_0[L^j](x,y) \\
&&- \frac{i}{2 \pi^2} \: m \: (m^\vee(\xi)-m^\wedge(\xi)) \; \xi\slsh \int_x^y
	L_j \: \xi^j \\
&&- \frac{i}{4 \pi^2} \: m \: (l^\vee(\xi)-l^\wedge(\xi)) \int_x^y (\alpha
	\: L^j_{\;,j} \; \xi\slsh + L \slsh ) \\
&&- \frac{i}{8 \pi^2} \: m \: (l^\vee(\xi)-l^\wedge(\xi)) \int_x^y (2\alpha-1)
	\; (\Pdd L_j) \: \xi^j \\
&&- \frac{i}{8 \pi^2} \: m \: (l^\vee(\xi)-l^\wedge(\xi)) \int_x^y (\alpha^2-
	\alpha) \; (\Box L_j) \: \xi^j \: \xi\slsh \\
&&+ \frac{1}{16 \pi^2} \: m^2 \: (l^\vee(\xi)-l^\wedge(\xi)) \; (L_j(y) +
	L_j(x)) \: \xi^j \\
&&+ \frac{1}{4 \pi^2} \: m^2 \; (l^\vee(\xi)-l^\wedge(\xi))
	\int_x^y L_j \: \xi^j \\
&&+ \frac{i}{16 \pi^2} \: m^2 \: (l^\vee(\xi)-l^\wedge(\xi)) \int_x^y
	(\partial_j L_m) \; \xi^m \: \xi_k \; \sigma^{jk} \\
&&+ {\cal{O}}(\xi^0) + {\cal{O}}(m^3) \\
\frac{i}{2} \: \Delta k_m[L^j_{\;.j}](x,y) &=& \frac{i}{2}
	\: \Delta k_0[L^j_{\;,j}](x,y)
\;+\; \frac{i}{8 \pi^2} \: m \: (l^\vee(\xi)-l^\wedge(\xi)) \int_x^y L^j_{\;,j}
	\; \xi\slsh \\
&& + {\cal{O}}(\xi^0) + {\cal{O}}(m^3) \spc .
\end{eqnarray*}
Setze nun in \Ref{a2_v0} ein.
\QED

Wir berechnen nun den Beitrag $\sim m^3$ bis zur Ordnung
${\cal{O}}(\xi^2)$:
\begin{Satz}
\label{zusatz_vkm}
In erster Ordnung St\"orungstheorie gilt
\begin{eqnarray}
\label{eq:a2_v8}
\lefteqn{ \Delta k^{(3)}(x,y) \;=\; \frac{i}{8 \pi^2} \:
	(l^\vee(\xi) - l^\wedge(\xi))
	\int_x^y L_j \xi^j \: \xi\slsh } \\
\label{eq:a2_v9}
&&+\frac{i}{16 \pi^2} \: (\Theta^\vee(\xi) - \Theta^\wedge(\xi))
	\int_x^y L\slsh \\
&&+\frac{i}{32 \pi^2} \: (\Theta^\vee(\xi) - \Theta^\wedge(\xi))
	\int_x^y (2\alpha-1) \: L^j_{\;,j} \: \xi\slsh \\
&&+\frac{i}{32 \pi^2} \: (\Theta^\vee(\xi) - \Theta^\wedge(\xi))
	\int_x^y (2\alpha-1) \: (\Pdd L_j) \; \xi^j \\
&&+\frac{i}{32 \pi^2} \: (\Theta^\vee(\xi) - \Theta^\wedge(\xi))
	\int_x^y (\alpha^2-\alpha) \: (\Box L_j) \: \xi^j \; \xi\slsh
	\;+\; {\cal{O}}(\xi^2) \;\;\; .
\end{eqnarray}
\end{Satz}
{\Beweis}
Nach Theorem \ref{a2_theorem_sm}, Lemma \ref{lemma_b3} und \Ref{a1_41a}
hat man f\"ur $\xi^0>0$
\begin{eqnarray*}
\lefteqn{ -i \frac{\partial}{\partial y^j} \Delta k^{(3)}[L^j](x,y)
\;=\; \frac{i}{8 \pi^3} \lint_x^y L\slsh \:+\:
	\frac{i}{8 \pi^3} \: \xi\slsh \: \frac{\partial}{\partial y^j}
	\lint L^j } \\
&&-\frac{i}{32 \pi^3} \left( 2 \lint_x^y dz \; (y-z)_j (\Pdd L^j)
	\:-\: 2 \lint_x^y dz \; (z-x)_j (\Pdd L^j) \right)
	\;+\; {\cal{O}}(\xi^2) \\
&=& \frac{i}{16 \pi^2} \: \Theta^\vee(\xi) \int_x^y L\slsh \:+\:
	\frac{i}{8 \pi^2} \: l^\vee(\xi) \int_x^y L_j \xi^j \:
	\xi\slsh \:+\:
	\frac{i}{16 \pi^2} \: \Theta^\vee(\xi) \int_x^y \alpha \:
	L^j_{\;,j} \: \xi\slsh \\
&&+\frac{i}{32 \pi^2} \: \Theta^\vee(\xi) \int_x^y (\alpha^2-\alpha)
	\: (\Box L_j) \xi^j \; \xi\slsh \\
&&+\: \frac{i}{32 \pi^2}
	\: \Theta^\vee(\xi) \int_x^y (2\alpha-1) \: (\Pdd L_j) \xi^j
	\;+\; {\cal{O}}(\xi^2) \\
\lefteqn{ \frac{i}{2} \: \Delta k^{(3)}[L^j_{\;,j}](x,y) \;=\;
	-\frac{i}{32 \pi^2} \: \Theta^\vee(\xi) \int_x^y L^j_{\;,j}
	\: \xi\slsh \;+\; {\cal{O}}(\xi^2) \spc . }
\end{eqnarray*}
\QED

\section{Differentialst\"orung durch axiales Potential}
Wir betrachten jetzt die St\"orung des Diracoperators
\Equ{a1_pva}
G \;=\; i \Pdd + \rho L^j \partial_j + \frac{1}{2} \: \rho L^j_{\;,j}
\EndEqu
mit einem reellen Vektorfeld $L$. 
F\"ur $\Delta k_m$ hat man in  erster Ordnung in $L$
\begin{eqnarray}
\label{eq:a2_pv0}
\Delta k_m(x,y)
&=& i \frac{\partial}{\partial y^j} \: \Delta k_m[i \rho L^j](x,y) - \frac{i}{2}
	\: \Delta k_m[i \rho L^j_{\;,j}] \spc ,
\end{eqnarray}
so da{\ss} wir die Rechnung teilweise auf das Ergebnis f\"ur pseudoskalare
St\"orungen, Theorem~\ref{a2_theorem_psm}, zur\"uckf\"uhren k\"onnen.

\begin{Thm}
\label{theorem_pvkm}
In erster Ordnung St\"orungstheorie gilt
\begin{eqnarray}
\Delta k_m(x,y) &=& i \rho \: \Delta k_0\left[i L^j \partial_j + \frac{i}{2}
	L^j_{\;,j}\right](x,y) \\
&&- \frac{1}{8 \pi^2} \: m \: \rho \: (l^\vee(\xi)-l^\wedge(\xi))
	\int_x^y (\Pdd L_j) \: \xi^j \\
\label{eq:a2_vr1}
&&+ \frac{i}{16 \pi^2} \: m^2 \: (l^\vee(\xi)-l^\wedge(\xi)) \; (L_j(y) +
	L_j(x)) \: \xi^j \; \rho \\
&&- \frac{1}{16 \pi^2} \: m^2 \: (l^\vee(\xi)-l^\wedge(\xi)) \int_x^y
	(\partial_j L_m) \; \xi^m \: \xi_k \; \rho \sigma^{jk} \\
&&+ {\cal{O}}(\xi^0) + {\cal{O}}(m^3) \spc \nonumber
\end{eqnarray}
\end{Thm}
{\Beweis}
Nach Theorem~\ref{a2_theorem_psm} hat man
\begin{eqnarray*}
i \frac{\partial}{\partial y^j} \: \Delta k_m[i \rho L^j](x,y) &=&
	\rho \frac{\partial}{\partial y^j} \: \Delta k_0[L^j](x,y) \\
&&- \frac{1}{8 \pi^2} \: m \: \rho \: (l^\vee(\xi)-l^\wedge(\xi))
	\int_x^y (\Pdd L_j) \: \xi^j \\
&&+ \frac{i}{16 \pi^2} \: m^2 \: (l^\vee(\xi)-l^\wedge(\xi)) \; (L_j(y) +
	L_j(x)) \: \xi^j \; \rho \\
&&- \frac{1}{16 \pi^2} \: m^2 \: (l^\vee(\xi)-l^\wedge(\xi)) \int_x^y
	(\partial_j L_m) \; \xi^m \: \xi_k \; \rho \sigma^{jk} \\
&&+ {\cal{O}}(\xi^0) + {\cal{O}}(m^3) \\
-\frac{i}{2} \: \Delta k_m[\rho L^j_{\;,j}](x,y) &=& -\frac{i}{2} \:
	\rho \: \Delta k_0[L^j_{\;,j}](x,y) \;+\; {\cal{O}}(\xi^0) \spc .
\end{eqnarray*}
\QED

Wir berechnen nun den Beitrag $\sim m^3$ bis zur Ordnung
${\cal{O}}(\xi^2)$:
\begin{Satz}
\label{zusatz_pvkm}
In erster Ordnung St\"orungstheorie gilt
\begin{eqnarray}
\Delta k^{(3)}(x,y) &=& \frac{1}{32 \pi^2} \: (\Theta^\vee(\xi)
	- \Theta^\wedge(\xi)) \int_x^y \rho (\Pdd L_j) \xi^j
	\;+\; {\cal{O}}(\xi^2) \spc .
\end{eqnarray}
\end{Satz}
{\Beweis}
Nach Theorem \ref{a2_theorem_psm} und \Ref{a1_41a}
hat man f\"ur $\xi^0>0$
\begin{eqnarray*}
\lefteqn{ i \frac{\partial}{\partial y^j}
	\Delta k^{(3)}[i \rho L^j](x,y) } \\
&=& \frac{\rho}{32 \pi^3} \left( 2 \lint_x^y dz \;
	(y-z)_j (\Pdd L^j) \:+\: 2 \lint_x^y dz \; (z-x)_j (\Pdd L^j)
	\right) \;+\; {\cal{O}}(\xi^2) \\
&=& \frac{\rho}{32 \pi^2} \: \Theta^\vee(\xi) \int_x^y (\Pdd L_j) \xi^j
	\;+\; {\cal{O}}(\xi^2) \\
\lefteqn{ -\frac{i}{2} \: \Delta k^{(3)}[i \rho L^j_{\;,j}](x,y) \;=\;
	{\cal{O}}(\xi^2) \spc . }
\end{eqnarray*}
\QED

\section{Bilineare Differentialst\"orung durch Vektorpotential}
Wir betrachten wieder die St\"orung des Diracoperators~\Ref{a1_bva}.
Mit Hilfe der Umformung
\begin{eqnarray*}
\Delta k_m[L\slsh (i \Pdd)] &=& - s_m \: L\slsh \: (i \Pdd) \: k_m \:-\:
	k_m \: L\slsh \: (i \Pdd) \: s_m \\
&=& m \: \Delta k_m[L\slsh] - k_m \: L\slsh
\end{eqnarray*}
und \Ref{a1_345a}
k\"onnen wir Gleichung \Ref{a1_345a} auf den Fall $m \neq 0$ \"ubertragen:
\begin{eqnarray}
\Delta k_m(x,y) &=& -i \Delta k_m \left[ i L^j \partial_j
	\right] (x,y) + i \Delta k_m \left[ L\slsh
	(i \Pdd) \right](x,y) + \frac{i}{2} \Delta k_m[L_{j,k} \sigma^{jk}](x,y)
	\nonumber \\
&=& -i \Delta k_m \left[ i L^j \partial_j + \frac{i}{2} \: L^j_{\;,j} \right](x,y)
	\:+\: \frac{i}{2} \: \Delta k_m[L_{j,k} \sigma^{jk}](x,y) \nonumber \\
\label{eq:a2_u}
&& + i m \: \Delta k_m[L\slsh](x,y) \:-\: i k_m(x,y) \: L\slsh(y) \:-\:
	\frac{1}{2} \: \Delta k_m[L^j_{\;,j}] \spc .
\end{eqnarray}
Damit l\"a{\ss}t sich die Rechnung teilweise auf diejenige f\"ur skalare, bilineare
und elektromagnetische St\"orungen sowie Differentialst\"orungen zur\"uckf\"uhren.

\begin{Thm}
\label{theorem_bvkm}
In erster Ordnung St\"orungstheorie gilt
\begin{eqnarray}
\Delta k_m(x,y) &=& \Delta k_0(x,y) \\
&&+\frac{1}{8 \pi^2} \:m\: (l^\vee(\xi) - l^\wedge(\xi)) \; (L\slsh(y) - L\slsh(x)) \\
\label{eq:a2_sv1}
&&-\frac{1}{16 \pi^2} \: m^2\: (l^\vee(\xi) - l^\wedge(\xi)) \; (L_j(y) + L_j(x))
	\: \xi_k \: \sigma^{jk} \\
&&-\frac{1}{8 \pi^2} \: m^2\: (l^\vee(\xi) - l^\wedge(\xi)) \int_x^y L^j_{\;,j}
	\: \xi\slsh \\
&&+\frac{i}{16 \pi^2} \: m^2\: (l^\vee(\xi) - l^\wedge(\xi)) \; (L_j(y) - L_j(x))
	\: \xi^j \\
&&+ {\cal{O}}(\xi^0) + {\cal{O}}(m^3) \spc . \nonumber
\end{eqnarray}
\end{Thm}
{\Beweis}
Wir erhalten aus Theorem~\ref{a2_thm_b} und Theorem~\ref{a1_theorem1}
\begin{eqnarray*}
\lefteqn{ \frac{i}{2} \: \Delta k_m[L_{i,j} \sigma^{ij}](x,y) \;=\;
	\frac{i}{2} \: \Delta k_0[L_{i,j} \sigma^{ij}](x,y) } \\
&&+\frac{i}{8 \pi^2} \:m\: (l^\vee(\xi) - l^\wedge(\xi)) \int_x^y
	\varepsilon^{ijkl} \: L_{i,j} \: \xi_k \; \rho \gamma_l \\
&&-\frac{1}{16 \pi^2} \: m^2\: (l^\vee(\xi) - l^\wedge(\xi)) \int_x^y L_{i,j} \: \xi^i
	\: \xi_k \: \sigma^{jk} \\
&&+\frac{1}{16 \pi^2} \: m^2\: (l^\vee(\xi) - l^\wedge(\xi)) \; (L_j(y) - L_j(x))
	\: \xi_k \: \sigma^{jk} \;+\; {\cal{O}}(\xi^0) + {\cal{O}}(m^3) \\
\lefteqn{ i m \: \Delta k_m[L\slsh](x,y) \;=\; \frac{1}{2 \pi^2} \:m \:
	(m^\vee(\xi) - m^\wedge(\xi)) \int_x^y L_j \xi^j \: \xi\slsh } \\
&&+\frac{1}{8 \pi^2} \:m\: (l^\vee(\xi) - l^\wedge(\xi)) \int_x^y (2\alpha-1) \:
	L^i_{\;,i} \; \xi\slsh \\
&&+\frac{1}{8 \pi^2} \: m\: (l^\vee(\xi) - l^\wedge(\xi))  \int_x^y (\alpha^2-\alpha)
	\: (\Box L_k) \: \xi^k \; \xi\slsh \\
&&+\frac{1}{8 \pi^2} \:m\: (l^\vee(\xi) - l^\wedge(\xi)) \int_x^y (2\alpha-1)
	\: (\Pdd L_j) \: \xi^j \\
&&-\frac{1}{8 \pi^2} \:m\: (l^\vee(\xi) - l^\wedge(\xi)) \; (L\slsh(y) + L\slsh(x)) \\
&&+\frac{1}{4 \pi^2} \:m\: (l^\vee(\xi) - l^\wedge(\xi)) \int_x^y L\slsh \\
&&-\frac{i}{8 \pi^2} \:m\: (l^\vee(\xi) - l^\wedge(\xi)) \int_x^y \varepsilon^{ijkl}
	\: L_{i,j} \: \xi_k \; \rho \gamma_l \\
&&+\frac{i}{4 \pi^2} \: m^2\: (l^\vee(\xi) - l^\wedge(\xi)) \int_x^y L_j \xi^j
	\;+\; {\cal{O}}(\xi^0) + {\cal{O}}(m^3) \\
\lefteqn{-i k_m(x,y) \: L\slsh(y) \;=\; -i k_0(x,y) \: L\slsh(y) } \\
&&+\frac{1}{4 \pi^2} \:m\: (l^\vee(\xi) - l^\wedge(\xi)) \; L\slsh(y) \\
&&+\frac{i}{8 \pi^2} \: m^2\: (l^\vee(\xi) - l^\wedge(\xi))  \; L_j(y) \: \xi^j \\
&&+\frac{1}{8 \pi^2} \: m^2\: (l^\vee(\xi) - l^\wedge(\xi)) \; \xi_j \: L_k(y) \;
	\sigma^{jk} \;+\; {\cal{O}}(m^3) \spc .
\end{eqnarray*}
Durch Vergleich mit Theorem~\ref{theorem_vkm} folgt
\begin{eqnarray*}
\lefteqn{\hspace*{-.5cm} -i \Delta (k_m-k_0) \left[ i L^j \partial_j + \frac{i}{2}
	\: L^j_{\;,j}
	\right](x,y) \:+\: \frac{i}{2} \: \Delta (k_m-k_0)[L_{i,j} \sigma^{ij}](x,y)
	} \\
\lefteqn{\hspace*{2cm}
	+\: i m \: \Delta k_m[L\slsh](x,y) \:-\: i (k_m-k_0)(x,y) \: L\slsh(y) } \\
&=&\frac{1}{8 \pi^2} \:m\: (l^\vee(\xi) - l^\wedge(\xi)) \; (L\slsh(y) - L\slsh(x)) \\
&&-\frac{1}{16 \pi^2} \: m^2\: (l^\vee(\xi) - l^\wedge(\xi)) \; (L_j(y) + L_j(x))
	\: \xi_k \: \sigma^{jk} \\
&&+\frac{i}{16 \pi^2} \: m^2\: (l^\vee(\xi) - l^\wedge(\xi)) \; (L_j(y) - L_j(x))
	\: \xi^j \;+\; {\cal{O}}(\xi^0) + {\cal{O}}(m^3) \spc .
\end{eqnarray*}
Man erh\"alt die Behauptung mit \Ref{a2_u} und Theorem~\ref{a2_theorem_sm}.
\QED

Wir berechnen nun den Beitrag $\sim m^3$ bis zur Ordnung
${\cal{O}}(\xi^2)$:
\begin{Satz}
\label{zusatz_bv}
In erster Ordnung St\"orungstheorie gilt
\begin{eqnarray}
\Delta k^{(3)}(x,y) &=& -\frac{1}{32 \pi^2} \: (\Theta^\vee(\xi) -
	\Theta^\wedge(\xi)) \; (L\slsh(y) - L\slsh(x)) \\
&&+\frac{1}{32 \pi^2} \: (\Theta^\vee(\xi) - \Theta^\wedge(\xi))
	\int_x^y L^j_{\;,j} \: \xi\slsh \;+\; {\cal{O}}(\xi^2) \spc .
\end{eqnarray}
\end{Satz}
{\Beweis}
Nach Theorem \ref{a2_thm_b} und Theorem \ref{a2_theorem2}
hat man f\"ur $\xi^0>0$
\begin{eqnarray*}
\frac{i}{2} \: \Delta k^{(3)}[L_{j,k} \sigma^{jk}](x,y)
	&=& -\frac{i}{32 \pi^2}
	\: (\Theta^\vee(\xi) - \Theta^\wedge(\xi)) \int_x^y
	\varepsilon^{ijkl} \: L_{i,j} \: \xi_k \; \rho \gamma_l
	\;+\: {\cal{O}}(\xi^2) \\
i \Delta k^{(2)}[L\slsh](x,y) &=& -\frac{1}{8 \pi^2}
	\: (l^\vee(\xi) - l^\wedge(\xi)) \int_x^y L_j \xi^j \:
	\xi\slsh \\
&&-\frac{1}{32 \pi^2} \: (\Theta^\vee(\xi) - \Theta^\wedge(\xi))
	\int_x^y (2\alpha-1) \: (\Pdd L_j) \xi^j \\
&&+\frac{1}{32 \pi^2} \: (\Theta^\vee(\xi) - \Theta^\wedge(\xi))
	\; (L\slsh(y) + L\slsh(x)) \\
&&-\frac{1}{16 \pi^2} \: (\Theta^\vee(\xi) - \Theta^\wedge(\xi))
	\int_x^y L\slsh \\
&&+\frac{i}{32 \pi^2} \: (\Theta^\vee(\xi) - \Theta^\wedge(\xi))
	\int_x^y \varepsilon^{ijkl} \: L_{i,j} \: \xi_k \;
	\rho \gamma_l \\
&&-\frac{1}{32 \pi^2} \: (\Theta^\vee(\xi) - \Theta^\wedge(\xi))
	\int_x^y (2\alpha-1) \: L^j_{\;,j} \: \xi\slsh \\
&&-\frac{1}{32 \pi^2} \: (\Theta^\vee(\xi) - \Theta^\wedge(\xi))
	\int_x^y (\alpha^2-\alpha) \: (\Box L_j) \xi^j \: \xi\slsh
	\;+\; {\cal{O}}(\xi^2) \\
-i k^{(3)}(x,y) \; L\slsh(y) &=& -\frac{1}{16 \pi^2}
	\: (\Theta^\vee(\xi) - \Theta^\wedge(\xi)) \; L\slsh(y)
	\;+\; {\cal{O}}(\xi^2)
\end{eqnarray*}
Durch Vergleich mit Satz~\ref{zusatz_vkm} folgt
\begin{eqnarray*}
\lefteqn{\hspace*{-.5cm} -i \Delta k^{(3)}
	\left[ i L^j \partial_j + \frac{i}{2} \: L^j_{\;,j} \right](x,y)
	\:+\: \frac{i}{2} \: \Delta k^{(3)}[L_{i,j} \sigma^{ij}](x,y) } \\
\lefteqn{\hspace*{2cm}
	+\: i \Delta k^{(2)}[L\slsh](x,y) \:-\: i k^{(3)}(x,y) \:
	L\slsh(y) } \\
&=& -\frac{1}{32 \pi^2} \: (\Theta^\vee(\xi) - \Theta^\wedge(\xi))
	\; (L\slsh(y) - L\slsh(x)) \;+\; {\cal{O}}(\xi^2) \spc .
\end{eqnarray*}
Man erh\"alt die Behauptung mit \Ref{a2_u} und Theorem~\ref{a2_theorem_sm}.
\QED

\section{Bilineare Differentialst\"orung durch axiales Potential}
Wir betrachten jetzt die St\"orung des Diracoperators
\Equ{a2_ba1}
G \;=\; i \Pdd \:+\: \rho L_j \: \sigma^{jk} \: \frac{\partial}{\partial x^k}
	\:+\: \frac{\rho}{2} \: L_{j,k} \sigma^{jk}
\EndEqu
mit einem reellen Vektorfeld $L$. Mit Hilfe der Umformung
\begin{eqnarray*}
\Delta k_m[\rho L\slsh (i \Pdd)] &=& - s_m \: \rho L\slsh \: (i \Pdd) \: k_m \:-\:
	k_m \: \rho L\slsh \: (i \Pdd) \: s_m \\
&=& m \: \Delta k_m[\rho L\slsh] - k_m \: \rho L\slsh
\end{eqnarray*}
und \Ref{a1_345a}
k\"onnen wir die Rechnung zum Teil auf diejenige f\"ur pseudoskalare, bilineare
und axiale St\"orungen sowie auf Differentialst\"orungen zur\"uckf\"uhren:
\begin{eqnarray}
\Delta k_m(x,y) &=& -i \Delta k_m \left[ \rho L^j \partial_j
	\right] (x,y) + \Delta k_m \left[ \rho L\slsh
	(i \Pdd) \right](x,y) + \frac{1}{2} \Delta k_m[\rho L_{j,k} \sigma^{jk}](x,y)
	\nonumber \\
&=& -i \Delta k_m \left[ \rho L^j \partial_j + \frac{\rho}{2} \: L^j_{\;,j} \right](x,y)
	\:+\: \frac{1}{2} \: \Delta k_m[\rho L_{j,k} \sigma^{jk}](x,y) \nonumber \\
\label{eq:a2_v}
&& - k_m(x,y) \: \rho L\slsh(y) \:+\: m \: \Delta k_m[\rho L\slsh](x,y) \:+\:
	\frac{1}{2} \: \Delta k_m[i \rho L^j_{\;,j}] \;\;\; .
\end{eqnarray}

\begin{Thm}
\label{theorem_bakm}
In erster Ordnung St\"orungstheorie gilt
\begin{eqnarray}
\Delta k_m(x,y) &=& i \rho \: \Delta k_0 \left[ i L_j \sigma^{jk} \partial_k
	+ \frac{i}{2} \: L_{j,k} \sigma^{jk} \right](x,y) \\
\label{eq:a2_bvr1}
&&+\frac{i}{2 \pi^2} \:m \:
	(m^\vee(\xi) - m^\wedge(\xi)) \int_x^y L_j \xi^j \; \rho \xi\slsh \\
\label{eq:a2_bvr2}
&&-\frac{i}{4 \pi^2} \:m\: (l^\vee(\xi) - l^\wedge(\xi)) \; \rho
	(L\slsh(y) + L\slsh(x)) \\
\label{eq:a2_bvr3}
&&+\frac{i}{4 \pi^2} \:m\: (l^\vee(\xi) - l^\wedge(\xi)) \int_x^y \rho L\slsh \\
&&+\frac{i}{8 \pi^2} \:m\: (l^\vee(\xi) - l^\wedge(\xi)) \int_x^y (2\alpha-1) \:
	L^i_{\;,i} \; \rho \xi\slsh \\
&&+\frac{i}{8 \pi^2} \: m\: (l^\vee(\xi) - l^\wedge(\xi))  \int_x^y (\alpha^2-\alpha)
	\: (\Box L_k) \: \xi^k \; \rho \xi\slsh \\
&&+\frac{i}{8 \pi^2} \:m\: (l^\vee(\xi) - l^\wedge(\xi)) \int_x^y (2\alpha-1)
	\: \rho (\Pdd L_j) \: \xi^j \\
&&+\frac{1}{8 \pi^2} \:m\: (l^\vee(\xi) - l^\wedge(\xi)) \int_x^y \varepsilon^{ijkl}
	\: L_{i,j} \: \xi_k \; \gamma_l \\
\label{eq:a2_bvr4}
&&-\frac{i}{4 \pi^2} \: m^2 \: (l^\vee(\xi) - l^\wedge(\xi)) \int_x^y L_j \: \xi_k
	\; \rho \sigma^{jk} \\
\label{eq:a2_bvr5}
&&-\frac{i}{16 \pi^2} \: m^2 \: (l^\vee(\xi) - l^\wedge(\xi)) \; (L_j(y) + L_j(x))
	\: \xi_k \; \rho \sigma^{jk} \\
&&-\frac{1}{16 \pi^2} \: m^2 \: (l^\vee(\xi) - l^\wedge(\xi)) \; (L_j(y) - L_j(x))
	\xi^j \; \rho \\
&&+ {\cal{O}}(\xi^0) + {\cal{O}}(m^3) \spc . \nonumber
\end{eqnarray}
\end{Thm}
{\Beweis}
Mit den Relationen
\begin{eqnarray*}
\rho \: \sigma^{ij} &=& \frac{i}{2} \: \varepsilon^{ijkl} \: \sigma_{kl} \\
\varepsilon^{ijkl} \: L_{i,j} \: \varepsilon_{klmn} &=& -2 (L_{m,n} - L_{n,m}) \\
\varepsilon^{ijkl} \: L_{i,j} \: \xi_k \: \xi^m \: \varepsilon_{lmno} \:
	\sigma^{no} &=& 2 (L_{i,j} - L_{j,i}) \: \xi^i \: \xi_k \: \sigma^{jk}
	\:+\: (L_{i,j}-L_{j,i}) \: \sigma^{ij} \: \xi^2
\end{eqnarray*}
erh\"alt man aus Theorem~\ref{a2_thm_b}
\begin{eqnarray*}
\lefteqn{ \frac{1}{2} \: \Delta (k_m-k_0)[\rho L_{i,j} \sigma^{ij}](x,y)
	\;=\; \frac{i}{4} \: \Delta (k_m-k_0) \left[ \varepsilon^{ijkl} \:
	L_{i,j} \: \sigma_{kl} \right](x,y) } \\
&=&\frac{i}{16 \pi^2} \:m\: (l^\vee(\xi) - l^\wedge(\xi)) \int_x^y
	\varepsilon^{ijkl} \: L_{i,j} \: \varepsilon_{klmn} \: \xi^m \;
	\rho \gamma^n \\
&&-\frac{1}{16 \pi^2} \: m^2 \: (l^\vee(\xi) - l^\wedge(\xi)) \int_x^y
	\varepsilon^{ijkl} \: L_{i,j} \: \xi_k \: \xi^m \: \sigma_{lm}
	\;+\; {\cal{O}}(\xi^0) + {\cal{O}}(m^3) \\
&=&-\frac{i}{8 \pi^2} \:m\: (l^\vee(\xi) - l^\wedge(\xi)) \int_x^y
	(L_{i,j}-L_{j,i}) \xi^i \: \rho \gamma^j \\
&&-\frac{i}{32 \pi^2} \: m^2 \: (l^\vee(\xi) - l^\wedge(\xi)) \int_x^y
	\varepsilon^{ijkl} \: L_{i,j} \: \xi_k \: \xi^m \: \varepsilon_{lmno}
	\: \rho \sigma^{no}
	\;+\; {\cal{O}}(\xi^0) + {\cal{O}}(m^3) \\
&=&-\frac{i}{8 \pi^2} \:m\: (l^\vee(\xi) - l^\wedge(\xi)) \int_x^y
	(L_{i,j}-L_{j,i}) \xi^i \: \rho \gamma^j \\
&&-\frac{i}{16 \pi^2} \: m^2 \: (l^\vee(\xi) - l^\wedge(\xi)) \int_x^y
	(L_{i,j}-L_{j,i}) \: \xi^i \: \xi_k \: \rho \sigma^{jk}
	\;+\; {\cal{O}}(\xi^0) + {\cal{O}}(m^3) \\
&=&-\frac{i}{8 \pi^2} \:m\: (l^\vee(\xi) - l^\wedge(\xi)) \int_x^y
	\rho \: (\Pdd L_i) \xi^i \: \rho \gamma^j \\
&&+\frac{i}{8 \pi^2} \:m\: (l^\vee(\xi) - l^\wedge(\xi)) \; \rho \:
	(L\slsh(y)-L\slsh(x)) \\
&&-\frac{i}{16 \pi^2} \: m^2 \: (l^\vee(\xi) - l^\wedge(\xi)) \int_x^y
	L_{i,j} \: \xi^i \: \xi_k \: \rho \sigma^{jk} \\
&&+\frac{i}{16 \pi^2} \: m^2 \: (l^\vee(\xi) - l^\wedge(\xi)) \;
	(L_j(y) - L_j(x)) \: \xi_k \; \rho \sigma^{jk}
	\;+\; {\cal{O}}(\xi^0) + {\cal{O}}(m^3) \;\;\; .
\end{eqnarray*}
Zus\"atzlich haben wir
\begin{eqnarray*}
- (k_m - k_0)(x,y) \: \rho L\slsh(y) &=& -\frac{i}{4 \pi^2}
	\: m\: (l^\vee(\xi) - l^\wedge(\xi)) \; \rho L\slsh(y) \\
&&-\frac{1}{8 \pi^2} \: m^2 \: (l^\vee(\xi) - l^\wedge(\xi)) \; L_j(y) \xi^j \; \rho \\
&&+\frac{i}{8 \pi^2} \: m^2 \: (l^\vee(\xi) - l^\wedge(\xi)) \; \xi_j \:
	L_k(y) \; \rho \sigma^{jk} \;+\; {\cal{O}}(m^3) \;\;\;.
\end{eqnarray*}
Durch Vergleich mit Theorem~\ref{theorem_pvkm} folgt
\begin{eqnarray*}
\lefteqn{ -i \Delta (k_m-k_0) \left[ \rho L^j \partial_j + \frac{\rho}{2} \:
	L^j_{\;,j} \right]
	(x,y) \:+\: \frac{1}{2} \: \Delta (k_m-k_0)[\rho L_{j,k} \sigma^{jk}](x,y) } \\
&&\hspace*{2cm} -\: (k_m-k_0)(x,y) \: \rho L\slsh(y) \\
&=&-\frac{i}{8 \pi^2} \:m\: (l^\vee(\xi) - l^\wedge(\xi)) \; \rho
	(L\slsh(y) + L\slsh(x)) \\
&&-\frac{1}{16 \pi^2} \: m^2 \: (l^\vee(\xi) - l^\wedge(\xi)) \; (L_j(y) - L_j(x))
	\xi^j \; \rho \\
&&-\frac{i}{16 \pi^2} \: m^2 \: (l^\vee(\xi) - l^\wedge(\xi)) \; (L_j(y) + L_j(x))
	\: \xi_k \; \rho \sigma^{jk} \;+\; {\cal{O}}(\xi^0) + {\cal{O}}(m^3) \;\;\; .
\end{eqnarray*}
Au{\ss}erdem haben wir nach Theorem~\ref{a1_theorem1} und Theorem \ref{a2_theorem2}
\begin{eqnarray*}
\lefteqn{ m \: \Delta k_m[\rho L\slsh](x,y) \;=\; \frac{i}{2 \pi^2} \:m \:
	(m^\vee(\xi) - m^\wedge(\xi)) \int_x^y L_j \xi^j \; \rho \xi\slsh } \\
&&+\frac{i}{8 \pi^2} \:m\: (l^\vee(\xi) - l^\wedge(\xi)) \int_x^y (2\alpha-1) \:
	L^i_{\;,i} \; \rho \xi\slsh \\
&&+\frac{i}{8 \pi^2} \: m\: (l^\vee(\xi) - l^\wedge(\xi))  \int_x^y (\alpha^2-\alpha)
	\: (\Box L_k) \: \xi^k \; \rho \xi\slsh \\
&&+\frac{i}{8 \pi^2} \:m\: (l^\vee(\xi) - l^\wedge(\xi)) \int_x^y
	(2\alpha-1) \: \rho (\Pdd L_j) \: \xi^j \\
&&-\frac{i}{8 \pi^2} \:m\: (l^\vee(\xi) - l^\wedge(\xi)) \; \rho
	(L\slsh(y) + L\slsh(x)) \\
&&+\frac{i}{4 \pi^2} \:m\: (l^\vee(\xi) - l^\wedge(\xi)) \int_x^y \rho L\slsh \\
&&+\frac{1}{8 \pi^2} \:m\: (l^\vee(\xi) - l^\wedge(\xi)) \int_x^y \varepsilon^{ijkl}
	\: L_{i,j} \: \xi_k \; \gamma_l \\
&&-\frac{i}{4 \pi^2} \: m^2 \: (l^\vee(\xi) - l^\wedge(\xi)) \int_x^y L_j \: \xi_k
	\; \rho \sigma^{jk} \;+\; {\cal{O}}(\xi^0) + {\cal{O}}(m^3) \spc .
\end{eqnarray*}
Die Behauptung folgt nun mit \Ref{a2_v} und Theorem~\ref{a2_theorem_psm}.
\QED

Wir berechnen nun den Beitrag $\sim m^3$ bis zur Ordnung
${\cal{O}}(\xi^2)$:
\begin{Satz}
In erster Ordnung St\"orungstheorie gilt
\begin{eqnarray}
\label{eq:a2_bvr6}
\lefteqn{ \Delta k^{(3)}(x,y) \;=\; -\frac{i}{8 \pi^2} \: (l^\vee(\xi) - l^\wedge(\xi))
	\int_x^y L_j \xi^j \; \rho \xi\slsh } \\
\label{eq:a2_bvr7}
&&+\frac{i}{16 \pi^2} \: (\Theta^\vee(\xi) - \Theta^\wedge(\xi))
	\; \rho (L\slsh(y) + L\slsh(x)) \\
\label{eq:a2_bvr8}
&&+\frac{i}{16 \pi^2} \: (\Theta^\vee(\xi) - \Theta^\wedge(\xi))
	\int_x^y \rho L\slsh \\
&&-\frac{i}{32 \pi^2} \: (\Theta^\vee(\xi) - \Theta^\wedge(\xi))
	\int_x^y (2\alpha-1) \: \rho (\Pdd L_j) \xi^j \\
&&-\frac{i}{32 \pi^2} \: (\Theta^\vee(\xi) - \Theta^\wedge(\xi))
	\int_x^y (2\alpha-1) \: L^j_{\;,j} \: \rho \xi\slsh \\
&&-\frac{i}{32 \pi^2} \: (\Theta^\vee(\xi) - \Theta^\wedge(\xi))
	\int_x^y (\alpha^2-\alpha) \: (\Box L_j) \xi^j \:
	\rho \xi\slsh \\
&&-\frac{1}{32 \pi^2} \: (\Theta^\vee(\xi) - \Theta^\wedge(\xi))
	\int_x^y \varepsilon^{ijkl} \: L_{i,j} \: \xi_k \; \gamma_l
	\:+\: {\cal{O}}(\xi^2) \spc .
\end{eqnarray}
\end{Satz}
{\Beweis}
Nach Theorem \ref{a2_thm_b} und Theorem \ref{a2_thm13}
hat man f\"ur $\xi^0>0$
\begin{eqnarray*}
\lefteqn{ \frac{1}{2} \: \Delta k^{(3)}[\rho L_{i,j} \sigma^{ij}](x,y)
	\;=\; \frac{i}{4} \: \Delta k^{(3)}[ \varepsilon^{ijkl} \:
	L_{i,j} \: \sigma_{kl}] } \\
&=&-\frac{i}{64 \pi^2} \: (\Theta^\vee(\xi) - \Theta^\wedge(\xi))
	\int_x^y \varepsilon^{ijkl} \: L_{i,j} \: \varepsilon_{klmn}
	\: \xi^m \; \rho \gamma^n \;+\; {\cal{O}}(\xi^2) \\
&=& \frac{i}{32 \pi^2} \: (\Theta^\vee(\xi) - \Theta^\wedge(\xi))
	\int_x^y (L_{i,j} - L_{j,i}) \: \xi^i \; \rho \gamma^j
	\;+\; {\cal{O}}(\xi^2) \\
&=& \frac{i}{32 \pi^2} \: (\Theta^\vee(\xi) - \Theta^\wedge(\xi))
	\int_x^y \rho (\Pdd L_j) \xi^j \\
&&-\frac{i}{32 \pi^2} \: (\Theta^\vee(\xi) - \Theta^\wedge(\xi))
	\; \rho (L\slsh(y)-L\slsh(x)) \;+\; {\cal{O}}(\xi^2) \\
\lefteqn{ -k^{(3)}(x,y) \; \rho L\slsh(y) \;=\; \frac{i}{16 \pi^2}
	\: (\Theta^\vee(\xi) - \Theta^\wedge(\xi)) \;
	\rho L\slsh(y) \;+\; {\cal{O}}(\xi^2) \spc . }
\end{eqnarray*}
Durch Vergleich mit Satz~\ref{zusatz_pvkm} folgt
\begin{eqnarray*}
\lefteqn{ -i \Delta k^{(3)} \left[ \rho L^j \partial_j + \frac{\rho}{2} \:
	L^j_{\;,j} \right](x,y) \:+\: \frac{1}{2} \: \Delta k^{(3)}[\rho L_{j,k} \sigma^{jk}](x,y) } \\
&&\hspace*{2cm} -\: k^{(3)}(x,y) \: \rho L\slsh(y) \\
&=&\frac{i}{32 \pi^2} \: (\Theta^\vee(\xi) - \Theta^\wedge(\xi))
	\; \rho (L\slsh(y) + L\slsh(x)) \;+\; {\cal{O}}(\xi^2) \spc .
\end{eqnarray*}
Au{\ss}erdem haben wir nach Theorem~\ref{a2_thm13} und Theorem
\ref{a2_theorem_psm}
\begin{eqnarray*}
\lefteqn{ \Delta k^{(2)}[\rho L\slsh](x,y) \;=\;
-\frac{i}{8 \pi^2} \: (l^\vee(\xi) - l^\wedge(\xi)) \int_x^y L_j \xi^j
	\; \rho \xi\slsh } \\
&&+\frac{i}{8 \pi^2} \: (\Theta^\vee(\xi) - \Theta^\wedge(\xi))
	\int_x^y \rho L\slsh \\
&&-\frac{i}{32 \pi^2} \: (\Theta^\vee(\xi) - \Theta^\wedge(\xi))
	\int_x^y (2\alpha-1) \: \rho (\Pdd L_j) \xi^j \\
&&+\frac{i}{32 \pi^2} \: (\Theta^\vee(\xi) - \Theta^\wedge(\xi))
	\; \rho (L\slsh(y) + L\slsh(x)) \\
&&-\frac{i}{16 \pi^2} \: (\Theta^\vee(\xi) - \Theta^\wedge(\xi))
	\int_x^y \rho L\slsh \\
&&-\frac{1}{32 \pi^2} \: (\Theta^\vee(\xi) - \Theta^\wedge(\xi))
	\int_x^y \varepsilon^{ijkl} \: L_{i,j} \: \xi_k \; \gamma_l \\
&&-\frac{i}{32 \pi^2} \: (\Theta^\vee(\xi) - \Theta^\wedge(\xi))
	\int_x^y (2\alpha-1) \: L^j_{\;,j} \: \rho \xi\slsh \\
&&-\frac{i}{32 \pi^2} \: (\Theta^\vee(\xi) - \Theta^\wedge(\xi))
	\int_x^y (\alpha^2-\alpha) \: (\Box L_j) \xi^j \:
	\rho \xi\slsh \:+\: {\cal{O}}(\xi^2) \\
\lefteqn{ \frac{1}{2} \: \Delta k^{(3)}[i \rho L^j_{\;,j}](x,y) \;=\;
	{\cal{O}}(\xi^2) \spc . }
\end{eqnarray*}
\QED

\chapter{St\"orungsrechnung f\"ur $p_0$ im Ortsraum}
\label{anh3}
In diesem Kapitel werden wir die Distribution $\tilde{p}_0$ f\"ur
verschiedene St\"orungen des Di\-rac\-ope\-ra\-tors im Ortsraum in erster
Ordnung berechnen.

Nach \Ref{2_28} m\"ussen wir die Gleichung
\[ \tilde{p}_0 \;=\; p_0 - p_0 {\cal{B}} s_0 - s_0 {\cal{B}} p_0
	\]
auswerten. Wir haben
\[  (p_0 \: {\cal{B}} \: s_0 + s_0 \: {\cal{B}} \: p_0)(x,y) \;=\; (i \Pdd_x) \; \left(
	P_0 \: {\cal{B}} \: S_0 + S_0 \: {\cal{B}} \: P_0 \right)(x,y) \;
	(i \Pdd_y)  , \]
dabei ist $S_0$ die Distribution~\Ref{a1_76}, $P_0$ ist im Ortsraum durch
\begin{eqnarray}
\label{eq:a3_170a}
P_0(x) &=& \int \frac{d^4k}{(2 \pi)^4} \; \delta(k^2) \: e^{-i k x}
    \;=\;  - \frac{1}{4 \pi^3} \: \frac{1}{x^2}
\end{eqnarray}
gegeben.
Der Pol auf dem Lichtkegel ist als Hauptwert zu behandeln, also f\"ur
$g \in C^\infty_c(M)$
\begin{eqnarray}
P_0 (g) &=& - \frac{1}{4 \pi^3} \hint d^4x \; \frac{1}{x^2} \: g(x) \nonumber \\
\label{eq:haupt1}
  &:=& - \frac{1}{4 \pi^3} \; \lim_{0 \neq \varepsilon \rightarrow 0}
	\frac{1}{2} \int d^4x \left( \frac{1}{x^2 + i \varepsilon}
	+ \frac{1}{x^2 - i \varepsilon} \right) \: g(x) \spc .
\end{eqnarray}
F\"ur die St\"orungsrechnung mu{\ss} man Ausdr\"ucke
der Form
\[      (P_0 f S_0 + S_0 f P_0) (x,y) \]
und deren partielle Ableitungen nach $x$, $y$ berechnen. Als technisches
Hilfsmittel werden dazu
allgemeinere Lichtkegelintegrale ben\"otigt:

\section{Verallgemeinerte Lichtkegelintegrale}
Zun\"achst m\"ussen wir Funktionen eingef\"uhren, die ein bestimmtes
Abfallverhalten im Unendlichen zeigen.
\begin{Def}
\label{a3_defa}
\begin{description}
\item[(i)]
  Eine Funktion $f \in C^\infty(M)$ hei{\ss}t {\bf zur $p$-ten Potenz
abfallend}, $p>0$, falls es zu jedem Punkt $0 \neq y \in M$ eine Umgebung $U
\subset M$ und Konstanten $\lambda_0, c$ gibt mit
\[  \left| \lambda^p \: f(\lambda z) \right| \;\leq\; c \spc 
	\forall \lambda > \lambda_0, \: z \in U	\spc . \]
Die zur ersten und zweiten Potenz abfallenden Funktionen werden auch
linear und quadratisch abfallend genannt.
\item[(ii)]
  Eine Funktion $f \in C^\infty(M)$ hei{\ss}t {\bf gleichm\"a{\ss}ig zur $p$-ten
Potenz abfallend}, $p>0$, falls die partiellen Ableitungen der Ordnung
$q$ von $f$, $q \geq 0$, zur $(p+q)$-ten Potenz abfallen.

Die Menge der gleichm\"a{\ss}ig zur $p$-ten Potenz abfallenden Funktionen wird mit
$C^\infty_p(M)$ bezeichnet.
\end{description}
\end{Def}
Offensichtlich folgt aus $f \in C^\infty_p$ und $g \in C^\infty_q$, da{\ss}
$f g \in C^\infty_{p+q}$.
\\[1em]
F\"ur $f \in C^\infty_2(M)$ soll nun das Integral
\Equ{a3_z}
  (\xint f)(y) \;:=\; - \int d^4z \: (l^\vee_y(z) - l^\wedge_y(z))
	\: (l^\vee(z) - l^\wedge(z)) \; f(z)
\EndEqu
als Distribution definiert werden.
Dazu wird folgendes kleine Lemma ben\"otigt:
\begin{Lemma}
\label{a3_lemma1}
\begin{description}
\item[(i)] F\"ur $g \in C^\infty_c(M)$ ist
\[  h(z) \;=\; \int d^4y \: (l^\vee(z-y) - l^\wedge(z-y)) \; g(y) \]
eine linear abfallende Funktion.
\item[(ii)] F\"allt $f$ zur dritten Potenz ab, so sind die Integrale
\[ \int d^4z \: l^\vee(z) \; f(z)  \spc,\spc
	 \int d^4z \: l^\wedge(z) \; f(z)  \]
wohldefiniert und endlich.
\end{description}
\end{Lemma}
{\Beweis}
\begin{description}
\item[(i)] Sei $0 \neq y \in M$. W\"ahle ein Bezugssystem mit $\vec{y} \neq 0$,
setze $\kappa = |\vec{y}|/2$. Da $g$ kompakten Tr\"ager besitzt, gibt es
$R>0$ mit $\mbox{supp } g \subset \R \times B_R(\vec{0})$.
Au{\ss}erdem ist $g$ beschr\"ankt, $|g|\leq c$.

W\"ahle $U = \R \times B_\kappa(\vec{y})$. F\"ur $z \in U$, $\lambda > 2 R /
\kappa$ hat man
\begin{eqnarray*}
\lambda \: |h(\lambda z)| &\leq& \lambda \int d^4y \: (l^\vee(\lambda z - y) +
	l^\wedge(\lambda z - y)) \; |g(y)| \\
&=& \int d^3\vec{x} \; \frac{\lambda}{2 |\vec{x}|} \left[ |g(\: (-|\vec{x}|,
	\vec{x}) +
	\lambda z)| + |g( \: (|\vec{x}|, \vec{x}) + \lambda z)| \right] \\
&\leq& \int d^3\vec{x} \; \frac{\lambda c}{|\vec{x}|} \: \chi_{B_R(\lambda
	\vec{z})}(\vec{x}) \spc .
\end{eqnarray*}
F\"ur $\vec{x} \in B_R(\lambda \vec{z})$ gilt $|\vec{x}| \geq |\lambda
\vec{z}| - R \geq \lambda \kappa - R \geq \frac{1}{2} \lambda \kappa$ und
somit
\[  \lambda \: |h(\lambda z)| \;\leq\; 2 \: \frac{c}{\kappa} \int d^3\vec{x} \;
	\chi_{B_R(\lambda \vec{z})}(\vec{x}) \;=\; \frac{8}{3} \: \pi
	\; \frac{c}{\kappa} \: R^3	\spc .		\]
\item[(ii)] In Polarkoordinaten $(t,r,\omega)$ hat man
\begin{eqnarray}
\label{eq:a3_d}
\int d^4z \: l^\vee(z) \; f(z)
&=& \frac{1}{2} \int_0^\infty r \: dr \int_{S^2} d\omega \; f(r,r,\omega)
	\spc .
\end{eqnarray}
Nach Voraussetzung f\"allt $f$ zur dritten Potenz ab, also
\[  \left| f(r,r,\omega) \right| \; \leq \; \frac{c}{r^3} \spc ,\]
und das Integral~\Ref{a3_d} ist endlich.
\end{description}
{\QED}
\begin{Def}
  Definiere f\"ur $f \in C^\infty_2(M)$ die Distribution $\sxint f$ durch
\begin{eqnarray}
\label{eq:a3_e}
  (\xint f)(g) &=& - \int d^4z \: (l^\vee(z) - l^\wedge(z)) \; h(z) \spc
	{\mbox{mit}} \\
  h(z) &=& f(z) \: \int d^4y \: (l^\vee(z-y) - l^\wedge(z-y)) \; g(y) \spc
	, \; g \in C^\infty_c(M). \nonumber
\end{eqnarray}
\end{Def}

Nach Lemma~\ref{a3_lemma1}, (i) ist $h \in C^\infty_3(M)$, somit
ist~\Ref{a3_e} nach Lemma~\ref{a3_lemma1}, (ii) sinnvoll definiert.

Die Bezeichnungen von Definition~\ref{def_a7} werden in analoger Weise
auch f\"ur $\sxint f$ verwendet, also
\begin{eqnarray*}
\xint_y f &=& \xint^y f \;=\; (\xint f)(y) \\
\xint_x^y f &=& \xint_y^x f \\
	&=& - \int d^4z \; (l^\vee(z-y) - l^\wedge(z-y))
	\: (l^\vee(z-x) - l^\wedge(z-x)) \; f(z) \spc .
\end{eqnarray*}
F\"ur alle noch folgenden Lichtkegelintegrale werden wir ebenfalls die
allgemeinere Notation von Definition \ref{def_a7}
verwenden, ohne darauf jedesmal im Einzelnen hinzuweisen.

F\"ur $y \in \I^\vee$ stimmt $\sxint f$ offensichtlich mit $\slint f$
\"uberein, f\"ur $y \in \Ra$ reduziert sich~\Ref{a3_z} auf das Integral von
$f$ \"uber ein zweidimensionales Hyperboloid.

Viele der Resultate, die in Anhang~\ref{anh1} f\"ur das Lichtkegelintegral
$\slint$ bewiesen wurden, lassen sich auf $\sxint f$ \"ubertragen.
Insbesondere haben wir in Analogie zu Satz \ref{a1_lemma1} das folgende Ergebnis:
\begin{Satz}
\label{a3_satz1}
$\sxint f$ ist harmonisch (im Distributionssinne).

Au{\ss}erdem ist $\sxint f$ eine Funktion, die auf $\I \cup \Ra$ glatt
ist. F\"ur $z \in \Li$ gilt
\begin{eqnarray}
\label{eq:a3_8}
\lim_{\I \ni y \rightarrow z} (\xint f)(y) &=& \frac{\pi}{2} \int_0^1
	f(\lambda z) \: d\lambda \\
\label{eq:a3_9}
\lim_{\Ra \ni y \rightarrow z} (\xint f)(y) &=& -\frac{\pi}{2} \int_{\sR
\setminus [0,1]} f(\lambda z) \: d\lambda
\end{eqnarray}
\end{Satz}
{\Beweis} Verl\"auft ganz \"ahnlich wie der Beweis von Satz~\ref{a1_lemma1}.
\QED

\begin{Satz}
\label{a3_satz2}
F\"ur die Ableitung von $\sxint f$ gilt (im Distributionssinne)
\Equ{a3_11}
  \partial_j (\xint f)(y) \;=\; (\xint h_j^{(y)})(y) + \pi \: y_j \: l(y)
	\; \int_{-\infty}^\infty f(\lambda y) \: d\lambda \spc .
\EndEqu
F\"ur $h_j^{(y)}$ kann man eine der Funktionen
\begin{eqnarray}
\label{eq:a3_13}
h_j^{(y)} &=& \frac{1}{2} \partial_j f(z) - \frac{z_j \: y^l - y_j
	\: z^l}{y^2} \; \partial_l f \\
\label{eq:a3_15}
h_j^{(y)}(z) &=& \left\{ \displaystyle \begin{array}{ll} \displaystyle
	\partial_j f(z) - \frac{1}{2} \Box_z (z_j \int_0^1 f(\alpha z) \:
	d\alpha) & \mbox{falls $y \in \I$}  \\[1em]
	\displaystyle \partial_j f(z) + \frac{1}{2} \Box_z (z_j \int_1^\infty
	f(\alpha z) \:
	d\alpha) & \mbox{falls $y \in \Ra$}
    \end{array}	\right.
\end{eqnarray}
setzen.

Falls $f$ sogar in $C^\infty_3(M)$, kann man f\"ur $h_j^{(y)}$ auch
\begin{eqnarray}
\label{eq:a3_12}
h_j^{(y)}(z) &=& \partial_j f(z) + 2 \frac{z_j \:(z-y)^l}{y^2} \: \partial_l f
	- 2 \frac{(y-2z)_j}{y^2} \: f  \spc {\mbox{oder}} \\
\label{eq:a3_14}
h_j^{(y)}(z) &=& 2 \frac{(y-2z)_j}{y^2} \: f(z) + 2 \frac{(y-z)_j \: z^l}{y^2}
	\: \partial_l f(z)
\end{eqnarray}
setzen.
\end{Satz}
{\Beweis}
Man kann ganz \"ahnlich wie in Anhang~\ref{anh1} vorgehen.
Daher soll der Beweis nur skizziert werden:

Zun\"achst betrachtet man den Fall, da{\ss} $f$ kompakten Tr\"ager besitzt.

Lemma~\ref{a1_lemma5} \"ubertr\"agt sich direkt auf das Lichtkegelintegral
$\sxint f$ und f\"uhrt auf~\Ref{a3_13}, \Ref{a3_12} und~\Ref{a3_14}.
Die Lemmata~\ref{a1_lemma6}, \ref{a1_lemma8} und \ref{a1_lemma9} gelten
w\"ortlich auch f\"ur $\sxint f$. Zum Beweis w\"ahlt man
an Stelle von Polarkoordinaten jeweils hyperbolische Koordinaten.

Bei Satz~\ref{a1_parabl} mu{\ss} man etwas aufpassen. Ersetzt man das
Lichtkegelintegral $\slint f$ durch $\sxint f$ und betrachtet den
Fall $y \in \Ra$, so verschwindet in~\Ref{a1_42aa} bei partieller
Integration in $\alpha$ der Randwert bei $\alpha=0$ nicht.
(Gleichung~\Ref{a1_44b} gilt nicht mehr, weil $\sxint_x^{\tilde{y}}$ f\"ur
$\Ra \ni \tilde{y} \rightarrow x$ divergiert.)
Ersetzt man in~\Ref{a1_55} jedoch den Term\ $\int_0^1 \alpha \;\: \Box_z
(z_j f(z))_{|z = \alpha y} \; d \alpha$ durch
$-\int_1^\infty \alpha \;\: \Box_z (z_j f(z))_{|z = \alpha y} \; d \alpha$, so
erh\"alt man bei der Umformung von~\Ref{a1_42aa} nur den gew\"unschten
Randterm bei $\alpha=1$. Auf diese Weise kann man~\Ref{a3_15} herleiten.

Den Beitrag $\sim l(y)$ in~\Ref{a3_11} erh\"alt man aus~\Ref{a3_8} und
\Ref{a3_9} genau wie im Beweis von Satz~\ref{a1_dis_abl}.

Die Voraussetzung $f \in C^\infty_c(M)$ kann man durch Approximation
abschw\"achen.
Dabei verwendet man, da{\ss} f\"ur $f \in C^\infty_2(M)$ auch $h_j^{(y)}$
gem\"a{\ss} \Ref{a3_13}, \Ref{a3_15} gleichm\"a{\ss}ig quadratisch abf\"allt.
Damit $h_j^{(y)}$ nach \Ref{a3_12}, \Ref{a3_14} gleichm\"a{\ss}ig quadratisch abf\"allt,
braucht man die st\"arkere Voraussetzung $f \in C^\infty_3(M)$.
\QED
Beachte, da{\ss} die Funktion $h_j^{(y)}$ in \Ref{a3_15} f\"ur $y \in \Ra$
einen Pol am Ursprung besitzt. Dies bereitet
f\"ur die Definition von $(\slint h_j^{(y)})(y)$ jedoch keine Probleme, weil das
Hyperboloid, \"uber das in~\Ref{a3_z} f\"ur festes
$y \in \I \cup \Ra$ integriert werden mu{\ss}, den
Ursprung nicht schneidet. Daher kann man $h_j^{(y)}$ in einer kleinen
Umgebung des Ursprungs beliebig ab\"andern und den Pol beseitigen.
\\[1em]
Wir werden nun weitere Lichtkegelintegrale einf\"uhren:
Sei $f \in C^\infty_2(M)$, $g \in C^\infty_c(M)$, setze f\"ur $z \neq 0$
\Equ{a3_20}
  h(z) \;=\; \int d^4y \; \delta(\bra z,y \ket) \: g(y) \spc .
\EndEqu
Das Integral ist endlich und h\"angt stetig von $z$ ab, wegen
\begin{eqnarray*}
  h(\lambda z) &=& \int d^4y \; \delta(\bra \lambda z,y \ket) \: g(y) \\
&=& \frac{1}{\lambda} \int d^4y \; \delta(\bra z,y \ket) \: g(y) \spc
	\mbox{, $\lambda > 0$}
\end{eqnarray*}
ist $h$ eine eine linear abfallende Funktion.
Nach Lemma~\ref{a3_lemma1}, (ii) ist daher das Integral
\[  \int d^4z \; l(z) \; f(z) \: h(z)  \]
wohldefiniert und endlich. Man kann also setzen:
\begin{Def}
Definiere f\"ur $f \in C^\infty_2(M)$ die Distribution $\sdint f$ durch
\[  (\dint f)(g) \;=\; \frac{1}{2} \int d^4z \; l(z) \: f(z) \: h(z)  \]
mit $h$ gem\"a{\ss}~\Ref{a3_20}.
\end{Def}
Formal ist $\sdint f$ durch das Integral
\Equ{a3_29}
  (\dint f)(y) \;=\; \frac{1}{2} \int d^4z \; l(z) \: \delta(\bra z,y \ket)
	\; f(z)
\EndEqu
gegeben. Offensichtlich hat man
\Equ{a3_181a}
   (\dint f)(\lambda y) \;=\; \frac{1}{|\lambda|} \: (\dint f)(y) \spc .
\EndEqu
Es gilt ferner
\begin{Satz}
\label{a3_satz5}
$\sdint f$ ist eine Funktion, die auf $\I$ identisch verschwindet und auf
$\Ra$ glatt ist. F\"ur $z \in \Li$ hat man
\[ \lim_{\Ra \ni y \rightarrow z} (\dint f)(y) \;=\; \frac{\pi}{2}
\int_{-\infty}^\infty f(\lambda z) \: d\lambda	\]
\end{Satz}
{\Beweis}
Verl\"auft \"ahnlich wie der Beweis von Satz~\ref{a1_lemma1}.
\QED
Das Lichtkegelintegral $\sdint f$ kann man auch als Grenzfall des
Lichtkegelintegrals $\sxint f$ darstellen:
\begin{Satz}
\label{a3_satz6}
Es gilt im Distributionssinne
\Equ{a3_27}
    (\dint f)(y) \;=\; - \lim_{0 \neq \lambda \rightarrow 0} |\lambda|
	(\xint f)(\lambda y) \spc .
\EndEqu
\end{Satz}
{\Beweis}
Man hat f\"ur $\lambda \neq 0$
\begin{eqnarray*}
(\xint f)(\lambda y) &=& - \int d^4z \; \delta(z^2) \: \delta((z-\lambda y)^2)
\; \epsilon(z^0) \: \epsilon(z^0-\lambda y^0) \; f(z) \\
&=& - \int d^4z \; \delta(z^2) \: \delta(- 2 \lambda \bra z,y \ket + \lambda^2
	y^2) \; \epsilon(z^0) \: \epsilon(z^0-\lambda y^0) \; f(z) \spc ,
\end{eqnarray*}
also
\begin{eqnarray*}
\lefteqn{- |\lambda| \int d^4y \; (\xint f)(\lambda y) \; g(y)}\\
&=& \frac{1}{2} \int d^4z \; \delta(z^2) \epsilon(z^0) \; f(z)
\int d^4y \; \delta(\bra z,y \ket - \frac{\lambda}{2} y^2) \; \epsilon(z^0 -
\lambda y^0) \; g(y) \spc .
\end{eqnarray*}
Im Grenzfall $\lambda \rightarrow 0$ konvergiert das innere Integral
gleichm\"a{\ss}ig
in $z$ gegen
\[	\int d^4y \; \delta(\bra z,y \ket) \: \epsilon(z^0) \; g(y) \spc , \]
woraus die Behauptung folgt.
\QED
F\"ur $y \in \Ra$ mu{\ss} die Funktion $f$ in~\Ref{a3_29} \"uber einen
zweidimensionalen Kegel integriert werden. Die Ableitungen von $f$ in
radialer Richtung liegen tangential zu diesem Kegel und k\"onnen partiell
integriert werden:
\begin{Lemma}
\label{a3_lemma3.9}
Es gilt f\"ur $y \in \Ra$ und die Funktion $g(z) = z^l \partial_l f(z)$
\Equ{a3_aaa}
    (\dint g)(y) \;=\; - (\dint f)(y)  \spc .
\EndEqu
\end{Lemma}
{\Beweis}
W\"ahle Zylinderkoordinaten $(t,x,r,\varphi)$\footnote{Also $y=r \cos \varphi,
z = r \sin \varphi$}, ohne Einschr\"ankung der Allgemeinheit kann man
$y=(0,x_0,0,0)$ annehmen.
\begin{eqnarray*}
( \dint g)(y) &=& \frac{1}{2} \int_{-\infty}^\infty dt \int_{-\infty}^\infty
dx \int_0^\infty r\; dr \int_0^{2 \pi} d \varphi \; \delta(t^2-r^2-x^2)
\; \delta(x \: x_0) \\
&&  \spc \times \left( t \frac{\partial}{\partial t} + x
\frac{\partial}{\partial x} + r \frac{\partial}{\partial r} \right) \:
f(t,r,r,\varphi)  \\
&=& \frac{1}{4 |x_0|} \int_0^{2 \pi} d\varphi \int_{-\infty}^\infty dt \;
t \frac{d}{d t} \: f(t,0, |t|, \varphi) \\
&=& - \frac{1}{4 |x^0|} \int_0^{2 \pi} d\varphi \int_{-\infty}^\infty dt \;
f(t,0,|t|, \varphi) \;=\; - (\dint f)(y)
\end{eqnarray*}
\QED
F\"ur die partiellen Ableitungen von $\sdint f$ hat man
\begin{Lemma}
\label{a3_lemma1a}
Es gilt im Distributionssinne
\Equ{a3_25}
\partial_j (\dint f)(y) \;=\; (\dint l_j^{(y)})(y)
\EndEqu
mit
\[  l_j^{(y)}(z) \;=\; - \frac{z_j \: y^l}{y^2} \; \partial_l f(z) - \frac{y_j}{y^2}
f(z) \spc . \]
Dies kann man auch schreiben als
\Equ{a3_26}
    \partial_j (\dint f)(y) \;=\; - \lim_{0 \neq \lambda \rightarrow 0}
	|\lambda| \; \frac{\partial}{\partial y^j} \: (\xint f)(\lambda y)
	\spc .
\EndEqu
\end{Lemma}
{\Beweis}
F\"ur eine Funktion $g \in C^\infty_c(M)$ hat man
\begin{eqnarray*}
\partial_j (\dint f)(g) &=& - \frac{1}{2} \int d^4z \; l(z) \: f(z) 
\int d^4y \; \delta(\bra z,y \ket) \; \partial_j g(y) \\
&=& \frac{1}{2} \int d^4z \; l(z) \: f(z) 
\int d^4y \; \delta^\prime(\bra z,y \ket) \; z_j \; g(y) \\
&=& \frac{1}{2} \int d^4z \; l(z) \: f(z) \; z_j \: \frac{\partial}{\partial
z^l} \int d^4y \; \frac{y^l}{y^2} \; \delta(\bra z,y \ket) \; g(y) \spc .
\end{eqnarray*}
Integriere nun partiell
\begin{eqnarray*}
&=& - \frac{1}{2} \int d^4z \; \frac{\partial}{\partial z^l} \left( z_j \:
	l(z) \; f(z) \right) \: \int d^4y \; \frac{y^l}{y^2} \;
	\delta(\bra z,y \ket) \; g(y) \\
&=& \frac{1}{2} \int d^4z \; \left( - g_{jl} \: f(z) - z_j \: \partial_l f
	\right)
	\int d^4y \; \frac{y^l}{y^2} \; \delta(\bra z,y \ket) \; g(y) \\
&& - \int d^4z \; l^\prime(z) \: f(z) \int d^4y \; \frac{1}{y^2}
	\bra z,y \ket \: \delta(\bra z,y \ket) \; g(y) \spc .
\end{eqnarray*}
Im letzten Summanden verschwindet das innere Integral, und es ergibt sich
\[  \;=\; \int d^4y \; (\dint l_j^{(y)})(y) \; g(y) \spc .   \]
Um Gleichung~\Ref{a3_26} zu beweisen, wendet man Satz~\ref{a3_satz2} mit
$h_j^{(y)}$ gem\"a{\ss} \Ref{a3_13} an
\begin{eqnarray*}
\lefteqn{ |\lambda| \: \frac{\partial}{\partial y^j} (\xint f)(\lambda y)
\;=\; |\lambda| \; \lambda (\xint h_j^{(\lambda y)})(\lambda y) } \\
&=& |\lambda| \: \lambda \xint_{\lambda y} dz \left( \frac{1}{2} \partial_j
f(z) - \frac{ z_j \: \lambda y_l - \lambda y_j \: z^l}{\lambda^2 y^2}
\: \partial_l f(z) \right) \\
&=& |\lambda| \xint_{\lambda y} dz \left( \frac{1}{2} \lambda \: \partial_j
f(z) - \frac{z_j \: y^l - y_j \: z^l}{y^2} \: \partial_l f(z) \right) \spc .
\end{eqnarray*}
Satz~\ref{a3_satz6} liefert
\[ - \lim_{0 \neq \lambda \rightarrow 0} |\lambda| \:
\frac{\partial}{\partial y^j} (\xint f)(\lambda y) \;=\; \dint_y -
\frac{z_j\: y^l - y_j \: z^l}{y^2} \: \partial_l f   \]
und nach Einsetzen von~\Ref{a3_aaa} die Behauptung.
\QED
Gleichung~\Ref{a3_26} ist leicht einsichtig, man erh\"alt sie aus~\Ref{a3_27}
durch formales Ableiten.
\\[1em]
Nun k\"onnen wir das Lichtkegelintegral $(\sxdint f)(y)$ definieren, das f\"ur
$y \in \I^\vee$ mit $(\slint f)(y)$ \"ubereinstimmt und auch auf dem
Lichtkegel (mit Ausnahme des Ursprungs) stetig ist:
\begin{Def}
Definiere f\"ur $f \in C^\infty_2(M)$ die Distribution $\sxdint f$ durch
\[   \xdint f \;=\; \xint f + \dint f	\spc .	\]
\end{Def}
\begin{Satz}
\label{a3_satz0}
$\sxdint f$ ist eine Funktion, die auf $M \setminus \{ 0 \}$ stetig ist.
F\"ur $0 \neq z \in \Li$ gilt
\Equ{a3_50}
   \lim_{y \rightarrow z} \xdint f \;=\; \frac{\pi}{2} \int_0^1 f(\alpha z)
	\; d\alpha \spc .
\EndEqu
Am Ursprung ist $\sxdint f$ jedoch im allgemeinen nicht stetig.
F\"ur $y \in \I$
hat man
\Equ{a3_51}
  \lim_{0 < \lambda \rightarrow 0} (\xdint f)(\lambda y) \;=\; \frac{\pi}{2}
\: f(0) \spc ,
\EndEqu
f\"ur $y \in \Ra$ dagegen
\Equ{a3_52}
  \lim_{0 < \lambda \rightarrow 0} (\xdint f)(\lambda y) \;=\; \frac{\pi}{2}
\: f(0) - \frac{1}{2} (\dint y^j \partial_j f)(y) \spc .
\EndEqu
\end{Satz}
{\Beweis}
Aus Satz~\ref{a3_satz1} und~\ref{a3_satz5} folgt unmittelbar, da{\ss} $\sxdint
f$ eine auf $M \setminus \{0\}$ stetige Funktion ist, die~\Ref{a3_50}
erf\"ullt.
Au{\ss}erdem ist klar, da{\ss} f\"ur $y \in \I^\vee$ die Lichtkegelintegrale
$(\sxdint f)(y)$ und $(\slint f)(y)$ \"ubereinstimmen. Gleichung~\Ref{a3_51}
folgt f\"ur $y \in \I^\vee$ somit aus Satz~\ref{a1_lemma1}, den Fall $y \in
\I^\wedge$ erh\"alt man durch Punktspiegelung am Ursprung.

Es bleibt also noch~\Ref{a3_52} abzuleiten. Man kann ein
Bezugssystem w\"ahlen mit $y=(0,x_0,0,0)$. In Zylinderkoordinaten
$(t,x,r,\varphi)$ erh\"alt man
\begin{eqnarray*}
(\dint f)(y) &=& \frac{1}{2\:|x_0|} \int_{-\infty}^\infty dt \int_0^\infty r\;
dr \int_0^{2 \pi} d\varphi \; \delta(t^2-r^2) \; f(t,0,r,\varphi) \\
&=& \frac{1}{4\:|x_0|} \int_0^{2 \pi} d\varphi \int_{-\infty}^\infty dt \;
	f(t,0,|t|, \varphi) \\
(\dint f)(\lambda y) &=& \frac{1}{4\:|\lambda x_0|} \int_0^{2 \pi} d\varphi
	\int_{-\infty}^\infty dt \; f(t,0,|t|, \varphi) \\
(\xint f)(\lambda y) &=& - \frac{1}{2\:|\lambda x_0|} \int_{-\infty}^\infty
dt \int_0^\infty r \; dr \int_0^{2 \pi} d\varphi \;
\delta(t^2-r^2-\frac{1}{4} \lambda^2 x_0^2) \: f(t,\frac{1}{2} \lambda x_0,
r,\varphi) \\
&=& -\frac{1}{4\:|\lambda x_0|} \int_0^{2\pi}
\int_{\sR \setminus [-\lambda \alpha, \lambda
\alpha ]} dt \; f(t, \frac{1}{2} \lambda x_0, r(t), \varphi)
\end{eqnarray*}
mit $r(t)=(t^2 - \alpha^2 \lambda^2)^{\frac{1}{2}}$, $\alpha=|x_0|/2$.
Somit gilt
\begin{eqnarray*}
(\xdint f)(\lambda y) &=& \frac{1}{4\:|\lambda x_0|} \int_0^{2 \pi} d\varphi
\int_{\sR \setminus [-\lambda \alpha, \lambda \alpha]} dt \; \left(
f(t,0,|t|,\varphi) - f(t,\frac{1}{2} \lambda x_0, r(t), \varphi) \right) \\
&& + \frac{1}{4\:|\lambda x_0|} \int_0^{2 \pi} d\varphi \int_{-\lambda
\alpha}^{\lambda \alpha} dt \; f(t,0,|t|, \varphi) \spc .
\end{eqnarray*}
Eine Taylorentwicklung in $\lambda$ liefert
\begin{eqnarray*}
&=& - \frac{x_0}{2} \: \frac{1}{4\:|x_0|} \int_0^{2 \pi} d\varphi
\int_{-\infty}^\infty dt \; \frac{\partial}{\partial x}
	f(t,0,|t|,\varphi) + \frac{1}{4\:|\lambda x_0|} \: \lambda |x_0|
	\: 2 \pi f(0) \;+\; {\cal{O}}(\lambda)
\end{eqnarray*}
und somit
\[ \lim_{0 < \lambda \rightarrow 0} (\xdint f)(\lambda y) \;=\; -
\frac{x_0}{2} \left( \dint \frac{\partial}{\partial x} f \right)(y) +
\frac{\pi}{2} \: f(0) \spc ,	\]
was, koordinateninvariant geschrieben, mit~\Ref{a3_52} \"ubereinstimmt.
\QED
F\"ur die partiellen Ableitungen von $\sxdint f$ erh\"alt man eine
einfache Formel:
\begin{Satz}
\label{a3_satz01}
Es gilt im Distributionssinne
\Equ{a3_30}
  \partial_j (\xdint f)(y) \;=\; (\xint h_j)(y) 
\EndEqu
mit
\[  h_j(z) \;=\; \partial_j f(z) - \frac{1}{2} \Box_z \left( z_j \: \int_0^1
f(\alpha z) \: d\alpha \right) \spc . \]
\end{Satz}
{\Beweis}
F\"ur $y \in \I$ hat man nach Satz~\ref{a3_satz5} und Satz~\ref{a3_satz2}
\[ \partial_j (\xdint f)(y) \;=\; \partial_j (\xint f)(y) \;=\; (\xint h_j)(y)
\spc . \]
F\"ur $y  \in \Ra$ gilt nach~\Ref{a3_26} und~\Ref{a3_11}
\begin{eqnarray}
\partial_j (\xdint f)(y) &=& \partial_j (\xint f)(y) + \partial_j (\dint
f)(y) \nonumber \\
&=& \partial_j (\xint f)(y) - \lim_{0 \neq \lambda \rightarrow 0} |\lambda| \:
\partial_j (\xint f)(\lambda y) \nonumber \\
\label{eq:a3_291}
&=& \xint_y (h_j + g_j) - \lim_{0 \neq \lambda \rightarrow 0} |\lambda| \:
\lambda \xint_{\lambda y} (h_j + g_j)
\end{eqnarray}
mit
\begin{eqnarray*}
g_j(z) &=& \frac{1}{2} \: \Box_z \left( z_j \int_0^\infty f(\alpha z) \:
	d\alpha \right) \\
&=& \int_0^\infty \alpha \: \partial_j f(\alpha z) d\alpha + \frac{1}{2} \:
	z_j \int_0^\infty \alpha^2 \: (\Box f)(\alpha z) \: d\alpha \spc .
\end{eqnarray*}
Mit Hilfe der Relation
\begin{eqnarray*}
g_j(\lambda z) &=& \int_0^\infty \alpha \: \partial_j f(\lambda \alpha z)
\: d\alpha + \frac{1}{2} \: \lambda z_j \int_0^\infty \alpha^2 \: (\Box
f)(\lambda \alpha z) \: d\alpha \\
&=& \frac{1}{|\lambda| \: \lambda} \: g_j(z)
\end{eqnarray*}
erh\"alt man
\[ |\lambda| \: \lambda \xint_{\lambda y} g_j(z) \: dz \;=\; |\lambda| \:
\lambda \xint_y g_j(\lambda z) \: dz \;=\; \xint_y g_j  \]
und nach Einsetzen in~\Ref{a3_291}
\[ \partial_j (\xdint f)(y) \;=\; (\xint h_j)(y) - \lim_{0 \neq \lambda
\rightarrow 0} |\lambda| \: \lambda (\xint h_j)(\lambda y) \spc . \]
Die Funktion $h_j$ ist in $C^\infty_2(M)$. Nach
Satz~\ref{a3_satz6} verschwindet der zweite Summand und man
erh\"alt~\Ref{a3_30}.

Es bleibt zu zeigen, da{\ss}~\Ref{a3_30} im Distributionssinne gilt. F\"ur
$g \in C^\infty_c(M)$ hat man
\begin{eqnarray*}
\partial_j (\xdint f)(g) &=& - \int d^4y \; (\xdint f)(y) \; \partial_j g(y)
	\\
&=& - \int_{\I} d^4y \; (\xdint f)(y) \; \partial_j g(y)
- \int_\Ra d^4y \; (\xdint f)(y) \; \partial_j g(y) \spc .
\end{eqnarray*}
Da $\sxdint f$ auf $\I \cup \Ra$ glatt ist, kann man partiell integrieren.
Die Randterme fallen weg, weil $\sxdint f$ auf dem Lichtkegel stetig ist
\begin{eqnarray*}
&=& \int_{\I \cup \Ra} d^4y \; \partial_j (\xdint f)(y) \; g(y) \\
&=& \int_{\I \cup \Ra} d^4y \; (\xint h_j)(y) \; g(y) \;=\; (\xint h_j)(g)
\spc ,
\end{eqnarray*}
da $\sxint h_j$ eine Funktion ist.
\QED
Man beachte, da{\ss} $\partial_j \sxdint f$ nach~\Ref{a3_30} auf dem
Lichtkegel im allgemeinen nicht stetig ist.
\\[1em]
Als letztes Lichtkegelintegral soll nun die Distribution
\[  (\veeint f)(y) \;=\; \int d^4z \; \frac{1}{(y-z)^2} \: l(z) \; f(z)  \]
eingef\"uhrt werden. Dazu braucht man folgendes Lemma:
\begin{Lemma}
\label{a3_lemma10}
F\"ur $g \in C^\infty_c(M)$ ist
\Equ{a3_22}
  h(z) \;=\; \int d^4y \; \frac{1}{(y-z)^2} \; g(y)
\EndEqu
eine linear abfallende Funktion. (Das Integral ist als Hauptwert definiert.)
\end{Lemma}
{\Beweis}
Sei $0 \neq y \in M$. W\"ahle ein Bezugssystem mit $\vec{y}\neq 0$, setze
$\kappa = |\vec{y}|/2$. Da $g$ kompakten Tr\"ager besitzt, gibt es $R > 0$
mit $\mbox{supp }g \subset \R \times B_R(\vec{0})$.

W\"ahle $U = \R \times B_\kappa(\vec{y})$. F\"ur $z \in U$, $\lambda > 2 R /
\kappa$ hat man
\begin{eqnarray*}
\lambda \: h(\lambda z) &=& \lambda \int d^4y \; \frac{1}{(y - \lambda z)^2}
	\; g(y) \\
&=&	\lambda \hint_{-\infty}^\infty dt\int_{B_R(\vec{0})} d^3 \vec{x} \;
\frac{1}{(t-\lambda z^0)^2 - (\vec{x}-\lambda \vec{z})^2} \; g(t,\vec{x})
\end{eqnarray*}
Eine Partialbruchzerlegung liefert
\Equ{a3_37}
   \;=\; \int_{B_R(\vec{0})} d^3 \vec{x} \hint_{-\infty}^\infty dt \; \frac{1}{2a}
  \left(
  \frac{1}{t-\lambda z^0-\lambda a} -
	\frac{1}{t-\lambda z^0+\lambda a} \right) \; g(t,\vec{x})
\EndEqu
mit $a=a(\lambda, \vec{x})=|\vec{z}-\vec{x}/\lambda|$.
Wegen
\[   a(\lambda, \vec{x}) \geq |\vec{z}| - \frac{|\vec{x}|}{\lambda}
\geq \kappa - \frac{R}{\lambda} \geq \frac{\kappa}{2}  \]
hat man
\Equ{a3_194b}
\lambda \; |h(\lambda z)| \;\leq\; \frac{1}{\kappa} \int_{B_R(\vec{0})}
	d^3\vec{x}
	\left| \hinti \left( \frac{1}{t-\lambda z^0-\lambda a} -
	\frac{1}{t-\lambda z^0+\lambda a} \right) \; g(t,\vec{x}) \right|
	\;\;\; .
\EndEqu
In~\Ref{a3_194b} ist das innere Integral, und somit auch der ganze Ausdruck,
in $\lambda$ und $z$ glm. beschr\"ankt.
\QED
\begin{Def}
Definiere f\"ur $f \in C^\infty_2(M)$ die Distribution $\sveeint f$ durch
\[ (\veeint f)(g) \;=\; \int d^4z \; l(z) \: h(z) \; f(z) \spc , g \in
C^\infty_c(M)  \]
mit $h$ gem\"a{\ss}~\Ref{a3_22}.
\end{Def}
Diese Definition ist sinnvoll, weil $h$ nach Lemma~\ref{a3_lemma10} eine
linear abfallende Funktion ist.
\begin{Satz}
\label{a3_satz99}
$\sveeint f$ ist harmonisch (im Distributionssinne).

Au{\ss}erdem ist $\sveeint f$ eine Funktion, die auf $\I \cup \Ra$ glatt ist.
In der N\"ahe des Lichtkegels besitzt $\sveeint f$ eine logarithmische
Divergenz, also f\"ur $z \in \Li$
\[  \lim_{y \rightarrow z} \left\{ (\veeint f)(y) - \frac{\pi}{2} \ln(|y^2|)
	\: \hint_{-\infty}^\infty d\lambda \: \epsilon(\lambda) \;
f(\lambda y) \right\} \;<\; \infty  \spc . \]
\end{Satz}
{\Beweis} Verl\"auft \"ahnlich wie der Beweis von Satz~\ref{a1_lemma1}.
\QED
Es ist g\"unstig, das Lichtkegelintegral $\sveeint f$ mit Hilfe von
$\sxint f$, $\sxdint f$ auszudr\"ucken. Dadurch wird es m\"oglich sein,
ohne gro{\ss}en Aufwand auch Formeln f\"ur die partiellen Ableitungen von
$\sveeint f$ abzuleiten.
\begin{Lemma}
\label{a3_lemma11}
Es gilt im Distributionssinne
\begin{eqnarray}
\label{eq:a3_32}
(\veeint f)(y) &=& \hinti d\lambda \; \frac{|\lambda|}{1-\lambda} \:
	(\xint f)(\lambda y)  \spc .
\end{eqnarray}
\end{Lemma}
{\Beweis}
Rechne f\"ur $g \in C^\infty_c(M)$:
\begin{eqnarray*}
\lefteqn{ \hinti d\lambda \; \frac{|\lambda|}{1-\lambda} \int d^4y \; 
	(\xint f)(\lambda y) \; g(y) } \\
&=& - \hinti d\lambda \; \frac{|\lambda|}{1-\lambda} \int d^4z \; 
(l^\vee(z) - l^\wedge(z)) \; f(z) \int d^4y \; (\l^\vee(z-\lambda y) -
l^\wedge(z-\lambda y)) \; g(y)
\end{eqnarray*}
F\"ur $y \in \I$ tragen die Terme $l^\vee(z) \: l^\wedge(z-\lambda y)$,
$l^\wedge(z) \: l^\vee(z-\lambda y)$ bei, f\"ur $y \in \Ra$ dagegen die
Terme $l^\vee(z) \: l^\vee(z-\lambda y)$,
$l^\wedge(z) \: l^\wedge(z-\lambda y)$. Dies f\"uhrt zu einem relativen
Minuszeichen f\"ur die F\"alle $y \in \I$, $y \in \Ra$ und man kann schreiben
\[ \;=\; \hinti d\lambda \; \frac{|\lambda|}{1-\lambda} \int d^4z \;
	l(z) \: f(z) \int d^4y \; l(z-\lambda y) \: \epsilon(y^2) \: g(y)
\spc . \]
Da die inneren Integrale stetig von $\lambda$ abh\"angen, kann man die
$\lambda$-Integration ausf\"uhren und erh\"alt einen Beitrag bei $\lambda = 2
\bra z,y \ket / y^2$
\begin{eqnarray*}
&=& \int d^4z \; l(z) \: f(z) \int d^4y \; \left| \frac{2 \bra z,y \ket}{y^2}
\right| \; \frac{y^2}{(y^2-2 \bra z,y \ket)} \frac{1}{|2 \bra z,y \ket|} \;
\epsilon(y^2) \: g(y) \\
&=& \int d^4z \; l(z) \: f(z) \int d^4y \; \frac{1}{(y-z)^2} \; g(y) \spc ,
\end{eqnarray*}
wobei im Integral nach $y$ wie gew\"unscht der Hauptwert auftritt.
\QED

Wir werden nun der Einfachheit halber annehmen, da{\ss} $f$ gleichm\"a{\ss}ig
zur dritten
Potenz abf\"allt.
\begin{Satz}
\label{a3_satz13}
F\"ur $f \in C^\infty_3(M)$ gilt im Distributionssinne
\Equ{a3_194a}
\partial_j (\veeint f)(y) \;=\; (\veeint h_j)(y) + \frac{2}{y^2}
\hint_{-\infty}^\infty
\frac{d\lambda}{|\lambda|} \; \xdint_{\lambda y} dz \; z_j \: f(z)
\EndEqu
mit
\[  h_j(z) \;=\; \partial_j f(z) - \frac{1}{2} \: \Box_z \left( z_j \;
\int_0^1 f(\alpha z) \: d\alpha \right) \spc . \]
\end{Satz}
{\Beweis}
Nach~\Ref{a3_32} gilt
\begin{eqnarray}
(\veeint f)(y) &=& \hinti d\lambda \; \frac{|\lambda|}{1-\lambda}
	\xint_{\lambda y} f \nonumber \\
\label{eq:a3_33a}
&=& \hinti d\lambda \; \frac{\epsilon(\lambda)}{1-\lambda} \xint_{\lambda y}
	f - \hinti d\lambda \; \epsilon(\lambda) \xint_{\lambda y} f \spc .
\end{eqnarray}
Aufgrund der Beziehung
\[  \hinti d\lambda \; \frac{\epsilon(\lambda)}{1-\lambda} \dint_{\lambda y} f
	\;=\; \hinti d\lambda \; \frac{1}{\lambda (1-\lambda)} \dint_y f \;=\;
	 0	\]
kann man im ersten Summanden von~\Ref{a3_33a} das Lichtkegelintegral
$\sxint f$ durch $\sxdint f$ ersetzen
\[   (\veeint f)(y) \;=\; \hinti d\lambda \; \frac{\epsilon(\lambda)}{1-\lambda}
	 \xdint_{\lambda y} f
	- \hinti d\lambda \; \epsilon(\lambda) \xint_{\lambda y} f \spc . \]
Durch Differentiation erh\"alt man unter Verwendung
von~\Ref{a3_30} und~\Ref{a3_32} die Gleichung
\Equ{a3_35}
  \partial_j (\veeint f)(y) \;=\; (\veeint h_j)(y) - \hinti d\lambda \;
  \epsilon(\lambda) \: \frac{\partial}{\partial y^j} (\xint f)(\lambda y)
  \spc. \EndEqu
Im zweiten Sumanden von~\Ref{a3_35} kann man durch geschickte partielle
Integration die Ableitung beseitigen:

Da $f$ gleichm\"a{\ss}ig zur dritten Potenz abf\"allt, kann man
Satz~\ref{a3_satz2} mit
$h_j^{(y)}$ in der Form~\Ref{a3_14} anwenden.
\begin{eqnarray}
\lefteqn{ \hint_{-\infty}^\infty d\lambda \; \epsilon(\lambda) \:
	\frac{\partial}{\partial y^j}
	\xint_{\lambda y} f \;=\; \hint_{-\infty}^\infty d\lambda \;
	|\lambda| \: \xint_{\lambda y} h_j^{(\lambda y)} } \nonumber \\
&=& \hint_{-\infty}^\infty d\lambda \; |\lambda| \xint_{\lambda y} dz \;
	\left( \frac{ 2 (\lambda y-2z)_j}{\lambda^2 y^2} \: f(z) +
	\frac{2 (\lambda y-z)_j}{\lambda^2 y^2} \: z^l \partial_l f(z)
	\right) \nonumber \\
&=& \frac{2}{y^2} \: \hint_{-\infty}^\infty d\lambda \; \epsilon(\lambda)
	\xint_{\lambda y} dz \; \left\{ z^l \frac{\partial}{\partial z^l}
	\left( (y-\frac{1}{\lambda} z)_j \: f(z) \right) +
	(y-\frac{1}{\lambda} z)_j \: f(z)
	\right\} \nonumber \\
&=& \frac{2}{y^2} \: \hint_{-\infty}^\infty d\lambda \; \epsilon(\lambda) \;
	\left\{ \xint_{\lambda y} dz \; (y-\frac{2}{\lambda} z)_j \: f(z) +
	\lambda \frac{d}{d \lambda}  \xint_{\lambda y} dz \;
	(y-\frac{1}{\lambda} z)_j \: f(z) \right\} \nonumber \\
\label{eq:a3_36}
&=& \frac{2}{y^2} \: \hint_{-\infty}^\infty d\lambda \; \epsilon(\lambda) \;
	\left\{ -\frac{1}{\lambda} \xint_{\lambda y} dz \; z_j \: f(z) +
	\frac{d}{d\lambda} \xint_{\lambda y} dz \; (\lambda y - z)_j \: f(z)
	 \right\}
\end{eqnarray}
An dieser Stelle tritt ein kleines Problem auf. Man kann n\"amlich
in~\Ref{a3_36} nicht ohne weiteres partiell integrieren, weil die Randterme
bei $\lambda=0$ f\"ur raumartiges $y$ nicht verschwinden.
Denn nach Satz~\ref{a3_satz1} besitzt das Lichtintegral $\sxint_{\lambda y}
dz \; (\lambda y - z)_j \: f(z)$ f\"ur $y \in \Ra$ bei $\lambda=0$ einen Pol.

Wir werden nun zun\"achst zeigen, da{\ss} man die Lichtkegelintegrale $\sxint$ in
\Ref{a3_36} durch $\sxdint$ ersetzen kann. Nach~\Ref{a3_181a} folgt
\begin{eqnarray*}
\lefteqn{ \hint_{-\infty}^\infty d\lambda \; \epsilon(\lambda)
	\; \frac{d}{d\lambda} \dint_{\lambda
	y} dz \; (\lambda y - z)_j \: f(z) } \\
&=& \hint_{-\infty}^\infty d\lambda \; \epsilon(\lambda) \frac{d}{d \lambda}
	\left\{ \epsilon(\lambda) \dint_{y} dz \; y_j \: f(z) -
	\frac{1}{|\lambda|} \dint_y dz \; z_j \: f(z) \right\} \\
&=& - \hint_{-\infty}^\infty d\lambda \; \epsilon(\lambda) \;
	\frac{1}{\lambda}
	\frac{1}{|\lambda|} \dint_y dz \; z_j \: f(z)
\;=\; - \hint_{-\infty}^\infty \frac{d \lambda}{|\lambda|} \: \dint_{\lambda y}
	dz \; z_j \: f(z)
\end{eqnarray*}
und nach Addieren zu~\Ref{a3_36}
\begin{eqnarray*}
\lefteqn{ \hint_{-\infty}^\infty d\lambda \; \epsilon(\lambda) \:
	\frac{\partial}{\partial y^j}
	(\xint f)(\lambda y) } \\
&=& \frac{2}{y^2} \hint_{-\infty}^\infty d\lambda \; \epsilon(\lambda) \left\{
- \frac{1}{\lambda} \xdint_{\lambda y} dz \; z_j \: f(z) + \frac{d}{d
\lambda} \xdint_{\lambda y} dz \; (\lambda y - z)_j \; f(z) \right\} \spc .
\end{eqnarray*}
Nun kann man problemlos partiell integrieren. Die Randwerte im Unendlichen
verschwinden, weil $f$ gleichm\"a{\ss}ig zur dritten Potenz abf\"allt. 
F\"ur $y\in \I$ fallen auch die Randwerte f\"ur $\lambda=0$ weg.
Falls $y \in \Ra$, hat man nach~\Ref{a3_52} 
\begin{eqnarray*}
\lefteqn{
\lim_{0 < \lambda \rightarrow 0} \xdint_{\lambda y} dz \; (\lambda y-z)_j \:
f(z) \;=\; - \lim_{0<\lambda \rightarrow 0} \xdint_{\lambda y} dz \; z_j
	\: f(z) } \\
&=& \frac{1}{2} \dint_y dz \; (y^l \frac{\partial}{\partial z^l}) (z^j f(z))
	\;=\; - \lim_{0>\lambda \rightarrow 0} \xdint_{\lambda y} dz \;
	(\lambda y - z)_j \: f(z) \spc ,
\end{eqnarray*}
so da{\ss} sich die Randwerte f\"ur $\lambda = \pm 0$ wegheben.

Man erh\"alt also
\[  \hint_{-\infty}^\infty d\lambda \;
\epsilon(\lambda) \; \frac{\partial}{\partial y^j} (\xint
f)(\lambda y) \;=\; -\frac{2}{y^2} \hint_{-\infty}^\infty
\frac{d\lambda}{|\lambda|} \: \xdint_{\lambda y} dz \; z_j \: f(z)  \]
und nach Einsetzen in~\Ref{a3_35} die Behauptung.
\QED
F\"ur Funktionen, die gleichm\"a{\ss}ig zur vierten Potenz abfallen, kann man
auch eine Formel f\"ur die zweiten Ableitungen herleiten:
\begin{Satz}
\label{a3_satz14}
F\"ur $f \in C^\infty_4(M)$ gilt im Distributionssinne
\begin{eqnarray}
\lefteqn{ \partial_{jk} (\veeint f)(y) \;=\; (\veeint h_{jk})(y) -
	\frac{1}{y^2} \hint_{-\infty}^\infty \frac{d\lambda}{|\lambda|} \:
	\xdint_{\lambda y} dz \; z_j \: z_k \; \Box f(z) } \nonumber \\
\label{eq:a3_39}
&& - \frac{8}{y^4} \hint_{-\infty}^\infty \frac{d\lambda}{|\lambda|}
	\left\{ \frac{1}{\lambda} \xdint_{\lambda y} dz \; z_j \: z_k \;
	f(z) + \frac{1}{2} \dint_{\lambda y} dz \; (y^l \partial_l)
	(z_j\: z_k \: f(z)) \right\}
\end{eqnarray}
mit $h_{jk}$ gem\"a{\ss}~\Ref{a1_60}.

Dabei ist mit $1/y^2$ wiederum der Hauptwert gemeint. Die Distribution $1/y^4$
bezeichnet die Ableitung des Hauptwertes, also f\"ur $g \in C^\infty_c(M)$
\begin{eqnarray}
\int d^4z \; \frac{1}{y^4} \: g(y) &:=& - \lim_{0 \neq \varepsilon
\rightarrow 0} \frac{1}{2} \: \int d^4y \; \frac{d}{d(y^2)} \left(
\frac{1}{y^2+i \varepsilon} + \frac{1}{y^2-i \varepsilon} \right) \; g(y) \nonumber \\
\label{eq:haupt2}
&=& \lim_{0 \neq \varepsilon \rightarrow 0} \frac{1}{2} \: \int d^4y
\left( \frac{1}{(y^2 + i \varepsilon)^2} + \frac{1}{(y^2-i \varepsilon)^2}
\right) \; g(y) \spc .
\end{eqnarray}
\end{Satz}
{\Beweis}
Aus Satz~\ref{a3_satz13} folgt
\begin{eqnarray}
\partial_{jk} (\veeint f)(y) &=& \frac{\partial}{\partial y^k} \left\{
(\veeint h_j)(y) + \frac{2}{y^2} \hint_{-\infty}^\infty
\frac{d\lambda}{|\lambda|}
\xdint_{\lambda y} dz \; z_j \: f(z) \right\} \nonumber \\
\label{eq:a3_38}
&=& (\veeint h_{jk})(y) + \frac{2}{y^2} \hint_{-\infty}^\infty
\frac{d \lambda}{|\lambda|} \xdint_{\lambda y} dz \; z_k \: h_j(z) \\
&& + \frac{\partial}{\partial y^k} \left\{ \frac{2}{y^2}
\hint_{-\infty}^\infty \frac{d\lambda}{|\lambda|} \xdint_{\lambda y} dz \;
z_j \: f(z) \right\} \spc .
\end{eqnarray}
Den zweiten Summanden in~\Ref{a3_38} kann man folgenderma{\ss}en umformen:
\begin{eqnarray*}
\lefteqn{ \hint_{-\infty}^\infty \frac{d\lambda}{|\lambda|}
	\xdint_{\lambda y} dz \; z_k \: h_j(z) } \\
&=& \hinti \frac{d\lambda}{|\lambda|} \xdint_{\lambda y}
	dz \; z_k \left( \partial_j f(z) - \int_0^1 \alpha \: \partial_j
	f(\alpha z) \: d\alpha - \frac{1}{2} \: z_j \int_0^1 \alpha^2 \:
	\Box f(\alpha z) \; d\alpha \right) \\
&=& \int \frac{d\lambda}{|\lambda|} \left( \xdint_{\lambda y} dz \; z_k \:
	\partial_j f(z) - \int_0^1 d\alpha \xdint_{\alpha \lambda y} dz \left(
	z_k \: \partial_j f(z) + \frac{1}{2} z_j \: z_k \: \Box f(z)
	\right) \right)
\end{eqnarray*}
Mit den Integralumformungen~\ref{integr_umf} erh\"alt man
\begin{eqnarray*}
&=& - \frac{1}{2} \hint_{-\infty}^\infty \frac{d\lambda}{|\lambda|} \;
	\xdint_{\lambda y} dz \; z_j \: z_k \; \Box f(z) \spc .
\end{eqnarray*}
Durch Einsetzen in~\Ref{a3_38} ergibt sich gerade der zweiten Summand
in~\Ref{a3_39}.

Es bleibt somit zu zeigen, da{\ss}
\begin{eqnarray*}
\lefteqn{ \frac{\partial}{\partial y^k} \left\{ \frac{2}{y^2} \:
	\hint_{-\infty}^\infty \frac{d\lambda}{|\lambda|}
	\xdint_{\lambda y} dz \; z_j \: f(z) \right\}  } \\
&=& - \frac{8}{y^4} \hint_{-\infty}^\infty \frac{d\lambda}{|\lambda|} \;
\left\{	\frac{1}{\lambda} \xdint_{\lambda y} dz \; z_j \: z_k \: f(z) +
\frac{1}{2} \dint_{\lambda y} dz \; (y^l \partial_l) (z_j \: z_k \: f(z))
	\right\} \spc .
\end{eqnarray*}
Hierzu sind zun\"achst eine Reihe von Umformungen n\"otig:
Nach Lemma~\ref{a3_lemma1a} und~\Ref{a3_181a} hat man
\begin{eqnarray}
\frac{1}{\lambda} \frac{\partial}{\partial y^k} \dint_{\lambda y}
dz \; z_j \: f(z) &=& - \frac{\epsilon(\lambda)}{\lambda^2} \dint_y dz \;
\left( \frac{z_k}{y^2} \: y^l \partial_l + \frac{y_k}{y^2} \right) (z_j f(z)
	\nonumber \\
&=& - \frac{\epsilon(\lambda)}{\lambda^2} \dint_y dz \; (y^l \partial_l)
	(z_j \: z_k \; f(z)) \nonumber \\
&=& - \frac{1}{\lambda y^2} \dint_{\lambda y} dz \; (y^l \partial_l)
	(z_j \: z_k \; f(z)) \nonumber \\
&=& - \frac{2}{y^2} \frac{1}{\lambda} \dint_{\lambda y} dz \; (y^l \partial_l)
	(z_j \: z_k \; f(z)) \nonumber \\
\label{eq:a3_41}
&& - \frac{1}{y^2} \: \frac{d}{d\lambda} \dint_{\lambda y} dz \; y^l
	\partial_l (z_j\: z_k \: f(z)) \spc .
\end{eqnarray}

Da $f \in C^\infty_4(M)$, f\"allt die Funktion $z_j f(z)$ gleichm\"a{\ss}ig zur
dritten Potenz ab. Anwendung von Satz~\ref{a3_satz2} mit $h_j^{(y)}$ in der
Form~\Ref{a3_14} liefert
\begin{eqnarray}
\lefteqn{ \frac{1}{\lambda} \frac{\partial}{\partial y^k} \xint_{\lambda y}
dz \; z_j \: f(z) } \nonumber \\
&=& \frac{2}{\lambda y^2} \xint_{\lambda y} dz \; \left( (z-\frac{2}{\lambda}
z)_k + (y-\frac{1}{\lambda} z)_k \; z^l \partial_l \right) (z_j f(z))
	\nonumber \\
&=& \frac{2}{\lambda y^2} \xint_{\lambda y} dz \; \left\{ z^l \partial_l
	\left( (y-\frac{1}{\lambda}z)_k \: z_j \: f(z) \right) +
	(y-\frac{1}{\lambda} z)_k \: z_j \: f(z) \right\} \nonumber \\
\label{eq:a3_40}
&=& \frac{2}{y^2} \; \frac{d}{d\lambda} \xint_{\lambda y} dz \;
	(y-\frac{1}{\lambda} z)_k \: z_j\: f(z)
	+ \frac{2}{y^2} \: \frac{1}{\lambda} \xint_{\lambda y} dz \;
	(y-\frac{2}{\lambda} z)_k \: z_j \: f(z) \spc .
\end{eqnarray}
\"Ahnlich wie im Beweis von Satz~\ref{a3_satz13} kann man in~\Ref{a3_40} das
Lichtkegelintegral $\sxint$ durch $\sxdint$ ersetzen, weil
\begin{eqnarray}
\lefteqn{ \frac{d}{d \lambda} \dint_{\lambda y} dz \; (y-\frac{1}{\lambda}
z)_k \: z_j \; f(z)  } \nonumber \\
&=& \frac{d}{d\lambda} \left\{ \frac{1}{|\lambda|} \dint_y dz \; y_k \: z_j
	\: f(z) - \frac{\epsilon(\lambda)}{\lambda^2} \dint_y dz \; z_k
	\: z_j \: f(z) \right\} \nonumber \\
&=& - \frac{\epsilon(\lambda)}{\lambda^2} \dint_y dz \; y_k \: z_j \: f(z) +
	2 \frac{\epsilon(\lambda)}{\lambda^3} \dint_y dz \; z_k \: z_j \: f(z)
	\nonumber \\
\label{eq:a3_42}
&=& - \frac{1}{\lambda} \dint_{\lambda y} dz \; (y-\frac{2}{\lambda}z)_k \:
	z_j \; f(z) \spc .
\end{eqnarray}
Durch Addition von~\Ref{a3_41}, \Ref{a3_40}, \Ref{a3_42}
und Integration \"uber $\lambda$ erh\"alt man
\begin{eqnarray*}
\lefteqn{ \hint_{-\infty}^\infty \frac{d\lambda}{|\lambda|} \:
	\frac{\partial}{\partial y^k} \xdint_{\lambda y} dz \; z_j \: f(z)
	} \\
&=& \frac{2}{y^2} \hint_{-\infty}^\infty \frac{d\lambda}{|\lambda|} \left\{
	\xdint_{\lambda y} dz \; (y-\frac{2}{\lambda} z)_k \: z_j \; f(z)
	- \dint_{\lambda y} dz \; (y^l \partial_l) (z_j \: z_k \: f(z))
	\right\} \\
&& + \frac{2}{y^2} \hint_{-\infty}^\infty d\lambda \; \epsilon(\lambda) \;
	\frac{d}{d \lambda} \left\{ \xdint_{\lambda y} dz \;
	(y - \frac{1}{\lambda} z)_k \: z_j \; f(z)
	- \frac{1}{2} \dint_{\lambda y} dz \; (y^l \partial_l)(z_j \:
	z_k \: f(z)) \right\} \; .
\end{eqnarray*}
Bei der partiellen Integration mu{\ss} man etwas aufpassen. Die Randwerte
im Unendlichen verschwinden, weil $f$ gleichm\"a{\ss}ig zur vierten Potenz
abf\"allt.
F\"ur $y \in \I$ fallen die Randwerte bei $\lambda=0$ wegen~\Ref{a3_51}
weg.
Der Fall $y \in \Ra$ mu{\ss} sorgf\"altig untersucht werden. Setze zur
Abk\"urzung
\[  g(\lambda) \;=\; \xdint_{\lambda y}
	dz \; (y- \frac{1}{\lambda} z)_k \: z_j \; f(z) - \frac{1}{2}
	\dint_{\lambda y} dz \; (y^l \partial_l)(z_j \:
	z_k \: f(z)) 	\spc . \]
Mit Hilfe von~\Ref{a3_52} folgt
\[	\lim_{\lambda \rightarrow 0} \lambda \: g(\lambda) \;=\; 0 \spc . \]
Man erh\"alt mit den Regeln von de l'Hopital
\begin{eqnarray*}
\lim_{0<\lambda \rightarrow 0} g(\lambda) &=& \lim_{0<\lambda \rightarrow
	0} \frac{d}{d\lambda} \left( \lambda \: g(\lambda) \right) \\
&=& - \lim_{0<\lambda \rightarrow 0} \left( \xdint_{\lambda y} dz \; y_k \:
z_j \: f(z) - \frac{d}{d\lambda} \xdint_{\lambda y} dz \; z_j \: z_k \: f(z)
	\right) \spc ,
\end{eqnarray*}
und mit Hilfe von Satz~\ref{a3_satz01}, \Ref{a3_52} und~\Ref{a3_27}
\begin{eqnarray*}
&=& - \lim_{0<\lambda \rightarrow 0} \left( \xdint dz \; y_k \: z_j \: f(z)
	- \lambda \xint_{\lambda y} h_j[z_j z_k f(z)] \right) \\
&=& \frac{1}{2} \dint_y dz \; (y^l \partial_l) (y_k \: z_j \: f(z))
	+ \dint_y h_j[z_j z_k f(z)] \spc .
\end{eqnarray*}
Der Grenzwert ist also endlich. Ganz analog erh\"alt man
\[ \lim_{0>\lambda \rightarrow 0} g(\lambda) \;=\;
- \frac{1}{2} \dint_y dz \; (y^l \partial_l) (y_k \: z_j \: f(z))
	- \dint_y h_j[z_j z_k f(z)] \spc ,	\]
so da{\ss} sich die Randwerte bei $\lambda = \pm 0$ wegheben.

Wir erhalten also die Gleichung
\begin{eqnarray*}
\lefteqn{ \hint_{-\infty}^\infty \frac{d\lambda}{|\lambda|} \;
\frac{\partial}{\partial y^k} \xdint_{\lambda y} dz \; z_j \: f(z) } \\
&=& \frac{2}{y^2} \hint_{-\infty}^\infty \frac{d\lambda}{|\lambda|} \:
	\left\{ \xdint_{\lambda y} dz \;
(y-\frac{2}{\lambda} z)_k \: z_j \: f(z) - \dint_{\lambda y} dz \;
(y^l \partial_l)(z_j \: z_k \: f(z)) \right\}
\end{eqnarray*}
und somit insgesamt
\begin{eqnarray*}
\lefteqn{ \frac{\partial}{\partial y^k} \left\{ \frac{2}{y^2} \:
	\hint_{-\infty}^\infty \frac{d\lambda}{|\lambda|}
	\xdint_{\lambda y} dz \; z_j \: f(z) \right\}  } \\
&=& - \frac{4}{y^4} \: y_k \hint_{-\infty}^\infty \frac{d\lambda}{|\lambda|}
	\xdint_{\lambda y} dz \; z_j \: f(z) 
	+ \frac{2}{y^2} \hint_{-\infty}^\infty \frac{d\lambda}{|\lambda|} \:
	\frac{\partial}{\partial y^k} \xdint_{\lambda y} dz \;
	z_j \: f(z) \\
&=& - \frac{8}{y^4} \hint_{-\infty}^\infty \frac{d\lambda}{|\lambda|} \;
\left\{	\frac{1}{\lambda} \xdint_{\lambda y} dz \; z_j \: z_k \: f(z) +
\frac{1}{2} \dint_{\lambda y} dz \; (y^l \partial_l) (z_j \: z_k \: f(z))
	\right\} \spc ,
\end{eqnarray*}
was den Beweis abschlie{\ss}t.
\QED
F\"ur $y \in \I$ kann man die erhaltenen Formeln f\"ur die partiellen
Ableitungen von $(\sveeint f)(y)$ noch vereinfachen:
\begin{eqnarray*}
\partial_j (\veeint f)(y) &=& (\veeint h_j)(y) + \frac{2}{y^2}
	\hint_{-\infty}^\infty \frac{d\lambda}{|\lambda|}
	\xint_{\lambda y} dz \; z_j \: f(z) \\
&=& (\veeint h_j)(y) + \frac{2}{y^2} \xint_y dz \; z_j \int_{-\infty}^\infty
	d\lambda \; \epsilon(\lambda) \: f(\lambda z) \\
\partial_{jk} (\veeint f)(y) &=& (\veeint h_{jk})(y) - \frac{1}{y^2}
	\hint_{-\infty}^\infty \frac{d\lambda}{|\lambda|}
	\xint_{\lambda y} dz \; z_j \: z_k \; \Box f(z) \\
&& - \frac{8}{y^4} \hint_{-\infty}^\infty d\lambda \;
	\frac{\epsilon(\lambda)}{\lambda^2} \xint_{\lambda y} dz
	\; z_j \: z_k \: f(z)	\\
&=& (\veeint h_{jk})(y) - \frac{1}{y^2} \xint_y dz \; z_j \: z_k
	\int_{-\infty}^\infty d\lambda \; |\lambda| \; \Box f(\lambda z) \\
&& - \frac{8}{y^4} \xint_y dz \; z_j \: z_k \int_{-\infty}^\infty d\lambda \;
	\epsilon(\lambda) \: f(\lambda z) \spc .
\end{eqnarray*}
Die Ableitungen von $(\sveeint f)(y)$ lassen sich also
mit Lichtkegelintegralen geeigneter, von $y$ unabh\"angiger Funktionen
darstellen. Beachte, da{\ss} die dabei auftretenden Integranden
\begin{eqnarray}
\label{eq:a3_60}
z_j \int_{-\infty}^\infty d\lambda \; \epsilon(\lambda) \: f(\lambda z) \\
z_j \: z_j \int_{-\infty}^\infty d\lambda \; |\lambda| \;
	\Box f(\lambda z)
\end{eqnarray}
homogen vom Grade $0$ sind; die Funktion
\Equ{a3_62}
   z_j \: z_k \int_{-\infty}^\infty d\lambda \; \epsilon(\lambda) \:
	f(\lambda z)
\EndEqu
ist homogen vom Grade $1$.
Diese Homogenit\"atseigenschaften \"ubertragen sich auch auf die
Lichtkegelintegrale, was die sp\"ateren Rechnungen vereinfachen wird.

Im Fall $y \in \Ra$ tritt das Problem auf, da{\ss} in
$\sxint_y$ \"uber ein Hyperboloid integriert werden mu{\ss}, das
Integrationsgebiet also nicht mehr kompakt ist.
Insbesondere divergieren die Lichtkegelintegrale \"uber die Funktionen
\Ref{a3_60} bis~\Ref{a3_62}.
Wegen dieser zus\"atzlichen Schwierigkeit sind die allgemeinen Formeln
etwas komplizierter, man kann die Integrale nur im Sinne eines Hauptwertes
definieren.
\\[1em]
Es g\"unstig, die Notation zu vereinfachen:
\begin{Def}
Definiere f\"ur $f \in C^\infty_2(M)$ die Distribution $\sdotint f$ durch
\[  (\dotint f)(y) \;=\; (\xint f)(y) + (\dint f)(y) + \frac{1}{2}
	(\dint y^j \partial_j f)(y) \spc .  \]
\end{Def}
\begin{Satz}
$\sdotint f$ ist eine stetige Funktion auf $M$, f\"ur $z \in \Li$ gilt
\[  \lim_{y \rightarrow z} (\dotint f)(y) \;=\; \frac{\pi}{2} \int_0^1
f(\alpha z) \; d\alpha \spc . \]
\end{Satz}
{\Beweis}
Nach Satz~\ref{a3_satz1}, Satz~\ref{a3_satz5} ist $\sdotint f$ eine
Funktion.
F\"ur $0 \neq z \in \Li$ folgt nach Satz~\ref{a3_satz0} und Satz~\ref{a3_satz5}
\begin{eqnarray}
\lim_{\I \ni y \rightarrow z} (\dotint f)(y) &=& \frac{\pi}{2} \int_0^1
f(\alpha z) \; d\alpha \nonumber \\
\label{eq:a3_150}
\lim_{\Ra \ni y \rightarrow z} (\dotint f)(y) &=& \frac{\pi}{2} \int_0^1
f(\alpha z) \; d\alpha + \frac{\pi}{4} \int_{-\infty}^\infty z^j \partial_j
f(\alpha z) \; d\alpha \\
&=& \frac{\pi}{2} \int_0^1 f(\alpha z) \; d\alpha \nonumber \spc ,
\end{eqnarray}
denn der zweite Summand in~\Ref{a3_150} f\"allt durch partielle Integration
weg.

F\"ur $y \in \I$ hat man nach~\Ref{a3_51}
\[ \lim_{0 < \lambda \rightarrow 0} (\dotint f)(\lambda y) \;=\;
\frac{\pi}{2} \: f(0)  \]
und f\"ur $y \in \Ra$ nach~\Ref{a3_52}
\begin{eqnarray*}
  \lim_{0 < \lambda \rightarrow 0} (\dotint f)(\lambda y) &=&
\lim_{0<\lambda \rightarrow 0} (\xdint f)(\lambda y) + \frac{1}{2} (\dint y^j
\partial_j f)(y) \\
&=& \frac{\pi}{2} \: f(0) \spc ,
\end{eqnarray*}
so da{\ss} $\sdotint f$ auch am Ursprung stetig ist.
\QED
Mit dem Lichtkegelintegral $\sdotint$ kann man die Ergebnisse von
Satz~\ref{a3_satz13} und Satz~\ref{a3_satz14} in der Form
\begin{eqnarray}
\label{eq:a3_151}
\partial_j (\veeint f)(y) &=& (\veeint h_j)(y) + \frac{2}{y^2}
\hint_{-\infty}^\infty \frac{d\lambda}{|\lambda|} \xdint_{\lambda y} dz \;
	z_j \: f(z) \\
\partial_{jk} (\veeint f)(y) &=& (\veeint h_{jk})(y) - \frac{1}{y^2}
\hint_{-\infty}^\infty \frac{d\lambda}{|\lambda|} \xdint_{\lambda y} dz \;
	z_j \: z_k \; \Box f(z) \nonumber \\
\label{eq:a3_152}
&& - \frac{8}{y^4} \hint_{-\infty}^\infty d\lambda \;
\frac{\epsilon(\lambda)}{\lambda^2} \dotint_{\lambda y} dz \; z_j \: z_k \;
	f(z)
\end{eqnarray}
schreiben.

F\"ur die Anwendungen im n\"achsten Abschnitt m\"ussen wir noch asymptotische
N\"aherungsformeln f\"ur die einzelnen darin auftretenden
Terme ableiten. F\"ur diese Entwicklung um den Lichtkegel werden drei
Lemmata ben\"otigt:
\begin{Lemma}
\label{a3_lemma25}
Es gilt im Distributionssinne
\Equ{a3_250}
  \dint_y dz \; z_k \: (\Box f)(z) \;=\; \dint_y dz \; \left( z_k \: \frac{y^l
y^m}{y^2} \: \partial_{lm} f + \frac{y_k y^l}{y^2} \: \partial_l f -
	\partial_k f \right) \spc .
\EndEqu
\end{Lemma}
{\Beweis}
Man hat (im Distributionssinne)
\begin{eqnarray*}
\lefteqn{ \dint_y dz \; z_k \left( \Box - \frac{y^l y^m}{y^2} \:
\partial_{lm} \right) \: f(z) } \\
&=& - \frac{1}{2} \int d^4z \: l(z) \: \delta(\bra z,y \ket) \; z_k \:
\left( g^{lm} - \frac{y^l y^m}{y^2} \right) \: \partial_{lm} f \spc .
\end{eqnarray*}
Integriere nun partiell
\begin{eqnarray}
\label{eq:a3_25A}
&=& \frac{1}{2} \int d^4z \; m(z) \: \delta(\bra z,y \ket) \; 2 z_l \: z_k \;
\left( g^{lm} - \frac{y^l y^m}{y^2} \right) \partial_m f \\
\label{eq:a3_25B}
&& + \frac{1}{2} \int d^4z \; l(z) \: \delta^\prime (\bra z,y \ket) \; y_l \: z_k \;
\left( g^{lm} - \frac{y^l y^m}{y^2} \right) \partial_m f \\
\label{eq:a3_25C}
&& + \frac{1}{2} \int d^4z \; l(z) \: \delta(\bra z,y \ket) \;
\left( g^{km} - \frac{y^k y^m}{y^2} \right) \partial_m f \spc .
\end{eqnarray}
Der Summand~\Ref{a3_25B} verschwindet, \Ref{a3_25A} l\"a{\ss}t sich schreiben
als
\[  \int d^4z \; m(z) \: \delta(\bra z,y \ket) \; z_k \: z^m
	\partial_m f \spc . \]
Die Ableitung von $f$ verl\"auft in radialer Richtung und kann wie in
Lemma~\ref{a3_lemma3.9} partiell integriert werden, dabei f\"allt~\Ref{a3_25A}
weg. Man erh\"alt somit
\[ \dint_y dz \; z_k \left( \Box - \frac{y^l y^m}{y^2} \:
\partial_{lm} \right) \: f \;=\; - \dint_y \left( g^{km} - \frac{y^k
y^m}{y^2} \right) \: \partial_m f \spc . \]
\QED

\begin{Lemma}
\label{a3_lemma26}
Sei $f \in C^\infty_2(M)$ und $f(0)=0$. Dann gilt
\Equ{a3_25G}
  \frac{\partial}{\partial y^j} \hinti \frac{d\lambda}{|\lambda|}
  \xdint_{\lambda y} f \;=\; - \frac{1}{2} \hinti d\lambda
  \; \epsilon(\lambda)
  \xint_{\lambda y} dz \; z_j \: \Box f(z) \spc .
\EndEqu
Unter der zus\"atzlichen Voraussetzung $\partial_k f(0)=0$ hat man
\Equ{a3_25H}
  \frac{\partial}{\partial y^j} \hinti d\lambda \;
  \frac{\epsilon(\lambda)}{\lambda^2} \dotint_{\lambda y} f \;=\;
	\hinti \frac{d\lambda}{|\lambda|} \xdint_{\lambda y} dz \: \left(
	\frac{1}{2} \partial_j f(z) - \frac{1}{4} z_j \: \Box f(z) \right)
	\spc .
\EndEqu
\end{Lemma}
{\Beweis}
Man beachte zun\"achst, da{\ss} die Bedingungen $f(0)=0$, $\partial_k f(0)=0$
notwendig sind, damit die Integrale auf der linken Seite von~\Ref{a3_25G} und
\Ref{a3_25H} existieren.

Nach~\Ref{a3_181a} und Satz~\ref{a3_satz01} gilt
\begin{eqnarray*}
  \lefteqn{ \frac{\partial}{\partial y^j} \hinti \frac{d\lambda}{|\lambda|}
  \xdint_{\lambda y} f
\;=\; \hinti d\lambda \; \epsilon(\lambda) \xint_{\lambda y} h_j[f] } \\
&=& \hinti d\lambda \; \epsilon(\lambda) \xint_{\lambda y} dz \; \left(
	\partial_j f - \int_0^1 d\alpha \; \alpha \: \partial_j f(\alpha z) -
	\frac{1}{2} \: z_j \int_0^1 d\alpha \; \alpha^2 \: (\Box f)(\alpha z)
	\right) \\
&=& \hinti d\lambda \; \epsilon(\lambda) \left( \xint_{\lambda y} \partial_j f
	- \int_0^1 d\alpha \; \alpha \xint_{\alpha \lambda y} dz \;
	(\partial_j f + \frac{1}{2} \: z_j \: \Box f(z) ) \right) \\
&=& \hinti d\lambda \; \epsilon(\lambda) \xint_{\lambda y} \partial_j f
	- \hinti d\lambda \; \epsilon(\lambda)
	\int_0^1 d\alpha \; \alpha \xint_{\alpha \lambda y} dz \;
	(\partial_j f + \frac{1}{2} \: z_j \: \Box f(z) ) \\
&=& - \frac{1}{2} \hinti d\lambda \; \epsilon(\lambda) \xint_{\lambda
	y} dz \; z_j \: \Box f(z) \spc .
\end{eqnarray*}
Beim Beweis von~\Ref{a3_25H} tritt die Schwierigkeit auf, da{\ss} von einzelnen,
in Zwischenrechnungen auftretenden Integralen der Hauptwert
sowie das uneigentliche Integral nicht existieren.
Fasse daher im Folgenden das Symbol $\hinti$ auf als
$\int_{[-\varepsilon^{-1}, -\varepsilon] \cup [\varepsilon,
\varepsilon^{-1}]}$ mit festem $\varepsilon > 0$. Dadurch sind
alle Ausdr\"ucke wohldefiniert, in der Endformel kann dann der Grenz\"ubergang
$\varepsilon \rightarrow 0$ durchgef\"uhrt werden.
\begin{eqnarray*}
\frac{\partial}{\partial y^j} \hinti d\lambda \:
\frac{\epsilon(\lambda)}{\lambda^2} \dotint_{\lambda y} f
&=& \hinti d\lambda \: \frac{\epsilon(\lambda)}{\lambda^2} \:
\frac{\partial}{\partial y^j} \xdint_{\lambda y} f + \frac{1}{2} \hinti
\frac{d\lambda}{|\lambda|} \frac{\partial}{\partial y^j} \dint_{\lambda y}
y^l \partial_l f
\end{eqnarray*}
Wende nun Satz~\ref{a3_satz01} und Lemma~\ref{a3_lemma1a} an
\begin{eqnarray*}
&=& \hinti \frac{d\lambda}{|\lambda|} \xint_{\lambda y} dz \; \left(
	\partial_j f(z) - \int_0^1 d\alpha \: \alpha \: \partial_j f(\alpha
	z) - \frac{1}{2} \: z_j \: \int_0^1 d\alpha \: \alpha^2 \: (\Box
	f)(\alpha z) \right)  \\
&& + \frac{1}{2} \hinti \frac{d\lambda}{|\lambda|} \dint_{\lambda y} dz
	\left( \partial_j f - z_j \: \frac{y^l y^m}{y^2} \: \partial_{lm}f -
	\frac{y_j y^l}{y^2} \: \partial_l f \right) \spc .
\end{eqnarray*}
Lemma~\ref{a3_lemma25} liefert
\begin{eqnarray}
\label{eq:a3_25K}
&=& \hinti \frac{d\lambda}{|\lambda|} \xint_{\lambda y} \partial_j f -
	\hinti \frac{d\lambda}{|\lambda|} \int_0^1 d\alpha \; \alpha
	\xint_{\alpha \lambda y} dz \: \left( \partial_j f(z) + \frac{1}{2}
	\: z_j \: \Box f(z) \right) \\
&& - \frac{1}{2} \hinti \frac{d\lambda}{|\lambda|} \dint_{\lambda y} dz \;
	z_j \; (\Box f)(z) \spc .
\end{eqnarray}
Beachte, da{\ss} an dieser Stelle im zweiten Summanden von~\Ref{a3_25K} 
die Variablentransformation $\lambda^\prime = \alpha \lambda$ mit
anschlie{\ss}ende Ausintegration von $\alpha$ nicht durchgef\"uhrt werden
darf. Dabei w\"urde sich n\"amlich der oben eingef\"uhrte Parameter
$\varepsilon$ im Hauptwertintegral ver\"andern.
Wir m\"ussen daher anders vorgehen:
\begin{eqnarray}
\label{eq:a3_25D}
&=& \hinti \frac{d\lambda}{|\lambda|} \xdint_{\lambda y} \partial_j f -
	\hinti \frac{d\lambda}{|\lambda|} \int_0^1 d\alpha \: \alpha
	\xdint_{\alpha \lambda y} dz \; \left( \partial_j f(z) + \frac{1}{2}
	\: z_j \: (\Box f)(z) \right) \\
\label{eq:a3_25E}
&& - \hinti \frac{d\lambda}{|\lambda|} \dint_{\lambda y} dz \left( \partial_j
f(z) + \frac{1}{2} \: z_j \; (\Box f)(z) \right) \\
\label{eq:a3_25F}
&& + \hinti \frac{d\lambda}{|\lambda|} \int_0^1 d\alpha \: \alpha
\dint_{\alpha \lambda y} dz \left( \partial_j
f(z) + \frac{1}{2} \: z_j \; (\Box f)(z) \right) \spc .
\end{eqnarray}
In~\Ref{a3_25F} kann man unter Verwendung von~\Ref{a3_181a} die Integration
\"uber $\alpha$ ausf\"uhren, \Ref{a3_25E} und \Ref{a3_25F} heben sich dabei weg.

In den beiden Summanden von~\Ref{a3_25D} ist der Hauptwert wegen~\Ref{a3_51},
\Ref{a3_52} und der Voraussetzung $f(0) = \partial_k f(0) =0$ jeweils
sinnvoll definiert\footnote{Die Lichtkegelintegrale $\sxdint_{\lambda y}$
verhalten sich f\"ur kleine $\lambda$ jeweils  $\sim \; {\mbox{const }}
\epsilon(\lambda) + {\cal{O}}(\lambda)$.}.
Daher k\"onnen wir jetzt im zweiten Summanden die Variablentransformation
durchf\"uhren und die Integration \"uber $\alpha$ ausf\"uhren:
\begin{eqnarray*}
&=& \hinti \frac{d\lambda}{|\lambda|} \xdint_{\lambda y} dz \: \left(
	\frac{1}{2} \partial_j f - \frac{1}{4} z_j \: \Box f(z) \right)
\end{eqnarray*}
\QED
\begin{Lemma}
\label{a3_lemma27}
Sei $0 \neq \bar{y} \in \Li$, $f \in C^\infty_2(M)$. Bei~\Ref{a3_25a} sei
die zus\"atzliche Voraussetzung $f(0)=\partial_k f(0)=0$, bei~\Ref{a3_25b} die
Bedingung $f(0)=0$ erf\"ullt. Dann gilt
\begin{eqnarray}
\label{eq:a3_25a}
\lim_{y \rightarrow \bar{y}} \hinti d\lambda \;
	\frac{\epsilon(\lambda)}{\lambda^2} \dotint_{\lambda y} f &=&
	\frac{\pi}{4} \hinti d\lambda \; \frac{\epsilon(\lambda)}{\lambda^2}
	\: f(\lambda \bar{y})  \\
\label{eq:a3_25b}
\lim_{y \rightarrow \bar{y}} \hinti \frac{d\lambda}{|\lambda|}
	\xdint_{\lambda y} f &=&
	\frac{\pi}{2} \hinti \frac{d\lambda}{|\lambda|}
	\: f(\lambda \bar{y})
\end{eqnarray}
und
\begin{eqnarray}
\label{eq:a3_25c}
\lim_{y \rightarrow \bar{y}} \left\{ \hinti d\lambda \;
	\epsilon(\lambda) \xint_{\lambda y} f
	+ \frac{\pi}{2} \ln(|y^2|)
	\: \hint_{-\infty}^\infty d\lambda \: \epsilon(\lambda) \;
	f(\lambda y) \right\} \;<\; \infty  \spc .
\end{eqnarray}
\end{Lemma}
{\Beweis}
Die Bedingungen $f(0)=0$, $\partial_k f(0)=0$ sind wiederum notwendig, damit
die Integrale existieren.

Aufgrund der in $\lambda$, $y$ gleichm\"a{\ss}igen Absch\"atzung
\[    \dotint_{\lambda y} f \leq {\mbox{const}}(f)  \]
und
\[  \lim_{y \rightarrow \bar{y}} \dotint_{\lambda y} f \;=\; \frac{\pi}{2}
\int_0^1 f(\alpha \lambda y) \: d\alpha	\spc \forall \lambda \neq 0  \]
folgt nach Lebesgues dominiertem Konvergenzsatz
\[  \lim_{y \rightarrow \bar{y}} \int_{\sR \setminus [-\varepsilon,
\varepsilon]} d\lambda \; \frac{\epsilon(\lambda)}{\lambda^2}
\dotint_{\lambda y} f \;=\; \frac{\pi}{2} \int_{\sR \setminus [-\varepsilon,
\varepsilon]} d\lambda \; \frac{\epsilon(\lambda)}{\lambda^2} \int_0^1
d\alpha \; f(\alpha \lambda \bar{y}) \spc .  \]
Im Grenzfall $\varepsilon \rightarrow 0$ erh\"alt man
\begin{eqnarray*}
\lim_{y \rightarrow \bar{y}} \hinti d\lambda \;
\frac{\epsilon(\lambda)}{\lambda^2} \dotint_{\lambda y} f &=& \frac{\pi}{2}
	\hinti d\lambda \; \frac{\epsilon(\lambda)}{\lambda^2} \int_0^1
	d\alpha \; f(\alpha \lambda \bar{y}) \\
&=& \frac{\pi}{2}
	\hinti d\lambda \; \frac{\epsilon(\lambda)}{\lambda^2} \:
	f(\lambda \bar{y}) \: \int_0^1 d\alpha \; \alpha \\
&=& \frac{\pi}{4} \hinti d\lambda \; \frac{d\lambda}{\lambda^2} \: f(\lambda
	\bar{y}) \spc .
\end{eqnarray*}
Gleichung~\Ref{a3_25b} folgt ganz analog.
Um~\Ref{a3_25c} abzuleiten, formt man folgenderma{\ss}en um:
\begin{eqnarray*}
\hinti d\lambda \; \epsilon(\lambda) \xint_{\lambda y} f &=&
	- \hinti d\lambda \; \frac{|\lambda|}{1-\lambda} \xint_{\lambda y} f +
	\hinti d\lambda \; \frac{\epsilon(\lambda)}{1-\lambda} \xint_{\lambda
	y} f \\
&=& - \veeint^y f +
	\hinti d\lambda \; \frac{\epsilon(\lambda)}{1-\lambda} \xint_{\lambda
	y} f
\end{eqnarray*}
Nach dominierter Konvergenz besitzt der zweite Summand f\"ur $y \rightarrow
\bar{y}$ einen endlichen Grenzwert. Das logarithmische Divergenzverhalten
auf dem Lichtkegel erh\"alt man nun aus Satz~\ref{a3_satz99}.
\QED
Jetzt k\"onnen wir die Entwicklung um den Lichtkegel durchf\"uhren:
\begin{Satz}
\label{a3_satz98}
F\"ur $f \in C^\infty_2(M)$, $f(0)=0$ gilt
\begin{eqnarray}
\hinti \frac{d\lambda}{|\lambda|} \xdint_{\lambda y} f &=& \frac{\pi}{2}
	\hinti \frac{d\lambda}{|\lambda|} \: f(\lambda y) \nonumber \\
&&+ \frac{\pi}{8} \: y^2 \: \ln (|y^2|) \; \hinti d\lambda \; |\lambda| \;
	(\Box f)(\lambda y)  \nonumber \\
\label{eq:a3_25e}
&&+ {\cal{O}}(y^2)  \spc .
\end{eqnarray}
Unter der zus\"atzlichen Bedingung $\partial_k f(0)=0$ hat man
\begin{eqnarray}
\hinti d\lambda \; \frac{\epsilon(\lambda)}{\lambda^2} \dotint_{\lambda y} f
&=& \frac{\pi}{4} \hinti d\lambda \; \frac{\epsilon(\lambda)}{\lambda^2}
	\: f(\lambda y)  \nonumber \\
&&- \frac{\pi}{16} \: y^2 \hinti d\lambda \; \epsilon(\lambda) \; (\Box
	f)(\lambda y) \nonumber \\
&&- \frac{\pi}{128} \: y^4 \: \ln(|y^2|) \; \hinti d\lambda \;
	\epsilon(\lambda) \: \lambda^2 \; (\Box^2 f)(\lambda y) \nonumber \\
\label{eq:a3_25f}
&&+ {\cal{O}}(y^4) \spc .
\end{eqnarray}
\end{Satz}
{\Beweis}
Man mu{\ss} zeigen, da{\ss}\footnote{Vergleiche auch Lemma~\ref{a2_lemma2}.}
\begin{description}
\item[a)] die Randwerte auf $\Li$ der Ableitungen nullter und erster Ordnung
der linken und rechten Seite von~\Ref{a3_25e} \"ubereinstimmen.
\item[b)] die Randwerte auf $\Li$ der Ableitungen nullter bis zweiter Ordnung
der linken und rechten Seite von~\Ref{a3_25f} \"ubereinstimmen.
\end{description}
Die partiellen Ableitungen kann man mit Lemma~\ref{a3_lemma26}, die Randwerte
auf dem Lichtkegel mit Lemma~\ref{a3_lemma27} direkt berechnen.
\QED

\section{St\"orungsrechnung f\"ur das elektromagnetische Feld}
\label{elek_p0}

Betrachte nun die St\"orung~\Ref{a1_54a} durch ein elektromagnetisches
Feld. In erster Ordnung St\"orungstheorie hat man
\[  \tilde{p}_0 \;=\; p_0 + \Delta p_0 \]
mit
\Equ{a3_24g}
  \Delta p_0 (x,y) \;=\; -e (p_0 \: \Aslsh \: s_0 \; + s_0 \: \Aslsh \:
p_0)(x,y) \spc .
\EndEqu
Wir werden im folgenden der Einfachheit halber annehmen, da{\ss} $A$
gleichm\"a{\ss}ig zur vierten Potenz abf\"allt.

Zun\"achst wollen wir eine Formel f\"ur $(s_0 \: \Aslsh \: p_0)(x,y)$ ableiten.
Zur besseren \"Ubersichtlichkeit werden wir den Fall $x=0$ betrachten,
den allgemeinen Fall erh\"alt man daraus durch eine
Parallelverschiebung.
\begin{Lemma}
\label{a3_lemma21}
Es gilt:
\begin{eqnarray}
\label{eq:a3_161}
(s_0 \: \Aslsh \: p_0)(0,y) &=& \frac{1}{\pi^4} \frac{1}{y^4}
\hint_{-\infty}^\infty d\lambda \: \frac{\epsilon(\lambda)}{\lambda^2}
	\dotint_{\lambda y} dz \; z\slsh \: z^k \: A_k(z) \\
\label{eq:a3_162}
&& - \frac{1}{8 \pi^4} \frac{1}{y^2} \hint_{-\infty}^\infty
\frac{d\lambda}{|\lambda|} \xdint_{\lambda y} dz \; z\slsh \: z_k \: j^k(z)
\\
\label{eq:a3_163}
&& + \frac{1}{16 \pi^4} \frac{1}{y^2} \hint_{-\infty}^\infty
\frac{d\lambda}{|\lambda|} \xdint_{\lambda y} dz \; \gamma^i \gamma^j \:
	F_{ij}(z) \; z\slsh \\
\label{eq:a3_164}
&& - \frac{1}{8 \pi^4} \veeint^y \left( \gamma^k h_{jk}[A^j] - \frac{1}{2}
\gamma^i \gamma^j \gamma^k \: h_k[\partial_i A_j] \right) 
\end{eqnarray}
\end{Lemma}
{\Beweis}
Unter Verwendung von \Ref{a1_76} und \Ref{a3_170a} kann man $s_0 \Aslsh
p_0$ mit dem Lichtkegelintegral $\sveeint$ ausdr\"ucken
\begin{eqnarray}
\lefteqn{ (s_0 \: \Aslsh \: p_0)(x,y) \;=\; - \Pdd_x (S_0 \Aslsh P_0)(x,y)
\Pdd_y } \nonumber \\
&=& - \frac{1}{16 \pi^4} \Pdd_x (\veeint_x^y \Aslsh) \Pdd_y \nonumber \\
&=& \frac{1}{16 \pi^4} \gamma^i \gamma^j \gamma^k \: \frac{\partial}{\partial
x^i} \frac{\partial}{\partial y^k} \veeint_x^y A_j \nonumber \\
\label{eq:a3_O}
&=& \frac{1}{16 \pi^4} \gamma^i \gamma^j \gamma^k
\left\{ \frac{\partial}{\partial y^k} \veeint_x^y \partial_i A_j
- \frac{\partial^2}{\partial y^i \partial y^k} \veeint_x^y A_j \right\} \\
&=& \frac{1}{16 \pi^4} \; \gamma^i \: \gamma^j \: \gamma^k
	\: \frac{\partial}{\partial
	y^k} \veeint_x^y \partial_i A_j - \frac{1}{8 \pi^4} \gamma^k \:
	\frac{\partial^2}{\partial y^j \partial y^k} \veeint_x^y A^j
	\nonumber \\
&& + \frac{1}{16 \pi^4} \; \Box_y \veeint_x^y \Aslsh \spc . \nonumber 
\end{eqnarray}
Setze nun $x=0$ und
wende Satz~\ref{a3_satz99} sowie \Ref{a3_151}, \Ref{a3_152} an:
\begin{eqnarray*}
(p_0 \: \Aslsh \: s_0)(0,y) &=& - \frac{1}{8 \pi^4} \veeint_x^y
\left( \gamma^k \:
h_{jk}[A^j] - \frac{1}{2} \gamma^i \gamma^j \gamma^k \: h_k[\partial_i A_j]
\right) \\
&& + \frac{1}{8 \pi^4} \frac{1}{y^2} \hint_{-\infty}^\infty
\frac{d\lambda}{|\lambda|} \xdint_{\lambda y} dz \; \gamma^i \gamma^j \;
	\partial_i A_j(z) \; z\slsh \\
&& + \frac{1}{8 \pi^4} \frac{1}{y^2} \hint_{-\infty}^\infty
\frac{d\lambda}{|\lambda|} \xdint_{\lambda y} dz \; z\slsh \: z_k \; \Box
A^k(z) \\
&& + \frac{1}{\pi^4} \frac{1}{y^4} \hint_{-\infty}^\infty d\lambda \;
\frac{\epsilon(\lambda)}{\lambda^2} \dotint_{\lambda y} dz \; z\slsh \: z_k \:
	A^k(z)
\end{eqnarray*}
F\"uhrt man nun die Ersetzung $\Box A_k = - j_k + \partial_{kl} A^l$ durch
und setzt die Umformung
\begin{eqnarray*}
\hint_{-\infty}^\infty \frac{d\lambda}{|\lambda|}
	\xdint_{\lambda y} dz \; z\slsh \: z^k \partial_{jk} A^j(z)
&=& \hint_{-\infty}^\infty \frac{d\lambda}{|\lambda|} \: \lambda^2
\frac{d}{d\lambda} \frac{1}{\lambda} \xdint_{\lambda y} dz \; z\slsh \;
	\partial_j A^j(z) \\
&=& - \hint_{-\infty}^\infty \frac{d\lambda}{|\lambda|} \xdint_{\lambda y}
	dz \; z\slsh \; \partial_j A^j(z)
\end{eqnarray*}
ein, ergibt sich die Behauptung.
\QED
Beachte, da{\ss} man~\Ref{a3_164} durch Einsetzen von~\Ref{a3_80a} mit Hilfe
des Feldst\"arketensors und dessen partiellen Ableitungen eichinvariant
schreiben kann.

Wir wollen nun \"uberpr\"ufen, da{\ss} die in Lemma~\ref{a3_lemma21}
abgeleitete Gleichung
das richtige Verhalten unter Eichtransformationen zeigt. F\"ur ein
elektromagnetisches Potential der Form
\[  e \Aslsh \;=\; \Pdd \Lambda  \spc,\;\;\; \Lambda \in C^\infty_3(M)  \]
fallen die Summanden~\Ref{a3_162} bis~\Ref{a3_164} weg. Man erh\"alt
\begin{eqnarray*}
(s_0 \: \Aslsh \: p_0)(0,y) &=& \frac{1}{e \pi^4} \frac{1}{y^4}
	\hint_{-\infty}^\infty d\lambda \; \frac{\epsilon(\lambda)}{\lambda^2}
	\dotint_{\lambda y} dz \; z\slsh \; z^k \: \partial_k \Lambda(z) \\
&=& \frac{1}{e \pi^4} \frac{1}{y^4}
	\hint_{-\infty}^\infty d\lambda \; \frac{\epsilon(\lambda)}{\lambda^2}
	\: \lambda^2 \frac{d}{d \lambda} \frac{1}{\lambda}
	\dotint_{\lambda y} dz \; z\slsh \; \Lambda(z) \\
&=& - \frac{2}{e \pi^4} \frac{1}{y^4} \; \lim_{\lambda \rightarrow 0}
	\frac{1}{\lambda} \dotint_{\lambda y} dz \; z\slsh \: \Lambda(z)
	\spc ,
\end{eqnarray*}
und unter Verwendung der Regeln von de l'Hopital\footnote{Beachte wegen
der Vorzeichen, da{\ss}
\[  p_0(0,y) \;=\; (i \Pdd \: P_0)(0,y) \;=\; - i \Pdd_y \: P_0(0,y)
	\;=\; - \frac{2i}{4 \pi^3} \frac{y \slsh}{y^4} \spc .  \]	}
\begin{eqnarray*}
&=& - \frac{2}{e \pi^4} \frac{1}{y^4} \; \frac{\pi}{4} \: y\slsh \; \Lambda(0)
\;=\; \frac{1}{ie} \: \Lambda(0) \; p_0(0,y) \spc .
\end{eqnarray*}
Dies ist das richtige Resultat, wie man an der Rechnung
\begin{eqnarray*}
\frac{1}{ie} \; s_0 \: [i \Pdd, \Lambda] \: p_0
&=& \frac{1}{ie} \: (i \Pdd \: s_0) \: \Lambda \: p_0 - \frac{1}{ie} \: s_0
	\Lambda \: (i \Pdd \; p_0) \;=\; \frac{1}{ie} \: \Lambda \: p_0
\end{eqnarray*}
sieht.

Der Summand~\Ref{a3_161} ist also f\"ur das richtige Transformationsverhalten
bei Eich\-trans\-for\-ma\-tio\-nen verantwortlich.
Man beachte jedoch folgende Schwierigkeit:
Im Gegensatz zu den Eichtermen, die in den Anh\"angen~\ref{anh1} und~\ref{anh2}
auftreten, f\"uhrt~\Ref{a3_161} bei allgemeinen St\"orungen nicht nur zu einer
Phasenverschiebung von $p_0(0,y)$.
In einem eichinvarianten Funktional, das $\tilde{p}_0$ enth\"alt, f\"allt
daher der Beitrag von~\Ref{a3_161} im allgemeinen nicht weg, wir m\"ussen
ihn noch berechnen.

Dazu f\"uhren wir eine Entwicklung um den Lichtkegel durch:
\begin{Satz}
\label{a3_satz97}
Es gilt:
\begin{eqnarray*}
\lefteqn{ (s_0 \: \Aslsh \: p_0)(0,y) } \\
&=& \frac{1}{4 \pi^3} \frac{y \slsh}{y^4} \; \hinti d\lambda \;
\epsilon(\lambda) \: A_k(\lambda y) \; y^k \\
&&+ \frac{1}{16 \pi^3} \frac{y \slsh}{y^2} \; \hinti d\lambda \;
\epsilon(\lambda) \: (\lambda^2 - \lambda) \; j_k \: y^k \\
&&- \frac{1}{16 \pi^3} \frac{1}{y^2} \; \hinti d\lambda \; \epsilon(\lambda)
\: (2 \lambda-1) \: F_{jk} \: \gamma^j \: y^k \\
&&- \frac{i}{32 \pi^3} \frac{1}{y^2} \; \hinti d\lambda \; \epsilon(\lambda)
\: \varepsilon^{ijkl} \: F_{ij} \: y_k \; \rho \gamma_l \\
&&+ \frac{1}{128 \pi^3} \: \ln (|y^2|) \; \hinti d\lambda \; \epsilon(\lambda)
\: (\lambda^4 - 2 \lambda^3 + \lambda^2) \; \Box j_k \: y^k \; y \slsh \\
&&- \frac{1}{128 \pi^3} \: \ln(|y^2|) \; \hinti d\lambda \; \epsilon(\lambda)
\: (4 \lambda^3 - 6 \lambda^2 + 2 \lambda) \; \Box F_{jk} \; \gamma^j \: y^k
\\
&& - \frac{i}{128 \pi^3} \: \ln(|y^2|) \; \hinti d\lambda \; \epsilon(\lambda)
\: (\lambda^2 - \lambda) \; \varepsilon^{ijkl} \; \Box F_{ij} \: y_k \;
\rho \gamma_l \\
&& - \frac{1}{16 \pi^3} \: \ln(|y^2|) \; \hinti d\lambda \; \epsilon(\lambda)
\: (\lambda^2 - \lambda) \; \gamma^k \: j_k \\
&& + {\cal{O}}(y^0) \spc .
\end{eqnarray*}
\end{Satz}
{\Beweis}
Die asymptotischen Entwicklungen von Satz~\ref{a3_satz98} liefern:
\begin{eqnarray*}
\lefteqn{\frac{1}{\pi^4} \frac{1}{y^4} \hinti d\lambda \;
\frac{\epsilon(\lambda)}{\lambda^2} \dotint_{\lambda y} dz \; z\slsh \: z^k
\: A_k(z) }\\
\;&=& \frac{1}{4 \pi^3} \frac{y\slsh}{y^4} \hinti d\lambda \; \epsilon(\lambda)
	\: A_k \: y^k \\
&&- \frac{1}{16 \pi^3} \frac{1}{y^2} \hinti d\lambda \; \epsilon(\lambda)
	\left(-\lambda^2 \: y\slsh \: y^k j_k + 2 \lambda F_{jk} \:
	\gamma^j y^k \right) \\
&&- \frac{1}{128 \pi^3} \: \ln(|y^2|) \hinti d\lambda \; \epsilon(\lambda)
	\left( 2 \lambda^2 \: \gamma^k j_k - 2 \lambda^3 y^k \: \Pdd j_k
	\right. \\
&& \spc \left. + 2 \lambda^3 \; \Box F_{jk} \; \gamma^j y^k - \lambda^4 \:
	y\slsh \: y^k \; \Box j_k \right) \\
&&+ {\cal{O}}(y^0) \\
\lefteqn{- \frac{1}{8 \pi^4} \frac{1}{y^2} \hinti \frac{d\lambda}{|\lambda|}
	\xdint_{\lambda y} dz \; z\slsh \; z_k \: j^k(z) } \\
&=& - \frac{1}{16 \pi^3}\frac{y \slsh}{y^2} \hinti d\lambda \; |\lambda| \:
	j_k \: y^k \\
&&- \frac{1}{64 \pi^3} \: \ln(|y^2|) \hinti d\lambda \; \epsilon(\lambda)
	\left( 2 \lambda \: \gamma^k j_k + 2 \lambda^2 \: y^k \: \Pdd j_k +
	\lambda^3 \: y \slsh \: y^k \; \Box j_k \right) \\
&&+ {\cal{O}}(y^0) \\
\lefteqn{\frac{1}{16 \pi^4} \frac{1}{y^2} \hinti \frac{d\lambda}{|\lambda|}
	\xdint_{\lambda y} dz \; \gamma^i \gamma^j \: F_{ij}(z) \; z\slsh }\\
&=& \frac{1}{32 \pi^3} \frac{1}{y^2} \hinti d\lambda \; \epsilon(\lambda) \;
	\gamma^i \gamma^j \: F_{ij} \; y\slsh \\
&&+ \frac{1}{128 \pi^3} \: \ln(|y^2|) \hinti d\lambda \; \epsilon(\lambda)
	\left( \lambda^2 \: \gamma^i \gamma^j \; \Box F_{ij} \; y \slsh +
	4 \lambda \: \gamma^k j_k \right)
	\;+\; {\cal{O}}(y^0) \\
\lefteqn{- \frac{1}{8 \pi^4} \veeint^y \left(\gamma^k h_{jk}[A^j] -\frac{1}{2}
\gamma^i \gamma^j \gamma^k \: h_k[\partial_i A_j] \right) } \\
&=& - \frac{1}{128 \pi^3} \: \ln(|y^2|) \hinti d\lambda \; \epsilon(\lambda)
	\left( 2 \lambda \; \Box F_{jk} \; \gamma^j \: y^k \right. \\
&& \spc\;\;\; \left. - \lambda^2 \; \Box j_k \; y^k \: y\slsh - i \lambda \:
	\varepsilon^{ijkl} \; \Box F_{ij} \: y_k \; \rho \gamma_l 
	\right)
	\;+\; {\cal{O}}(y^0)
\end{eqnarray*}
Die Integranden sind dabei jeweils an der Stelle $\lambda y$ auszuwerten.

St\"orend sind die Terme der Form $y^k \; \Pdd j_k$. Beachte dazu, da{\ss}
\begin{eqnarray*}
y^k \: \Pdd j_k &=& y^k \: \gamma^j \: (\partial_{jkl} A^l - \Box \partial_j
	A_k)  \\
\Box F_{jk} \: \gamma^j y^k &=& y^k \: \gamma^j \: (\Box \partial_j A_k -
	\Box \partial_k A_j) \spc ,
\end{eqnarray*}
also
\[   y^k \: \Pdd j_k \;=\; - \Box F_{jk} \; \gamma^j y^k + \gamma^l \:
	y^k \partial_k j_l \spc .	\]
Setzt man diese Beziehung in die asymptotischen Formeln ein, ersetzt die
radiale Ableitung $y^j \partial_k j_l$ durch die Ableitung nach $\lambda$
und integriert partiell, so folgt die Behauptung.
\QED
Nun k\"onnen wir nach~\Ref{a3_24g} auch direkt die Formel f\"ur $\tilde{p}_0$
angeben:
\begin{Thm}
\label{a3_theorem1}
In erster Ordnung St\"orungstheorie gilt
\begin{eqnarray}
\label{eq:a3_111b}
\Delta p_0(x,y) &=& - i e \left( \int_x^y A_j \right) \: \xi^j \; p_0(x,y) \\
\label{eq:a3_112a}
   && - \frac{e}{8 \pi^3} \left( \int_x^y (\alpha^2-\alpha) \; \xi \slsh
	\: \xi^k \: j_k \right) \;\; \frac{1}{\xi^2} \\
\label{eq:a3_113a}
   && + \frac{e}{8 \pi^3} \left( \int_x^y (2 \alpha - 1) \; \xi^j \:
	\gamma^k \: F_{kj} \right) \;\; \frac{1}{\xi^2} \\
\label{eq:a3_114a}
   && + \frac{i e}{16 \pi^3} \left( \int_x^y \varepsilon^{ijkl} \; F_{ij} \:
	\xi_k \; \rho \gamma_l \right) \;\; \frac{1}{\xi^2} \\
\label{eq:a3_115a}
  && - \frac{e}{64 \pi^3} \: \int_x^y (\alpha^4 - 2 \alpha^3
	+ \alpha^2) \: \xi \slsh \: \xi_k \; \Box j^k \;\; \ln(|\xi^2|) \\
 \label{eq:a3_116a}
 && + \frac{e}{64 \pi^3} \: \int_x^y (4 \alpha^3
	- 6 \alpha^2 + 2 \alpha) \: \xi^j \: \gamma^k \: (\Box F_{kj})
	\;\; \ln(|\xi^2|) \\
  \label{eq:a3_117a}
&& + \frac{i e}{64 \pi^3} \: \int_x^y (\alpha^2 - \alpha) \;
	\varepsilon^{ijkl} \: (\Box F_{ij}) \: \xi_k \; \rho \gamma_l
	\;\; \ln(|\xi^2|) \\
  \label{eq:a3_118a}
&& + \frac{e}{8 \pi^3} \: \int_x^y (\alpha^2 - \alpha)
	\: \gamma^k \: j_k \;\; \ln(|\xi^2|) \\
  && + {\cal{O}}(\xi^0) \spc . \nonumber
\end{eqnarray}
\end{Thm}
{\Beweis}
Folgt direkt aus Satz~\ref{a3_satz97} durch Translation und Vertauschung
von $x$ und $y$.
\QED
Auffallend ist die \"Ahnlichkeit zu Theorem~\ref{a1_theorem1} und
Satz~\ref{a1_randwert}. Ein Beitrag $\sim l(\xi)$ in $\tilde{k}_0$ entspricht
in $\tilde{p}_0$ einem Pol $\sim \xi^{-2}$.
An Stelle der in der N\"ahe des Lichtkegels beschr\"ankten Terme in
$\tilde{k}_0$ treten nun logarithmisch divergente Beitr\"age auf.
\\[1em]
Beachte, da{\ss} wir bei der Entwicklung von $\Delta p_0$ in
Theorem~\ref{a3_theorem1} die Terme der Ordnung $\xi^0$ nicht berechnet haben.
Es zeigt sich, da{\ss} sie nicht mit einem Linienintegral ausgedr\"uckt werden
k\"onnen, sondern da{\ss} dabei auch die St\"orung au{\ss}erhalb der
Verbindungsgeraden durch $x$ und $y$ beitr\"agt.
Daher sind diese Terme schwieriger zu behandeln.
Eine Entwicklung um den Lichtkegel ist nicht ausreichend,
wir m\"ussen geeignete Umformungen f\"ur die Lichtkegelintegrale finden.
Unser Ziel ist ein einfacher Ausdruck f\"ur $\frac{1}{2} \{ s_0 \Aslsh
p_0(0,y), y\slsh \}$. Zun\"achst wird ein Lemma ben\"otigt:

\begin{Lemma}
\label{a3_lemma4}
Es gilt
\begin{eqnarray}
\label{eq:a3_M}
y^j \frac{\partial}{\partial y^k} \xint_{\lambda y} F_j^{\;\:k} &=&
	\xint_{\lambda y} dz \; z^k j_k(z) \\
\label{eq:a3_N}
y^j \frac{\partial}{\partial y^k} \veeint^{y} F_j^{\;\:k} &=&
	\veeint^{y} dz \; z^k j_k(z) \spc .
\end{eqnarray}
\end{Lemma}
{\Beweis}
Nach Lemma \ref{a3_lemma1a} hat man
\begin{eqnarray}
\frac{\partial}{\partial y^j} \dint_{\lambda y} dz \; z^k \:
	F^j_{\;\:k}(z)
&=& \dint_{\lambda y} dz \; \left( -\frac{z^j}{y^2} \: y^l \partial_l \:
	(z^k \: F_{jk}(z)) - \frac{y^j}{y^2} \: z^k \: F_{jk}(z) \right)
	\nonumber \\
\label{eq:a3_G}
&=& - \frac{1}{y^2} \dint_{\lambda y} dz \; \left( z^j y^k + y^j z^k
	\right) F_{jk}(z) \;=\; 0 \spc .
\end{eqnarray}
Weiterhin gilt
\begin{eqnarray}
y^j \frac{\partial}{\partial y^k} \xint_{\lambda y} F_j^{\;\:k}
&=& \lambda \xint_{\lambda y} dz \; \left( 2 \: \frac{y^j \: (\lambda
	y-2z)^k}{\lambda^2 y^2} + 2 \: \frac{y^j \: (\lambda y -z)^k}{\lambda^2
	y^2} \: z^l \partial_l \right) \; F_{jk}(z) \nonumber \\
&=& -\frac{4}{\lambda y^2} \xint_{\lambda y} dz \; y^j z^k \: F_{jk}(z)
	- \frac{2}{\lambda y^2} \xint_{\lambda y} dz \; y^j z^k \;
	z^l \partial_l F_{jk}(z) \nonumber \\
&=& -\frac{2}{\lambda y^2} \xint_{\lambda y} dz \; y^j  \: (1+z^l \partial_l)
	(z^k \: F_{jk}(z)) \nonumber \\
\label{eq:a3_H}
&=& -\frac{1}{\lambda} \frac{\partial}{\partial y^j} \xint_{\lambda y}
	dz \; z^k \: F^j_{\;\:k}(z) \spc ,
\end{eqnarray}
wobei in der ersten und letzten Zeile jeweils Satz~\ref{a3_satz2} mit
$h_j$ gem\"a{\ss} \Ref{a3_14} angewendet wurde.
Durch Kombination von \Ref{a3_G} und \Ref{a3_H} erh\"alt man nach
Satz~\ref{a3_satz01}
\begin{eqnarray*}
y^j \frac{\partial}{\partial y^k} \xint_{\lambda y} F_j^{\;\:k} &=&
	-\frac{1}{\lambda} \frac{\partial}{\partial y^j} \xdint_{\lambda y}
	dz \; z^k \; F^j_{\;\:k}(z) \\
&=& -\xint_{\lambda y} dz \; \left( \partial_j(z^k \: F^j_{\;\:k}(z)) -
	\frac{1}{2} \Box_z \left( z^j \: z^k \int_0^1 \alpha \: F_{jk}(\alpha
	z) \; d\alpha \right) \right) \\
&=& \xint_{\lambda y} dz \; z^k \: j_k(z) \spc .
\end{eqnarray*}
Gleichung~\Ref{a3_N} folgt hieraus mit Hilfe von Lemma~\ref{a3_lemma11}:
\begin{eqnarray*}
y^j \frac{\partial}{\partial y^k} \veeint^y F_j^{\;\:k} &=& y^j 
	\frac{\partial}{\partial y^k} \hinti d\lambda \;
	\frac{|\lambda|}{1-\lambda} \xint_{\lambda y} F_j^{\;\:k} \\
&=& \hinti d\lambda \; \frac{|\lambda|}{1-\lambda} \xint_{\lambda y}
	dz \; z^k \: j_k(z) \\
&=& \veeint^y dz \; z^k \: j_k(z)
\end{eqnarray*}
\QED

\begin{Satz}
\label{a3_satz55}
Es gilt
\begin{eqnarray}
\label{eq:a3_X}
\frac{1}{2} \left\{ (s_0 \Aslsh p_0)(0,y) , y \slsh \right\} &=&
\frac{1}{2 \pi^4} \frac{1}{y^2} \hinti \frac{d\lambda}{|\lambda|}
	\xdint_{\lambda y} dz \; z^k \: A_k(z) \\
\label{eq:a3_Y}
&& + \frac{1}{8 \pi^4} \veeint^y dz \; z^k \int_0^1 d\alpha \; \alpha^2
	\: j_k(\alpha z) \\
\label{eq:a3_Z}
&& - \frac{i}{16 \pi^4} \left( \frac{1}{2} \veeint^y dz \; z^k \partial_k
	F_{ij} + y^k \: \partial_i \veeint^y F_{jk} \right) \sigma^{ij}
	\;\;\; .
\end{eqnarray}
\end{Satz}
{\Beweis}
Wir gehen zur\"uck zu Gleichung~\Ref{a3_O}.
\begin{eqnarray}
\lefteqn{ \frac{1}{2} \left\{ (s_0 \Aslsh p_0)(0,y) , y\slsh \right\} }
	\nonumber \\
&=& \frac{1}{32 \pi^4} \left\{ \gamma^i \gamma^j \gamma^k, y\slsh \right\}
	\left( \partial_k \veeint^y \partial_i A_j - \partial_{ik}
	\veeint^y A_j \right) \nonumber \\
&=& \frac{1}{16 \pi^4} \left( \gamma^i \gamma^j y^k - \gamma^i y^j \gamma^k
 + y^i \gamma^j  \gamma^k \right) \left( \partial_k \veeint^y \partial_i A_j
	- \partial_{ik} \veeint^y A_j \right) \nonumber \\
&=& \frac{1}{16 \pi^4} \left( g^{ij} \: y^k - g^{ik} \: y^j
 + g^{jk} \: y^i \right) \left( \partial_k \veeint^y \partial_i A_j
	- \partial_{ik} \veeint^y A_j \right) \nonumber \\
&& - \frac{i}{16 \pi^4} \left( \sigma^{ij} \: y^k - \sigma^{ik} \: y^j
 + \sigma^{jk} \: y^i \right) \left( \partial_k \veeint^y \partial_i A_j
	- \partial_{ik} \veeint^y A_j \right) \nonumber \\
\label{eq:a3_U}
&=& \frac{1}{16 \pi^4} \left( y^j \partial_j \veeint^y \partial_k A^k +
	y^j \: \partial_k \veeint^y F_j^{\;\:k} - 2 \: y^j \partial_{jk}
	\veeint^y A^k \right) \\
\label{eq:a3_V}
&& - \frac{i}{16 \pi^4} \left( \frac{1}{2} \: y^k \partial_k \veeint^y F_{ij}
	- y^k \: \partial_i \veeint^y F_{kj} \right) \: \sigma^{ij} \spc .
\end{eqnarray}
\Ref{a3_V} stimmt mit~\Ref{a3_Z} \"uberein,
wir m\"ussen also noch die einzelnen Terme in~\Ref{a3_U} umformen. Dabei nutzen
wir aus, da{\ss} die partielle Ableitung in Richtung von $y$ leicht ausgef\"uhrt
werden kann:
\begin{eqnarray*}
y^j \partial_j \veeint^y \partial_k A^k &=& \veeint^y dz \; z^j \partial_{jk}
	A^k(z) \\
y^j \partial_{jk} \veeint^y A^k &=& \frac{\partial}{\partial y^k} \left(
	y^j \partial_j \veeint^y A^k \right) - \partial_k \veeint^y A^k
\;=\; \partial_k \veeint^y dz \; (z^j \partial_j -1) \: A^k(z)
\end{eqnarray*}
Wende nun~\Ref{a3_151} an:
\begin{eqnarray}
\label{eq:a3_P}
&=& \frac{2}{y^2} \hinti \frac{d\lambda}{|\lambda|} \xdint_{\lambda y} dz \;
	z_k \: (z^j \partial_j -1) A^k(z)
	+ \veeint^y dz \; h_k[ (z^j \partial_j-1) A^k(z) ]
\end{eqnarray}
Der erste Summand l\"a{\ss}t sich noch vereinfachen
\begin{eqnarray}
\lefteqn{ \hspace*{-3cm} \hinti \frac{d\lambda}{|\lambda|} \xdint_{\lambda y} dz
	\; z_k \: (z^j \partial_j - 1) \: A^k(z)
\;=\; \hinti \frac{d\lambda}{|\lambda|} \xdint_{\lambda y} dz \;
	(z^j \partial_j -2) \: z_k \: A^k(z) } \nonumber \\
\label{eq:a3_377a}
&=& \hinti \frac{d\lambda}{|\lambda|} \left( \lambda \frac{d}{d\lambda} -2
	\right) \xdint_{\lambda y} dz \; z_k \: A^k(z) \\
\label{eq:a3_260a}
&=& -2 \hinti \frac{d\lambda}{|\lambda|} \xdint_{\lambda y} dz \;
	z^k \: A_k(z) \;\;\; .
\end{eqnarray}
Beachte, da{\ss} bei der partiellen Integration von \Ref{a3_377a} die Randwerte
bei $\lambda=\pm0$ nach \Ref{a3_51}, \Ref{a3_52} und \Ref{a3_181a} wegfallen.

Wir berechnen nun die Funktion $h_k$ in~\Ref{a3_P}:
\begin{eqnarray*}
\lefteqn{ h_k[ (z^j \partial_j -1) A^k(z) ] }\\
&=& z^j \partial_{jk} A^k - \int_0^1 d\alpha \; \alpha^2 \: z^j \partial_{jk}
	A_k(\alpha z) - \frac{1}{2} \: z_k \int_0^1 d\alpha \; \alpha^2
	\: (\alpha \: z^j \partial_j + 1) \: \Box A^k(\alpha z) \\
&=& z^j \partial_{jk} A^k(z) - \partial_k A^k(z) - \frac{1}{2} z_k \:
	\Box A^k(z) + 2 \int_0^1 d\alpha \; \alpha \: \partial_k A^k(\alpha z)
		\\
&& \spc + z_k \int_0^1 d\alpha \; \alpha^2 \: \Box A^k(\alpha z)
\end{eqnarray*}
Ersetze $\Box A_k = - j_k + \partial_{kl}A^l$ und integriere partiell:
\Equ{a3_260b}
 =\; - z_k \int_0^1 d\alpha \; \alpha^2 \: j^k(\alpha z) + \frac{1}{2} \:
	z_k \: j^k(z) + \frac{1}{2} z^j \partial_{jk} A^k(z)
\EndEqu
Nun setzen wir~\Ref{a3_260a} und \Ref{a3_260b} in~\Ref{a3_P} ein und erhalten
\begin{eqnarray*}
y^j \partial_{jk} \veeint^y A^k &=& - \frac{4}{y^2} \hinti
	\frac{d\lambda}{|\lambda|} \xdint_{\lambda y} dz \; z^k \:  A_k(z) \\
&& + \frac{1}{2} \veeint^y dz \; z_k \: \left( j^k(z) - 2 \int_0^1 d\alpha
	\; \alpha^2 \: j^k(\alpha z) \right) \\
&& + \frac{1}{2} \veeint^y dz \; z^j \partial_{jk} A^k(z) \spc .
\end{eqnarray*}
Den zweiten Summanden in~\Ref{a3_U} kann man schlie{\ss}lich mit Hilfe von
Lemma~\ref{a3_lemma4} berechnen.
\QED

\section{St\"orungsrechnung f\"ur das Gravitationsfeld}
Wie in Abschnitt~\ref{grav_k0} arbeiten wir mit der linearisierten
Gravitationstheorie in symmetrischer Eichung und der Koordinatenbedingung
\Ref{a1_210}. Der Diracoperator ist durch \Ref{a1_208} gegeben, f\"ur
$\Delta p_0$ erh\"alt man ganz analog wie in Abschnitt \ref{grav_k0}
\Equ{a3_g1}
\Delta p_0(x,y) \;=\; \left( \frac{1}{4} \:h(x) + \frac{3}{4} \: h(y) \right)
	\: p_0(x,y) - \frac{i}{e} \frac{\partial}{\partial y^k}
	\: \Delta p_0[\gamma^j h_j^k](x,y) \spc .
\EndEqu
Mit dieser Gleichung k\"onnen wir die Rechnung zum Teil auf diejenige
f\"ur das elektromagnetische Feld, Theorem \ref{a3_theorem1}, zur\"uckf\"uhren.
\begin{Thm}
\label{thm_gp0}
In erster Ordnung St\"orungstheorie gilt in symmetrischer Eichung
\begin{eqnarray}
\label{eq:a3_g3}
\Delta{p}_0(x,y) &=&
	- \left( \int_x^y h^k_j \right) \xi^j \frac{\partial}{\partial y^k}
	\;\; p_0(x,y) \\
\label{eq:a3_g4}
  && + \frac{i}{4 \pi^3} \; \frac{1}{\xi^4}
	\left( \int_x^y (2 \alpha -1) \; \gamma^i \: \xi^j \:\xi^k
	\; (h_{jk},_i - h_{ik},_j)  \right) \\
\label{eq:a3_g5}
  && + \frac{1}{8 \pi^3} \; \frac{1}{\xi^4}
	\left( \int_x^y \varepsilon^{ijlm} \; (h_{jk},_i
	- h_{ik},_j) \: \xi^k \; \xi_l \: \rho \gamma_m \right) \\
\label{eq:a3_g6}
  && + \frac{1}{2} \left( \int_x^y (\alpha^2 - \alpha) \; \xi^j \; \xi^k \;
	R_{jk} \right) \;\; p_0(x,y) \\
\label{eq:a3_g7}
  && + \frac{i}{32 \pi^3}  \; \frac{1}{\xi^2}
	\left( \int_x^y (\alpha^4 - 2 \alpha^3 + \alpha^2)
	\; \xi \slsh \; \xi^j \; \xi^k \; \Box R_{jk} \right) \\
\label{eq:a3_g8}
  && - \frac{i}{32 \pi^3}  \; \frac{1}{\xi^2}
	\left( \int_x^y (6 \alpha^2 - 6 \alpha + 1) \;
	\xi \slsh \; R \right) \\
\label{eq:a3_g9}
  && + \frac{i}{32 \pi^3} \; \frac{1}{\xi^2}
	\left( \int_x^y (4 \alpha^3 - 6 \alpha^2 + 2 \alpha)
	\; \xi^j \: \xi^k \: \gamma^l \; R_{j[k},_{l]} \right) \\
\label{eq:a3_g10}
  && + \frac{1}{16 \pi^3} \; \frac{1}{\xi^2}
	\left( \int_x^y (\alpha^2 - \alpha) \;
	\varepsilon^{ijlm} \: R_{ki},_j \: \xi^k \: \xi_l \: \rho \gamma_m
	\right) \\
\label{eq:a3_g11}
  && - \frac{i}{8 \pi^3} \; \frac{1}{\xi^2}
	\left( \int_x^y (\alpha^2 - \alpha) \; \xi^j \:
	\gamma^k \: G_{jk} \right) \\
\label{eq:a3_g12}
  && + {\cal{O}}(\ln(|\xi^2|)) \nonumber \spc  .
\end{eqnarray}
\end{Thm}
{\Beweis}
Bei der St\"orungsrechnung f\"ur das elektromagnetische Feld in Abschnitt
\ref{elek_p0} und \ref{elek_k0} hatten wir gesehen, da{\ss} die
asymptotischen Entwicklungsformeln von $\Delta k_0$ und $\Delta p_0$
bis zur Ordnung $\xi^0$ bzw. $\ln(|\xi^2|)$ durch die Ersetzungen
\begin{eqnarray}
k_0 &\longrightarrow& p_0 \nonumber \\
i \pi \; (m^\vee(\xi)-m^\wedge(\xi)) &\longrightarrow&
	\frac{1}{\xi^4} \nonumber \\
-i \pi \; (l^\vee(\xi)-l^\wedge(\xi)) &\longrightarrow& 
	\frac{1}{\xi^2} \nonumber \\
\label{eq:a3_393a}
- i \pi \; \Theta(\xi^2) \: \epsilon(\xi^0) &\longrightarrow&
	\ln(|\xi^2|)
\end{eqnarray}
ineinander \"ubergehen. Diese Ersetzungen werden durch das Berechnen
partieller Ableitungen nach $y$ nicht zerst\"ort (es spielt also keine
Rolle, ob man zuerst die Ersetzungen durchf\"uhrt und dann nach $y$
ableitet oder umgekehrt).
Durch Vergleich von \Ref{a1_223} und \Ref{a3_g1} sieht man daher, da{\ss} sich
die Ersetzungen auch f\"ur das Gravitationsfeld \"ubertragen.
Damit folgt die Behauptung direkt aus Theorem \ref{a1_theorem2}.
\QED

\section{Skalare St\"orung}
Wir betrachten wieder die skalare St\"orung \Ref{a1_s1}. F\"ur $\Delta p_0$
hat man in erster Ordnung in $\Xi$
\Equ{a3_s1}
   \Delta p_0 \;=\; - (s_0 \: \Xi \: p_0 + p_0 \: \Xi \: s_0) \spc .
\EndEqu
\begin{Thm}
\label{theorem_sp0}
In erster Ordnung St\"orungstheorie gilt
\begin{eqnarray}
\label{eq:a3_s2}
\Delta p_0(x,y) &=& - \frac{1}{2} \: (\Xi(y)+\Xi(x)) \; p^{(1)}(x,y) \\
&&+ \frac{i}{8 \pi^3} \frac{1}{\xi^2} \int_x^y
	\; (\partial_j \Xi) \; \xi_k \; \sigma^{jk} \\
&&+ \frac{i}{32 \pi^3} \: \ln(|\xi^2|) \int_x^y
	(\alpha^2 - \alpha) \; (\partial_j \Box \Xi) \;
	\xi_k \; \sigma^{jk} \\
&&- \frac{1}{32 \pi^3} \: \ln(|\xi^2|) \int_x^y
	\Box \Xi \\
&& +\; {\cal{O}}(\xi^0) \spc . \nonumber
\end{eqnarray}
\end{Thm}
{\Beweis}
Es gilt
\begin{eqnarray*}
\lefteqn{ (s_0 \: \Xi \: p_0)(0,y) \;=\; -\frac{1}{16 \pi^4} \: \Pdd_x
	\left( \veeint_x^y \Xi \right) \Pdd_{y \; |y=0}
\;=\; -\frac{1}{16 \pi^4} \left( \veeint^y (\Pdd \Xi) \: \Pdd_y +
	\Box_y \veeint^y \Xi \right) } \\
&=& \frac{1}{8 \pi^4} \frac{1}{y^2} \hinti \frac{d\lambda}{|\lambda|}
	\xdint_{\lambda y} dz \; (\Pdd \Xi) \: z\slsh \\
&& + \frac{1}{16 \pi^4} \veeint^y dz \left( (\Box \Xi)(z) - \int_0^1 d\alpha
	\; \alpha \: (\Box \Xi)(\alpha z) - \frac{1}{2} \int_0^1 d\alpha
	\; \alpha^2 \: (\Pdd \Box \Xi)(\alpha z) \; z \slsh \right) \\
&=& \frac{1}{8 \pi^4} \frac{1}{y^2} \hinti \frac{d\lambda}{|\lambda|}
	\xdint_{\lambda y} dz \left( z^j \partial_j \Xi - i
	\: (\partial_j \Xi) \: z_k \; \sigma^{jk} \right) \\
&& + \frac{1}{32 \pi^4} \veeint^y dz \left( (\Box \Xi)(z) + i \int_0^1 d\alpha
	\; \alpha^2 \: (\partial_j \Box \Xi)(\alpha z) \; z_k \; \sigma^{jk}
	\right) \spc .
\end{eqnarray*}
Nach Satz \ref{a3_satz99} und Satz \ref{a3_satz98} hat man die asymptotischen
Entwicklungen
\begin{eqnarray*}
\lefteqn{ \hinti \frac{d\lambda}{|\lambda|} \xdint_{\lambda y} dz \left(
	z^j \partial_j \Xi - i \: (\partial_j \Xi) \; z_k \; \sigma^{jk}
	\right) } \\
&=& - \pi \: \Xi(0) - \frac{i\pi}{2} \hinti d\lambda \; \epsilon(\lambda)
	\: (\partial_j \Xi)(\lambda y) \: y_k \; \sigma^{jk} \\
&& - \frac{i \pi}{8} \: y^2 \: \ln(|y^2|) \hinti d\lambda \;
	\epsilon(\lambda) \; \lambda^2 \;(\partial_j \Box \Xi)(\lambda
	y) \: y_k \; \sigma^{jk} \;+\; {\cal{O}}(y^2) \\
\lefteqn{ \veeint^y dz \left( (\Box \Xi)(z) + i \int_0^1 d\alpha \; \alpha^2
	\; (\partial_j \Box \Xi)(\alpha z) \; z_k \; \sigma^{jk} \right) } \\
&=& \frac{\pi}{2} \: \ln(|y^2|) \hinti d\lambda \; \epsilon(\lambda)
	\left( (\Box \Xi)(\lambda y) + i \lambda \; (\partial_j \Box
	\Xi)(\lambda y) \; y_k \; \sigma^{jk} \right) \;+\; {\cal{O}}(y^0)
\end{eqnarray*}
und somit
\begin{eqnarray}
\lefteqn{ (s_0 \: \Xi \: p_0)(0,y)
\;=\; -\frac{1}{8 \pi^3} \frac{1}{y^2} \; \Xi(0) } \nonumber \\
&&- \frac{i}{16 \pi^3} \frac{1}{y^2} \hinti d\lambda \; \epsilon(\lambda)
	\; (\partial_j \Xi)(\lambda y) \; y_k \; \sigma^{jk} \nonumber \\
&&- \frac{i}{64 \pi^3} \: \ln(|y^2|) \hinti d\lambda \; \epsilon(\lambda)
	\; (\lambda^2 - \lambda) \; (\partial_j \Box \Xi)(\lambda y) \;
	y_k \; \sigma^{jk} \nonumber \\
\label{eq:a3_395a}
&&+ \frac{1}{64 \pi^3} \ln(|y^2|) \hinti d\lambda \; \epsilon(\lambda)
	\; (\Box \Xi)(\lambda y) \;+\; {\cal{O}}(y^0) \spc .
\end{eqnarray}
Durch Translation um $x$ und Vertauschung von $x$, $y$ erh\"alt man
hieraus Formeln f\"ur $(s_0 \:\Xi\: p_0)(x,y)$ und $(p_0 \:\Xi\: s_0)(x,y)$.
Setze nun in \Ref{a3_s1} ein.
\QED

\section{Bilineare St\"orung}
Wir betrachten wieder die bilineare St\"orung \Ref{a1_b00}. In erster Ordnung
in $B$ hat man
\Equ{a3_b10}
\Delta p_0 \;=\; - ( s_0 \: B_{jk} \: \sigma^{jk} \: p_0 + p_0 \: B_{jk} \:
	\sigma^{jk} \: s_0 ) \spc .
\EndEqu

\begin{Thm}
In erster Ordnung St\"orungstheorie gilt
\begin{eqnarray}
\label{eq:a1_b11}
\Delta p_0(x,y) &=& -\frac{1}{\pi^3} \; \frac{1}{\xi^4}
	\int_x^y B_{ij} \; \xi^i \: \xi_k \; \sigma^{jk} \\
\label{eq:a1_b13}
&&+ \frac{1}{8 \pi^3} \; \frac{1}{\xi^2} \; (B_{jk}(y) +
	B_{jk}(x)) \; \sigma^{jk} \\
\label{eq:a1_b14}
&&- \frac{1}{2 \pi^3} \; \frac{1}{\xi^2} \int_x^y B_{jk} \;
	\sigma^{jk} \\
\label{eq:a1_b12}
&&+ \frac{i}{4 \pi^3} \; \frac{1}{\xi^2} \int_x^y \xi_j \;
	B^{jk}_{\;\;,k} \\
\label{eq:a1_b15}
&&+ \frac{1}{4 \pi^3} \; \frac{1}{\xi^2} \int_x^y (2\alpha-1) \;
	(\xi^k \: B_{jk,i} + \xi_i \: B_{jk,}^{\;\;\;k}) \;
	\sigma^{ij} \\
\label{eq:a1_b16}
&&+ \frac{1}{4 \pi^3} \; \frac{1}{\xi^2} \int_x^y
	(\alpha^2-\alpha) \; (\Box B_{ij}) \: \xi^i \: \xi_k \; \sigma^{jk} \\
\label{eq:a1_b17}
&&- \frac{1}{8 \pi^3} \; \frac{1}{\xi^2} \int_x^y
	\varepsilon^{ijkl} \; B_{ij,k} \; \xi_l \; \rho \\
&&+ {\cal{O}}(\ln(|\xi^2|) \spc . \nonumber
\end{eqnarray}
\end{Thm}
{\Beweis}
Unter Verwendung von \Ref{a1_b8} erh\"alt man
\begin{eqnarray}
\lefteqn{ (s_0 \: B_{jk} \: \sigma^{jk} \: p_0)(x,y) \;=\; \frac{1}{16 \pi^4}
	\: \Pdd_x \: \sigma^{jk} \: \Pdd_y \veeint_x^y B_{jk}
\;=\; \frac{i}{16 \pi^4} \: \frac{\partial}{\partial x^a}
	\frac{\partial}{\partial y^b} \veeint_x^y \gamma^a \: \gamma^j \gamma^k
	\: \gamma^b \; B_{jk} } \nonumber \\
&=& \frac{i}{8 \pi^4} \: \frac{\partial}{\partial x^j}
	\frac{\partial}{\partial y^k} \veeint_x^y B^{jk}
\;+\; \frac{1}{16 \pi^4} \: \frac{\partial}{\partial x^a} \frac{\partial}
	{\partial y_a} \veeint_x^y B_{jk} \; \sigma^{jk} \nonumber \\
\label{eq:a3_b20}
&&+ \frac{1}{8 \pi^4} \: \frac{\partial}{\partial x^a}
	\frac{\partial}{\partial y^b} \veeint_x^y (B^a_{\;\:j} \: \sigma^{jb}
	+ B^b_{\;\:j} \: \sigma^{ja} )
\;+\; \frac{1}{16 \pi^4} \; \rho \: \varepsilon^{ijab} \:
	\frac{\partial}{\partial x^a} \frac{\partial}{\partial y^b}
	\veeint_x^y B_{ij} \;\; . \spc
\end{eqnarray}
Berechne nun mit Hilfe von \Ref{a3_151} und \Ref{a3_152} die auftretenden
Terme nacheinander. Unter Verwendung von Satz \ref{a3_satz99} und Satz
\ref{a3_satz98} f\"uhren wir eine Entwicklung um den Lichtkegel durch,
dabei lassen wir alle Beitr\"age der Ordnung $\ln(|\xi^2|)$ weg.
Um die Notation zu vereinfachen, setzen wir $x=0$.
\begin{eqnarray*}
\lefteqn{ \frac{\partial}{\partial x^j} \frac{\partial}{\partial y^k}
	\veeint_x^y B^{jk}_{\;\;\;\;|x=0} \;=\;
	\frac{\partial}{\partial y^k} \veeint^y
	\partial_j B^{jk} - \frac{\partial^2}{\partial y^j \: \partial y^k}
	\veeint^y B^{jk}
\;=\; - \frac{\partial}{\partial y^j} \veeint^y B^{jk}_{\;\;,k} } \\
&=& - \frac{2}{y^2} \hinti \frac{d\lambda}{|\lambda|} \xdint_{\lambda y}
	dz \; z_j \: B^{jk}_{\;\;,k} - \veeint^y h_j[B^{jk}_{\;\;,k}] \\
&=& - \frac{\pi}{y^2} \hinti d\lambda \; \epsilon(\lambda) \: y_j \;
	B^{jk}_{\;\;,k}(\lambda y) + {\cal{O}}(\ln(|y^2|)) \\
\lefteqn{\frac{\partial}{\partial x^a} \frac{\partial}{\partial y_a}
	\veeint_x^y B_{jk \; | x=0} \;=\; \frac{\partial}{\partial y^a}
	\veeint^y \partial^a B_{jk} - \Box_y \veeint^y B_{jk}
\;=\; \frac{\partial}{\partial y^a} \veeint^y \partial^a B_{jk} } \\
&=& \frac{2}{y^2} \hinti \frac{d\lambda}{|\lambda|} \xdint_{\lambda y}
	dz \; z^a \partial_a B_{jk} + \veeint^y h_a[\partial^a B_{jk}] \\
&=& \frac{\pi}{y^2} \hinti d\lambda \; \epsilon(\lambda) \: y^a \partial_a
	B_{jk}(\lambda y) + {\cal{O}}(\ln(|y^2|)
	\;=\; -\frac{2\pi}{y^2} \: B_{jk}(0) + {\cal{O}}(\ln(|y^2|)) \\
\lefteqn{ \frac{\partial}{\partial x^a} \frac{\partial}{\partial y^b}
	\veeint_x^y (B^a_{\;\:j} \: \sigma^{jb} + B^b_{\;\:j} \:
	\sigma^{ja})_{|x=0} } \\
&=& \frac{\partial}{\partial y^b} \veeint^y (\partial_a B^a_{\;\:j}
	\: \sigma^{jb} + \partial_a B^b_{\;\:j} \: \sigma^{ja} )
	- 2 \: \frac{\partial^2}{\partial y^a \: \partial y^b} \veeint^y
	B^a_{\;\:j} \; \sigma^{jb} \\
&=& \frac{16}{y^4} \hinti d\lambda \; \frac{\epsilon(\lambda)}{\lambda^2}
	\dotint_{\lambda y} dz \; B^a_{\;\:j} \; z_a \: z_b \; \sigma^{jb} \\
&&+ \frac{2}{y^2} \hinti \frac{d\lambda}{|\lambda|} \xdint_{\lambda y} dz \;
	(\Box B^a_{\;\:j}) \; z_a \: z_b \: \sigma^{jb} - 2 \veeint^y
	h_{ab}[B^a_{\;\:j}] \; \sigma^{jb} \\
&&+ \frac{2}{y^2} \hinti \frac{d\lambda}{|\lambda|} \xdint_{\lambda y}
	dz \; z_b \;
	(\partial_a B^a_{\;\:j} \: \sigma^{jb} + \partial_a B^b_{\;\:j} \:
	\sigma^{ja}) \\
&&+ \veeint^y h_b[\partial_a B^a_{\;\:j} \: \sigma^{jb} + \partial_a
	B^b_{\;\:j}j \: \sigma^{ja}] \\
&=& \frac{4\pi}{y^4} \hinti d\lambda \; \epsilon(\lambda) \; B_{aj}(\lambda
	y) \: y^a \: y_b \; \sigma^{jb} \\
&&+ \frac{2\pi}{y^2} \hinti d\lambda \; \epsilon(\lambda) \; B_{jk}(\lambda
	y) \: \sigma^{jk} \\
&&- \frac{\pi}{y^2} \hinti d\lambda \; \epsilon(\lambda) \:
	(\lambda^2-\lambda) \; (\Box B_{aj})(\lambda y)
	\: y^a \: y_b \; \sigma^{jb} \\
&&- \frac{\pi}{y^2} \hinti d\lambda \; \epsilon(\lambda) \: (2\lambda -1)
	\; (y^a \: B_{aj,b}(\lambda y) + y_b \: B_{aj,}^{\;\;\;a}(\lambda y))
	\; \sigma^{jb}
	\;+\; {\cal{O}}(\ln(|y^2|)) \\
\lefteqn{ \varepsilon^{ijab} \: \frac{\partial}{\partial x^a}
\frac{\partial}{\partial y^b} \veeint_x^y B_{ij \; |x=0}
\;=\; \varepsilon^{ijab} \: \left( \frac{\partial}{\partial y^b}
	\veeint^y \partial_a B_{ij} - \frac{\partial^2}{\partial y^a \:
	\partial y^b} \veeint^y B_{ij} \right) } \\
&=& \varepsilon^{ijab} \: \frac{\partial}{\partial y^b} \veeint^y \partial_a
	B_{ij}
\;=\; \frac{2}{y^2} \: \varepsilon^{ijab} \hinti \frac{d\lambda}{|\lambda|}
	\xdint_{\lambda y} dz \; z_b \: \partial_a B_{ij}
	+ \varepsilon^{ijab} \veeint^y h_b[\partial_a B_{ij}] \\
&=& \frac{\pi}{y^2} \hinti d\lambda \; \epsilon(\lambda) \:
	\epsilon^{ijkl} \: B_{ij,k}(\lambda y) \; y_l +
	{\cal{O}}(\ln(|y^2|))
\end{eqnarray*}
Wir setzen in \Ref{a3_b20} ein:
\begin{eqnarray*}
\lefteqn{ (s_0 \: B_{jk} \: \sigma^{jk} \: p_o)(0,y)
\;=\; \frac{1}{2 \pi^3} \frac{1}{y^4} \hinti d\lambda \; \epsilon(\lambda)
	\; B_{ij}(\lambda y) \; y^i \: y_k \; \sigma^{jk} } \\
&&- \frac{i}{8 \pi^3} \frac{1}{y^2} \hinti d\lambda \; \epsilon(\lambda)
	\: y_j \; B^{jk}_{\;\;,k}(\lambda y) \\
&&- \frac{1}{8 \pi^3} \frac{1}{y^2} \; B_{jk}(0) \; \sigma^{jk} \\
&&+ \frac{1}{4 \pi^3} \frac{1}{y^2} \hinti d\lambda \; \epsilon(\lambda)
	\: B_{jk}(\lambda y) \; \sigma^{jk} \\
&&- \frac{1}{8 \pi^3} \frac{1}{y^2} \hinti d\lambda \; \epsilon(\lambda)
	\: (2\lambda-1) \; (y^k \: B_{jk,i}(\lambda y) +
	y_i \: B_{jk,}^{\;\;\;k}(\lambda y)) \; \sigma^{ij} \\
&&- \frac{1}{8 \pi^3} \frac{1}{y^2} \hinti d\lambda \; \epsilon(\lambda)
	\: (\lambda^2-\lambda) \; (\Box B_{ij})(\lambda y) \;
	y^i \: y_k \;  \sigma^{jk} \\
&&+ \frac{1}{16 \pi^3} \frac{1}{y^2} \hinti d\lambda \; \epsilon(\lambda)
	\; \varepsilon^{ijkl} \: B_{ij,k}(\lambda y) \;
	y_l \; \rho
\;+\; {\cal{O}}(\ln(|\xi^2|))
\end{eqnarray*}
Durch Translation um $x$ und Vertauschung von $x$, $y$ erh\"alt man eine
entsprechende Formel f\"ur $(p_0 \: B_{ij} \: \sigma^{ij} \: s_0)(x,y)$.
Setze nun in \Ref{a3_b10} ein.
\QED

\section{Differentialst\"orung durch Vektorpotential}
Wir betrachten wieder die St\"orung des Diracoperators \Ref{a1_va}. F\"ur
$\Delta p_0$ hat man in erster Ordnung in $L$
\[ \Delta p_0(x,y) \;=\; -i \frac{\partial}{\partial y^j} \:
	\Delta p_0[L^j](x,y) + \frac{i}{2} \: \Delta p_0[L^j_{\;,j}](x,y) \spc . \]
\begin{Thm}
In erster Ordnung St\"orungstheorie gilt
\begin{eqnarray}
\label{eq:a3_v1}
\Delta p_0(x,y) &=& + \frac{i}{4 \pi^3} \frac{1}{\xi^4} \;
	(L^j(y) + L^j(x)) \: \xi_j \\
\label{eq:a3_v2}
&&- \frac{1}{4 \pi^3} \frac{1}{\xi^4} \int_x^y L_{m,j} \;
	\xi^m \: \xi_k \; \sigma^{jk} \\
\label{eq:a3_v3}
&&- \frac{i}{16 \pi^3} \frac{1}{\xi^2} \; (L^j_{\;,j}(y)
	- L^j_{\;,j}(x)) \\
\label{eq:a3_v4}
&&+ \frac{1}{16 \pi^3} \frac{1}{\xi^2} \int_x^y
	(2\alpha-1) \: L^m_{\;,mj} \: \xi_k \; \sigma^{jk} \\
\label{eq:a3_v5}
&&+ \frac{1}{8 \pi^3} \frac{1}{\xi^2} \int_x^y L_{k,j} \;
	\sigma^{jk} \\
\label{eq:a3_v6}
&&+ \frac{1}{16 \pi^3} \frac{1}{\xi^2} \int_x^y 
	(\alpha^2-\alpha) \: (\Box L_{m,j}) \; \xi^m \: \xi_k \; \sigma^{jk}\\
\label{eq:a3_v7}
&&+ \frac{i}{16 \pi^3} \frac{1}{\xi^2} \int_x^y (\Box L_j) \:
	\xi^j \\
&&+ {\cal{O}}(\ln(|\xi^2|)) \spc . \nonumber
\end{eqnarray}
\end{Thm}
{\Beweis}
Man kann genau wie im Beweis von Theorem~\ref{thm_gp0} argumentieren.
Durch Vergleich von Theorem~\ref{theorem_sk0} und Theorem~\ref{theorem_sp0}
folgt die Behauptung aus Theorem~\ref{theorem_vk0} durch die formalen
Ersetzungen~\Ref{a3_393a}.
\QED

\section{Bilineare Differentialst\"orung durch Vektorpotential}
Wir betrachten wieder die St\"orung des Diracoperators \Ref{a1_bva}. F\"ur
$\Delta p_0$ hat man in erster Ordnung in $L$
\begin{eqnarray*}
\Delta p_0(x,y) &=& -i \Delta p_0 \left[ i L^j \partial_j + \frac{i}{2} \: L^j_{\;,j}
	\right](x,y) \:-\: i p_0(x,y) \: L \slsh(y) \\
&&+ \frac{i}{2} \: \Delta p_0[L_{j,k} \sigma^{jk}](x,y) \:-\: \frac{1}{2} \:
	\Delta p_0[L^j_{\;,j}](x,y) \spc .
\end{eqnarray*}

\begin{Thm}
In erster Ordnung St\"orungstheorie gilt
\begin{eqnarray}
\label{eq:a3_bv1}
\Delta k_0(x,y) &=& -\frac{i}{4 \pi^3} \; \frac{1}{\xi^4} \;
	(L_j(y) + L_j(x)) \: \xi_k \: \sigma^{jk} \\
&&- \frac{1}{4 \pi^3} \; \frac{1}{\xi^4} \; (L_j(y) - L_j(x)) \: \xi^j \\
&&-\frac{1}{16 \pi^3} \;\frac{1}{\xi^2} \; (L^k_{\;,k}(y)+L^k_{\;,k}(x)) \\
&&-\frac{i}{16 \pi^3} \;\frac{1}{\xi^2} \;
	\int_x^y L^m_{\;,mj} \; \xi_k \; \sigma^{jk} \;+\; {\cal{O}}(\ln(|\xi^2|))
	\;\;\;.
\end{eqnarray}
\end{Thm}
{\Beweis}
Man kann genau wie im Beweis von Theorem~\ref{thm_gp0} argumentieren.
Durch Vergleich von Theorem~\ref{theorem_sk0} und Theorem~\ref{theorem_sp0}
folgt die Behauptung aus Theorem~\ref{theorem_bvk0} durch die formalen
Ersetzungen~\Ref{a3_393a}.
\QED

\chapter{St\"orungsrechnung f\"ur $p_m$ im Ortsraum}
\label{anh4}
Wie in Anhang~\ref{anh2} werden wir $m>0$ annehmen, au{\ss}erdem werden wir die
Rechnung nur bis zur Ordnung $m^2$ durchf\"uhren.
Man hat
\Equ{a4_1}
  p_m(x) \;=\; (i \Pdd + m) \: P_{m^2}(x)
\EndEqu
mit
\begin{eqnarray*}
P_{m^2}(x) &=& \int \frac{d^4k}{(2 \pi)^4} \; \delta(k^2-m^2) \; e^{-i k x} \\
&=& \left\{ \begin{array}{ll}
	\displaystyle \frac{m^2}{8 \pi^2} \frac{Y_1(m \tau)}{m \tau}
		& , x \in \I \\[1em]
	\displaystyle \frac{m^2}{4 \pi^3} \frac{K_1(m \tau)}{m \tau}
		& , x \in \Ra
   \end{array} \right. \spc ,
\end{eqnarray*}
wobei $\tau = \sqrt{|x^2|}$ gesetzt wurde; $Y_1$, $K_1$ sind Besselfunktionen
zweiter Art.

Die Entwicklung von~\Ref{a4_1} und~\Ref{a2_6} nach $m$ liefert
\begin{eqnarray*}
p_m(x) &=& \left( i \Pdd + m - \frac{i}{2} \: m^2 \: x \slsh \right) P_0(x) +
{\cal{O}}(m^3) \\
s_m(x) &=& \left( i \Pdd + m - \frac{i}{2} \: m^2 \: x \slsh \right) S_0(x) +
{\cal{O}}(m^3) \spc .
\end{eqnarray*}
\section{St\"orungsrechnung f\"ur das elektromagnetische Feld}
Betrachte nun speziell die St\"orung~\Ref{a1_54a} durch ein
elektromagnetisches Feld. In erster Ordnung St\"orungstheorie hat man
\[  \tilde{p}_m \;=\; p_m + \Delta p_m \]
mit
\Equ{a4_a10}
  \Delta p_m \;=\; - e \left( p_m \:\Aslsh\: s_m + s_m \:\Aslsh\: p_m \right)
   \spc .
\EndEqu
Zun\"achst wollen wir eine Gleichung f\"ur $(s_m \: \Aslsh \: p_m)(x,y)$
ableiten; der Einfachheit halber setzen wir $x=0$.
\begin{Lemma}
\label{a4_lemma1}
Es gilt
\begin{eqnarray*}
(s_m \: \Aslsh \: p_m)(0,y) &=& (s_0 \: \Aslsh \: p_0)(0,y) \\
&&- m \: \frac{i}{4 \pi^4} \frac{1}{y^2} \hinti \frac{d\lambda}{|\lambda|}
	\xdint_{\lambda y} dz \; A_j(z) \: z^j \\
&&- m \: \frac{i}{16 \pi^4} \veeint^y dz \; z^k\int_0^1 d\alpha \; \alpha^2 \:
	j_k(\alpha z) \\
&&+ m \: \frac{1}{32 \pi^4} \veeint^y F_{ij} \: \sigma^{ij} \\
&&+ m^2 \: \frac{1}{8 \pi^4} \frac{y \slsh}{y^2} \hinti
	\frac{d\lambda}{|\lambda|} \xdint_{\lambda y} dz \; A_j(z) \: z^j \\
&&+ m^2 \: \frac{i}{64 \pi^4} \veeint^y dz \; \varepsilon^{ijkl} \; F_{ij}(z)
	\: y_k \; \rho \gamma_l \\
&&+ m^2 \: \frac{1}{32 \pi^4} \veeint^y dz \; F_{ij} \; \gamma^i \: (2z-y)^j \\
&&+ m^2 \: \frac{1}{32 \pi^4} \veeint^y dz  \; z^k \int_0^1 d\alpha \;
	\alpha^2 \; j_k(\alpha z) \; y\slsh \\
&&+ {\cal{O}}(m^3) \spc .
\end{eqnarray*}
\end{Lemma}
{\Beweis}
Betrachte die Beitr\"age verschiedener Ordnung in $m$ nacheinander:
\begin{enumerate}
\item Terme $\sim m^0$:\\
F\"uhrt auf die St\"orungsrechnung f\"ur $m=0$, vergleiche
Lemma~\ref{a3_lemma21}.
\item Terme $\sim m$:\\ Man hat den Beitrag
\begin{eqnarray}
\lefteqn{ m \left( (i \Pdd) \: S_0 \Aslsh P_0 + S_0 \Aslsh P_0 \: (i \Pdd)
	\right)(x,y)  } \nonumber \\
\label{eq:a4_243a}
&=& \frac{im}{16 \pi^4} \left( \Pdd_x (\veeint_x^y \Aslsh) + (\veeint_x^y
\Aslsh) \Pdd_y \right) \\
&=& \frac{im}{16 \pi^4} \left( \veeint_x^y \Pdd \Aslsh
	- \Pdd_y (\veeint_x^y \Aslsh) - \frac{\partial}{\partial y^j}
	(\veeint_x^y \Aslsh) \gamma^j \right) \\
\label{eq:a4_243b}
&=& \frac{im}{16 \pi^4} \left( \veeint_x^y \Pdd \Aslsh
	- 2 \frac{\partial}{\partial y^j}
	(\veeint_x^y A^j) \right) \spc . \nonumber
\end{eqnarray}
Setze nun $x=0$ und wende \Ref{a3_151} an:
\begin{eqnarray*}
\lefteqn{ m \left( (i \Pdd) \: S_0 \Aslsh P_0 + S_0 \Aslsh P_0 \: (i \Pdd)
	\right)(0,y)  } \\
&=& - \frac{im}{4 \pi^4} \frac{1}{y^2} \hinti \frac{d\lambda}{|\lambda|}
	\xdint_{\lambda_y} dz \; A^j(\lambda z) \: z_j \\
&&+ \frac{im}{16 \pi^4} \veeint^y \left( -2 h_j[A^j] + \Pdd \Aslsh \right)
	\spc .
\end{eqnarray*}
F\"uhre nun noch die Ersetzungen
\begin{eqnarray}
\Pdd \Aslsh &=& \partial_j A^j - \frac{i}{2} \: F_{ij} \: \sigma^{ij} \spc ,
	\\
h_j[A^j](z) &=& \partial_j A^j(z) - \int_0^1 d\alpha \; \alpha \: \partial_j
	A^j(\alpha z) - \frac{1}{2} \: z_j \int_0^1 d\alpha \; \alpha^2 \;
	(\Box A^j)(\alpha z) \nonumber \\
\label{eq:a4_4}
&=& \frac{1}{2} \int_0^1 d\alpha \; \alpha^2 \; j_k(\alpha z) \: z^k
	+ \frac{1}{2} \partial_j A^j(z)
\end{eqnarray}
durch.
\item Terme $\sim m^2$:\\ Man hat den Beitrag
\begin{eqnarray}
\lefteqn{ \frac{m^2}{16 \pi^4} \left( \veeint_x^y \Aslsh + \frac{1}{2} \Pdd_x
\veeint_x^y dz \; \Aslsh(z) \; (z-y)_j \gamma^j \right. \nonumber } \\
\label{eq:a4_3}
&&\spc\spc \left. + \frac{1}{2} \veeint_x^y
dz \; (x-z)_j \gamma^j \: \Aslsh(z) \; \Pdd_y \right) \spc .
\end{eqnarray}
Beachte f\"ur $f \in C^\infty_3$ zun\"achst die Beziehungen
\begin{eqnarray}
\lefteqn{ \Pdd_y \veeint_x^y dz \; f(z) \; (z-y)_j \gamma^j \;=\; \Pdd_y \int
d^4z \; \frac{(z-y)_j \gamma^j}{(y-z)^2} \; l(z) \; f(z) } \nonumber \\
\label{eq:a4_246z}
&=& -2 \int d^4z \; \frac{1}{(y-z)^2} \; l(z) \; f(z) \;=\; -2 \veeint_x^y f
	\\
\lefteqn{ \Pdd_y \veeint_x^y dz \; (x-z)^j \gamma_j \; f(z) } \nonumber \\
\label{eq:a4_246y}
&=& \Pdd_x \int d^4z \; (x-z)^j \gamma_j \: l(z) \; \frac{1}{(y-z)^2} \:
	f(z) \;=\; 2 \veeint_x^y f \spc .
\end{eqnarray}
Mit \Ref{a4_246z} kann man den zweiten Summanden in~\Ref{a4_3} umformen
\begin{eqnarray}
\lefteqn{ \Pdd_x \veeint_x^y dz \; \Aslsh(z) \; (z-y)_j \gamma^j }
	\nonumber \\
&=& \veeint_x^y dz \; (\Pdd \Aslsh)(z) \; (z-y)_j \gamma^j -
	\Pdd_y \veeint_x^y dz \; \Aslsh(z) \; (z-y)_j \gamma^j \nonumber \\
\label{eq:a4_246a}
&=& \veeint_x^y dz \; (\Pdd \Aslsh)(z) \; (z-y)_j \gamma^j - 2
	\frac{\partial}{\partial y^k} \veeint_x^y dz \; A^k(z) \; (z-y)_j
	\gamma^j - 2 \veeint_x^y \Aslsh \; , \spc
\end{eqnarray}
aus \Ref{a4_246y} erh\"alt man ganz analog
\begin{eqnarray}
\lefteqn{ \veeint_x^y dz \; (x-z)_j \gamma^j \; \Aslsh(z) \; \Pdd_y }
	\nonumber \\
\label{eq:a4_246b}
&=& \veeint_x^y dz \; (x-z)_j \gamma^j \; (\Pdd \Aslsh)(z) - 2 
	\frac{\partial}{\partial y^k} \veeint_x^y dz \; (x-z)_j \gamma^j
	\; A^k(z) - 2 \veeint_x^y \Aslsh \spc .
\end{eqnarray}
Setzt man dies in~\Ref{a4_3} ein und setzt $x=0$, so erh\"alt man f\"ur den
Beitrag der Ordnung $m^2$ den Ausdruck
\begin{eqnarray*}
\frac{m^2}{16 \pi^4} \left( \frac{1}{2} \veeint^y dz \; \left( (\Pdd
\Aslsh)(z) \; (z-y)_j \gamma^j - z_j \gamma^j \; (\Pdd \Aslsh)(z) \right)
	\right. \\
\spc \left. + \frac{\partial}{\partial y^k} \veeint^y dz \; A^k(z) \; y \slsh
	- \veeint^y \Aslsh \right) \spc .
\end{eqnarray*}
Wende nun \Ref{a3_151} an und setze die Gleichungen~\Ref{a4_4}
sowie
\begin{eqnarray*}
\lefteqn{ (\Pdd \Aslsh) \; (z-y)_j \gamma^j - z_jy^j \; (\Pdd \Aslsh) }\\
&=& - \partial_k A^k \: y \slsh + F_{ij} \: \gamma^i \: (2z-y)_j +
	\frac{i}{2} \: \varepsilon^{ijkl} \: F_{ij} \: y_k \; \rho \gamma_l
\end{eqnarray*}
ein.
\end{enumerate}
\QED
Nun f\"uhren wir wieder eine Entwicklung um den Lichtkegel durch:
\begin{Satz}
\label{a4_lemma2}
Es gilt
\begin{eqnarray*}
\lefteqn{ (s_m \: \Aslsh \: p_m)(0,y) \;=\; (s_0 \: \Aslsh \: p_0)(0,y) } \\
&&- m \; \frac{i}{8 \pi^3}\frac{1}{y^2} \hinti d\lambda \; \epsilon(\lambda)
	\: A_j(\lambda y) \: y^j \\
&&+ m \; \frac{1}{64 \pi^3} \; \ln(|y^2|) \hinti d\lambda \; \epsilon(\lambda)
	\; F_{ij}(\lambda y) \; \sigma^{ij} \\
&&+ m \; \frac{i}{32 \pi^3} \; \ln(|y^2|) \hinti d\lambda \; \epsilon(\lambda)
	\: (\lambda^2 - \lambda) \; j_k(\lambda y) \: y^k \\
&&+ m^2 \; \frac{1}{16 \pi^3} \frac{y\slsh}{y^2} \hinti d\lambda \;
	\epsilon(\lambda) \: A_j(\lambda y) \: y^j \\
&&+ m^2 \; \frac{1}{64 \pi^3} \; \ln(|y^2)| \hinti d\lambda \;
	\epsilon(\lambda)
	\: (2\lambda -1) \: F_{jk}(\lambda y) \: \gamma^j \: y^k \\
&&+ m^2 \; \frac{i}{128 \pi^3} \; \ln(|y^2|) \hinti d\lambda \;
	\epsilon(\lambda)
	\: \varepsilon^{ijkl} \: F_{ij}(\lambda y) \: y_k \; \rho \gamma_l \\
&&- m^2 \; \frac{1}{64 \pi^3} \; \ln(|y^2|) \hinti d\lambda \;
	\epsilon(\lambda)
	\:(\lambda^2 - \lambda) \; j_k(\lambda y) \: y^k \; y\slsh \\
&&+ {\cal{O}}(m^3) + {\cal{O}}(y^0) \spc .
\end{eqnarray*}
\end{Satz}
{\Beweis}
Folgt unmittelbar aus Lemma~\ref{a4_lemma1} und Satz~\ref{a3_satz98}.
\QED
\begin{Thm}
\label{a4_thm1}
In erster Ordnung St\"orungstheorie gilt
\begin{eqnarray}
\lefteqn{\Delta p_m(x,y) \;=\; \Delta p_0(x,y) } \nonumber \\
\label{eq:a4_111f}
&& - i e \left( \int_x^y A_j \right) \xi^j \; (p_m-p_0)(x,y) \\
&& - \frac{e}{32 \pi^3} \; m \; \ln(|\xi^2|) \; \int_x^y F_{ij} \:
	\sigma^{ij} \\
\label{eq:a4_113f}
&& - \frac{ie}{16 \pi^3} \; m \: \ln(|\xi^2|) \; \int_x^y (\alpha^2-\alpha) \:
	j_k \: \xi^k \\
&& - \frac{e}{32 \pi^3} \; m^2 \: \ln(|\xi^2|) \; \int_x^y (2 \alpha -1) \:
	\gamma^i \: F_{ij} \: \xi^j \\
&& - \frac{ie}{64 \pi^3} \; m^2 \: \ln(|\xi^2|) \; \int_x^y \varepsilon^{ijkl}
	\; F_{ij} \:  \xi_k \; \rho \gamma_l \\
\label{eq:a4_116f}
&& + \frac{e}{32 \pi^3} \; m^2 \: \ln(|\xi^2|) \; \int_x^y (\alpha^2-\alpha)
	\; j_k \: \xi^k \; \xi \slsh \\
&& + {\cal{O}}(m^3) + {\cal{O}}(\xi^0) \spc . \nonumber
\end{eqnarray}
\end{Thm}
{\Beweis}
Folgt direkt aus Satz~\ref{a4_lemma2}.
\QED

\section{St\"orungsrechnung f\"ur das Gravitationsfeld}
Wie in Abschnitt~\ref{grav_k0} arbeiten wir mit der linearisierten
Gravitationstheorie in symmetrischer Eichung und der Koordinatenbedingung
\Ref{a1_210}. Die St\"orung des Diracoperators \Ref{a1_208} f\"uhrt
analog wie in Abschnitt \ref{grav_km} auf
\Equ{a4_g29}
\Delta p_m \;=\; \left(\frac{1}{4} \: h(x) + \frac{3}{4} \: h(y) \right) \:
	p_m(x,y) - \frac{i}{e} \frac{\partial}{\partial y^k} \:
	\Delta p_m[\gamma^j h_j^k](x,y) \spc .
\EndEqu
Damit k\"onnen wir die Rechnung zum Teil auf Lemma~\ref{a4_lemma1} und
Theorem~\ref{a4_thm1} zur\"uck\-f\"uh\-ren.

\begin{Thm}
In erster Ordnung St\"orungstheorie gilt
\begin{eqnarray}
\Delta p_m(x,y) &=& \Delta p_0(x,y) \\
&&- \left( \int_x^y h^k_j \right) \: \xi^j \; \frac{\partial}{\partial y^k}
	\left( p_m(x,y)-p_0(x,y) \right) \\
&&+ \frac{i}{2} \: m \left( \int_x^y h_{ki,j} \right) \: \xi^k \;
	\sigma^{ij} \; \pe(x,y) \\
\label{eq:a4_gra}
&&+ \frac{1}{2} \: m \int_x^y (\alpha^2-\alpha) \;
	R_{jk} \; \xi^j \: \xi^k \; p^{(1)}(x,y) \\
&&+ \frac{1}{16 \pi^3} \: m \: \ln(|\xi^2|) \int_x^y
	(\alpha^2-\alpha+\frac{1}{4}) \; R \\
&&- \frac{1}{64 \pi^3} \: m \: \ln(|\xi^2|) \int_x^y (\alpha^4 - 2 \alpha^3 +
	\alpha^2) \; (\Box R_{jk}) \; \xi^j \: \xi^k \\
&&+ \frac{i}{32 \pi^3} \: m \: \ln(|\xi^2|) \int_x^y (\alpha^2-\alpha) \;
	R_{ki,j} \; \xi^k \; \sigma^{ij} \\
&&+ \frac{i}{16 \pi^3} \: m^2 \: \frac{1}{\xi^2} \int_x^y (2\alpha-1)
	\; (h_{jk,i} - h_{ik,j}) \: \gamma^i \; \xi^j \: \xi^k \\
&&- \frac{1}{16 \pi^3} \: m^2 \: \frac{1}{\xi^2} \int_x^y \varepsilon^{ijlm}
	h_{jk,i} \; \xi^k \: \xi_l \; \rho \gamma_m \\
\label{eq:a4_grb}
&&- \frac{i}{16 \pi^3} \: m^2 \: \frac{1}{\xi^2} \int_x^y (\alpha^2-\alpha)
	\; R_{jk} \; \xi^j \: \xi^k \: \xi \slsh \\
&&+ {\cal{O}}(\xi^0) + m^2 \: {\cal{O}}(\ln(|\xi^2|)) + {\cal{O}}(m^3)
	\spc . \nonumber
\end{eqnarray}
\end{Thm}
{\Beweis}
F\"ur diejenigen Summanden, die auf dem Lichtkegel wie $\xi^{-2}$ oder
$\xi^{-4}$ divergieren, k\"onnen wir genau wie im Beweis von
Theorem~\ref{thm_gp0} argumentieren: Die asymptotischen Entwicklungsformeln
f\"ur $\Delta p_m[\Aslsh]$ und $\Delta k_m[\Aslsh]$, Theorem~\ref{a4_thm1} und
Theorem~\ref{a2_theorem2}, gehen durch die formalen Ersetzungen
\Ref{a3_393a} ineinander \"uber. Durch Vergleich von \Ref{a4_g29} und
\Ref{a2_g10} \"ubertr\"agt sich daher auch das Ergebnis von Theorem \ref{thm_gpm}
durch diese Ersetzungen auf die St\"orungsrechnung f\"ur $\Delta p_m$.

Beachte, da{\ss} man bei den Termen der Ordnung $\ln(|\xi^2|)$ anders vorgehen
mu{\ss}. F\"ur ihre Berechnung ist nach \Ref{a4_g29} n\"amlich auch der Beitrag
von $\Delta p_m[\Aslsh]$ der Ordnung $\xi^2 \: \ln(|\xi^2|)$ wichtig, den
wir in Theorem~\ref{a4_thm1} jedoch nicht ber\"ucksichtigt haben.

Da wir die in $\xi^2$ logarithmisch divergenten Terme nur in erster Ordnung
in $m$ bestimmen wollen, m\"ussen wir nur noch $\Delta \pe$ untersuchen.
Nach \Ref{a4_243b}, \Ref{a3_151} und \Ref{a3_152} hat man
\begin{eqnarray}
i \: \frac{\partial}{\partial y^k} (s \: \gamma^j h^k_j \:
	p)^{(1)}(0,y)
&=& -\frac{1}{16 \pi^4} \: \frac{\partial}{\partial y^k} \left(
	\veeint^y (\Pdd h_j^k) \: \gamma^j - 2 \frac{\partial}{\partial y^j}
	\veeint^y h^{jk} \right) \nonumber \\
\label{eq:a4_g30}
&=& -\frac{1}{\pi^4} \frac{1}{y^4} \hinti d\lambda \;
	\frac{\epsilon(\lambda)}{\lambda^2} \dotint_{\lambda y} dz \;
	z^j \: z^k \; h_{jk} \\
\label{eq:a4_g31}
&&- \frac{1}{8 \pi^4} \frac{1}{y^2} \hinti \frac{d\lambda}{|\lambda|}
	\xdint_{\lambda y} dz \; z^j \: z^k \; \Box h_{jk} \\
\label{eq:a4_g32}
&&- \frac{1}{8 \pi^4} \frac{1}{y^2} \hinti \frac{d\lambda}{|\lambda|}
	\xdint_{\lambda y} dz \; z_k \: (\Pdd h^k_j) \; \gamma^j \\
\label{eq:a4_g33}
&&- \frac{1}{16 \pi^4} \veeint^y h_k[(\Pdd h^k_j)] \; \gamma^j \\
\label{eq:a4_g34}
&&+ \frac{1}{8 \pi^4} \veeint^y h_{jk}[h^{jk}] \spc .
\end{eqnarray}
Nun f\"uhren wir mit Hilfe von Satz~\ref{a3_satz99} und Satz~\ref{a3_satz98}
eine Entwicklung um den Lichtkegel durch. F\"ur die einzelnen Summanden
\Ref{a4_g30} bis \Ref{a4_g34} erh\"alt man:
\begin{eqnarray*}
\Ref{a4_g30} &=& -\frac{1}{4 \pi^3} \frac{1}{y^4} \hinti d\lambda \;
	\epsilon(\lambda) \; h_{jk}(\lambda y) \; y^j \: y^k \\
&&- \frac{1}{16 \pi^3} \frac{1}{y^2} \hinti d\lambda \; \epsilon(\lambda)
	\: \lambda^2 \; R_{jk}(\lambda y) \; y^j \: y^k \\
&&+ \frac{1}{32 \pi^3} \: \ln(|y^2|) \hinti d\lambda \; \epsilon(\lambda)
	\: \lambda^2 \; R(\lambda y) \\
&&- \frac{1}{128 \pi^3} \: \ln(|y^2|) \hinti d\lambda \; \epsilon(\lambda)
	\: \lambda^4 \; (\Box R_{jk})(\lambda y) \; y^j \: y^k
	\;+\; {\cal{O}}(y^0) \\
\Ref{a4_g31} &=& \frac{1}{16 \pi^3} \frac{1}{y^2} \hinti d\lambda \; |\lambda|
	\; R_{jk}(\lambda y) \; y^j \: y^k \\
&&-\frac{1}{32 \pi^3} \: \ln(|y^2|) \hinti d\lambda \; |\lambda| \;
	R(\lambda y) \\
&&+\frac{1}{64 \pi^3} \: \ln(|y^2|) \hinti d\lambda \; \epsilon(\lambda) \:
	\lambda^3 \; (\Box R_{jk})(\lambda y) \; y^j \: y^k \;+\;
	{\cal{O}}(y^0) \\
\Ref{a4_g32} &=& \frac{1}{16 \pi^3} \frac{1}{y^2} \; h(0)
\;+\; \frac{i}{16 \pi^3} \frac{1}{y^2} \hinti d\lambda \; \epsilon(\lambda)
	\; y^k \: h_{jk,i}(\lambda y) \; \sigma^{ij} \\
&&-\frac{i}{64 \pi^3} \: \ln(|y^2|) \hinti d\lambda \; \epsilon(\lambda)
	\: \lambda^2 \; R_{jk,i}(\lambda y) \; \sigma^{ij} \;+\;
	{\cal{O}}(y^0) \\
\Ref{a4_g33} &=& \frac{1}{128 \pi^3} \: \ln(|y^2|) \hinti d\lambda \;
	\epsilon(\lambda) \; R(\lambda y) \\
&&+ \frac{i}{64 \pi^3} \: \ln(|y^2|) \hinti d\lambda \; |\lambda| \;
	y^k \: R_{jk,i}(\lambda y) \; \sigma^{ij} \;+\; {\cal{O}}(y^0) \\
\Ref{a4_g34} &=& -\frac{1}{128 \pi^2} \: \ln(|y^2|) \hinti d\lambda \;
	\epsilon(\lambda) \: \lambda^2 \; (\Box R_{jk})(\lambda y) \; y^j
	\: y^k \;+\; {\cal{O}}(y^0)
\end{eqnarray*}
Aufsummieren dieser Terme liefert
\begin{eqnarray*}
\lefteqn{ i \: \frac{\partial}{\partial y^k} \: (s \: \gamma^j h_j^k \:
	p)^{(1)}(0,y)
\;=\; -\frac{1}{4 \pi^3} \frac{1}{y^4} \hinti d\lambda \;
	\epsilon(\lambda) \; h_{jk}(\lambda y) \; y^j \: y^k } \\
\spc &&+ \frac{i}{16 \pi^3} \frac{1}{y^2} \hinti d\lambda \; \epsilon(\lambda)
	\; y^k \: h_{jk,i}(\lambda y) \; \sigma^{ij} \\
&&- \frac{1}{16 \pi^3} \frac{1}{y^2} \hinti d\lambda \; \epsilon(\lambda)
	\: (\lambda^2-\lambda) \; R_{jk}(\lambda y) \; y^j \: y^k \\
&&+ \frac{1}{32 \pi^3} \: \ln(|y^2|) \hinti d\lambda \; \epsilon(\lambda)
	\: (\lambda^2-\lambda + \frac{1}{4}) \; R(\lambda y) \\
&&- \frac{1}{128 \pi^3} \: \ln(|y^2|) \hinti d\lambda \; \epsilon(\lambda)
	\: (\lambda^4-2\lambda^3+\lambda^2) \; (\Box R_{jk})(\lambda y)
	\; y^j \: y^k \\
&&- \frac{i}{64 \pi^3} \: \ln(|y^2|) \hinti d\lambda \; \epsilon(\lambda)
	\: (\lambda^2-\lambda) \; R_{jk,i}(\lambda y) \; \sigma^{ij} \\
&&- \frac{1}{4} \: h(0) \: p^{(1)}(0,y) + {\cal{O}}(y^0) \spc .
\end{eqnarray*}
Durch Translation erh\"alt man eine Gleichung f\"ur $i \partial_k
(s \: \gamma^j h_j^k \: p)^{(1)}(x,y)$. Die entsprechende Formel, bei
der $s_m$ und $p_m$ vertauscht sind, 
leitet man daraus wieder ab, indem man $x$ durch $y$ ersetzt und
umgekehrt\footnote{
$A^\ast$ bezeichnet die Adjungierte der
$4\times4$-Matrix $A$ bzgl. des Spinskalarproduktes.}:
\begin{eqnarray*}
\lefteqn{ i \: \frac{\partial}{\partial y^k} \: (p_m \: \gamma^j h_j^k \:
	s_m)(x,y) \;=\; i \: \frac{\partial}{\partial y^k} \: (s_m \:
	\gamma^j  h_j^k \: p_m)(y,x)^\ast }\\
&=& i (s_m \; \gamma^j h_{jk,}^{\;\;\;k} \; p_m)(y,x)^\ast -
	 i \: \frac{\partial}{\partial x^k} \: (s_m \: \gamma^j h_j^k \:
	p_m)(y,x)^\ast \\
&=& \frac{1}{2} (s_m \; (i \Pdd h) \; p_m)(y,x)^\ast -
	 i \: \frac{\partial}{\partial x^k} \: (s_m \: \gamma^j h_j^k \:
	p_m)(y,x)^\ast \\
&=& \frac{1}{2} \: h(y) \; p_m(x,y) - i \: \frac{\partial}{\partial x^k} \:
	(s_m \: \gamma^j h_j^k \: p_m)(y,x)^\ast \spc .
\end{eqnarray*}
Durch Einsetzen in \Ref{a4_a10} und \Ref{a4_g29} folgt die Behauptung.
\QED

\section{Axiale St\"orung}
Wir betrachten die axiale St\"orung \Ref{a2_x0}. F\"ur die Auswirkung der
St\"orung auf $p_m$ hat man in erster Ordnung
\Equ{a4_y0}
\Delta p_m \;=\; -e \: ( p_m \: \rho \Aslsh \: s_m + s_m \: \rho \Aslsh \: p_m
	) \spc .
\EndEqu
Wir berechnen zun\"achst $s_m \: \rho \Aslsh \: p_m$:
\begin{Lemma}
Es gilt
\begin{eqnarray}
\label{eq:a4_y1}
(s_m \; \rho \Aslsh \: p_m)(0,y) &=& -\rho \; (s_m \: \Aslsh \: p_m)(0,y) \\
\label{eq:a4_y2}
&& - \frac{im}{4 \pi^4} \frac{1}{y^2} \hinti \frac{d\lambda}{|\lambda|}
	\xdint_{\lambda y} dz \; \rho \Aslsh \; z\slsh \\
\label{eq:a4_y3}
&& - \frac{im}{16 \pi^4} \veeint^y dz \int_0^1 d\alpha \; \alpha^2 \;
	j_k(\alpha z) \; z^k \; \rho \\
\label{eq:a4_y4}
&& - \frac{im}{16 \pi^4} \veeint^y \partial_j A^j \; \rho \\
\label{eq:a4_y5}
&& + \frac{m}{8 \pi^4} \veeint^y h_j[A_k] \; \rho \sigma^{jk} \\
\label{eq:a4_y6}
&& + \frac{m^2}{8 \pi^4} \veeint^y \rho \Aslsh \\
&& + {\cal{O}}(m^3) \spc . \nonumber
\end{eqnarray}
\end{Lemma}
{\Beweis}
Wir betrachten die Beitr\"age verschiedener Ordnung in $m$ nacheinander:
\begin{enumerate}
\item Terme $\sim m^0$:
\[  (s_0 \: \rho \Aslsh \: p_0)(0,y) \;=\; -\rho \; (s_0 \: \Aslsh \:
	p_0)(0,y) \]
\item Terme $\sim m$:
\begin{eqnarray}
\label{eq:a4_305a}
\lefteqn{ (s_m \: \rho \Aslsh \: p_m)(0,y) \;=\; \frac{im}{16 \pi^4} \left(
	\Pdd_x (\veeint_x^y \rho \Aslsh) + (\veeint_x^y \rho \Aslsh) \Pdd_y
	\right)_{|x=0} } \\
&=& -\rho \: (s_m \: \Aslsh \: p_m)(0,y) - \frac{im}{8 \pi^4} \:
	\frac{\partial}{\partial y^j} \veeint^y \rho \Aslsh \: \gamma^j
	\nonumber \\
&=& -\rho \: (s_m \: \Aslsh \: p_m)(0,y) - \frac{im}{4 \pi^4} \frac{1}{y^2}
	\hinti \frac{d\lambda}{|\lambda|} \xdint_{\lambda y} dz \;
	\rho \Aslsh \: z\slsh - \frac{im}{8 \pi^4} \: \rho \veeint^y
	h_j[\Aslsh] \: \gamma^j \nonumber
\end{eqnarray}
Setze nun noch \Ref{a2_x10}, \Ref{a2_x11} ein.
\item Terme $\sim m^2$:
\begin{eqnarray*}
(s_m \: \rho \Aslsh \: p_m)(x,y) &=& 
\frac{m^2}{16 \pi^4} \left( \veeint_x^y \rho \Aslsh + \frac{1}{2} \Pdd_x
\left( \veeint_x^y dz \; \rho \Aslsh(z) \; (z-y)_j \gamma^j \right) \right. \\
&&\spc\spc \left. + \frac{1}{2} \left( \veeint_x^y
dz \; (x-z)_j \gamma^j \: \rho \Aslsh(z) \right) \: \Pdd_y \right) \\
&=& - \rho \: (s_m \: \Aslsh \: p_m)(x,y) + \frac{m^2}{8 \pi^4} \veeint_x^y
	\rho \Aslsh
\end{eqnarray*}
\end{enumerate}
\QED
Nun f\"uhren wir wieder eine Entwicklung um den Lichtkegel durch:
\begin{Lemma}
\label{a4_lemma20}
Es gilt
\begin{eqnarray*}
\lefteqn{ (s_m \: \rho \Aslsh \: p_m)(0,y) \;=\; - \rho \: (s_m \: \Aslsh \:
p_m)(0,y) } \\
&& - \frac{im}{8 \pi^3} \frac{1}{y^2} \hinti d\lambda \; \epsilon(\lambda) \;
	\rho \Aslsh(\lambda y) \: y\slsh \\
&& + \frac{im}{32 \pi^3} \: \ln(|y^2|) \hinti d\lambda \; \epsilon(\lambda) \;
	(\lambda^2-\lambda) \: j_k(\lambda y) \: y^k \; \rho \\
&& - \frac{im}{32 \pi^3} \: \ln(|y^2|) \hinti d\lambda \; \epsilon(\lambda) \;
	\partial_j A^j(\lambda y)  \; \rho \\
&& + \frac{m}{32 \pi^3} \: \ln(|y^2|) \hinti d\lambda \; |\lambda| \;
	F_{jk}(\lambda y) \; \rho \sigma^{jk} \\
&& - \frac{m}{32 \pi^3} \: \ln(|y^2|) \hinti d\lambda \; \epsilon(\lambda)
	\; (\lambda^2-\lambda) \; \Box A_j(\lambda y) \: y_k \; \rho
	\sigma^{jk} \\
&& + \frac{m^2}{16 \pi^3} \: \ln(|y^2|) \hinti d\lambda \; \epsilon(\lambda)
	\; \rho \Aslsh(\lambda y) \\
&& + {\cal{O}}(m^3) + {\cal{O}}(y^0) \spc .
\end{eqnarray*}
\end{Lemma}
{\Beweis}
Die asymptotischen Enwicklungen von Satz \ref{a3_satz98} und \ref{a3_satz99}
liefern:
\begin{eqnarray*}
\lefteqn{ \frac{1}{y^2} \hinti \frac{d\lambda}{|\lambda|} \xdint_{\lambda y}
	dz \; \Aslsh \: z\slsh
\;=\; \frac{\pi}{2} \frac{1}{y^2} \hinti d\lambda \; \epsilon(\lambda) \;
	\Aslsh(\lambda y) \; y\slsh } \\
&& + \frac{\pi}{8} \: \ln(|y^2|) \hinti d\lambda \; |\lambda| \; \Box_z
	\left( A_j \: z^j \right)_{|z=\lambda y} \\
&& - \frac{i \pi}{8} \: \ln(|y^2|) \hinti d\lambda \; |\lambda| \; \Box_z
	\left( A_j \: z_k \right)_{|z=\lambda y} \; \sigma^{jk} 
	+ {\cal{O}}(y^0)\\
&=& \frac{\pi}{2} \frac{1}{y^2} \hinti d\lambda \; \epsilon(\lambda) \;
	\Aslsh(\lambda y) \: y\slsh \\
&& - \frac{\pi}{8} \: \ln(|y^2|) \hinti d\lambda \; \epsilon(\lambda) \:
	\lambda^2 \: j_k(\lambda y) \: y^k \\
&& + \frac{i \pi}{8} \: \ln(|y^2|) \hinti d\lambda \; |\lambda| \;
	F_{jk}(\lambda y) \; \sigma^{jk} \\
&& - \frac{i \pi}{8} \: \ln(|y^2|) \hinti d\lambda \;
	\epsilon(\lambda) \: \lambda^2 \;
	\Box A_j(\lambda y) \: y_k \; \sigma^{jk} + {\cal{O}}(y^0) \\
\lefteqn{ \veeint^y dz \int_0^1 d\alpha \; \alpha^2 \; j_k(\alpha z) \; z^k
	} \\
&=& \frac{\pi}{2} \: \ln(|y^2|) \hinti d\lambda \; \epsilon(\lambda)
	\int_0^1 d\alpha \; \alpha^2 \; j_k(\alpha \lambda y) \;
	\lambda y^k  + {\cal{O}}(y^0) \\
&=& \frac{\pi}{2} \: \ln(|y^2|) \hinti d\lambda \; |\lambda| \;
	j_k(\lambda y) \; y^k + {\cal{O}}(y^0) \\
\lefteqn{ \veeint^y h_j[A_k]
\;=\; \frac{\pi}{2} \ln(|y^2|) \hinti d\lambda \; \epsilon(\lambda) \left(
	\partial_j A_k(\lambda y) - \int_0^1 d\alpha \; \alpha \: \partial_j
	A_k(\alpha \lambda y) \right. } \\
&& \left. \spc - \frac{1}{2} \lambda y_j \int_0^1 d\alpha \;
	\alpha^2 \; \Box A_k(\alpha \lambda y) \right) + {\cal{O}}(y^0) \\
&=& - \frac{\pi}{4} \: \ln(|y^2|) \hinti d\lambda \; |\lambda| \;
	y_j \: (\Box A_k)(\lambda y) + {\cal{O}}(y^0)
\end{eqnarray*}
\QED

\begin{Thm}
In erster Ordnung St\"orungstheorie gilt
\begin{eqnarray}
\label{eq:a4_z1}
\Delta p_m(x,y) &=& - \rho \: \Delta p_0[\Aslsh](x,y)  \\
\label{eq:a4_z2}
&& - \frac{ie}{4 \pi^3} \: m \: \frac{1}{\xi^2} \int_x^y \rho \:
	\frac{1}{2} [\xi\slsh, \Aslsh] \\
\label{eq:a4_z3}
&& -\frac{e}{32 \pi^3} \: m \: \ln(|\xi^2|) \int_x^y (2\alpha -1) \;
	F_{jk} \; \rho \sigma^{jk} \\
\label{eq:a4_z4}
&& +\frac{ie}{16 \pi^3} \: m \:  \ln(|\xi^2|) \int_x^y \partial_j A^j \;
	\rho \\
\label{eq:a4_z5}
&& +\frac{e}{16 \pi^3} \: m \:  \ln(|\xi^2|) \int_x^y (\alpha^2-\alpha) \;
	\Box A_j \; \xi_k \; \rho \sigma^{jk} \\
\label{eq:a4_z6}
&& +\frac{e}{8 \pi^3} \: m^2 \: \frac{1}{\xi^2} \int_x^y A_j \: \xi^j \;
	\rho \xi\slsh \\
\label{eq:a4_z7}
&& -\frac{e}{8\pi^3} \: m^2 \: \ln(|\xi^2|) \int_x^y \rho \Aslsh \\
\label{eq:a4_z8}
&& +\frac{e}{32 \pi^3} \: m^2 \: \ln(|\xi^2|) \int_x^y (2\alpha-1) \;
	F_{jk} \: \xi^k \; \rho \gamma^j \\
\label{eq:a4_z9}
&& +\frac{ie}{64 \pi^3} \: m^2 \: \ln(|\xi^2|) \int_x^y \varepsilon^{ijkl}
	\; F_{ij} \: \xi_k \; \gamma_l \\
\label{eq:a4_z10}
&& -\frac{e}{32 \pi^3} \: m^2 \: \ln(|\xi^2|) \int_x^y (\alpha^2-\alpha) \;
	j_k \: \xi^k \; \rho \xi\slsh \\
\label{eq:a4_z11}
&& +\frac{ie}{16 \pi^3} \: m^3 \: \ln(|\xi^2|) \int_x^y \rho\:
	\frac{1}{2} [\xi\slsh, \Aslsh] \\
&& + {\cal{O}}(m^4) + {\cal{O}}(\xi^0) \spc . \nonumber
\end{eqnarray}
\end{Thm}
{\Beweis}
Wir setzen das Ergebnis von Satz~\ref{a4_lemma2} in \ref{a4_lemma20} ein.
Die Beitr\"age $\sim m^0$, $m^2$ erh\"alt man unter Verwendung von \Ref{a4_y0}
direkt durch Translation und Vertauschung von $x$ und $y$.
F\"ur die Ordnung $\sim m$ beachte man das Zwischenresultat
\begin{eqnarray*}
(s_m \: \rho \Aslsh \: p_m)(0,y) &=& -\rho \: (s_0 \: \Aslsh \: p_0)(0,y) \\
&& +\frac{im}{8 \pi^3} \frac{1}{y^2} \hinti d\lambda \; \epsilon(\lambda)
	\: \rho \; \frac{1}{2} [\xi\slsh,\Aslsh(\lambda y)] \\
&& +\frac{m}{64 \pi^3} \: \ln(|y^2|) \hinti d\lambda \; \epsilon(\lambda)
	\: (2\lambda-1) \; F_{jk}(\lambda y) \; \rho \sigma^{jk} \\
&& -\frac{im}{32 \pi^3} \: \ln(|y^2|) \hinti d\lambda \; \epsilon(\lambda)
	\: \partial_j A^j(\lambda y) \; \rho \\
&& -\frac{m}{32 \pi^3} \: \ln(|y^2|) \hinti d\lambda \; \epsilon(\lambda)
	\: (\lambda^2-\lambda) \; (\Box A_j)(\lambda y) \: y_k \;
	\rho \sigma^{jk} \\
&& + {\cal{O}}(m^2) + {\cal{O}}(\xi^0) \spc .
\end{eqnarray*}
Der Beitrag $\sim m^3$ folgt direkt aus der Form der Pseudoeichterme
\Ref{a2_x20} (die bei $p_m$ die gleiche Gestalt haben wie bei $k_m$);
die Stomterme $\sim m^3$ sind bereits von der Ordnung $\xi^0$.
\QED
Die Summanden \Ref{a4_z2}, \Ref{a4_z6} und \Ref{a4_z11} sind Pseudoeichterme,
\Ref{a4_z7} ist der Massenterm.

Beachte, da{\ss} der Stromterm \Ref{a4_z5} anstelle des Maxwellstromes
den Ausdruck $\Box A_j$ enth\"alt. Man sieht daran, da{\ss} bei axialen
St\"orungen keine lokale Eichinvarianz vorhanden ist.

\section{Skalare St\"orung}
Wir betrachten wieder die skalare St\"orung~\Ref{a1_s1}. F\"ur
$\Delta p_m$ hat man in erster Ordnung in $\Xi$
\Equ{a4_spm}
	\Delta p_m \;=\; - \left( s_m \: \Xi \: p_m + p_m \: \Xi \: s_m
	\right) \spc .
\EndEqu

\begin{Thm}
\label{a4_theorem_spm}
In erster Ordnung St\"orungstheorie gilt
\begin{eqnarray}
\Delta p_m(x,y) &=& \Delta p_0(x,y) \\
\label{eq:a4_s2}
&&- 2 m \left( \int_x^y \Xi \right) \; p^{(2)}(x,y) \\
&&+ \frac{i}{16 \pi^3} \: m \: \ln(|\xi^2|) \int_x^y (2\alpha-1) \;
	(\Pdd \Xi) \\
&&+ \frac{i}{16 \pi^3} \: m \: \ln(|\xi^2|) \int_x^y (\alpha^2-\alpha) \;
	(\Box \Xi) \; \xi\slsh \\
&&- \frac{1}{32 \pi^3} \: m^2 \: \ln(|\xi^2|) \; (\Xi(y) + \Xi(x)) \\
&&- \frac{1}{8 \pi^3} \: m^2 \: \ln(|\xi^2|) \int_x^y \Xi \\
&&- \frac{i}{32 \pi^3} \: m^2 \: \ln(|\xi^2|) \int_x^y 
	(\partial_j \Xi) \: \xi_k \; \sigma^{jk} \\
&&+ {\cal{O}}(m^3) + {\cal{O}}(\xi^0) \spc . \nonumber
\end{eqnarray}
\end{Thm}
{\Beweis}
Untersuche die Terme verschiedener Ordnung in $m$ nacheinander.
\begin{enumerate}
\item Terme $\sim m$:
Unter Verwendung von \Ref{a3_151} und den asymptotischen Entwicklungen von 
Satz~\ref{a3_satz98} hat man
\begin{eqnarray*}
\lefteqn{ (s \: \Xi \: p)^{(1)}(0,y) \;=\; \frac{i}{16 \pi^4} \left(
	\Pdd_x (\veeint_x^y \Xi) + (\veeint_x^y \Xi) \Pdd_y \right)_{|x=0} }
	\\
&=& \frac{i}{16 \pi^4} \left( \veeint^y (\Pdd \Xi) - 2 \: \Pdd_y
	\veeint^y \Xi \right) \\
&=& -\frac{i}{4 \pi^4} \frac{1}{y^2} \hinti \frac{d\lambda}{|\lambda|}
	\xdint_{\lambda y} dz \; z \slsh \: \Xi(z)
	\;+\; \frac{i}{16 \pi^4} \veeint^y \left( (\Pdd \Xi) - 2 \: \gamma^j
	\: h_j[\Xi] \right) \\
&=& -\frac{i}{8 \pi^3} \frac{1}{y^2} \hinti d\lambda \; y\slsh \;
	\Xi(\lambda y) \\
&&-\frac{i}{32 \pi^3} \: \ln(|\xi^2|) \hinti d\lambda \; \epsilon(\lambda)
	\: (2\lambda-1) \; (\Pdd \Xi)(\lambda y) \\
&&- \frac{i}{32 \pi^3} \: \ln(|\xi^2|) \hinti d\lambda \; \epsilon(\lambda)
	\: (\lambda^2-\lambda) \; (\Box \Xi)(\lambda y) \;
	\xi\slsh \;+\; {\cal{O}}(y^0)
\spc .
\end{eqnarray*}
\item Terme $\sim m^2$:
\begin{eqnarray*}
(s \: \Xi \: p)^{(2)}(x,y) &=& \frac{1}{16 \pi^4} \left\{ \veeint_x^y \Xi
	+ \frac{1}{2} \: \Pdd_x \left( \veeint_x^y dz \; \Xi(z) \; (z-y)^j
	\gamma_j \right) \right. \\
&& \left. \spc + \frac{1}{2} \left( \veeint_x^y dz \; (x-z)^j \gamma_j
	\; \Xi(z) \right) \: \Pdd_y \right\}
\end{eqnarray*}
Mit Hilfe von~\Ref{a4_246z} und~\Ref{a4_246y} erh\"alt man die Umformungen
\begin{eqnarray*}
\Pdd_x \left( \veeint_x^y dz \: \Xi(z) \; (z-y)^j \gamma_j \right)
	&=& \veeint_x^y dz \;
	(\Pdd \Xi)(z) \; (z-y)^j \gamma_j \;+\; 2 \veeint_x^y \Xi \\
\left( \veeint_x^y dz \; (x-z)^j \gamma_j \; \Xi(z) \right) \: \Pdd_y
	&=& - \veeint_x^y dz
	\; (x-z)^j \gamma_j \; (\Pdd \Xi)(z) \;+\; 2 \veeint_x^y \Xi
\end{eqnarray*}
und somit
\begin{eqnarray*}
\lefteqn{ (s \: \Xi \: p)^{(2)}(0,y) \;=\; \frac{3}{16 \pi^4} \veeint^y \Xi
	} \\
&&+ \frac{1}{32 \pi^4} \veeint^ydz \: \left( (\Pdd \Xi)(z) \; (z-y)^j \gamma_j
	- (x-z)^j \gamma_j \; (\Pdd \Xi)(z) \right) \\
&=& \frac{1}{32 \pi^3} \: \ln(|y^2|) \hinti d\lambda \; \epsilon(\lambda)
	\; \Xi(\lambda y) \\
&&+ \frac{1}{64 \pi^3} \: \ln(|y^2|) \hinti d\lambda \; \epsilon(\lambda)
	\: (2\lambda -1) \; (\partial_j \Xi)(\lambda z) \; y^j \\
&&+ \frac{i}{64 \pi^3} \: \ln(|y^2|) \hinti d\lambda \; \epsilon(\lambda)
	\; (\partial_j \Xi)(\lambda y) \; y_k \; \sigma^{jk}
	\;+\; {\cal{O}}(y^0) \\
&=& \frac{1}{32 \pi^3} \: \ln(|y^2|) \; \Xi(0) \\
&&+ \frac{1}{16 \pi^3} \: \ln(|y^2|) \hinti d\lambda \; \epsilon(\lambda)
	\; \Xi(\lambda y) \\
&&+ \frac{i}{64 \pi^3} \: \ln(|y^2|) \hinti d\lambda \; \epsilon(\lambda)
	\; (\partial_j \Xi)(\lambda y) \; y_k \; \sigma^{jk}
	\;+\; {\cal{O}}(y^0) \spc .
\end{eqnarray*}
\end{enumerate}
Durch Translation um $x$ und Vertauschung von $x$, $y$ erh\"alt man
Gleichungen f\"ur $(s_m \: \Xi \: p_m)(x,y)$ und $(p_m \: \Xi \: s_m)(x,y)$.
Setze nun in \Ref{a4_spm} ein.
\QED

\section{Pseudoskalare St\"orung}
Wir betrachten wieder die pseudoskalare St\"orung~\Ref{a1_ps1}. F\"ur
$\Delta p_m$ hat man in erster Ordnung in $\Xi$
\Equ{a4_pspm}
	\Delta p_m \;=\; - i \left( s_m \: \rho \Xi \: p_m + p_m \: \rho \Xi
	\: s_m \right) \spc .
\EndEqu

\begin{Thm}
\label{a4_theorem_pspm}
In erster Ordnung St\"orungstheorie gilt
\begin{eqnarray}
\Delta p_m(x,y) &=& -i \rho \: \Delta p_0[\Xi](x,y) \\
&&- \frac{1}{16 \pi^3} \: m \: \rho \: \ln(|\xi^2|) \int_x^y (\Pdd \Xi) \\
&&+ \frac{i}{32 \pi^3} \: m^2 \: \ln(|\xi^2|) \; (\Xi(y) + \Xi(x)) \; \rho \\
&&- \frac{1}{32 \pi^3} \: m^2 \: \ln(|\xi^2|) \int_x^y 
	(\partial_j \Xi) \: \xi_k \; \rho \sigma^{jk} \\
&&+ {\cal{O}}(m^3) + {\cal{O}}(\xi^0) \spc . \nonumber
\end{eqnarray}
\end{Thm}
{\Beweis}
Untersuche die Terme verschiedener Ordnung in $m$ nacheinander.
\begin{enumerate}
\item Terme $\sim m$:
Unter Verwendung von \Ref{a3_151} und den asymptotischen Entwicklungen von 
Satz~\ref{a3_satz98} hat man
\begin{eqnarray*}
\lefteqn{ (s \: \rho \Xi \: p)^{(1)}(0,y) \;=\; - \frac{1}{16 \pi^4} \left(
	\Pdd_x (\veeint_x^y \rho \Xi) + (\veeint_x^y \rho \Xi) \Pdd_y
	\right)_{|x=0} } \\
&=& \frac{1}{16 \pi^4} \: \rho \left(
	\Pdd_x (\veeint_x^y \rho \Xi) + \Pdd_y (\veeint_x^y \rho \Xi)
	\right)_{|x=0} \\
&=& \frac{1}{16 \pi^4} \: \rho \veeint^y (\Pdd \Xi) \spc .
\end{eqnarray*}
\item Terme $\sim m^2$:
\begin{eqnarray*}
(s \: \rho \Xi \: p)^{(2)}(x,y) &=& \frac{i}{16 \pi^4} \left\{ \veeint_x^y
	\rho \Xi
	+ \frac{1}{2} \: \Pdd_x \left( \veeint_x^y dz \; \rho \Xi(z) \; (z-y)^j
	\gamma_j \right) \right. \\
&& \left. \spc + \frac{1}{2} \left( \veeint_x^y dz \; (x-z)^j \gamma_j
	\; \rho \Xi(z) \right) \: \Pdd_y \right\} \\
&=& - i \rho \: (s \: \Xi \: p)^{(2)}(x,y) + \frac{i}{8 \pi^4} \:
	\rho \veeint_x^y \Xi 
\end{eqnarray*}
\end{enumerate}
Durch Translation um $x$ und Vertauschung von $x$, $y$ erh\"alt man
Gleichungen f\"ur $(s_m \: \rho \Xi \: p_m)(x,y)$ und $(p_m \: \rho \Xi \: s_m)(x,y)$.
Setze nun in \Ref{a4_pspm} ein.
\QED

\section{Bilineare St\"orung}
Wir betrachten wieder die bilineare St\"orung~\Ref{a1_b00}.
F\"ur $\Delta p_m$ hat man in erster Ordnung in $B$
\Equ{a4_bpm}
  \Delta p_m \;=\; - \left( s_m \: B_{jk} \: \sigma^{jk} \; p_m +
	p_m \: B_{jk} \: \sigma^{jk} \: s_m \right) \spc .
\EndEqu
\begin{Thm}
In erster Ordnung St\"orungstheorie gilt
\begin{eqnarray}
\Delta p_m(x,y) &=& \Delta p_0(x,y) \\
&&+ \frac{i}{4 \pi^3} \: m \: \frac{1}{\xi^2} \int_x^y
	\varepsilon^{ijkl} \; B_{ij} \: \xi_k \; \rho \gamma_l \\
&&- \frac{1}{8 \pi^3} \: m \; \ln(|\xi^2|) \int_x^y
	B_{jk,}^{\;\;\;k} \; \gamma^j \\
&&+ \frac{i}{16 \pi^3} \: m \; \ln(|\xi^2|) \int_x^y (2\alpha-1) \;
	\varepsilon^{ijkl} \: B_{ij,k} \; \rho \gamma_l \\
&&+ \frac{i}{16 \pi^3} \: m \; \ln(|\xi^2|) \int_x^y (\alpha^2-\alpha) \;
	\varepsilon^{ijkl} \; (\Box B_{ij}) \: \xi_k
	\; \rho \gamma_l \\
&&- \frac{1}{4 \pi^3} \: m^2 \: \frac{1}{\xi^2} \int_x^y
	B_{ij} \: \xi^i \: \xi_k \; \sigma^{jk} \\
&&- \frac{1}{32 \pi^3} \: m^2 \; \ln(|\xi^2|) \;
	(B_{jk}(y) + B_{jk}(x)) \; \sigma^{jk} \\
&&- \frac{i}{16 \pi^3} \: m^2 \; \ln(|\xi^2|) \int_x^y
	\xi_j \; B^{jk}_{\;\;\:,k} \\
&&- \frac{1}{16 \pi^3} \: m^2 \; \ln(|\xi^2|)
	\int_x^y (2\alpha-1) \; (\xi^k \: B_{jk,i} + \xi_i \:
	B_{jk,}^{\;\;\;k}) \; \sigma^{ij} \\
&&- \frac{1}{16 \pi^3} \: m^2 \; \ln(|\xi^2|)
	\int_x^y (\alpha^2-\alpha) \; (\Box B_{ij}) \: \xi^i \: \xi_k \;
	\sigma^{jk} \\
&&+ \frac{1}{32 \pi^3} \: m^2 \; \ln(|\xi^2|)
	\int_x^y \varepsilon^{ijkl} \; B_{ij,k} \: \xi_l \; \rho \\
&&+ {\cal{O}}(\xi^0) + {\cal{O}}(m^3) \spc . \nonumber
\end{eqnarray}
\end{Thm}
{\Beweis}
Wir betrachten die Beitr\"age verschiedener Ordnung in $m$ nacheinander:
\begin{enumerate}
\item Terme $\sim m$:
\begin{eqnarray*}
(s_m \: B_{jk} \: \sigma^{jk} \: p_m)(x,y) &=& \frac{im}{16 \pi^4}
	\left( \Pdd_x (\veeint_x^y B_{jk} \: \sigma^{jk}) +
	(\veeint_x^y B_{jk} \: \sigma^{jk}) \Pdd_y \right)
\end{eqnarray*}
Setze~\Ref{a2_gam1} und~\Ref{a2_gam2} ein
\begin{eqnarray*}
&=& -\frac{m}{8 \pi^4} \: \frac{\partial}{\partial x^j} \veeint_x^y B^{jk} \:
	\gamma_k + \frac{im}{16 \pi^4} \: \varepsilon^{ijkl} \:
	\frac{\partial}{\partial x^i} \veeint_x^y B_{jk} \; \rho \gamma_l \\
&& -\frac{m}{8 \pi^4} \: \frac{\partial}{\partial y^j} \veeint_x^y B^{jk} \:
	\gamma_k - \frac{im}{16 \pi^4} \: \varepsilon^{ijkl} \:
	\frac{\partial}{\partial y^i} \veeint_x^y B_{jk} \; \rho \gamma_l\\
&=& -\frac{m}{8 \pi^4} \veeint_x^y B^{jk}_{\;\;\;,j} \: \gamma_k +
	\frac{im}{16 \pi^4} \: \varepsilon^{ijkl} \veeint_x^y B_{ij,k} \;
	\rho \gamma_l - \frac{im}{8 \pi^4} \: \varepsilon^{ijkl} \:
	\frac{\partial}{\partial y^i} \veeint_x^y B_{jk} \; \rho \gamma_l
	\;\; .
\end{eqnarray*}
Jetzt wenden wir \Ref{a3_151} an und f\"uhren mit Satz~\ref{a3_satz99} und
\ref{a3_satz98} eine Entwicklung um den Lichtkegel durch. Setze der
Einfachheit halber $x=0$.
\begin{eqnarray*}
\lefteqn{ (s_m \: B_{jk} \: \sigma^{jk} \: p_m)(x,y) \;=\; \frac{m}{8 \pi^4}
	\veeint^y B^{jk}_{\;\;\;,k} \: \gamma_j } \\
&&- \frac{im}{4 \pi^4} \frac{1}{y^2} \hinti \frac{d\lambda}{|\lambda|}
	\xdint_{\lambda y} dz \; \varepsilon^{ijkl} \; z_i \: B_{jk} \;
	\rho \gamma_l \\
&&- \frac{im}{16 \pi^4} \: \varepsilon^{ijkl} \veeint^y dz \: \left(
	B_{ij,k}(z) - 2 \int^z \alpha \; B_{ij,k} \right) \: \rho \gamma_l \\
&&+ \frac{im}{16 \pi^4} \: \varepsilon^{ijkl} \veeint^y dz \; \zeta_i
	\int^z \alpha^2 \; (\Box B_{jk}) \; \rho \gamma_l \\
&=& -\frac{im}{8 \pi^3} \frac{1}{y^2} \hinti d\lambda \; \epsilon(\lambda)
	\: \varepsilon^{ijkl} \; B_{ij}(\lambda y) \: \xi_k \; \rho \gamma_l \\
&&+ \frac{m}{16 \pi^3} \: \ln(|y^2|) \hinti d\lambda \; \epsilon(\lambda)
	\: B_{jk,}^{\;\;\;k}(\lambda y) \; \gamma^j \\
&&- \frac{im}{32 \pi^3} \: \ln(|y^2|) \hinti d\lambda \; \epsilon(\lambda)
	\: (2\lambda-1) \: \varepsilon^{ijkl} \; B_{ij,k}(\lambda y) \;
	\rho \gamma_l \\
&&- \frac{im}{32 \pi^3} \: \ln(|y^2|) \hinti d\lambda \; \epsilon(\lambda)
	\: (\lambda^2-\lambda) \: \varepsilon^{ijkl} \; (\Box B_{ij}(\lambda
	y) \: \xi_k \; \rho \gamma_l \;+\; {\cal{O}}(\xi^0)
\end{eqnarray*}
\item Terme $\sim m^2$:
\begin{eqnarray}
\lefteqn{ (s_m \: B_{jk} \: \sigma^{jk} \: p_m)(0,y) \;=\;
\frac{1}{16 \pi^4} \: m^2 \lint_x^y
	B_{jk} \: \sigma^{jk} } \nonumber \\
&&\spc\spc + \frac{m^2}{32 \pi^4} \: \Pdd_x \left( \lint_x^y dz
	\; B_{jk}(z) \: \sigma^{jk} \; (z-y)^m \gamma_m \right) \nonumber \\
\label{eq:a4_bi0}
&&\spc\spc + \frac{m^2}{32 \pi^4} \left( \lint_x^y dz \;
	(x-z)^m \gamma_m \; B_{jk}(z) \: \sigma^{jk} \right) \Pdd_y
\end{eqnarray}
Wir berechnen nun den zweiten und dritten Summanden. Dazu l\"osen wir das
Produkt der Diracmatrizen mit Hilfe von~\Ref{a1_b8} auf, berechnen mit
\Ref{a3_151} die Ableitungen und f\"uhren mit Satz~\ref{a3_satz99}
und~\ref{a3_satz98} eine Entwicklung um den Lichtkegel durch.
Setze wieder $x=0$ und beachte die Identit\"aten
\begin{eqnarray*}
\lefteqn{ \frac{\partial}{\partial y^m} \veeint_x^y dz \; B_{jk} \:
	\sigma^{jk} \; (z-y)^m } \\
&=& \frac{\partial}{\partial y^m} \int d^4z \;
	\frac{(z-y)^m}{(y-z)^2} \: l(z-x) \; B_{jk}(z) \: \sigma^{jk}
\;=\; -2 \veeint_x^y B_{jk} \: \sigma^{jk} \\
\lefteqn{ \frac{\partial}{\partial x^m} \veeint_x^y dz \; (x-z)^m \;
	B_{jk}(z) \: \sigma^{jk} } \\
&=& \frac{\partial}{\partial x^m} \int d^4z \;
	(x-z)^m \: l(z-x) \; \frac{1}{(y-z)^2} \; B_{jk}(z) \: \sigma^{jk}
\;=\; 2 \veeint_x^y B_{jk} \: \sigma^{jk} \spc .
\end{eqnarray*}
Man erh\"alt auf diese Weise
\begin{eqnarray*}
\lefteqn{ \Pdd_x \left(\veeint_x^y dz \; B_{jk}(z) \: \sigma^{jk} \;
	(z-y)^m \gamma_m \right)_{|x=0}
\;=\; 2 i \: \frac{\partial}{\partial x^j} \veeint_x^y dz \; B^{jk}(z) \:
	(z-y)_k } \\
&&+ \frac{\partial}{\partial x^m} \veeint_x^y dz \; B_{jk}(z) \:
	\sigma^{jk} \; (z-y)^m
\;+\; 2 \: \frac{\partial}{\partial x^j} \veeint_x^y dz \; B^{jk}(z) \:
	\sigma_{km} \; (z-y)^m \\
&&+ 2 \: \frac{\partial}{\partial x^j} \veeint_x^y dz \; B_{mk}(z) \:
	\sigma^{kj} \; (z-y)^m
\;+\; \frac{\partial}{\partial x^k} \veeint_x^y \varepsilon^{ijkl} \;
	B_{ij}(z) \; (z-y)_l \; \rho \; _{|x=0} \\
&=& -\frac{4 \pi}{y^2} \hinti d\lambda \; \varepsilon(\lambda) \:
	(\lambda-1) \; B^{jk}(\lambda y) \; y_j \: y^m \;
	\sigma_{km} \\
&&+ \frac{\pi}{2} \: \ln(|y^2|) \hinti d\lambda \; \epsilon(\lambda) \:
	(\lambda-1) \; \frac{d}{d\lambda} B_{jk}(\lambda y) \; \sigma^{jk} \\
&&+ \pi \: \ln(|y^2|) \hinti d\lambda \; \epsilon(\lambda) \:
	(2\lambda-1) \; B_{jk}(\lambda y) \: \sigma^{jk} \\
&&- i \pi \: \ln(|y^2|) \hinti d\lambda \; \epsilon(\lambda) \:
	(\lambda-1) \; B^{jk}_{\;\;,k}(\lambda y) \; y_j \\
&&+ \frac{\pi}{2} \: \ln(|y^2|) \hinti d\lambda \; \epsilon(\lambda) \:
	(\lambda-1) \; \varepsilon^{ijkl} \; B_{ij,k}(\lambda y)
	\: y_l \; \rho \\
&&- \pi \: \ln(|y^2|) \hinti d\lambda \; \epsilon(\lambda) \:
	(2\lambda^2-3\lambda+1) \; (B_{jk,}^{\;\;\;j} \: y_m
	+ B_{jk,m} \: y^j)
	\; \sigma^{km} \\
&&- \pi \: \ln(|y^2|) \hinti d\lambda \; \epsilon(\lambda) \:
	(\lambda^3-2\lambda^2+\lambda) \; (\Box B^{jk})(\lambda y)
	\: y_j \: y^m \; \sigma_{km} \;+\; {\cal{O}}(y^0) \\
\lefteqn{\left(\veeint_x^y dz \; (x-z)^m \gamma_m \; B_{jk}(z) \: \sigma^{jk}
	 \right) \: \Pdd_y \; _{|x=0}
\;=\; - 2 i \: \frac{\partial}{\partial y^k} \veeint^y dz \;
	(x-z)_j \; B^{jk}(z) } \\
&&- \frac{\partial}{\partial y^m} \veeint^y dz \; (x-z)^m \; B_{jk}(z) \:
	\sigma^{jk}
\;-\; 2 \: \frac{\partial}{\partial y^j} \veeint^y dz \; (x-z)^m \; B^{jk}(z)
	 \: \sigma_{km} \\
&&- 2 \: \frac{\partial}{\partial y^j} \veeint^y dz \; (x-z)^m \; B_{mk}(z)
	 \: \sigma^{kj}
\;-\; \frac{\partial}{\partial y^l} \veeint^y \varepsilon^{ijkl} \;
	B_{ij}(z) \; (x-z)_k \; \rho \\
&=& \frac{4 \pi}{y^2} \hinti d\lambda \; |\lambda|
	\; B^{jk}(\lambda y) \; y_j \: y^m \;
	\sigma_{km} \\
&&+ \frac{\pi}{2} \: \ln(|y^2|) \hinti d\lambda \; |\lambda| \:
	\frac{d}{d\lambda} B_{jk}(\lambda y) \; \sigma^{jk} \\
&&- \pi \: \ln(|y^2|) \hinti d\lambda \; \epsilon(\lambda) \:
	(2\lambda-1) \; B_{jk}(\lambda y) \: \sigma^{jk} \\
&&+ i \pi \: \ln(|y^2|) \hinti d\lambda \; |\lambda| \:
	 B^{jk}_{\;\;,k}(\lambda y) \; y_j \\
&&- \frac{\pi}{2} \: \ln(|y^2|) \hinti d\lambda \; |\lambda| \:
	\varepsilon^{ijkl} \; B_{ij,k}(\lambda y)
	\: y_l \; \rho \\
&&+ \pi \: \ln(|y^2|) \hinti d\lambda \; \epsilon(\lambda) \:
	(2\lambda^2-\lambda) \; (B_{jk,}^{\;\;\;j} \: y_m
	+ B_{jk,m} \: y^j)
	\; \sigma^{km} \\
&&+ \pi \: \ln(|y^2|) \hinti d\lambda \; \epsilon(\lambda) \:
	(\lambda^3-\lambda^2) \; (\Box B^{jk})(\lambda y)
	\: y_j \: y^m \; \sigma_{km} \;+\; {\cal{O}}(y^0) \\
\end{eqnarray*}
und somit
\begin{eqnarray*}
\lefteqn{ (s_m \: B_{jk} \: \sigma^{jk} \: p_m)(0,y) } \\
&=& \frac{m^2}{8 \pi^3} \frac{1}{y^2} \hinti d\lambda \; \epsilon(\lambda)
	\: B_{jk}(\lambda y) \; y^j \: y_m \: \sigma^{km} \\
&&+ \frac{m^2}{32 \pi^3} \: \ln(|y^2|) \; B_{jk}(0) \: \sigma^{jk} \\
&&+ \frac{i m^2}{32 \pi^3} \: \ln(|y^2|) \hinti d\lambda \; \epsilon(\lambda)
	\: B^{jk}_{\;\;,k}(\lambda y) \; y_j \\
&&- \frac{m^2}{64 \pi^3} \: \ln(|y^2|) \hinti d\lambda \; \epsilon(\lambda) \:
	\varepsilon^{ijkl} \; B_{ij,k}(\lambda y) \; y_l \: \rho \\
&&+ \frac{m^2}{32 \pi^3} \: \ln(|y^2|) \hinti d\lambda \; \epsilon(\lambda) \:
	(2\lambda-1) \; (B_{jk,}^{\;\;\;k} \: y^m
	+ B_{jk,m} \: y^j) \; \sigma^{km} \\
&&+ \frac{m^2}{32 \pi^3} \: \ln(|y^2|) \hinti d\lambda \; \epsilon(\lambda) \:
	(\lambda^2-\lambda) \; (\Box B^{jk})(\lambda y)
	\: y_j \: y^m \; \sigma_{km} \\
&&+ {\cal{O}}(y^0) \spc .
\end{eqnarray*}
\end{enumerate}
Durch Translation um $x$ und Vertauschung von $x$, $y$ erh\"alt man
hieraus Formeln f\"ur $(s_m \:B_{jk} \: \sigma^{jk}\: p_m)(x,y)$ und
$(p_m \:B_{jk} \: \sigma^{jk}\: s_m)(x,y)$.
Setze nun in \Ref{a4_bpm} ein.
\QED

\section{Differentialst\"orung durch Vektorpotential}
Wir betrachten wieder die St\"orung des Diracoperators \Ref{a1_va}. F\"ur
$\Delta p_m$ hat man in erster Ordnung in $L$
\[ \Delta p_m(x,y) \;=\; -i \frac{\partial}{\partial y^j} \:
	\Delta p_m[L^j](x,y) + \frac{i}{2} \: \Delta p_m[L^j_{\;,j}](x,y) \spc . \]
\begin{Thm}
In erster Ordnung St\"orungstheorie gilt
\begin{eqnarray}
\Delta p_m(x,y) &=& \Delta p_0(x,y) \\
\label{eq:a4_v1}
&&- im \left(\int_x^y L_j \right) \: \xi^j \; p_0(x,y) \\
&&+ \frac{1}{8 \pi^3} \frac{1}{\xi^2} \: m \int_x^y (2\alpha-1)
	\; (L^j_{\;,j} \: \xi\slsh + (\Pdd L_j) \: \xi^j ) \\
\label{eq:a4_v2}
&&+ \frac{1}{4 \pi^3} \frac{1}{\xi^2} \: m  \int_x^y L \slsh \\
&&+ \frac{1}{8 \pi^3} \frac{1}{\xi^2} \: m \int_x^y
	(\alpha^2-\alpha) \; (\Box L_j) \; \xi^j \: \xi\slsh \\
\label{eq:a4_v3}
&&+ \frac{i}{16 \pi^3} \frac{1}{\xi^2} \: m^2 \; (L_j(y) +
	L_j(x)) \; \xi^j \\
\label{eq:a4_v4}
&&+ \frac{i}{4 \pi^3} \frac{1}{\xi^2} \: m^2 \; \int_x^y L_j \: \xi^j \\
&&- \frac{1}{16 \pi^3} \frac{1}{\xi^2} \: m^2 \int_x^y
	\partial_j L_m \; \xi^m \: \xi_k \; \sigma^{jk} \\
&&+ {\cal{O}}(\ln(|\xi^2|)) + {\cal{O}}(m^3) \spc . \nonumber
\end{eqnarray}
\end{Thm}
{\Beweis}
Wir k\"onnen genau wie im Beweis von Theorem~\ref{thm_gp0} argumentieren.
Durch Vergleich von Theorem~\ref{a2_theorem_sm} und
Theorem~\ref{a4_theorem_spm}
folgt die Behauptung aus Theorem~\ref{theorem_vkm} durch die formalen
Ersetzungen~\Ref{a3_393a}.
\QED

\section{Differentialst\"orung durch axiales Potential}
Wir betrachten wieder die St\"orung des Diracoperators \Ref{a1_pva}. F\"ur
$\Delta p_m$ hat man in erster Ordnung in $L$
\[ \Delta p_m(x,y) \;=\; -i \frac{\partial}{\partial y^j} \:
	\Delta p_m[i \rho L^j](x,y) + \frac{i}{2} \: \Delta p_m[i \rho L^j_{\;,j}] \spc . \]
\begin{Thm}
\label{a4_thm_pvkm}
In erster Ordnung St\"orungstheorie gilt
\begin{eqnarray}
\Delta p_m(x,y) &=& i \rho \: \Delta p_0\left[i L^j \partial_j + \frac{i}{2}
	L^j_{\;,j}\right](x,y) \\
&&- \frac{i}{8 \pi^3} \frac{1}{\xi^2} \: m \: \rho \int_x^y (2\alpha-1) (\Pdd L_j)
	\xi^j \\
\label{eq:a4_vr1}
&&- \frac{1}{16 \pi^3} \frac{1}{\xi^2} \: m^2 \; (L_j(y) +
	L_j(x)) \; \xi^j \; \rho \\
&&- \frac{i}{16 \pi^3} \frac{1}{\xi^2} \: m^2 \int_x^y
	(\partial_j L_m) \; \xi^m \: \xi_k \; \rho \sigma^{jk} \\
&&+ {\cal{O}}(\ln(|\xi^2|)) + {\cal{O}}(m^3) \spc . \nonumber
\end{eqnarray}
\end{Thm}
{\Beweis}
Wir k\"onnen genau wie im Beweis von Theorem~\ref{thm_gp0} argumentieren.
Durch Vergleich von Theorem~\ref{a2_theorem_psm} und
Theorem~\ref{a4_theorem_pspm}
folgt die Behauptung aus Theorem~\ref{theorem_pvkm} durch die formalen
Ersetzungen~\Ref{a3_393a}.
\QED

\section{Bilineare Differentialst\"orung durch Vektorpotential}
Wir betrachten wieder die St\"orung des Diracoperators \Ref{a1_bva}. F\"ur
$\Delta p_m$ hat man in erster Ordnung in $L$
\begin{eqnarray*}
\Delta p_m(x,y)
&=& -i \Delta p_m \left[ i L^j \partial_j + \frac{i}{2} \: L^j_{\;,j} \right](x,y)
	\:+\: \frac{i}{2} \: \Delta p_m[L_{i,j} \sigma^{ij}](x,y) \nonumber \\
&& + i m \: \Delta p_m[L\slsh](x,y) \:-\: i p_m(x,y) \: L\slsh(y) \:-\:
	\frac{1}{2} \: \Delta p_m[L^j_{\;,j}] \spc .
\end{eqnarray*}

\begin{Thm}
In erster Ordnung St\"orungstheorie gilt
\begin{eqnarray}
\Delta p_m(x,y) &=& \Delta p_0(x,y) \\
&&+\frac{i}{8 \pi^3} \:m\: \frac{1}{\xi^2} \; (L\slsh(y) - L\slsh(x)) \\
&&-\frac{i}{8 \pi^3} \: m^2\: \frac{1}{\xi^2} \int_x^y L^j_{\;,j}
	\: \xi\slsh \\
\label{eq:a4_bv1}
&&-\frac{i}{16 \pi^3} \: m^2\: \frac{1}{\xi^2} \; (L_j(y) + L_j(x))
	\: \xi_k \: \sigma^{jk} \\
&&-\frac{1}{16 \pi^3} \: m^2\: \frac{1}{\xi^2} \; (L_j(y) - L_j(x))
	\: \xi^j \\
&&+ {\cal{O}}(\ln(|\xi^2|)) + {\cal{O}}(m^3) \spc . \nonumber
\end{eqnarray}
\end{Thm}
{\Beweis}
Wir k\"onnen genau wie im Beweis von Theorem~\ref{thm_gp0} argumentieren.
Durch Vergleich von Theorem~\ref{a2_theorem_sm} und
Theorem~\ref{a4_theorem_spm}
folgt die Behauptung aus Theorem~\ref{theorem_bvkm} durch die formalen
Ersetzungen~\Ref{a3_393a}.
\QED

\section{Bilineare Differentialst\"orung durch axiales Potential}
Wir betrachten wieder die St\"orung des Diracoperators \Ref{a2_ba1}. F\"ur
$\Delta p_m$ hat man in erster Ordnung in $L$
\begin{eqnarray*}
\Delta k_m(x,y)
&=& -i \Delta p_m \left[ \rho L^j \partial_j + \frac{\rho}{2} \: L^j_{\;,j} \right](x,y)
	\:+\: \frac{i}{2} \: \Delta p_m[\rho L_{i,j} \sigma^{ij}](x,y) \\
&& - p_m(x,y) \: \rho L\slsh(y) \:+\: m \: \Delta p_m[\rho L\slsh](x,y) \:+\:
	\frac{1}{2} \: \Delta p_m[i \rho L^j_{\;,j}] \spc .
\end{eqnarray*}

\begin{Thm}
In erster Ordnung St\"orungstheorie gilt
\begin{eqnarray}
\Delta p_m(x,y) &=& \Delta p_0(x,y) \\
\label{eq:a4_bvr1}
&&+\frac{1}{2 \pi^3} \:m \:
	\frac{1}{\xi^4} \int_x^y L_j \xi^j \; \rho \xi\slsh \\
&&-\frac{1}{8 \pi^3} \:m\: \frac{1}{\xi^2} \int_x^y (2\alpha-1) \:
	L^i_{\;,i} \; \rho \xi\slsh \\
&&-\frac{1}{8 \pi^3} \: m\: \frac{1}{\xi^2}  \int_x^y (\alpha^2-\alpha)
	\: (\Box L_k) \: \xi^k \; \rho \xi\slsh \\
&&-\frac{1}{8 \pi^3} \:m\: \frac{1}{\xi^2} \int_x^y (2\alpha-1) \: \rho (\Pdd L_j) \:
	\xi^j \\
\label{eq:a4_bvr2}
&&+\frac{1}{4 \pi^3} \:m\: \frac{1}{\xi^2} \; \rho (L\slsh(y) + L\slsh(x)) \\
\label{eq:a4_bvr3}
&&-\frac{1}{4 \pi^3} \:m\: \frac{1}{\xi^2} \int_x^y \rho L\slsh \\
&&+\frac{i}{8 \pi^3} \:m\: \frac{1}{\xi^2} \int_x^y \varepsilon^{ijkl}
	\: L_{i,j} \: \xi_k \; \gamma_l \\
\label{eq:a4_bvr4}
&&+\frac{1}{4 \pi^3} \: m^2 \: \frac{1}{\xi^2} \int_x^y L_j \: \xi_k
	\; \rho \sigma^{jk} \\
&&-\frac{i}{16 \pi^3} \: m^2 \: \frac{1}{\xi^2} \; (L_j(y) - L_j(x))
	\xi^j \; \rho \\
\label{eq:a4_bvr5}
&&+\frac{1}{16 \pi^3} \: m^2 \: \frac{1}{\xi^2} \; (L_j(y) + L_j(x))
	\: \xi_k \; \rho \sigma^{jk} \\
&&+ {\cal{O}}(\ln(|\xi^2|)) + {\cal{O}}(m^3) \spc . \nonumber
\end{eqnarray}
\end{Thm}
{\Beweis}
Wir k\"onnen genau wie im Beweis von Theorem~\ref{thm_gp0} argumentieren.
Durch Vergleich von Theorem~\ref{a2_theorem_sm} und
Theorem~\ref{a4_theorem_spm}
folgt die Behauptung aus Theorem~\ref{theorem_bakm} durch die formalen
Ersetzungen~\Ref{a3_393a}.
\QED

\chapter{St\"orungsrechnung h\"oherer Ordnung}
\label{anh6}
In diesem Kapitel werden wir die Distribution $P(x,y)$ f\"ur endliche
St\"orungen des Diracoperators untersuchen. Dazu werden wir die
Beitr\"age jeder Ordnung St\"orungstheorie asymptotisch um den Lichtkegel
entwickeln und anschlie{\ss}end alle Beitr\"age explizit aufsummieren.
Als Ergebnis erhalten wir asymptotische Entwicklungsformeln f\"ur $P(x,y)$,
die nicht-perturbativ g\"ultig sind.

Als Vorbereitung betrachten wir eine allgemeine St\"orung des Diracoperators durch
ein lokales Potential, also
\Equ{a6_1Z}
G \;=\; i \Pdd \:+\: B
\EndEqu
mit einer matrixwertigen Funktion $B(x)$. Zur Einfachheit nehmen wir an,
da{\ss} $B$ im Unendlichen glm. zur vierten Potenz abf\"allt. Die Spindimension
sei $4n$.

Bei der formalen St\"orungsentwicklung treten
Operatorprodukte der Form
\Equ{a6_1A}
A_1 \:B\: A_2 \:B\: \cdots \:B\: A_{l-1} \:B\: A_l
\EndEqu
auf, wobei $A_j$ f\"ur Faktoren $p_m$, $k_m$ oder $s_m$ steht.
\begin{Lemma}
\label{a6_lemmaa}
Das Operatorprodukt \Ref{a6_1A} ist als Operator auf den Schwartzfunktionen
wohldefiniert.
\end{Lemma}
{\Beweis}
Wir f\"uhren den Beweis im Impulsraum durch.
Beachte zun\"achst, da{\ss} die Fouriertransformierte von $B \in C^\infty_4$
stetig ist und im Unendlichen st\"arker als polynomial abf\"allt
\[ \sup_p |p^\alpha \: \tilde{B}(p)| \;\leq\; \infty \spc
	{\mbox{f\"ur alle Multi-Indizes $\alpha$}}. \]
F\"ur eine gegebene Schwartzfunktion ist $f$ ebenfalls stetig und
f\"allt im Unendlichen st\"arker als polynomial ab.
Folglich ist das Faltungsprodukt
\[ \widetilde{(B \:A_l\: f)} \;=\; \tilde{B} \:*\: \tilde{A}_l \: \tilde{f} \]
stetig und f\"allt im Unendlichen st\"arker als polynomial ab.
Nun k\"onnen wir induktiv mit $A_{l-1}, \ldots, A_1$ multiplizieren und
nach jedem Schritt mit $\tilde{B}$ falten.
Man erh\"alt als Ergebnis im
Impulsraum eine Distribution, die im Unendlichen st\"arker als polynomial
abf\"allt. Im Ortsraum ist dies eine glatte Funktion.
\QED
Damit sind die Beitr\"age jeder Ordnung St\"orungstheorie wohldefiniert und
endlich.

Wir entwickeln nun \Ref{a6_1A} nach Potenzen von $m$
und ordnen die Beitr\"age au{\ss}erdem nach der Anzahl der auftretenden
partiellen Ableitungen des Potentials $B$ an.
F\"ur diese Terme kann man das singul\"are Verhalten auf dem Lichtkegel
genau beschreiben:
\begin{Lemma}
\label{a6_lemma0}
Der Beitrag $\sim m^p$ von \Ref{a6_1A}, bei dem $q$ Faktoren des Potentials
$B$ partiell abgeleitet sind, verh\"alt sich auf dem Lichtkegel wie
\Equ{a6_1B}
\left\{ \begin{array}{cl}
	{\cal{O}}(\xi\slsh \: \xi^{-4+p+q}) & {\mbox{f\"ur $p+q=0,2$}} \\
	{\cal{O}}(\xi^{-3+p+q}) & {\mbox{f\"ur $p+q=1$}} \\
	{\cal{O}}(\ln(|\xi^2|)) & {\mbox{f\"ur $p+q=3$}} \\
	{\cal{O}}(\xi\slsh \: \ln(|\xi^2|)) & {\mbox{f\"ur $p+q=4$}} \\
	{\cal{O}}(\xi^0) & {\mbox{f\"ur $p+q>4$}}
		\end{array} \right. \spc .
\EndEqu
\end{Lemma}
Falls das St\"orpotential einen skalaren/pseudoskalaren Anteil enth\"alt,
k\"onnen wir diese Absch\"atzung noch verbessern:
\begin{Lemma}
\label{a6_lemma01}
F\"ur das Potential
\[ B \;=\; B_0 \:-\: m \Phi \:-\: im \rho \Psi \]
mit $B_0, \Phi, \Psi \in C^\infty_4(M, \R^{4n \times 4n})$ und
\[ \Phi(x) \:=\: \1_4 \otimes (\Phi_{ij})_{i,j=1,\cdots,4n} \;\;\;, \spc
	\Psi(x) \:=\: \1_4 \otimes (\Psi_{ij})_{i,j=1,\cdots,4n} \]
verh\"alt sich der Beitrag $\sim m^p$ von \Ref{a6_1A}, bei dem
$q$ Faktoren des Potentials $B$ partiell abgeleitet sind, auf dem Lichtkegel
wie \Ref{a6_1B}.
\end{Lemma}
Man beachte, da{\ss} in Lemma \ref{a6_lemma01} im Gegensatz zu Lemma \ref{a6_lemma1}
auch das Potential $B$ nach $m$ entwickelt wird.

Nach Lemma \ref{a6_lemma1} und Lemma \ref{a6_lemma01} wird die Singularit\"at
auf dem Lichtkegel bei den Beitr\"agen h\"oherer Ordnung in der Masse und
h\"oherer Ordnung in den Ableitungen der Potentiale schw\"acher.
Zu vorgegebener Ordnung ${\cal{O}}(\ln^\alpha(|\xi^2|) \: \xi^{2 \beta})$
auf dem Lichtkegel tragen damit nur ganz bestimmte Terme der St\"orungsentwicklung
bei, was die Rechnungen in diesem Kapitel wesentlich erleichtert.

\section{Eichterme/Pseudoeichterme}
\label{a6_ab1}
Wir wollen nun f\"ur den Diracoperator
\Equ{a6_100}
i \Pdd \;+\; \chi_L \: \Aslsh_R \;+\; \chi_R \: \Aslsh_L
\EndEqu
die Eichterme und Pseudoeichterme in beliebiger Ordnung in $\Aslsh_{L\!/\!R}$
zu berechnen, dabei bezeichnet $\chi_{L\!/\!R} = \frac{1}{2} (1 \mp \gamma^5)$
die chiralen Projektoren.
Im allgemeinen Fall ist die Spindimension $4n$, $n\geq1$ und
\[ \Aslsh_{L\!/\!R} = \gamma_j \: (A^j_{L\!/\!R\:ab})_{a,b=1,\cdots,n} \spc . \]
Wir f\"uhren eine Entwicklung sowohl nach $m$ als auch nach dem singul\"aren
Verhalten auf dem Lichtkegel durch. In diesem Abschnitt ber\"ucksichtigen wir
in $P(x,y)$
die Beitr\"age der Ordnung ${\cal{O}}(m^2) + {\cal{O}}(\ln(|\xi^2|))$, in
Abschnitt \ref{a6_ab2} werden die Beitr\"age $\sim m^2$ bis zur
Ordnung ${\cal{O}}(\xi^0)$ behandelt.
Da alle Beitr\"age der Ordnung ${\cal{O}}(m^3)$ auf dem Lichtkegel h\"ochstens
logarithmisch singul\"ar sind, haben wir dann die Eich- und
Pseudoeichterme f\"ur endliche St\"orungen bis zur Ordnung ${\cal{O}}(\ln(|\xi^2|))$
vollst\"andig berechnet.

Es treten zwei Schwierigkeiten auf: Zum einen m\"ussen wir nichtlokale
Linienintegrale studieren, au{\ss}erdem brauchen die Potentiale nicht
miteinander zu kommutieren:
\[ \left[ A^j_{L\!/\!R}(x) , \: A^k_{L\!/\!R}(y) \right] \;\neq\; 0 \]
Wir gehen das Problem schrittweise an und beginnen mit einem Spezialfall,
bei dem diese Komplikationen noch nicht auftreten:
\subsection{Kommutative Eichpotentiale}
\label{a6_ab1.1}
Der Diracoperator habe die Form
\Equ{a6_91}
G \;=\; i \Pdd \;+\; \chi_L \: (\Pdd \Lambda_R) \;+\; \chi_R \: (\Pdd \Lambda_L)
\EndEqu
mit reellen Funktionen $\Lambda_L, \Lambda_R$.
Die St\"orung des Diracoperators hat \"Ahnlichkeit mit den bei
$U(1)$-Eichtransformationen auftretenden Potentialen. Dadurch, da{\ss}
wir f\"ur die links- und rechtsh\"andige Komponente unterschiedliche
Potentiale zulassen, k\"onnen wir die St\"orung \Ref{a6_91} aber nicht
global wegeichen.

Wir berechnen zun\"achst einzelne Operatorprodukte:
\begin{Lemma}
\label{a6_lemma1}
F\"ur $\Aslsh=(\Pdd \Lambda)$, $\Lambda \in C^\infty(\R^4)$ gilt
\begin{eqnarray}
\label{eq:a6_2}
k_0 \:\Aslsh\: p_0 &=& \ke \:\Aslsh\: p_0 + k_0 \:\Aslsh\:\pe \;=\; 0 \\
\label{eq:a6_2a}
s_0 \:\Aslsh\: s_0 &=& -i \Lambda \: s_0 + i s_0 \: \Lambda \\
\label{eq:a6_3}
s_0 \:\Aslsh\: p_0 &=& -i \Lambda \: p_0 \\
\label{eq:a6_6}
s_0 \:\Aslsh\: k_0 &=& -i \Lambda \: k_0 \spc .
\end{eqnarray}
\end{Lemma}
{\Beweis}
Durch Entwicklung von
\[ k_m \:\Aslsh\: p_m \;=\; -i k_m \:[i\Pdd, \Lambda]\: p_m \;=\; 0 \]
nach $m$ erh\"alt man \Ref{a6_2}.
Die Gleichungen \Ref{a6_2a}, \Ref{a6_3} und \Ref{a6_6} folgen aus
\begin{eqnarray*}
s_0 \:\Aslsh\: s_0 &=& -i s_0 \:[i\Pdd,\Lambda]\: s_0 \;=\; -i \Lambda
	\: s_0 + i s_0 \: \Lambda \\
s_0 \:\Aslsh\: p_0 &=& -i s_0 \:[i\Pdd,\Lambda]\: p_0 \;=\; i \Lambda
	\: p_0 \\
s_0 \:\Aslsh\: k_0 &=& -i s_0 \:[i\Pdd,\Lambda]\: k_0 \;=\; i \Lambda
	\: k_0 \spc .
\end{eqnarray*}
\QED
Bei einigen Operatorprodukten l\"a{\ss}t sich die Anzahl der Faktoren $s_0$, $p_0$ stark
reduzieren:
\begin{Lemma}
Mit $\Aslsh=(\Pdd \Lambda)$, $B \!\slsh = (\Pdd \phi)$ gilt f\"ur $m,n \geq1$
\begin{eqnarray}
\label{eq:a6_10}
(s_0 \: \Aslsh)^n \: p_0 &=& \frac{(-i\Lambda)^n}{n!} \: p_0 \\
\label{eq:a6_10a}
\se \: \Aslsh \:(s_0 \: \Aslsh)^{n-1} \: p_0 &=& \se \: \Aslsh \:
	\frac{(-i\Lambda)^{n-1}}{(n-1)!} \: p_0 \\
\label{eq:a6_11}
(s_0 \: \Aslsh)^n \: \pe &=& \sum_{p=0}^{n-1} \; \frac{(-i\Lambda)^p}{p!}
	\: s_0 \: \frac{(i \Lambda)^{n-p-1}}{(n-p-1)!} \:\Aslsh\: \pe \\
\lefteqn{\hspace*{-3.2cm}(s_0 \: \Aslsh)^n\: \se \:B\!\slsh\: (s_0 \: B\!\slsh)^{m-1}
	\: p_0 } \nonumber \\
\label{eq:a6_12}
&=& \sum_{p=0}^{n-1} \; \frac{(-i\Lambda)^p}{p!}
	\: s_0 \:\frac{(i\Lambda)^{n-p-1}}{(n-p-1)!} \:\Aslsh\:
	\se \:B\!\slsh\: \frac{(-i\phi)^{m-1}}{(m-1)!} \: p_0 \; .
	\hspace*{1cm}
\end{eqnarray}
\end{Lemma}
{\Beweis}
Nach \Ref{a6_3} gilt
\begin{eqnarray*}
(s_0 \:\Aslsh)^{n-1} \: s_0 \:\Aslsh\: p_0 &=& -i (s_0 \: \Aslsh)^{n-2}
	\Aslsh \Lambda \: p_0 \\
&=& -\frac{i}{2} (s_0 \: \Aslsh)^{n-2} \: s_0 \: \Pdd(\Lambda^2) \: p_0
	\;=\; \cdots \;=\; \frac{(-i\Lambda)^n}{n!} \: p_0 \spc .
\end{eqnarray*}
Daraus erh\"alt man \Ref{a6_10} und \Ref{a6_10a}. Zur Berechnung von
\Ref{a6_11} wenden wir \Ref{a6_2a} iterativ von links nach rechts an, also
\[ (s_0 \: \Aslsh)^n \: \pe \;=\; -i \Lambda \:(s_0 \Aslsh)^{n-1} \pe \;+\;
	\frac{i}{2} \: s_0 \: \Pdd(\Lambda^2) \: (s_0 \: \Aslsh)^{n-2} \:
	\pe \;=\; \cdots \spc . \]
Nach $n-1$ Schritten erhalten wir Terme der Form
\Equ{a6_13}
c(p) \; \Lambda^p \: s_0 \: \Lambda^{n-p-1} \: \Aslsh \: \pe \spc .
\EndEqu
Wir m\"ussen noch die Konstanten $c(p)$ bestimmen. Dazu ist es
g\"unstig, zun\"achst bei dem $p$-ten Faktor $\Aslsh$ Relation \Ref{a6_2a}
anzuwenden
\begin{eqnarray*}
(s_0 \: \Aslsh)^n \: \pe &=& -i (s_0 \: \Aslsh)^{p-1} \: \Lambda \:
	(s_0 \: \Aslsh)^{n-p} \; \pe \\
&& \hspace*{1cm} +\; \frac{i}{2}
	(s_0 \: \Aslsh)^{p-1} \: s_0 \: \Pdd(\Lambda^2) \:
	(s_0 \: \Aslsh)^{n-p-1} \: \pe \spc .
\end{eqnarray*}
Wir brauchen dabei nur den ersten Summanden zu ber\"ucksichtigen, denn im
zweiten Summanden stehen der $p$-te und $(p+1)$-te Faktor $\Lambda$
direkt nebeneinander, was bei weiterer Anwendung von Lemma \ref{a6_lemma1}
nicht zu einem Term der Form \Ref{a6_13} f\"uhren kann.
Durch iterative Anwendung dieses Argumentes erhalten wir als einzigen
Beitrag zu \Ref{a6_13} den Ausdruck
\[ \frac{(-i\Lambda)^p}{p!} \: (s_0 \: \Aslsh)^{n-p} \: \pe \spc . \]
Jetzt wenden wir das gleiche Argument bei den verbleibenden Faktoren $\Aslsh$
von links nach rechts an und erhalten f\"ur \Ref{a6_13} den Ausdruck
\[ \frac{(-i\Lambda)^p}{p!} \:s_0 \:\frac{(i\Lambda)^{n-p-1}}{(n-p-1)!}
	\:\Aslsh\: \pe \]
und damit die Behauptung \Ref{a6_11}.

Zum Beweis von \Ref{a6_12} k\"onnen wir bei den rechten $m-1$ Faktoren
$B \! \slsh$
iterativ \Ref{a6_3} anwenden, bei den Faktoren $\Aslsh$ geht man genau
vor wie bei der Herleitung von \Ref{a6_11}.
\QED
Wir wollen nun den Operator $V p_m$ mit $V$ gem\"a{\ss} \Ref{2_62}
nach Potenzen von $m$ entwickeln und f\"ur den St\"orpoperator
\Equ{a6_20}
{\cal{B}} \;=\; \chi_L \: \Aslsh_R \;+\; \chi_R \: \Aslsh_L
\EndEqu
bis zur Ordnung ${\cal{O}}(m^2)$ berechnen.

Zur Ordnung $\sim m^0$ haben die Summanden in \Ref{2_62} die Form
\[ C_0(1,Q) \:{\cal{B}}\: C_0(2,Q) \cdots C_0(l-1,Q) \:{\cal{B}}\:
	C_0(l,Q) \:{\cal{B}}\: p_0 \spc . \]
Bei Einsetzen von \Ref{a6_20} k\"onnen wir die Projektoren $\chi_{L\!/\!R}$ unter
Verwendung der Antikommutatorrelationen $\{\chi_{L\!/\!R}, \: C_0\} = \{\chi_{L\!/\!R}, \:
	\Aslsh_{L\!/\!R}\}=0$ nach vorne bringen und erhalten
\begin{eqnarray}
&=& \chi_L \: C_0(1,Q) \:\Aslsh_L\: C_0(2,Q) \:\Aslsh_L \cdots \Aslsh_L
	\:C_0(l,Q)\: \Aslsh_L \:p_0 \nonumber \\
\label{eq:a6_21}
&&+ \chi_R \: C_0(1,Q) \:\Aslsh_R\: C_0(2,Q) \:\Aslsh_R \cdots \Aslsh_R
	\:C_0(l,Q)\: \Aslsh_R \:p_0 \spc .
\end{eqnarray}
Zur Ordnung $\sim m$ ist bei dem Operatorprodukt in \Ref{2_62} genau ein
Faktor der lineare Term in $m$, wie z.B. in
\[ C_0(1,Q) \:{\cal{B}}\: C_0(2,Q) \cdots {\cal{B}}\: C^{(1)}(p,Q)
	\:{\cal{B}} \cdots C_0(l,Q) \:{\cal{B}}\: p_0 \spc, 1<p<l \;\;\; . \]
Unter Verwendung von $[\chi_{L\!/\!R}, \: C^{(1)}]=0$ k\"onnen wir diesen Ausdruck
in der Form
\begin{eqnarray*}
&=& \chi_L \;C_0(1,Q) \:\Aslsh_L\: C_0(2,Q) \cdots \Aslsh_L \: C^{(1)}(p,Q) \:
	\Aslsh_R \cdots C_0(l,Q) \:\Aslsh_R \: p_0 \\
&&+ \chi_R \;C_0(1,Q) \:\Aslsh_R\: C_0(2,Q) \cdots \Aslsh_R \: C^{(1)}(p,Q) \:
	\Aslsh_L \cdots C_0(l,Q) \:\Aslsh_L \: p_0
\end{eqnarray*}
umschreiben.

Wir brauchen nur diejenigen Summanden zu ber\"ucksichtigen, die auch
bei der unit\"aren Transformation (\ref{eq:2_69}) auftreten:
\begin{Lemma}
\label{a6_lemma4}
F\"ur die St\"orung \Ref{a6_20} des Diracoperators gilt
\[ U p_m \;=\; V p_m \;+\; {\cal{O}}(m^2) \]
mit $U, V$ gem\"a{\ss} \Ref{2_69}, \Ref{2_62}.
\end{Lemma}
{\Beweis}
Wir m\"ussen zeigen, da{\ss} alle Terme in \Ref{2_62}, in denen ein Faktor
$k_m$ vorkommt, bis zur Ordnung ${\cal{O}}(m^2)$ verschwinden.

Zur Ordnung $\sim m^0$ haben diese Terme die Form
\[ \cdots k_0 \:{\cal{B}}\: (s_0 \:{\cal{B}})^n \: p_0 \spc , \;\; n\geq 0
	\spc , \]
wobei die P\"unktchen \"`$\cdots$\"' f\"ur einen beliebigen Vorfaktor stehen.
Nach Umschreiben gem\"a{\ss} \Ref{a6_21} enth\"alt jeder der beiden Summanden
einen Faktor
\Equ{a6_25}
k_0 \:\Aslsh\: (s_0 \: \Aslsh)^n \: p_0
\EndEqu
mit $\Aslsh=\Aslsh_{L\!/\!R}$, der nach mehrfacher Anwendung von \Ref{a6_3}
und \Ref{a6_2} verschwindet:
\[ \;=\; k_0 \:\Aslsh\: \frac{(-i\Lambda)^n}{n!} \: p_0 \;=\; k_0 \:
	\left[ i \Pdd, \: \frac{(-i\Lambda)^{n+1}}{(n+1)!} \right] \:
	p_0 \;=\; 0 \]
Wir kommen zu den Summanden $\sim m$ in \Ref{2_62}. Da die Anzahl der
Faktoren $k_m$ immer gerade ist, gen\"ugt es, den Fall $\#Q\geq2$ zu
betrachten.
Damit bleibt zu zeigen, da{\ss} die Produkte
\begin{eqnarray}
\label{eq:a6_22}
\cdots k_0 \:{\cal{B}}\: (s_0 \: {\cal{B}})^n \: p_0 \\
\label{eq:a6_23}
\cdots k_0 \:{\cal{B}} \: (s_0 \: {\cal{B}})^n \: k_0 \:(s_0 \: {\cal{B}})^p \: s^{(1)}
	\: {\cal{B}} \: (s_0\:{\cal{B}})^q \: p_0 \\
\label{eq:a6_24}
\cdots k_0 \:{\cal{B}}\: (s_0 \: {\cal{B}})^p \: \ke \: {\cal{B}} \: (s_0 \: {\cal{B}})^q
	\: p_0
\end{eqnarray}
null sind. Die Ausdr\"ucke \Ref{a6_22} und \Ref{a6_23} enthalten den Faktor
\Ref{a6_25} bzw.
\[ k_0 \: \Aslsh \: (s_0 \: \Aslsh)^n \: k_0 \;=\; k_0 \:
	\left[ i \Pdd, \: \frac{(-i\Lambda)^{n+1}}{(n+1)!} \right] \:
	k_0 \;=\; 0 \]
und fallen weg.
Bei \Ref{a6_24} verwendet man die Umformungen
\begin{eqnarray*}
k_0 \:{\cal{B}} \: (s_0 \: {\cal{B}})^p \: \ke \: {\cal{B}} \: (s_0 \: {\cal{B}})^q
	\: p_0
&=& k_0 \:\left[ i \Pdd, \: \frac{(i\Lambda)^{p+1}}{(p+1)!} \right] \: \ke \:
	\left[ i \Pdd, \: \frac{(-i\Lambda)^{q+1}}{(q+1)!} \right] \: p_0 \\
&=& - k_0 \:\frac{(i\Lambda)^{p+1}}{(p+1)!}\: (i \Pdd) \: \ke (i \Pdd) \:
	\frac{(-i\Lambda)^{q+1}}{(q+1)!} \: p_0
\end{eqnarray*}
und $(i \Pdd) \: \ke (i \Pdd) = i \Pdd \: k_0 = 0$.
\QED

\begin{Lemma}
\label{a6_lemma5}
F\"ur den St\"oroperator \Ref{a6_20} gilt
\begin{eqnarray}
\label{eq:a6_92}
\chi_{L} \: V \: p_m &=&  \chi_{L} \: e^{i \Lambda_{L}} \: p_0 \\
\label{eq:a6_93}
&&+\; m \: \chi_{L} \: \left( 1 - e^{i \Lambda_{L}} \: s_0 \:
	e^{-i \Lambda_{L}} \: \Aslsh_{L} \right) \: \pe \\
\label{eq:a6_94}
&&-\;m\: \chi_{L} \: \left(1-e^{i\Lambda_{L}} \: s_0 \:
	e^{-i \Lambda_{L}} \: \Aslsh_{L} \right) \: \se \: \Aslsh_{R} \:
	e^{i \Lambda_{R}} \: p_0 \;+\; {\cal{O}}(m^2) \;\;\;.\hspace*{1cm}
\end{eqnarray}
F\"ur die rechtsh\"andige Komponente gilt die analoge Gleichung, wenn man
die Indizes $L, R$ vertauscht.
\end{Lemma}
{\Beweis}
Wir betrachten nur die linksh\"andige Komponente, die rechtsh\"andige folgt
analog. Nach Lemma \ref{a6_lemma4} gen\"ugt es, den Operator $U$ gem\"a{\ss}
\Ref{2_69} zu berechnen.
Zur Ordnung $\sim m^0$ haben wir nach \Ref{a6_10}
\begin{eqnarray*}
\chi_L \: U \: p_0 &=& \chi_L \:\sum_{l=0}^\infty (-s_0 \: {\cal{B}})^l \: p_0 \;=\;
	\chi_L \: \sum_{l=0}^\infty (-s_0 \: \Aslsh_L)^l \: p_0 \\
&=& \chi_L \: \sum_{l=0}^\infty \: \frac{(i \Lambda_L)^l}{l!} \: p_0 \;=\;
	\chi_L \: e^{i \Lambda_L} \: p_0 \spc .
\end{eqnarray*}
Den Beitrag $\sim m$ k\"onnen wir in der Form
\begin{eqnarray*}
(\chi_L \: U \: p_m)^{(1)} &=& \chi_L \: \sum_{l=0}^\infty (-s_0 \: {\cal{B}})^l
	\: \pe \;+\; \chi_L \: \sum_{m,n=0}^\infty (-s_0 \:
	{\cal{B}})^n \: (-\se \: {\cal{B}}) \: (-s_0 \: {\cal{B}})^m
	\: p_0 \\
&=& \chi_L \: \pe \;+\; \chi_L \: \sum_{n=1}^\infty (-s_0 \: \Aslsh_L)^n \: \pe \\
&&+\: \chi_L \: \sum_{m=0}^\infty \:(-\se \: \Aslsh_R) \: (-s_0 \: \Aslsh_R)^m \:
	p_0 \\
&&+\: \chi_L \: \sum_{n,m=1}^\infty \:(-s_0 \: \Aslsh_L)^n \: (-\se \: \Aslsh_R)
	\: (-s_0 \: \Aslsh_R)^{m-1} \: p_0
\end{eqnarray*}
umschreiben und erhalten mit \Ref{a6_11}, \Ref{a6_10}
\begin{eqnarray*}
&=& \chi_L \: \pe \;-\; \chi_L \: \sum_{n=1}^\infty \sum_{p=0}^{n-1} \:
	\frac{(i \Lambda_L)^p}{p!} \: s_0 \: \frac{(-i
	\Lambda_L)^{n-p-1}}{(n-p-1)!}\:\Aslsh_L \: \pe \\
&&-\: \chi_L \: \sum_{n=0}^\infty \: \se \: \Aslsh_R  \: \frac{(i
	\Lambda_R)^n}{n!} \: p_0 \\
&&+\: \chi_L \: \sum_{n,m=1}^\infty \sum_{p=0}^{n-1} \: \frac{(i\Lambda_L)^p}{p!}
	\:s_0\: \frac{(-i\Lambda_L)^{n-p-1}}{(n-p-1)!} \: \Aslsh_L \:
	\se \: \Aslsh_R \: \frac{(i\Lambda_R)^{m-1}}{(m-1)!} \: p_0 \spc .
\end{eqnarray*}
Durch Umordnen der Reihen folgt die Behauptung.
\QED
Nach diesen Vorbereitungen k\"onnen wir auch $\tilde{p}_m, \tilde{k}_m$
berechnen:
\begin{Satz}
\label{a6_satz1}
F\"ur $\tilde{p}_m$, $\tilde{k}_m$ gilt mit der symbolischen Ersetzung
$C_m=p_m$ oder $C_m=k_m$
\begin{eqnarray}
\lefteqn{\chi_L \: \tilde{C}_m \;=\; \chi_L \: e^{i \Lambda_L} \: C_0 \:
	e^{-i \Lambda_L} } \nonumber \\
&&+ m \: \chi_L \left(1-e^{i\Lambda_L} \:s_0\: e^{-i\Lambda_L}\: \Aslsh_L 
	\right) \:C^{(1)}\: \left(1-\Aslsh_R\:e^{i\Lambda_R} \:s_0\:
	e^{-i\Lambda_R} \right) \nonumber \\
&&- m \: \chi_L \left(1-e^{i\Lambda_L} \:s_0\: e^{-i\Lambda_L}\: \Aslsh_L 
	\right) \:s^{(1)}\: \Aslsh_R \: e^{i \Lambda_R} \:C_0\:
	e^{-i\Lambda_R} \nonumber \\
\label{eq:a6_30a}
&&- m \: \chi_L \: e^{i \Lambda_L} \:C_0\: e^{-i \Lambda_L} \: \Aslsh_L \:
	\se \left(1-\Aslsh_R\:e^{i\Lambda_R} \:s_0\:
	e^{-i\Lambda_R} \right) \;+\; {\cal{O}}(m^2) \;\;\; .
\end{eqnarray}
F\"ur die rechtsh\"andige Komponente gilt die analoge Gleichung, wenn man die
Indizes $L$, $R$ vertauscht.
\end{Satz}
{\Beweis}
Nach Lemma \ref{a6_lemma4} und \Ref{2_69} gilt
\Equ{a6_95}
\tilde{C}_m \;=\; V \: C_m \: V^* \;=\; \sum_{p,q=0}^\infty
	(-s_m \: {\cal{B}})^p \:C_m \: (-{\cal{B}} \: s_m)^q \;+\; {\cal{O}}(m^2)
	\spc .
\EndEqu
Bei der Entwicklung dieser Gleichung nach $m$ k\"onnen wir den Beitrag
$\sim m^0$ direkt aus \Ref{a6_92} zusammensetzen.

Zur Ordnung $\sim m$ treten drei verschiedene Beitr\"age auf:
Wenn f\"ur einen der ersten $p$ Faktoren $s$ in \Ref{a6_95} der Operator
$\se$ eingesetzt wird, haben wir f\"ur $V$ den Ausdruck \Ref{a6_94},
f\"ur $V^*$ dagegen \Ref{a6_92} zu verwenden.
Falls einer der letzten $q$ Faktoren $s$ linear in $m$ ist, erhalten wir
entsprechend f\"ur $V$ \Ref{a6_92} und f\"ur $V^*$ \Ref{a6_94}.
Im Fall, da{\ss} in \Ref{a6_95} der Operator $C^{(1)}$ auftritt, m\"ussen wir
sowohl f\"ur $V$ als auch f\"ur $V^*$ den Term \Ref{a6_93} benutzen.
\QED
Um Gleichung \Ref{a6_30a} besser interpretieren zu k\"onnen, wollen wir eine
Entwicklung um den Lichtkegel durchf\"uhren.
\begin{Lemma}
\label{a6_lemma2}
Es sei $E(x)$ ein beliebiges Matrixfeld und $\Aslsh=(\Pdd \Lambda)$, 
$B \!\slsh = (\Pdd \phi)$. Mit der symbolischen Ersetzung $C_m=p_m$ oder
$C_m=k_m$ gilt
\begin{eqnarray}
\label{eq:a6_4}
\left(\se \:E\: C_0\right)(x,y) &=& \frac{i}{4} \inti d\lambda \;
	\epsilon(\lambda) \: E \; \xi\slsh \;
	C^{(1)} (x,y) \;+\; {\cal{O}}(\ln(|\xi^2|)) \\
\label{eq:a6_5}
\left(s_0 \:E\: C^{(1)} \right)(x,y) &=& \frac{i}{4} \inti d\lambda \;
	\epsilon(\lambda) \: \xi\slsh \; E(z) \;
	C^{(1)}(x,y) \;+\; {\cal{O}}(\ln(|\xi^2|)) \\
\left(s_0 \:B\!\slsh\: \se \Aslsh\: C_0 \right)(x,y) &=&
	\frac{1}{4} \: \phi(x) \; \inti d\lambda \; \epsilon(\lambda)
	\; \Aslsh(z) \: \xi\slsh \; C^{(1)}(x,y) \nonumber \\
\label{eq:a6_9}
&& \hspace*{-4cm} -\frac{1}{4} \: \inti d\lambda \; \epsilon(\lambda) \:
	(\phi(z) \: \Aslsh(z)+ \Lambda(z) B \!\slsh (z)) \: \xi\slsh
	\; C^{(1)}(x,y) \;+\;
	{\cal{O}}(\ln(|\xi^2|)) \\
\left(C_0 \:B\!\slsh\: \se \Aslsh\: s_0 \right)(x,y) &=&
	\frac{1}{4} \: \inti d\lambda \; \epsilon(1-\lambda)
	\; B \!\slsh (z) \: \xi\slsh \; C^{(1)}(x,y) \; \Lambda(y) \nonumber \\
\label{eq:a6_29}
&& \hspace*{-4.5cm} -\frac{1}{4} \: \inti d\lambda \; \epsilon(1-\lambda) \:
	(\phi(z) \: \Aslsh(z)+ \Lambda(z) B\!\slsh(z)) \: \xi\slsh
	\; C^{(1)}(x,y) \;+\;
	{\cal{O}}(\ln(|\xi^2|)) \\
\left(s_0 \:B\!\slsh\: C^{(1)} \Aslsh\: s_0 \right)(x,y) &=&
	\frac{1}{4} \: \phi(x) \; \inti d\lambda \; \epsilon(1-\lambda)
	\; \Aslsh(z) \: \xi\slsh \; C^{(1)}(x,y) \nonumber \\
\label{eq:a6_30}
&& \hspace*{-4cm} + \frac{1}{4} \: \inti d\lambda \;
	\epsilon(\lambda) \: B\!\slsh(z)) \: \xi\slsh
	\; \pe(x,y) \; \Lambda(y) \;+\;
	{\cal{O}}(\ln(|\xi^2|)) \spc , \hspace*{1cm}
\end{eqnarray}
wobei zur Abk\"urzung $z=\lambda y + (1-\lambda)x$ gesetzt wurde.
\end{Lemma}
{\Beweis}
Nach \Ref{a4_243a}, \Ref{a3_194a}, \Ref{a3_25e} und \Ref{a5_b} haben wir
\begin{eqnarray*}
\lefteqn{ \left(\se \:E\: p_0\right)(x,y) \;=\; \left(S_0 \:E\:
	\chi_0 \:(i \Pdd)\right)(x,y) } \\
&=& -\frac{i}{16 \pi^4} \: \frac{\partial}{\partial y^j} \veeint_x^y
	E \gamma^j \\
&=& \frac{i}{4} \: \inti d\lambda \; \epsilon(\lambda) \;
	E(z) \: \xi\slsh \; \pe(x,y) \;+\;
	{\cal{O}}(\ln(|\xi^2|)) \spc .
\end{eqnarray*}
Damit folgt \Ref{a6_4} f\"ur $C_m=p_m$, Gleichung \Ref{a6_5} folgt analog.

Nach \Ref{a1_75}, \Ref{a1_76} und \Ref{a3_z} gilt
\begin{eqnarray*}
(S_0 \:E\: K_0)(x,y) &=& \frac{i}{16 \pi^3} \: \int d^4z \;
	l(z) \; (l^\vee_y(z)-l^\wedge_y(z)) \; E(z)
\;=\; \frac{i}{16 \pi^3} \; \xint_x^y \hat{E}
\end{eqnarray*}
mit $\hat{E}(z)=\epsilon(z^0-x^0) \: A(z)$. Unter Verwendung von
\Ref{a3_11}, \Ref{a5_dd} erh\"alt man
\begin{eqnarray*}
\left(\se \:E\: k_0\right)(x,y) &=& \left(S_0 \:E\:
	K_0 \: (i \Pdd)\right)(x,y) \\
&=& -\frac{1}{16 \pi^2} \; l(\xi) \; \inti d\lambda \; \hat{E}(z) \;
	\xi\slsh \;+\; {\cal{O}}(\xi^0) \\
&=& -\frac{1}{16 \pi^2} \; l(\xi) \: \epsilon(\xi^0) \; \inti d\lambda
	\; \epsilon(\lambda) \; E(z) \;
	\xi\slsh \;+\; {\cal{O}}(\xi^0) \\
&=& \frac{i}{4} \; \inti d\lambda \; \epsilon(\lambda) \;  E(z)
	\; \xi\slsh \; \ke(x,y) \;+\; {\cal{O}}(\xi^0) \spc .
\end{eqnarray*}
Das ist \Ref{a6_4} f\"ur den Fall $C_m=k_m$,
Gleichung \Ref{a6_5} kann man genauso ableiten.

Zum Beweis von \Ref{a6_9} wenden wir nacheinander \Ref{a6_4} und \Ref{a6_5}
an:
\begin{eqnarray*}
\lefteqn{ \left(s_0 \:B\!\slsh\: \se \:\Aslsh\: C_0\right)(x,y) \;=\; \int
	d^4u \; s_0(x,u) \: B\!\slsh(u) \; (\se \:\Aslsh\: C_0)(u,y) } \\
&=& \frac{i}{4} \: \int d^4u \; s_0(x,u) \: B\!\slsh(u) \; \inti d\alpha \;
	\epsilon(\alpha) \; \Aslsh(\alpha y + (1-\alpha)u) \:
	(y\slsh-u\slsh) \; C^{(1)}(u,y) \\
&& \hspace*{1cm} +\; {\cal{O}}(\ln(|\xi^2|)) \\
&=& -\frac{1}{16} \: \inti d\beta \; \epsilon(\beta) \inti d\alpha \;
	\epsilon(\alpha) \; \xi\slsh \; B\!\slsh(\beta y+(1-\beta)x) \\
&& \times \; \Aslsh((\alpha+(1-\alpha)\beta) y + (1-\alpha)(1-\beta)x) \;
	(1-\beta) \; \xi\slsh \; C^{(1)}(x,y) \;+\; {\cal{O}}(\ln(|\xi^2|))
\end{eqnarray*}
Nach der Variablentransformation $\alpha \rightarrow \alpha+(1-\alpha)\beta$
erh\"alt man
\begin{eqnarray*}
&=& -\frac{1}{16} \: \inti d\beta \; \epsilon(\beta) \inti d\alpha \;
	\epsilon(\alpha-\beta) \; \xi\slsh\;B\!\slsh(\beta y+(1-\beta)x) \\
&& \hspace*{2cm}\times \; \Aslsh(\alpha y + (1-\alpha) x) \;
	\xi\slsh \; C^{(1)}(x,y) \;+\; {\cal{O}}(\ln(|\xi^2|)) \spc .
\end{eqnarray*}
Wir setzen die Relation
\Equ{a6_31}
\xi\slsh \: v\slsh w\slsh \: \xi\slsh \;=\; 2 \: v_j \xi^j \: w\slsh
	\xi\slsh - 2 \: w_j \xi^j \: v\slsh \xi\slsh \;+\;
	{\cal{O}}(\xi^2)
\EndEqu
ein und integrieren jeweils in einer der Variablen $\alpha, \beta$
partiell:
\begin{eqnarray*}
&=& -\frac{1}{8} \: \inti d\beta \inti d\alpha \; \epsilon(\beta) \:
	\epsilon(\alpha-\beta) \; \frac{d}{d\beta}
	\phi(\beta y+(1-\beta)x)) \\
&& \hspace*{1cm} \times \; \Aslsh(\alpha y + (1-\alpha)x) \; \xi\slsh \;
	C^{(1)}(x,y) \\
&& +\frac{1}{8} \: \inti d\beta \inti d\alpha \; \epsilon(\beta) \:
	\epsilon(\alpha-\beta) \; \frac{d}{d\alpha}
	\Lambda(\alpha y+(1-\alpha)x)) \\
&& \hspace*{1cm} \times \; B \! \slsh(\beta y + (1-\beta)x) \; \xi\slsh \;
	C^{(1)}(x,y)  \;+\; {\cal{O}}(\ln(|\xi^2|))\\
&=& -\frac{1}{4} \: \inti d\alpha \; \epsilon(\alpha) \; \phi(\alpha
	y+(1-\alpha)x) \; \Aslsh(\alpha y+(1-\alpha)x) \; \xi\slsh \;
	C^{(1)}(x,y) \\
&&+ \frac{1}{4} \: \inti d\alpha \; \epsilon(\alpha) \; \phi(x) \;
	\Aslsh(\alpha y+(1-\alpha)x) \; \xi\slsh \; C^{(1)}(x,y) \\
&&- \frac{1}{4} \: \inti d\beta \; \epsilon(\beta) \; \Lambda(\beta
	y+(1-\beta)x) \; B \! \slsh(\beta y+(1-\beta)x) \; \xi\slsh \;
	C^{(1)}(x,y)  \;+\; {\cal{O}}(\ln(|\xi^2|))
\end{eqnarray*}
Gleichung \Ref{a6_29} folgt aus \Ref{a6_9} durch Bildung der Adjungierten.
Der Beweis von \Ref{a6_30} verl\"auft \"ahnlich wie die Herleitung von
\Ref{a6_9}
\begin{eqnarray*}
\lefteqn{ \left(s_0 \:B\!\slsh\: C^{(1)} \:\Aslsh\: p_0\right)(x,y) \;=\; \int
	d^4u \; s_0(x,u) \: B\!\slsh(u) \; (C^{(1)} \:\Aslsh\: s_0)(u,y) } \\
&=& \frac{i}{4} \: \int d^4u \; s_0(x,u) \: B\!\slsh(u) \; \inti d\alpha \;
	\epsilon(\alpha) \; \Aslsh(\alpha u + (1-\alpha)y) \:
	(y\slsh-u\slsh) \; C^{(1)}(u,y) \\
&& \hspace*{1cm} +\; {\cal{O}}(\ln(|\xi^2|)) \\
&=& -\frac{1}{16} \: \inti d\beta \; \epsilon(\beta) \inti d\alpha \;
	\epsilon(\alpha) \; \xi\slsh \; B\!\slsh(\beta y+(1-\beta)x) \\
&& \times \; \Aslsh((1-\alpha+\alpha\beta)y+ (\alpha-\alpha \beta)x) \;
	(1-\beta) \; \xi\slsh \; C^{(1)}(x,y) \;+\; {\cal{O}}(\ln(|\xi^2|))
	\;\;\; .
\end{eqnarray*}
Wir f\"uhren die Variablentransformation $\alpha \rightarrow 1-\alpha+\alpha
\beta$ durch, setzen \Ref{a6_31} ein und integrieren partiell:
\begin{eqnarray*}
&=& -\frac{1}{16} \: \inti d\beta \; \epsilon(\beta) \inti d\alpha \;
	\epsilon(1-\alpha) \; \xi\slsh\;B\!\slsh(\beta y+(1-\beta)x) \\
&& \hspace*{2cm}\times \; \Aslsh(\alpha y + (1-\alpha) x) \;
	\xi\slsh \; C^{(1)}(x,y) \;+\; {\cal{O}}(\ln(|\xi^2|)) \\
&=& \frac{1}{4} \: \phi(x) \: \inti d\alpha \; \epsilon(1-\alpha) \;
	\Aslsh(\alpha y+(1-\alpha)x) \; \xi\slsh \; C^{(1)}(x,y) \\
&& +\frac{1}{4} \: \Lambda(y) \: \inti d\beta \; \epsilon(\beta) \;
	B \! \slsh(\beta y+(1-\beta)x) \; \xi\slsh \;
	C^{(1)}(x,y) \;+\; {\cal{O}}(\ln(|\xi^2|))\spc .
\end{eqnarray*}
\QED
\begin{Satz}
\label{a6_satz5}
F\"ur $\tilde{p}_m$, $\tilde{k}_m$ gilt mit der symbolischen Ersetzung
$C_m=k_m$ oder $C_m=p_m$
\begin{eqnarray}
\label{eq:a6_40}
\tilde{C}_m(x,y) &=& \chi_L \: e^{i \Lambda_L(x)-i \Lambda_L(y)}
	\: C_m(x,y) \;+\; \chi_R \: e^{i \Lambda_R(x)-i \Lambda_R(y)} \:
	C_m(x,y) \\
\label{eq:a6_41}
&& -\frac{m}{2} \: \chi_L \: e^{i \Lambda_L(x) - i \Lambda_R(y)} \int_x^y
	\left( \Pdd \: e^{-i\Lambda_L + i\Lambda_R} \right) \: \xi\slsh
	\; C^{(1)}(x,y) \\
\label{eq:a6_42}
&& -\frac{m}{2} \: \chi_R \: e^{i \Lambda_R(x) - i \Lambda_L(y)} \int_x^y
	\left( \Pdd \: e^{-i\Lambda_R + i\Lambda_L} \right) \: \xi\slsh
	\; C^{(1)}(x,y) \\
&& +\; {\cal{O}}(\ln(|\xi^2|)) \;+\; {\cal{O}}(m^2) \spc . \nonumber 
\end{eqnarray}
\end{Satz}
{\Beweis}
Es gen\"ugt wieder, die linksh\"andige Komponente zu betrachten, die
rechts\-h\"an\-di\-ge folgt durch Vertauschung der Indizes $L, R$.
Durch Ausmultiplizieren von \Ref{a6_30a} erh\"alt man
\begin{eqnarray}
\lefteqn{ \chi_L \: \tilde{C}_m \;=\; \chi_L \: e^{i \Lambda_L} \: C_0 \:
	e^{-i\Lambda_L} \;+\; m \: \chi_L \: C^{(1)} } \nonumber \\
&& -\:m \: \chi_L \: e^{i \Lambda_L} \: \left( s_0 \:\Aslsh_L\: e^{-i\Lambda_L}
	\: C^{(1)} + C_0 \:\Aslsh_L\: e^{-i\Lambda_L} \: s^{(1)} \right)
	\nonumber \\
&& -\:m \: \chi_L \: \left(C^{(1)} \: e^{i \Lambda_R} \: \Aslsh_R \: s_0 +
	s^{(1)} \: e^{i \Lambda_R} \: \Aslsh_R \: C_0 \right) \: e^{-i
	\Lambda_R} \nonumber \\
&& +\:m \: \chi_L \: e^{i \Lambda_L} \: \left(s_0 \:\Aslsh_L\: e^{-i \Lambda_L}
	\: C^{(1)} \: \Aslsh_R \: e^{i \Lambda_R} \: s_0
	+C_0 \:\Aslsh_L\: e^{-i \Lambda_L} \: s^{(1)} \:
	\Aslsh_R \: e^{i \Lambda_R} \: s_0 \right. \nonumber \\
\label{eq:a6_759a}
&&\hspace*{1cm}\left. +s_0 \:\Aslsh_L\: e^{-i \Lambda_L} \: s^{(1)} \:
	\Aslsh_R \: e^{i \Lambda_R} \: C_0 \right) \: e^{-i \Lambda_R}
	\;+\; {\cal{O}}(\ln(|\xi^2|)) \;+\; {\cal{O}}(m^2) \;\;\; . \spc
\end{eqnarray}
Nach Einsetzen der asymptotischen Entwicklungen von Lemma \ref{a6_lemma2}
folgt
\begin{eqnarray*}
\lefteqn{ \chi_L \: \tilde{C}_m(x,y) \;=\; \chi_L \: e^{i \Lambda_L(x)-
	i\Lambda_L(y)} \: C_0(x,y) \;+\; m \: \chi_L \: C^{(1)}(x,y) } \\
&&-\frac{im}{2} \:\chi_L\: e^{i \Lambda_L(x)} \left(\int_x^y \xi\slsh \:
	\Aslsh_L \: e^{-i\Lambda_L} \right) \; C^{(1)}(x,y) \\
&&-\frac{im}{2} \:\chi_L\: e^{-i \Lambda_R(y)} \left(\int_x^y \Aslsh_R \:
	\xi\slsh \: e^{i\Lambda_R} \right) \; C^{(1)}(x,y) \\
&&-\frac{im}{2} \:\chi_L\: e^{i \Lambda_L(x)-i \Lambda_R(y)} \left(\int_x^y
	e^{-i\Lambda_L +i\Lambda_R} \:(\Aslsh_R-\Aslsh_L)\: \xi\slsh
	\right) \; C^{(1)}(x,y) \\
&&+\frac{im}{2} \:\chi_L\: e^{-i \Lambda_R(x)} \left(\int_x^y \Aslsh_R \:
	\xi\slsh \: e^{i\Lambda_R} \right) \; C^{(1)}(x,y) \\
&&-\frac{im}{2} \:\chi_L\: e^{i \Lambda_L(x)} \left(\int_x^y \Aslsh_L \:
	\xi\slsh \: e^{-i\Lambda_L} \right) \; C^{(1)}(x,y)
	\;+\; {\cal{O}}(\ln(|\xi^2|))
\end{eqnarray*}
und somit die Behauptung.
\QED
Dieses Ergebnis kann man auch direkt einsehen: Der Beitrag \Ref{a6_40}
beschreibt eine Phasentransformation der Distribution
$C_m(x,y)$, die Summanden \Ref{a6_41}, \Ref{a6_42} modifizieren
das Transformationsverhalten von $C^{(1)}$.
Im Spezialfall $\Lambda_L=\Lambda_R$ verschwinden \Ref{a6_41}, \Ref{a6_42};
aus \Ref{a6_40} erhalten wir das \"ubliche Verhalten bei
$U(1)$-Pha\-sen\-trans\-for\-ma\-tio\-nen der Elektrodynamik.

Im allgemeineren Fall $\Lambda_L \neq \Lambda_R$ f\"uhrt \Ref{a6_40} zus\"atzlich zu einer
relativen Phasenverschiebung der links- und rechtsh\"andigen Komponente.
Um zu verstehen, warum nun die Beitr\"age \Ref{a6_41}, \Ref{a6_42} ben\"otigt
werden, berechnen wir die Adjungierte von \Ref{a6_40}:
\begin{eqnarray}
\lefteqn{ \left(\chi_L \: e^{i\Lambda_L} \: C_m \: e^{-i\Lambda_L} \:+\:
	\chi_R \: e^{i\Lambda_R} \: C_m \: e^{-i\Lambda_R} \right)^*(x,y) }
	\nonumber \\
&=& \left( e^{i\Lambda_L} \: C_m \: e^{-i\Lambda_L} \: \chi_R \:+\:
	e^{i\Lambda_R} \: C_m \: e^{-i\Lambda_R} \: \chi_L \right)(x,y)
	\nonumber \\
&=& \chi_L \: e^{i \Lambda_L(x)-i \Lambda_L(y)} \: C_0(x,y) \;+\;
	\chi_R \: e^{i \Lambda_R(x)-i \Lambda_R(y)} \: C_0(x,y) \nonumber \\
&& +\: m \: \chi_L \: e^{i \Lambda_R(x)-i \Lambda_R(y)}
	\: C^{(1)} (x,y) \;+\; m \: \chi_R \: e^{i \Lambda_L(x)-i
	\Lambda_L(y)} \: C^{(1)}(x,y) \nonumber \\
\label{eq:a6_47}
&&+\; {\cal{O}}(m^2)
\end{eqnarray}
Bei dem in $m$ linearen Beitrag hat sich die links- und rechtsh\"andige
Komponente gerade vertauscht. Wir sehen also, so da{\ss} \Ref{a6_40}
allein nicht hermitesch ist.
Mit Hilfe der Relationen
\begin{eqnarray}
\int_x^y \left(\Pdd \: e^{-i\Lambda_L + i\Lambda_R}\right) \xi\slsh
&=& -\int_x^y \xi\slsh \left(\Pdd \: e^{-i\Lambda_L + i\Lambda_R}\right)
	+ 2 \int_x^y \xi^j \partial_j e^{-i\Lambda_L + i\Lambda_R}
	\nonumber \\
\label{eq:a6_753a}
&& \hspace*{-2.5cm} =\; -\int_x^y \xi\slsh \left(\Pdd \: e^{-i\Lambda_L +
	i\Lambda_R}\right)
	+ 2e^{-i\Lambda_L(y) + i\Lambda_R(y)}
	- 2e^{-i\Lambda_L(x) + i\Lambda_R(x)} \\
\int_x^y \left(\Pdd \: e^{-i\Lambda_R + i\Lambda_L}\right) \xi\slsh
&=& -\int_x^y \xi\slsh \left(\Pdd \: e^{-i\Lambda_R + i\Lambda_L}\right)
	+ 2e^{-i\Lambda_R(y) + i\Lambda_L(y)}
	- 2e^{-i\Lambda_R(x) + i\Lambda_L(x)} \nonumber 
\end{eqnarray}
k\"onnen wir $\tilde{C}_m$ umformen:
\begin{eqnarray}
\label{eq:a6_43}
\lefteqn{ \tilde{C}_m(x,y) \;=\; \chi_L \: e^{i \Lambda_L(x)-i \Lambda_L(y)}
	\: C_0(x,y) \;+\; \chi_R \: e^{i \Lambda_R(x)-i \Lambda_R(y)} \:
	C_0(x,y) } \\
\label{eq:a6_44}
&& +\:m \: \chi_L \: e^{i \Lambda_R(x)-i \Lambda_R(y)}
	\: C^{(1)}(x,y) \;+\; m \: \chi_R \: e^{i \Lambda_L(x)-i \Lambda_L(y)} \:
	C^{(1)}(x,y) \\
\label{eq:a6_45}
&& +\:\frac{m}{2} \: \chi_L \: e^{i \Lambda_L(x) - i \Lambda_R(y)} \int_x^y
	\xi\slsh \left( \Pdd \: e^{-i\Lambda_L + i\Lambda_R} \right)
	\; C^{(1)}(x,y) \\
\label{eq:a6_46}
&& +\:\frac{m}{2} \: \chi_R \: e^{i \Lambda_R(x) - i \Lambda_L(y)} \int_x^y
	\xi\slsh \left( \Pdd \: e^{-i\Lambda_R + i\Lambda_L} \right)
	\; C^{(1)}(x,y) \\
&&+\; {\cal{O}}(\ln(|\xi^2|)) \;+\; {\cal{O}}(m^2) \nonumber
\end{eqnarray}
Der Beitrag \Ref{a6_47} stimmt mit \Ref{a6_43}+\Ref{a6_44} \"uberein,
die Adjungierten von \Ref{a6_45} und \Ref{a6_46} sind gerade \Ref{a6_42}
bzw. \Ref{a6_41}. Unser Ausdruck f\"ur $\tilde{C}_m$ wird also erst durch
die Beitr\"age \Ref{a6_41}, \Ref{a6_42} hermitesch. Mit
dieser \"Uberlegung k\"onnen wir auch alle Vorzeichen und Vorfaktoren
in \ref{a6_satz5} kontrollieren.

\subsection{Nichtabelsche Eichpotentiale}
Wir wollen nun die Ergebnisse von Satz \ref{a6_satz5} auf den allgemeineren
Operator \Ref{a6_100} erweitern.
Im Gegensatz zum vorigen Abschnitt arbeiten wir hier nicht mit der
Transformation $V$, sondern berechnen zun\"achst $\tilde{s}^\vee_m$
gem\"a{\ss} \Ref{2_t2}. Daraus erhalten wir mit Hilfe von \Ref{2_tm}
$\tilde{k}_m$ und mit einem Analogieargument auch $\tilde{p}_m$.

Der Grund f\"ur dieses Vorgehen liegt darin, da{\ss} im Operator $V$ nichtlokale
Linienintegrale auftreten. Auf die damit verbundenen Probleme werden wir erst
in Abschnitt \ref{a6_sect4} eingehen.

Im folgenden bezeichne $\Aslsh$ bei einer Spindimension von $4n$ stets
ein Matrixfeld der Form
\[ \Aslsh(x) \;=\; \gamma_j \: (A^j_{ab})_{a,b=1,\ldots,n} \spc . \]
\begin{Lemma}
\label{a6_lemma6}
Es gelten die asymptotischen Entwicklungen
\begin{eqnarray}
\label{eq:a6_59a}
(\svn \:\Aslsh\: \svn)(x,y) &=& i \int_x^y A_j \: \xi^j \;
	s_0^\vee(x,y) \\
\label{eq:a6_59b}
&&+\frac{1}{2} \left( \int_x^y (2 \alpha-1) \: \xi^j \gamma^k F_{kj} \right)
	\; \sve(x,y) \\
\label{eq:a6_59c}
&&-\frac{1}{2}  \left( \int_x^y (\alpha^2-\alpha) \: \xi \slsh \xi^k \: j_k
	 \right) \; \sve(x,y) \\
\label{eq:a6_59d}
&&+\frac{i}{4} \left( \int_x^y \varepsilon^{ijkl} \: F_{ij} \: \xi_k
      \; \rho \gamma_l \right) \; \sve(x,y) \;+\;
	{\cal{O}}(\xi^0) \\
\label{eq:a6_51}
(\svn \:\Aslsh\: \sve)(x,y) &=& \frac{i}{2} \int_x^y
	\xi\slsh \Aslsh \; \sve(x,y) \;+\; {\cal{O}}(\xi^0)\\
\label{eq:a6_52}
(\sve \:\Aslsh\: \svn)(x,y) &=& \frac{i}{2} \int_x^y
	\Aslsh \: \xi\slsh \; \sve(x,y) \;+\; {\cal{O}}(\xi^0)
\end{eqnarray}
mit $F_{jk}=\partial_j A_k-\partial_k A_j$, $j^k=F^{jk}_{\;,l}$.\\
(Zur \"ubersichtlicheren Notation schreiben wir bei $s^\vee$ den
Index $^{(1)}$ entgegen der sonstigen Konvention nach unten.)
\end{Lemma}
{\Beweis}
Nach \Ref{2_f7} und \Ref{2_21}, \Ref{2_22} gilt
\begin{eqnarray*}
\svn(x,y) &=& -\frac{1}{2\pi} \: i \Pdd_x \l^\vee(y-x) \;=\;
	\frac{i}{\pi} \: \xi\slsh \: m^\vee(\xi) \\
\sve(x,y) &=& -\frac{1}{2\pi} \: l^\vee(y-x)
\end{eqnarray*}
und somit
\begin{eqnarray*}
(\svn \:\Aslsh\: \svn)(x,y) &=& -\frac{1}{4 \pi^2} \: \Pdd_x \left(\int
	d^4z \; l^\vee(z-x) \: \Aslsh(z) \: l^\vee(y-z) \right) \Pdd_y \\
&=& \frac{1}{4\pi^2} \: \gamma^i \gamma^j \gamma^k
   \: \frac{\partial}{\partial x^i} \frac{\partial}{\partial y^k}
   \lint_x^y A_j \spc .
\end{eqnarray*}
Das stimmt bis auf einen Faktor $2 \pi\:/\:ie$ mit \Ref{a1_95} \"uberein.
Aus \Ref{a1_87}-\Ref{a1_89} folgt Gleichung \Ref{a6_59a}-\Ref{a6_59d}.
Zum Beweis von \Ref{a6_51} wenden wir Satz \ref{a1_dis_abl} an
\begin{eqnarray*}
(\svn \:\Aslsh\: \sve)(x,y) &=& \frac{i}{4 \pi^2} \: \Pdd_x \lint_x^y \Aslsh \\
&=& -\frac{i}{4 \pi} \: \int_x^y \xi\slsh \Aslsh \: l^\vee(\xi) \;+\;
	{\cal{O}}(\xi^0) \spc ,
\end{eqnarray*}
\Ref{a6_52} folgt analog.
\QED
Bei der iterativen Anwendung dieser Formeln treten {\bf{zeitgeordnete
Linienintegrale}} auf:
\begin{Def}
\label{a6_def1}
Wir definieren
\begin{eqnarray}
\T \left(\int_x^y A_j \: (y-x)^j\right)^n &=& \int_x^y dz_1 \int_x^{z_1}
	dz_2 \cdots \int_x^{z_{n-1}} dz_n \nonumber \\
\label{eq:a6_53}
&& \times \; A_{j_n}(z_n)\:(z_n-x)^{j_n} \cdots
	A_{j_1}(z_1)\:(z_1-x)^{j_1} \\
\label{eq:a6_54}
\Texp \left(\int_x^y A_j \: (y-x)^j \right) &=& \sum_{n=0}^\infty \;
	\T \left(\int_x^y A_j \: (y-x)^j \right)^n \spc .
\end{eqnarray}
\end{Def}
Wir stellen die wichtigsten Eigenschaften der zeitgeordneten Integrale zusammen:
\begin{Lemma}
Die Reihe \Ref{a6_54} konvergiert absolut. Es gelten die
Differentialgleichungen
\begin{eqnarray}
\label{eq:a6_55}
\xi^k \frac{\partial}{\partial y^k} \T \left(\int_x^y  A_j \: \xi^j
	\right)^n &=& \T\left(\int_x^y A_j \xi^j \right)^{n-1} \; A_k(y)
	\: \xi^k \;\;\;, n \geq 1 \\
\label{eq:a6_56}
\xi^k \frac{\partial}{\partial y^k} \Texp \left(\int_x^y A_j \: \xi^j \right)
	&=& \Texp \left(\int_x^y A_j \xi^j \right) \; A_k(y) \: \xi^j
\end{eqnarray}
mit den Randbedingungen
\begin{eqnarray}
\label{eq:a6_57}
\T \left(\int_x^x A_j \: (x-x)^j \right)^n &=& \left\{
	\begin{array}{cl} 1 & {\mbox{f\"ur $n=0$}} \\
			  0 & {\mbox{sonst}} \end{array} \right. \\
\label{eq:a6_58}
\Texp \left(\int_x^x A_j \: (x-x)^j \right) &=& 1 \spc .
\end{eqnarray}
Durch \Ref{a6_55}-\Ref{a6_58} sind die zeitgeordneten Integrale vollst\"andig
charakterisiert.
F\"ur drei Punkte $x, y, z$, die auf einer Geraden liegen, gilt
\Equ{a6_50b}
\Texp \left(\int_x^z A_j \: (z-x)^j \right) \cdot \Texp \left(\int_z^y
	A_k \: (y-z)^k\right) \;=\; \Texp \left(\int_x^y A_j (y-x)^j
	\right) \;\; .
\EndEqu
F\"ur die Adjungierten (der Matrizen) hat man
\begin{eqnarray}
\label{eq:a6_59}
\left[ \T (\int_x^y A_j \: (y-x)^j)^n \right]^* &=& \T \left(- \int_y^x A_j
	\: (x-y)^j \right)^n \\
\label{eq:a6_60}
\left[ \Texp (\int_x^y A_j \: (y-x)^j) \right]^* &=& \Texp \left(-
	\int_y^x A_j \: (x-y)^j \right) \spc .
\end{eqnarray}
\end{Lemma}
{\Beweis}
Nach Wahl einer Parametrisierung $z_j=\lambda_j \: y + (1-\lambda_j)x$
gilt
\[ \T \left(\int_x^y A_j \xi^j\right)^n \;=\;
	\int_0^1 d\lambda_1 \int_0^{\lambda_1}
	d\lambda_2 \cdots \int_0^{\lambda_{n-1}} d\lambda_n \;
	A_{j_n}(z_n) \: \xi^{j_n} \cdots A_{j_1}(z_1) \: \xi^{j_1} \spc . \]
Damit haben wir die zeitgeordneten Integrale auf die Dysonreihe
zur\"uckgef\"uhrt; der einzige Unterschied bei unserer Definition besteht
darin, da{\ss} die \"`sp\"ateste\"' Matrix $A(z_1)$ ganz rechts (und nicht
ganz links)
steht. Die Konvergenz von \Ref{a6_54} folgt mit der Absch\"atzung
\Equ{a6_106}
\left\| \: \T (\int_x^y A_j \xi^j)^n \: \right\| \;\leq\; \frac{1}{n!} \:
	\sup_{\lambda \in [0,1]}
	\left\| A_j(\lambda y + (1-\lambda) y) \: \xi^j \right\|^n \spc .
\EndEqu
Die Differentialgleichungen \Ref{a6_55}, \Ref{a6_56} sowie
\Ref{a6_57}, \Ref{a6_58} erh\"alt man nach der
Umformung
\[ \xi^k \frac{\partial}{\partial y^k} \T \left(\int_x^y  A_j \xi^j
	\right)^n \;=\; \frac{d}{ds}_{|s=1} \int_0^s \!d\lambda_1
	\int_0^{\lambda_1} \! d\lambda_2 \cdots \int_0^{\lambda_{n-1}}
	\! d\lambda_n \;
	A_{j_n}(z_n) \: \xi^{j_n} \cdots A_{j_1}(z_1) \: \xi^{j_1} \]
genau wie bei der Dysonreihe, \Ref{a6_50b} ist eine Konsequenz von
\Ref{a6_56}, \Ref{a6_58}.
Zum Beweis von \Ref{a6_59} mu{\ss} man eine Variablentransformation durchf\"uhren,
wir betrachten exemplarisch den Fall $n=2$. Wir haben
\begin{eqnarray*}
\left[ \T (\int_x^y A_j \: (y-x)^j)^2 \right]^* &=& \int_0^1 d\lambda_1
	\int_0^{\lambda_1} d\lambda_2 \; A^*_{j_1}(z_1) \: A^*_{j_2}(z_2) \;
	\xi^{j_1} \: \xi^{j_2} \\
&=& \int_0^1 d\bar{\lambda}_1 \int_{\bar{\lambda}_1}^1 d\bar{\lambda}_2
	\; A^*_{j_1}(\bar{z}_1) \: A^*_{j_2}(\bar{z}_2) \;
	\xi^{j_1} \: \xi^{j_2}
\end{eqnarray*}
mit $\bar{z}_j=\bar{\lambda}_j \: y + (1-\bar{\lambda}_j) x$ und somit
\[ \;=\; \int_0^1 d\bar{\lambda}_2 \int_0^{\bar{\lambda}_2}
	d\bar{\lambda}_1
	\; A^*_{j_1}(\bar{z}_1) \: A^*_{j_2}(\bar{z}_2) \;
	\xi^{j_1} \: \xi^{j_2} \;=\; \T \left( -\int_y^x A_j \:
	(x-y)^j \right)^2 \spc . \]
Gleichung \Ref{a6_60} folgt aus \Ref{a6_59} durch Summation \"uber $n$.
\QED
Durch Bildung der Adjungierten von \Ref{a6_55}, \Ref{a6_56} erh\"alt man
die Differentialgleichungen
\begin{eqnarray}
\label{eq:a6_55a}
\xi^k \frac{\partial}{\partial x^k} \T \left(\int_x^y  A_j \xi^j
	\right)^n &=& - A_k(x) \: \xi^k \;\T\left(\int_x^y A_j \xi^j
	\right)^{n-1} \;\;\;, n \geq 1 \\
\label{eq:a6_56a}
\xi^k \frac{\partial}{\partial x^k} \Texp \left(\int_x^y A_j \xi^j \right)
	&=& -A_j(x) \: \xi^j \;\Texp \left(\int_x^y A_j \xi^j \right)
	\spc .
\end{eqnarray}

F\"ur die Ableitung von $\Texp(\int_x^y A_j \: \xi^j)$ in Richtung des
Vektors $y-x$ gelten die einfachen Regeln \Ref{a6_56} und \Ref{a6_56a}.
F\"ur partielle Ableitungen haben wir i.a. keine entsprechenden einfachen
Gleichungen, also insbesondere
\begin{eqnarray*}
\frac{\partial}{\partial y^k} \: \Texp \left( \int_x^y A_j \: \xi^j \right)
	&\neq&  \Texp \left( \int_x^y A_j \: \xi^j \right) \: A_k(y) \\
\frac{\partial}{\partial x^k} \: \Texp \left( \int_x^y A_j \: \xi^j \right)
	&\neq&  - A_k(x) \: \Texp \left( \int_x^y A_j \: \xi^j \right) \spc .
\end{eqnarray*}
F\"ur eine \"ubersichtliche Notation ist ein solcher ``Ableitungsoperator'' dennoch
n\"utzlich:
\begin{Def}
Wir vereinbaren die Schreibweise
\begin{eqnarray*}
\hat{\Pdd}_y \:  \Texp \left( \int_x^y A_j \: \xi^j \right) &:=&
	\Texp \left( \int_x^y A_j \: \xi^j \right) \: \Aslsh(y) \\
\hat{\Pdd}_x \:  \Texp \left( \int_x^y A_j \: \xi^j \right) &:=& -\Aslsh(x) \:
	\Texp \left( \int_x^y A_j \: \xi^j \right) \spc .
\end{eqnarray*}
Auf alle anderen Funktionen wirkt $\hat{\Pdd}$ wie der Operator $\Pdd$, also
beispielsweise
\begin{eqnarray*}
\lefteqn{ \hat{\Pdd}_z \: \Texp \left( \int_x^z A_j \: (z-x)^j \right) \: f(z) \:
	\Texp \left( \int_z^y B_k \: (y-z)^j \right) } \\
&:=& \Texp \left( \int_x^z A_j \: (z-x)^j \right) \:(\Aslsh \:f + (\Pdd f) -
	f \: B \! \slsh)_{|z} \: \Texp \left( \int_z^y B_k \: (y-z)^k \right) \;\;\; .
\end{eqnarray*}
\end{Def}
Wir verwenden f\"ur das zeitgeordnete Exponential \Ref{a6_54} auch die
k\"urzere Schreibweise
\[ \T e^{\int_x^y A_j \xi^j} \spc . \]
\begin{Lemma}
\label{a6_lemma8}
Es gilt die asymptotische Entwicklung
\begin{eqnarray*}
\lefteqn{ \chi_L \: \tilde{s}^\vee_m(x,y) \;=\; \chi_L \: \Texp \left(-i \int_x^y
	A_L^j \: \xi_j \right) \; s^\vee_m(x,y) } \\
&&-\frac{1}{2}\:\chi_L \int_x^y dz \; (2 \alpha-1) \;\T e^{-i \int_x^z
	A_L^a \: (z-x)_a} \; \xi_j \: \gamma_k \: F_L^{kj}
	\;\T e^{-i \int_z^y A_L^b \: (y-z)_b} \; \sve(x,y) \\
&&+\frac{1}{2}\:\chi_L \int_x^y dz \;(\alpha^2-\alpha)\;\T e^{-i \int_x^z
	A_L^a \: (z-x)_a} \;\xi \slsh \: \xi_k \: j_L^k
	\;\T e^{-i \int_z^y A_L^b \: (y-z)_b} \; \sve(x,y) \\
&&-\frac{i}{4}\:\chi_L \int_x^y dz \;\T e^{-i \int_x^z
	A_L^a \: (z-x)_a} \;\varepsilon_{ijkl} \: F_L^{ij} \: \xi^k
	\; \rho \gamma^l
	\;\T e^{-i \int_z^y A_L^b \: (y-z)_b} \; \sve(x,y) \\
&&-\frac{m}{2}\:\chi_L \int_x^y dz \;\T e^{-i \int_x^z
	A_L^a \: (z-x)_a} \;(-i \Aslsh_L(z) + i \Aslsh_R(z))
	\:\xi\slsh \;\T e^{-i \int_z^y A_R^b \: (y-z)_b} \; \sve(x,y) \\
&&+ {\cal{O}}(\xi^0) \;+\; {\cal{O}}(m^2)
\end{eqnarray*}
mit
\begin{eqnarray}
\label{eq:a6_62}
F_L^{jk} &=& \partial^j A_L^k - \partial^k A_L^j - i [A_L^j,A_L^k] \\
\label{eq:a6_63}
j_L^k &=& g_{ml} \: \left[ \partial^m - i A_L^m, \: F_L^{kl} \right] \spc .
\end{eqnarray}
F\"ur die rechtsh\"andige Komponente hat man die analoge Gleichung, wenn man
die Indizes $L,R$ vertauscht.
\end{Lemma}
{\Beweis}
Wir entwickeln die einzelnen Summanden der Gleichung
\Equ{a6_60a}
\chi_L \: \tilde{s}_m^\vee(x,y) \;=\; \chi_L \: \sum_{n=0}^\infty \;
	\left(-s_m^\vee \: {\cal{B}}\right)^n \: s_m^\vee
\EndEqu
nach Potenzen von $m$ und wenden iterativ Lemma \ref{a6_lemma6} an.
Zur Ordnung ${\cal{O}}(\ln(|\xi^2|))$ brauchen wir nach Lemma
\ref{a6_lemma0} nur die drei folgenden F\"alle zu ber\"ucksichtigen:
\begin{description}
\item{a)} Alle Faktoren $s^\vee$ in \Ref{a6_60a} sind $s_0^\vee$, und bei
der asymptotischen Entwicklung von $s_0^\vee \:.\: s_0^\vee$ tritt immer der
f\"uhrende Beitrag \Ref{a6_59a} auf.
\item{b)} Alle Faktoren $s^\vee$ in \Ref{a6_60a} sind $s_0^\vee$, und bei
der asymptotischen Entwicklung von $s_0^\vee \:.\: s_0^\vee$ tritt genau einmal
der Beitrag \Ref{a6_59b}-\Ref{a6_59d} auf.
\item{c)} Genau ein Faktor $s^\vee$ in \Ref{a6_60a} ist $m \: \se$, und bei
der asymptotischen Entwicklung von $s_0^\vee \:.\: s_0^\vee$ tritt immer der
Beitrag \Ref{a6_59a} auf.
\end{description}
Wir betrachten diese F\"alle nacheinander, mit \"`$\asymp$\"' bezeichnen wir die
jeweils auftretenden Beitr\"age:
\begin{description}
\item{Zu a)} Wir haben
\begin{eqnarray*}
\lefteqn{ \chi_L \: \left((-\svn \:{\cal{B}})^n \: \svn \right)(x,y)
	\;=\; \chi_L \: \left((-\svn \: \Aslsh_L)^n \: \svn \right)(x,y) } \\
&=& -\chi_L \int d^4z \; \left(\svn \:\Aslsh_L\: \svn\right)(x,z) \:
	\left(-\Aslsh_L \: \svn \right)^{n-1}(z,y) \\
&\asymp& \chi_L \int d^4z \; \svn(x,z) \:\left(-i \int_x^z A_L^j\: (z-x)_j \right)\:
	\left(-\Aslsh_L \: \svn \right)^{n-1}(z,y) \\
&=& \cdots \;=\; \chi_L \int d^4z \; \svn(x,z) \: \T\left(-i\int_x^z
	A_L^j \: (z-x)_j\right)^p \: \left(-\Aslsh_L \:
	\svn\right)^{n-p}(z,y) \\
&=& \chi_L \: \T\left(-i\int_x^y A_L^j \: \xi_j\right)^n \; \svn(x,y) \spc .
\end{eqnarray*}
Nach Summation \"uber $n$ erh\"alt man
\Equ{a6_785a}
\chi_L \: \tilde{s}_m^\vee(x,y) \;\asymp\; \chi_L \; \Texp \left(-i \int_x^y
	A_L^j \: \xi_j \right) \; s_0^\vee(x,y) \spc .
\EndEqu
\item{Zu b)} Wir nehmen an, da{\ss} mit $n=p+q+1$, $p,q\geq0$ bei dem
$(p+1)$-ten Faktor
${\cal{B}}$ die niedrigeren Entwicklungsterme \Ref{a6_59b}-\Ref{a6_59d}
auftreten. Mit analogen Umformungen wie unter a) erh\"alt man den Ausdruck
\begin{eqnarray*}
\chi_L \: \left((-\svn \:{\cal{B}})^n \: \svn \right)(x,y)
&\asymp& \chi_L \int d^4z \; \svn(x,z) \; \T\left(-i\int_x^z A_L^j \:
	(z-x)_j\right)^p \; \Aslsh_L(z) \\
&& \hspace*{1cm} \times \; \T\left(-i\int_z^y A_L^k \: (y-z)_k\right)^q
	\; \svn(z,y) \spc .
\end{eqnarray*}
Wir summieren \"uber $p$, $q$
\begin{eqnarray}
\chi_L \: \tilde{s}_m^\vee(x,y) &\asymp& \chi_L \int d^4z \; \svn(x,z) \; \Texp
	\left(-i\int_x^z A_L^j \: (z-x)_j\right) \; \Aslsh_L(z) \nonumber \\
\label{eq:a6_61}
&& \hspace*{1.5cm} \times \; \Texp \left(-i\int_z^y A_L^k \: (y-z)_k\right)
	\; \svn(z,y)
\end{eqnarray}
und bestimmen die Entwicklungsterme \Ref{a6_59b}-\Ref{a6_59d}.
Bei der Berechnung von $F_{jk}$, $j^k$ mu{\ss}
man beachten, da{\ss} auch die zeitgeordneten Integrale in \Ref{a6_61}
differenziert werden m\"ussen. Diese Ableitungen k\"onnen zun\"achst nicht
ausgef\"uhrt werden, nach Kontraktion mit $\xi^j$ ist aber \Ref{a6_55}
anwendbar, und es ergeben sich f\"ur Feldst\"arke und Strom die Ausdr\"ucke
\Ref{a6_62}, \Ref{a6_63}.
\item{Zu c)}
F\"ur $p, q \geq 1$ gilt
\begin{eqnarray*}
\lefteqn{\chi_L \: \left( (-\svn \: {\cal{B}})^p \: \sve \:
	(-{\cal{B}} \: \svn)^q \:\right)(x,y) } \\
&=& \chi_L \: \left( (-\svn \: \Aslsh_L)^p \: \sve \: (-\Aslsh_R \: \svn)^q \:
	\right)(x,y) \\
&\asymp& \chi_L \int d^4z_1 \int d^4z_2 \; \svn(x,z_1) \: \T\left(-i\int_x^{z_1}
	A_L^j \: (z_1-x)_j\right)^{p-1} \; \Aslsh_L(z_1) \\
&& \hspace*{.5cm} \times \; \sve(z_1,z_2) \; \Aslsh_R(z_2) \;
	\T\left(-i\int_{z_2}^y A_L^j \: (y-z_2)_j \right)^{q-1} \;
	\svn(z_2,y) \spc .
\end{eqnarray*}
Wenn $p$ oder $q$ verschwinden, haben wir eine analoge Gleichung.
Summation \"uber $p$, $q$ liefert
\begin{eqnarray}
\lefteqn{\chi_L \: \tilde{s}_m^\vee(x,y) \;\asymp\; m \: \chi_L \:
	\sum_{p,q=0}^\infty \: \left( (-\svn \: \Aslsh_L)^p \:\sve\:
	(-\Aslsh_R \: \svn)^q \right)(x,y) } \nonumber \\
&=& \chi_L \: \sve(x,y) \;-\; \chi_L \: \int d^4z \; \sve(x,z) \;\Aslsh_R(z) \:
	\T e^{-i\int_z^y A_R^j \: (y-z)_j} \; \svn(z,y) \nonumber \\
&& - \chi_L \: \int d^4z \: \svn(x,z) \; \T e^{-i \int_x^z A_L^j \: (z-x)_j}
	\: \Aslsh_L(z) \; \sve(z,y) \nonumber \\
&&+ \chi_L \int d^4z_1 \int d^4z_2 \; \svn(x,z_1) \; \Texp
	\left(-i\int_x^{z_1} A_L^j \: (z_1-x)_j\right) \: \Aslsh_L(z_1) \nonumber \\
\label{eq:a6_787a}
&& \hspace*{.5cm} \times \; \sve(z_1,z_2) \; \Aslsh_R(z_2) \;
	\Texp \left(-i\int_{z_2}^y A_L^k \: (y-z_2)_k \right) \;
	\svn(z_2,y) \;\;\; . \hspace*{1cm}
\end{eqnarray}
Wir wenden nacheinander \Ref{a6_51}, \Ref{a6_52} an und vereinfachen die
Diracmatrizen mit \Ref{a6_31}. Mit Hilfe von \Ref{a6_56}, \Ref{a6_56a}
k\"onnen wir die Faktoren $\Aslsh^j_{L\!/\!R} \: \xi_j$ als Ableitung umschreiben
und partiell integrieren.
\end{description}
\QED
Man beachte, da{\ss} man f\"ur $F_{ij}$, $j^k$ jetzt der Ausdruck der
nichtabelschen Eichtheorien auftritt. Das ist nat\"urlich kein Zufall, sondern
war nach dem bekannten Eichtransformationsverhalten von $s^\vee_m$ zu erwarten.

Nach der Definitionsgleichung \Ref{2_tm} erhalten wir aus Lemma
\ref{a6_lemma8} unmittelbar auch eine asymptotische Entwicklung f\"ur
$\tilde{k}_m$. Um auch eine Gleichung f\"ur $\tilde{p}_m$ abzuleiten,
ben\"otigen wir folgendes Lemma:
\begin{Lemma}
\label{a6_lemma7}
Es sei $E(x)$ ein beliebiges Matrixfeld.
Die asymptotischen Entwicklungsformeln bis zur Ordnung
${\cal{O}}(\ln(|\xi^2|)$ der Distributionen in jeder Zeile von
\[      \begin{array}{ccc}
(k_0 \: E \: p_0)(x,y) &\;\;\;\;\;,\;\;\;\;\;& (k_0 \: E \: k_0)(x,y) \\[0.2em]
(\ke \: E \: p_0)(x,y) &\;\;\;\;\;,\;\;\;\;\;& (\ke \: E \: k_0)(x,y) \\[0.2em]
(k_0 \: E \: \pe)(x,y) &\;\;\;\;\;,\;\;\;\;\;& (k_0 \: E \: \ke)(x,y) \\[0.2em]
(s_0 \: E \: p_0)(x,y) &\;\;\;\;\;,\;\;\;\;\;& (s_0 \: E \: k_0)(x,y)
\end{array}     \]
unterscheiden sich jeweils nur durch die symbolische Ersetzung
$p \leftrightarrow k$.
\end{Lemma}
{\Beweis}
Nach \Ref{a1_75}, \Ref{a1_76} und \Ref{a3_170a} haben wir die Gleichungen
\[      \begin{array}{rclcrcl} \displaystyle
(K_0 \:E\: P_0)(x,y) &=& \displaystyle -\frac{i}{16 \pi^5} \: \veeint_x^y \hat{E}
	&\;\;\;\:,\:\;\;\;&
(K_0 \:E\: K_0)(x,y) &=& \displaystyle \frac{1}{16 \pi^4} \: \xint_x^y E
\\[.5em]
(S_0 \:E\: P_0)(x,y) &=& \displaystyle \frac{1}{16 \pi^4} \: \veeint_x^y E
	&\;\;\;\:,\:\;\;\;&
(S_0 \:E\: K_0)(x,y) &=& \displaystyle -\frac{i}{16 \pi^3} \: \xint_x^y \hat{E}
\end{array} \]
mit $\hat{E}(z)=\epsilon(z^0-x^0) \: E(z)$.
Die gesuchten Distributionen erh\"alt man daraus durch
Differentiation. Die Behauptung folgt aus der Analogie der
Ableitungsregeln \ref{a3_satz2}, \ref{a3_satz13} in Verbindung mit
den asymptotischen Entwicklungsformeln von Satz \ref{a3_satz98}.
\QED

\begin{Satz}
\label{a6_satz6}
F\"ur $\tilde{p}_m, \tilde{k}_m$ gilt mit der symbolischen Ersetzung
$C_m=k_m$ oder $C_m=p_m$
\begin{eqnarray}
\lefteqn{ \chi_L \: \tilde{C}_m(x,y) \;=\; \chi_L \: \Texp \left(-i \int_x^y
	A_L^j \: \xi_j \right) \; C_m(x,y) } \nonumber \\
&&-\frac{1}{2}\:\chi_L \int_x^y dz \; (2 \alpha-1) \;\T e^{-i \int_x^z
	A_L^a \: (z-x)_a} \; \xi_j \: \gamma_k \: F_L^{kj}
	\;\T e^{-i \int_z^y A_L^b \: (y-z)_b} \; C^{(1)}(x,y)
	\nonumber \\
&&+\frac{1}{2}\:\chi_L \int_x^y dz \;(\alpha^2-\alpha)\;\T e^{-i \int_x^z
	A_L^a \: (z-x)_a} \;\xi \slsh \: \xi_k \: j_L^k
	\;\T e^{-i \int_z^y A_L^b \: (y-z)_b} \; C^{(1)}(x,y)
	\nonumber \\
&&-\frac{i}{4}\:\chi_L \int_x^y dz \;\T e^{-i \int_x^z
	A_L^a \: (z-x)_a} \;\varepsilon_{ijkl} \: F_L^{ij} \: \xi^k
	\; \rho \gamma^l
	\;\T e^{-i \int_z^y A_L^b \: (y-z)_b} \; C^{(1)}(x,y)
	\nonumber \\
&&-\frac{m}{2}\:\chi_L \int_x^y dz \;\T e^{-i \int_x^z
	A_L^a \: (z-x)_a} \;(-i \Aslsh_L(z) + i \Aslsh_R(z))
	\:\xi\slsh \;\T e^{-i \int_z^y A_R^b \: (y-z)_b} \; C^{(1)}(x,y)
	\nonumber \\
\label{eq:a6_84}
&&+ {\cal{O}}(\ln(|\xi^2|)) \;+\; {\cal{O}}(m^2)
\end{eqnarray}
mit $F_L^{ij}$, $j_L^k$ gem\"a{\ss} \Ref{a6_62}, \Ref{a6_63}.
F\"ur die rechtsh\"andige Komponente hat man die analoge Gleichung, wenn
man die Indizes $L$, $R$ vertauscht.
\end{Satz}
{\Beweis}
F\"ur $C_m=k_m$ ist die Behauptung nach \Ref{2_tm} klar.

F\"ur den Fall $C_m=p_m$ vergleichen wir \Ref{2_62a} und \Ref{2_123a}:
Setzt man in diese Gleichungen $V$ gem\"a{\ss} \Ref{2_62} ein, ergibt sich
f\"ur $\tilde{k}_m$ der Ausdruck \Ref{2_65} und f\"ur $\tilde{p}_m$
entsprechend
\begin{eqnarray}
\tilde{p}_m &=& \sum_{l_1, l_2 =0}^\infty \; (-1)^{l_1+l_2}
\sum_{ \scriptsize
	\begin{array}{cc} \scriptsize Q_1 \in {\cal{P}}(l_1) , \\
		\scriptsize \# Q_1 \; {\mbox{gerade}}
	\end{array} }
\sum_{ \scriptsize
	\begin{array}{cc} \scriptsize Q_2 \in {\cal{P}}(l_2) , \\
		\scriptsize \# Q_2 \; {\mbox{gerade}}
	\end{array} } \!\!\!
	(i \pi)^{\#Q_1 + \#Q_2} \; c(\#Q_1/2) \: c(\#Q_2/2) \nonumber \\
\label{eq:a6_86}
&& \hspace*{-1cm} \times \; A_m(Q_1, 1) \: {\cal{B}} \cdots {\cal{B}} \:
	A_m(Q_1, l_1) \:
	{\cal{B}} \: p_m \: {\cal{B}} \: A_m(Q_2, 1) \: {\cal{B}} \cdots
	{\cal{B}} \: A_m(Q_2, l_2) \;\;\; .
\end{eqnarray}
Aus \Ref{2_65} erh\"alt man Gleichung \Ref{a6_84} f\"ur den Fall $C_m=k_m$, indem man
jeden Summanden iterativ bis zur Ordnung ${\cal{O}}(\ln(|\xi^2|))$
entwickelt und die Summe ausf\"uhrt.
Die Summanden in \Ref{a6_86} unterscheiden sich von denjenigen in
\Ref{2_65} lediglich dadurch, da{\ss} ein Operator $k_m$ durch $p_m$ ersetzt
ist. Nach Lemma \ref{a6_lemma7} und \Ref{a6_4}, \Ref{a6_5} entspricht das
gerade einer symbolischen Ersetzung von $k_m$ durch $p_m$
in der Endformel.
\QED
Die abgeleiteten Gleichungen f\"ur $\tilde{p}_m, \tilde{k}_m$ sind eine
Verallgemeinerung von Satz \ref{a6_satz5}.
Wir k\"onnen wieder explizit \"uberpr\"ufen, da{\ss} $\tilde{p}_m, \tilde{k}_m$
hermitesch sind; dabei verwendet man anstelle von
\Ref{a6_753a} die Umformung
\begin{eqnarray*}
\lefteqn{ \int_x^y dz \; \T e^{-i \int_x^z A_L^a \: (z-x)_a} \;
	(-i A_L^j + i A_R^j) \: \xi_j \; \T e^{-i \int_z^y A_R^b
	\: (y-z)_b} } \\
&=& \int_x^y dz \; \xi^j \frac{\partial}{\partial z_j}
	\left\{ \Texp\left(-i \int_x^z A_L^a \: (z-x)_a \right) \; \Texp\left(
	-i \int_z^y A_R^b \: (y-z)_b \right) \right\} \\
&=& \Texp\left(-i \int_x^y A^j_L \: \xi_j \right) \;-\;
	\Texp\left(-i \int_x^y A^j_R \: \xi_j \right) \spc .
\end{eqnarray*}

\subsection{St\"orungsrechnung mit Massenasymmetrie}
\label{a6_sect4}
Satz \ref{a6_satz6} l\"a{\ss}t sich direkt auf die St\"orungsrechnung
mit Massenasymmetrie \"uber\-tra\-gen.

\begin{Satz}
\label{a6_satz7}
Mit der symbolischen Ersetzung $C=k$ oder $C=p$ gilt
\begin{eqnarray}
\lefteqn{ \chi_L \: \tilde{C}(x,y) \;=\; \chi_L \: \Texp \left(-i \int_x^y
	A_L^j \: \xi_j \right) \; C(x,y) } \nonumber \\
&&-\frac{1}{2}\:\chi_L \int_x^y dz \; (2 \alpha-1) \;\T e^{-i \int_x^z
	A_L^a \: (z-x)_a} \; \xi_j \: \gamma_k \: F_L^{kj}
	\;\T e^{-i \int_z^y A_L^b \: (y-z)_b} \; C^{(1)}(x,y)
	\nonumber \\
&&+\frac{1}{2}\:\chi_L \int_x^y dz \;(\alpha^2-\alpha)\;\T e^{-i \int_x^z
	A_L^a \: (z-x)_a} \;\xi \slsh \: \xi_k \: j_L^k
	\;\T e^{-i \int_z^y A_L^b \: (y-z)_b} \; C^{(1)}(x,y)
	\nonumber \\
&&-\frac{i}{4}\:\chi_L \int_x^y dz \;\T e^{-i \int_x^z
	A_L^a \: (z-x)_a} \;\varepsilon_{ijkl} \: F_L^{ij} \: \xi^k
	\; \rho \gamma^l
	\;\T e^{-i \int_z^y A_L^b \: (y-z)_b} \; C^{(1)}(x,y)
	\nonumber \\
&&-\frac{m}{2}\:\chi_L \int_x^y dz \;\Texp \left(-i \int_x^z
	A_L^a \: (z-x)_a\right) \;(-i \Aslsh_L(z) \: Y + i Y \: \Aslsh_R(z))
	\:\xi\slsh \nonumber \\
&& \hspace*{4cm} \times \;\Texp \left(-i \int_z^y A_R^b \: (y-z)_b\right)
	\; C^{(1)}(x,y) \nonumber \\
\label{eq:a6_101}
&&+ {\cal{O}}(\ln(|\xi^2|)) \;+\; {\cal{O}}(m^2)
\end{eqnarray}
mit $F_L^{ij}$, $j_L^k$ gem\"a{\ss} \Ref{a6_62}, \Ref{a6_63}.
F\"ur die rechtsh\"andige Komponente hat man die analoge Gleichung, wenn
man die Indizes $L$, $R$ vertauscht.
\end{Satz}
{\Beweis}
Der Beweis von Lemma \ref{a6_lemma8} gilt w\"ortlich auch mit
Massensymmetrie, wenn wir $s^\vee_m$ durch
\[ s^\vee \;=\; \frac{1}{2 \pi i} \: (k+s) \]
ersetzen. Wir haben
\[ s^\vee \;=\; s^\vee_0 \;+\; m \: Y \:  s^\vee_{(1)} \;+\; {\cal{O}}(m^2)
	\spc , \]
so da{\ss} sich alle Rechnungen und damit auch das Ergebnis
von Lemma \ref{a6_lemma8} bis auf die zus\"atzlich auftretende
Massenmatrix $Y$ \"ubertragen.
Die Behauptung folgt daraus genau wie Satz \ref{a6_satz6} mit
Hilfe von Lemma \ref{a6_lemma7}.
\QED

\subsection{Zus\"atzliche freie Asymmetrie von $P$}
Die St\"orung von $P$ wird durch die Gleichung
\Equ{a6_102}
\tilde{P} \;=\; V \: P \: V^*
\EndEqu
mit $V$ gem\"a{\ss} \Ref{2_62} beschrieben.
Bisher haben wir au{\ss}er in Abschnitt \ref{a6_ab1.1} nicht mit dem Operator
$V$ gearbeitet, sondern nur die St\"orungsrechnung f\"ur $s^\vee$ untersucht.
Wenn $P$ eine freie Asymmetrie enth\"alt (also, wenn der freie Projektor $P$ die
chirale Symmetrie bricht oder wenn die Massenasymmetrie
der Fermionen nicht in der
St\"orungsrechnung ber\"ucksichtigt wird) ist das nicht ausreichend, wir
m\"ussen den Operator $V$ berechnen und $\tilde{P}$ mit Hilfe von
Gleichung \Ref{a6_102} bestimmen.

Um die dabei auftretenden Schwierigkeiten zu erl\"autern, betrachten wir
eine Entwicklung bis zur Ordnung ${\cal{O}}(\xi^{-2})$.
Mit der Abk\"urzung $z=\lambda y + (1-\lambda)x$ und
$\hat{A}(z)=\epsilon(z^0-x^0) \: A(z)$ gilt
\begin{eqnarray*}
(s_0 \:\Aslsh\: p_0)(x,y) &=& -\frac{1}{16 \pi^4} \: \Pdd_x \left(
	\veeint_x^y \Aslsh \right) \Pdd_y \\
&=& \frac{i}{2} \inti d\lambda \; \epsilon(\lambda) \; A_j(z) \: \xi^j
	\; p_0(x,y) \;+\; {\cal{O}}(\xi^{-2}) \\
(k_0 \:\Aslsh\: p_0)(x,y) &=& \frac{i}{16 \pi^5} \: \Pdd_x \left(
	\xint_x^y \hat{A} \! \slsh \right) \Pdd_y \\
&=& \frac{1}{2\pi} \inti d\lambda \; A_j(z) \: \xi^j
	\; p_0(x,y) \;+\; {\cal{O}}(\xi^{-2}) \spc .
\end{eqnarray*}
Durch Iteration erh\"alt man daraus nach einer Variablentransformation
\begin{eqnarray}
\left( (s_0 \: \Aslsh)^n \: p_0 \right) &=& \left(\frac{i}{2} \right)^n \:
	\inti \!\! d\lambda_1 \; \epsilon(\lambda_1) \inti \!\! d\lambda_2 \:
	\epsilon(\lambda_2-\lambda_1) \cdots
	\inti \!\! d\lambda_n \; \epsilon(\lambda_n-\lambda_{n-1}) \nonumber \\
\label{eq:a6_103}
&& \hspace*{2cm} \times \; A_{j_1}(z_1) \: \xi^{j_1} \cdots A_{j_n} \:
	\xi^{j_n} \;+\; {\cal{O}}(\xi^{-2}) \\
\left( k_0 \: \Aslsh \:(s_0 \: \Aslsh)^n \: p_0 \right) &=&
	\frac{1}{2\pi} \left(\frac{i}{2} \right)^n \:
	\inti \!\! d\lambda_0  \inti \!\! d\lambda_1 \: \epsilon(\lambda_1-\lambda_0)
	\cdots
	\inti \!\! d\lambda_n \; \epsilon(\lambda_n-\lambda_{n-1}) \nonumber \\
\label{eq:a6_104}
&& \hspace*{2cm} \times \; A_{j_0}(z_0) \: \xi^{j_0} \cdots A_{j_n} \:
	\xi^{j_n} \;+\; {\cal{O}}(\xi^{-2})
\end{eqnarray}
mit $z_j=\lambda_j \: y + (1-\lambda_j) x$.
Bei der Berechnung von $V p_0$ nach \Ref{2_62} treten also geschachtelte
Integrale auf. Diese {\bf{nichtlokalen Linienintegrale}} haben
\"Ahnlichkeit mit den zeitgeordneten Integralen von Definition \ref{a6_def1},
als wesentlicher Unterschied ist das Integrationsgebiet jetzt aber nicht
beschr\"ankt.

Da{\ss} in den bereits abgeleiteten Formeln f\"ur $\tilde{p}_m$, $\tilde{k}_m$
keine nichtlokalen Linienintegrale auftreten, hat folgenden Grund:
wenn wir beispielsweise die Gleichungen
\[ \tilde{p}_0(x,y) \;=\; (V \:p_0\: V^*)(x,y) \;\;\;,\spc
	\tilde{k}_0(x,y) \;=\; (V \:k_0\: V^*)(x,y) \]
nach Potenzen von ${\cal{B}}$ entwickeln und die asymptotischen Formeln
\Ref{a6_103}, \Ref{a6_104} einsetzen, heben sich alle Beitr\"age der
Integrale au{\ss}erhalb der Verbindungsstrecke $\overline{xy}$ weg,
und man erh\"alt die zeitgeordneten Integrale \Ref{a6_53}.

Eine allgemeine mathematische Behandlung der nichtlokalen Linienintegrale ist
aufwendig und schwierig. Vor allem deswegen, weil bei der Reihe
\Equ{a6_105}
\sum_{n=0}^\infty \; \inti \!\! d\lambda_1 \; \epsilon(\lambda_1)
	\inti \!\! d\lambda_2 \: \epsilon(\lambda_2-\lambda_1) \cdots \!
	\inti \!\! d\lambda_n \; \epsilon(\lambda_n-\lambda_{n-1}) \;
	 A_{j_1}(z_1) \: \xi^{j_1} \cdots A_{j_n}(z_n) \: \xi^{j_n}
\EndEqu
Konvergenzprobleme auftreten. Die Absch\"atzung \Ref{a6_106} f\"ur
zeitgeordnete Integrale l\"a{\ss}t sich nicht auf \Ref{a6_105} \"ubertragen.
Damit ist nicht klar, ob unsere St\"orungsentwicklung auch dann noch
im Distributionssinne konvergiert, wenn nichtlokale Linienintegrale
auftreten.

Wie in \cite{F3} \"uberlegt wurde, m\"ussen die mit der Matrix $X$ nicht
kommutierenden Eichpotentiale bei einer Spindimension
von $4n$ die Form
\Equ{a6_107}
\chi_L \: iU_R (\Pdd U_R^{-1}) \;+\; \chi_R \: iU_L (\Pdd U_L^{-1})
\EndEqu
mit Matrixfeldern $U_{L\!/\!R}(x) \in U(n)$ haben.
Diese Annahme ist auch eine wesentliche technische Vereinfachung.
Wir werden $\tilde{P}$ nur f\"ur diesen Spezialfall berechnen.

Im ersten Schritt wollen wir annehmen, da{\ss} \Ref{a6_107} die einzige
St\"orung des Diracoperators ist, also
\Equ{a6_108}
G \;=\; i \Pdd \;+\; \chi_L \: iU_R (\Pdd U_R^{-1}) \;+\; \chi_R \: iU_L
	(\Pdd U_L^{-1}) \spc .
\EndEqu
Bei vierkomponentigen Wellenfunktionen geht \Ref{a6_108} in den Operator \Ref{a6_91}
mit
\[ U_{L\!/\!R}(x) \;=\; e^{i \Lambda_{L\!/\!R}(x)} \]
\"uber.
\begin{Satz}
\label{a6_satz10}
Mit der symbolischen Ersetzung $C=p$ oder $C=k$ gilt
\begin{eqnarray*}
\lefteqn{ \chi_L \: \tilde{C}(x,y) \;=\; \chi_L \: U_L(x) \:X_L\: U_L^{-1}(y) \;
	C_0(x,y) } \\
&&+m \: \chi_L \: U_L(x) \:X_L\: U_L^{-1}(y) \: Y \; C^{(1)}(x,y) \\
&&-\frac{m}{4} \: \chi_L \: U_L(x) \inti d\lambda \; \left\{ \epsilon(\lambda) \:
	\Pdd(U_L^{-1}(z) \:Y\: U_R(z)) \: X_R \: \xi\slsh \right. \\
&&\left. \hspace*{3cm} + \: \epsilon(1-\lambda)
	\:X_L\: \Pdd((U_L^{-1}(z) \:Y\: U_R(z)) \: \xi\slsh \right\} \;
	U_R^{-1}(y) \; C^{(1)}(x,y) \\
&&+ {\cal{O}}(\ln(|\xi^2|)) \spc ,
\end{eqnarray*}
wobei wir die Abk\"urzung $z=\lambda y + (1-\lambda)x$ verwenden.
F\"ur die rechtsh\"andige Komponente hat man die analoge Gleichung, wenn
man die Indizes $L$, $R$ vertauscht.
\end{Satz}
{\Beweis}
F\"ur ein unit\"ares Matrixfeld $U(x) \in U(n)$ erf\"ullt der Ausdruck $U(x) \: U^{-1}(y)$
die Differentialgleichung
\[ \xi^j \: \frac{\partial}{\partial y^j} \: U(x) \: U^{-1}(y) \;=\;
	U(x) \: U^{-1}(y) \:(U(y) \: \partial_j U^{-1}(y)) \: \xi^j \]
mit Anfangsbedingung $U(x) \: U^{-1}(x)=1$. Durch Vergleich mit \Ref{a6_56},
\Ref{a6_58} folgt
\[ \Texp \left( \int_x^y U (\partial_j U^{-1}) \right) \: \xi^j
	\;=\; U(x) \: U^{-1}(y) \spc . \]
F\"ur den Diracoperator \Ref{a6_108} verschwinden au{\ss}erdem Feldst\"arke
und Noetherstrom \Ref{a6_62}, \Ref{a6_63}.
Damit tragen in Lemma \ref{a6_lemma8} nur die Beitr\"age \Ref{a6_785a},
\Ref{a6_787a} bei, also
\begin{eqnarray*}
\chi_L \: \tilde{s}^\vee_m &=& \chi_L \: U_L \:s^\vee_0\: U_L^{-1} \\
&&+ \: m \: \chi_L \: \left(1-iU_L \:\svn\: (\Pdd U_L^{-1})\right) \:\sve\:
	\left(1+i (\Pdd U_R) \:\svn\: U_R^{-1} \right)
	+ {\cal{O}}(m^2) \;\;\; .
\end{eqnarray*}
Mit Hilfe von \Ref{2_tm} und dem Analogieargument von Lemma \ref{a6_lemma7} folgt
f\"ur $\tilde{p}_m$, $\tilde{k}_m$
\begin{eqnarray}
\lefteqn{ \chi_L \: \tilde{C}_m \;=\; \chi_L \: U_L \: C_0 \: U_L^{-1}
	\;+\; m \: \chi_L \: C^{(1)} } \nonumber \\
&& -\:im \: \chi_L \: U_L \left( s_0 \: (\Pdd U_L^{-1})
	\: C^{(1)} + C_0 \:(\Pdd U_L^{-1}) \: s^{(1)} \right)
	\nonumber \\
&& +\:im \: \chi_L \: \left(C^{(1)} \:(\Pdd U_R)\: s_0 +
	s^{(1)} \:(\Pdd U_R) \: C_0 \right) \: U_R^{-1} \nonumber \\
&& +\:m \: \chi_L \: U_L \: \left(s_0 \:(\Pdd U_L^{-1})
	\: C^{(1)} \:(\Pdd U_R) \: s_0
	+C_0 \: (\Pdd U_L^{-1})\: s^{(1)} \:
	(\Pdd U_R)\: s_0 \right. \nonumber \\
\label{eq:a6_110}
&&\hspace*{1cm}\left. +s_0 \:(\Pdd U_L^{-1})\: s^{(1)} \:
	(\Pdd U_R) \: C_0 \right) \: U_R^{-1}
	\;+\; {\cal{O}}(\ln(|\xi^2|)) \spc .
\end{eqnarray}
Dabei mu{\ss} man beachten, da{\ss} die Ausdr\"ucke
\[ k_0 \:(\Pdd U_L^{-1}) \: C^{(1)} \: (\Pdd U_R) \: k_0 \;\;,\;\;\;
 C_0 \:(\Pdd U_L^{-1}) \: k^{(1)} \: (\Pdd U_R) \: k_0 \;\;,\;\;\;
 k_0 \:(\Pdd U_L^{-1}) \: k^{(1)} \: (\Pdd U_R) \: C_0 \]
verschwinden, denn es gilt beispielsweise
\begin{eqnarray*}
k_0 \:(\Pdd U_L^{-1}) \: C^{(1)} \: (\Pdd U_R) \: k_0 &=&
	-  k_0 \:[i \Pdd,\: U_L^{-1}] \: C^{(1)} \: [i \Pdd, \:U_R] \: k_0 \\
&=&  k_0 \: U_L^{-1} \: (i\Pdd) C^{(1)} (i \Pdd) \: U_R \: k_0 \;=\; 0
	\spc .
\end{eqnarray*}
Gleichung \Ref{a6_110} ist eine unmittelbare Verallgemeinerung von
\Ref{a6_759a} auf den nichtabelschen Fall.
Wir k\"onnen \Ref{a6_110} in der Form
\Equ{a6_112a}
\tilde{C}_m \;=\; V \: C_m \: V^*
\EndEqu
mit dem Operator $V$ schreiben, der analog zu Lemma \ref{a6_lemma5} durch
\begin{eqnarray}
\chi_{L} \: V \: p_m &=&  \chi_{L} \:U_L \: p_0 \;+\; m \: \chi_{L} \:
	\left( 1-i U_L \: s_0 \: (\Pdd U_L^{-1}) \right) \: \pe \nonumber \\
\label{eq:a6_111}
&&+\;im\: \chi_{L} \: \left(1-i U_L\: s_0 \:
	(\Pdd U_L^{-1}) \right) \: \se \: (\Pdd U_R) \: p_0
	\;+\; {\cal{O}}(m^2) 
\end{eqnarray}
gegeben ist, was man genau wie im Beweis von Satz \ref{a6_satz1}
direkt nachrechnen kann.

Der Operator $V$, (\ref{eq:2_62}), ist gerade so konstruiert, da{\ss}
die Bedingung \Ref{a6_112a} erf\"ullt ist. Deswegen ist klar, da{\ss} der durch
\Ref{a6_111} gegebene Operator tats\"achlich mit \Ref{2_62} \"ubereinstimmt.
Aus diesem Grund k\"onnen wir Gleichung \Ref{a6_110} auf die
St\"orungsrechnung mit
Massenasymmetrie f\"ur $p$, $k$ \"ubertragen: die Massenasymmetrie
ber\"ucksichtigen wir wie in Satz \ref{a6_satz7} durch die Ersetzung
\Equ{a6_144}
\se \longrightarrow Y \; \se \spc ,
\EndEqu
die Asymmetrie in $P$ nach Vergleich von \Ref{a6_112a}, \Ref{a6_102} durch die
Ersetzungen
\Equ{a6_145}
C_0 \longrightarrow X \: C_0 \;\;\;,\spc C^{(1)} \longrightarrow Y
	\; C^{(1)} \spc .
\EndEqu
Nun f\"uhren wir eine asymptotische Entwicklung um den Lichtkegel durch.
Dabei k\"onnen wir Lemma \ref{a6_lemma2} anwenden, denn die dort abgeleiteten
Formeln sind offensichtlich auch dann g\"ultig, wenn $A_j, B_k$ miteinander
nicht kommutierende Matrizen sind. Man erh\"alt
\begin{eqnarray}
\lefteqn{ \chi_L \: \tilde{C}(x,y) \;=\; \chi_L \: U_L(x) \: X_L \: U_L^{-1}(y) \; C_0(x,y)
	\;+\; m \: \chi_L \: Y \: C^{(1)}(x,y) } \\
&&+\frac{m}{4} \: \chi_L \: U_L(x) \: \inti d\lambda \; \xi\slsh \:
	\left( \epsilon(\lambda) \; (\Pdd U_L^{-1})(z) \: Y \right.
	\nonumber \\
\label{eq:a6_112}
&&\left. \hspace*{3.5cm} + \:\epsilon(1-\lambda) \: X_L \: (\Pdd U_L^{-1})(z) \:
	Y \right) \; C^{(1)}(x,y) \\
&&-\frac{m}{4} \: \chi_L \:\inti d\lambda \;
	\left( \epsilon(\lambda) \; Y \: (\Pdd U_R)(z) \: X_R \right.
	\nonumber \\
\label{eq:a6_113}
&&\left. \hspace*{3.5cm} + \: \epsilon(1-\lambda) \: Y \: (\Pdd U_R)(z)
	\right) \: \xi\slsh \; C^{(1)}(x,y) \\
&&-\frac{m}{4} \: \chi_L \: U_L(x) \: \inti d\lambda \; \epsilon(\lambda) \:
	\: \left( (U_L^{-1})(z) \:Y\: (\Pdd U_R)(z) \right. \nonumber \\
&&\left. \hspace*{3.5cm} + \: (\Pdd U_L^{-1})(z) \: Y \: U_R(z) \right) \:
	X_R \: \xi\slsh \: U_R^{-1}(y) \; C^{(1)}(x,y) \\
&&-\frac{m}{4} \: \chi_L \: U_L(x) \: \inti d\lambda \; \epsilon(1-\lambda)
	\:X_L \: \left( (U_L^{-1})(z) \:Y\: (\Pdd U_R)(z) \right. \nonumber \\
&&\left. \hspace*{3.5cm} + \: (\Pdd U_L^{-1})(z) \: Y \: U_R(z) \right) \:
	\xi\slsh \: U_R^{-1}(y) \; C^{(1)}(x,y) \\
&&+\frac{m}{4} \: \chi_L \:\inti d\lambda \;
	\left( \epsilon(\lambda) \; Y \: (\Pdd U_R)(z) \: X_R \right.
	\nonumber \\
\label{eq:a6_116}
&&\left. \hspace*{3.5cm} + \: \epsilon(1-\lambda) \: Y \: (\Pdd U_R)(z)
	\right) \: \xi\slsh \; C^{(1)}(x,y) \\
&&+\frac{m}{4} \: \chi_L \: U_L(x) \: \inti d\lambda \;
	\left( \epsilon(\lambda) \; (\Pdd U_L^{-1})(z) \: Y \right.
	\nonumber \\
\label{eq:a6_117}
&&\left. \hspace*{3.5cm} + \: \epsilon(1-\lambda) \: X_L \: (\Pdd U_L^{-1})(z) \:
	Y \right) \: \xi\slsh \; C^{(1)}(x,y) \\
&&+ {\cal{O}}(\ln(|\xi^2)) \spc . \nonumber
\end{eqnarray}
Die Terme \Ref{a6_113} und \Ref{a6_116} heben sich weg. Bei
\Ref{a6_112}+\Ref{a6_117} kann man partiell integrieren, dabei
fallen die Randwerte bei $\lambda=\pm \infty$ wegen
$U(z_{|\lambda=\pm \infty})=1$ weg.
\QED

Wir wollen nun Satz \ref{a6_satz10} auf den allgemeineren Diracoperator
\Equ{a6_150}
G \;=\; \chi_L \: U_R (i \Pdd + \Aslsh_R) U_R^{-1} \;+\;
	\chi_R \: U_L (i \Pdd + \Aslsh_L) U_L^{-1}
\EndEqu
erweitern. Die links- und rechtsh\"andigen Potentiale $\Aslsh_{L\!/\!R}$
sollen mit $X$ kommutieren, also
\Equ{a6_150a}
\left[ \Aslsh_{L\!/\!R}(x), \: X \right] \;=\; 0 \spc {\mbox{f\"ur alle
	$x \in M$.}}
\EndEqu
Im Spezialfall $\Aslsh_{L\!/\!R} \equiv 0$ geht \Ref{a6_150} in den
Diracoperator \Ref{a6_108} \"uber.
\begin{Lemma}
Mit der symbolischen Ersetzung $C_m=p_m$ oder $C_m=k_m$ und der
Abk\"urzung $z=\lambda y + (1-\lambda)x$ gilt die
asymptotische Entwicklung
\begin{eqnarray}
\label{eq:a6_140}
\lefteqn{(s_0 \: \Aslsh \: C_0)(x,y) \;=\; \frac{i}{2} \: \left(
	\inti d\lambda \; \epsilon(\lambda) \; A_j(z) \: \xi^j \right)
	\; C_0(x,y) } \\
\label{eq:a6_141}
&&+\frac{1}{4} \: \left(\inti d\lambda \; (2 \lambda-1) \;
	\xi^j \: \gamma^k \: F_{kj}(z) \right) \; C^{(1)}(x,y) \\
\label{eq:a6_142}
&&-\frac{1}{4}\: \left(\inti d\lambda \;(\lambda^2-\lambda)\; \xi \slsh \:
	\xi_k \: j^k(z) \right) \; C^{(1)}(x,y) \\
\label{eq:a6_143}
&&+\frac{i}{8} \: \left(\inti d\lambda \; \varepsilon_{ijkl} \:
	F^{ij}(z) \: \xi^k \; \rho \gamma^l \right)
	\; C^{(1)}(x,y) \;+\; {\cal{O}}(\ln(|\xi^2|)
	\spc .
\end{eqnarray}
\end{Lemma}
{\Beweis}
Nach Lemma \ref{a6_lemma7} gen\"ugt es, die asymptotische Entwicklung
f\"ur den Fall $C_m=p_m$ zu beweisen. Die Behauptung ist damit ein
Spezialfall von Satz \ref{a3_satz97}.
\QED

\begin{Satz}
\label{a6_thm1}
F\"ur den Diracoperator \Ref{a6_150} gilt mit der symbolischen Ersetzung
$C=p$ oder $C=k$
\begin{eqnarray}
\label{eq:a6_151}
\lefteqn{\chi_L \: \tilde{C}(x,y) \;=\; \chi_L \: U_L(x) \:
	\Texp \left(-i \int_x^y A_L^j \: \xi_j\right) \:X_L\:
	U_L^{-1}(y) \; C_0(x,y) } \\
&&-\frac{1}{2}\:\chi_L \: U_L(x) \: X_L
	\int_x^y dz \; (2 \alpha-1) \;\T e^{-i \int_x^z
	A_L^a \: (z-x)_a} \; \xi_j \: \gamma_k \: F_L^{kj}(z) \nonumber \\
\label{eq:a6_152}
&&\hspace*{6cm} \times \; \T e^{-i \int_z^y A_L^b \: (y-z)_b} \:U_L^{-1}(y)
	\; C^{(1)}(x,y) \\
&&+\frac{1}{2}\:\chi_L \: U_L(x) \: X_L
	\int_x^y dz \;(\alpha^2-\alpha)\;\T e^{-i \int_x^z
	A_L^a \: (z-x)_a} \;\xi \slsh \: \xi_k \: j_L^k(z) \nonumber \\
\label{eq:a6_153}
&&\hspace*{6cm} \times \;\T e^{-i \int_z^y A_L^b \: (y-z)_b} \: U_L^{-1}(y)
	\; C^{(1)}(x,y) \\
&&-\frac{i}{4}\:\chi_L \: U_L(x) \: X_L
	\int_x^y dz \;\T e^{-i \int_x^z
	A_L^a \: (z-x)_a} \;\varepsilon_{ijkl} \: F_L^{ij}(z) \: \xi^k
	\; \rho \gamma^l\nonumber \\
\label{eq:a6_154}
&&\hspace*{6cm} \times \;\T e^{-i \int_z^y A_L^b \: (y-z)_b} \: U_L^{-1}(y)
	\; C^{(1)}(x,y) \\
\label{eq:a6_155}
&&+ m \: \chi_L \: U_L(x) \: \Texp \left(-i\int_x^y A_L^j \: \xi_j \right)
	\: X_L \: U_L^{-1}(y) \:Y \; C^{(1)}(x,y) \\
&&-\frac{m}{4} \: \chi_L \: U_L(x) \inti d\lambda \nonumber \\
&&\hspace*{.5cm} \times \; \left\{ \epsilon(\lambda) \: \hat{\Pdd}_z \left(
	\T e^{-i\int_x^z A_L^j \: (z-x)_j} \: U_L^{-1}(z) \:Y\: U_R(z) \:
	\T e^{-i\int_z^y A_R^k \: (y-z)_k} \right) \:X_R \:\xi\slsh \right.
	\nonumber \\
&&\hspace*{1cm} \left. +\:\epsilon(1- \lambda) \: X_L \: \hat{\Pdd}_z \left(
	\T e^{-i\int_x^z A_L^j \: (z-x)_j} \: U_L^{-1}(z) \:Y\: U_R(z) \:
	\T e^{-i\int_z^y A_R^k \: (y-z)_k} \right) \:\xi\slsh \right\}
	\nonumber \\
\label{eq:a6_156}
&&\hspace*{0.5cm} \times \; U_R^{-1}(y) \; C^{(1)}(x,y) \;+\;
	{\cal{O}}(\ln(|\xi^2|)) \;+\; {\cal{O}}(m^2)
\end{eqnarray}
mit $F_L^{ij}$, $j_L^k$ gem\"a{\ss} \Ref{a6_62}, \Ref{a6_63}. Zur
Abk\"urzung wurde $z=\lambda y + (1-\lambda)x$ gesetzt.
F\"ur die rechtsh\"andige Komponente gilt die analoge Gleichung, wenn
man die Indizes $L$, $R$ vertauscht.
\end{Satz}
{\Beweis}
F\"ur das zeitgeordnete Integral der Potentiale hat man
\[ \Texp \left( \int_x^y (-i U A_j U^{-1} \:+\: U (\partial_j U^{-1})) \:
	\xi^j \right) \;=\; U(x) \; \Texp \left(-i \int_x^y A_j \: \xi^j
	\right) \: U^{-1}(y) \;\;\; , \]
wie man durch partielle Ableitung in Richtung $\xi$ und Vergleich
mit der Differentialgleichung \Ref{a6_56} sowie \Ref{a6_58} verifiziert.

Bei den Operatoren $p_0, k_0, s_0^\vee$ verwenden wir die Kurzschreibweise
\Equ{a6_823a}
\left( \bar{\:.\:} \right)_{L\!/\!R}(x,y) \;:=\; \Texp \left(-i \int_x^y A_{L\!/\!R}^j \:
	\xi^j \right) \; (\:.\:)(x,y) \spc .
\EndEqu
Wir betrachten nun die F\"alle a) bis c) in Lemma \ref{a6_lemma8} nacheinander
und untersuchen, welche Beitr\"age sich jeweils f\"ur $\tilde{P}$
ergeben. Dazu schreiben wir die Formeln f\"ur $\tilde{p}_m$, $\tilde{k}_m$
so um, da{\ss} man die Gleichung f\"ur $p, k$ wie beim Beweis von Satz
\ref{a6_satz10} durch die Ersetzungen \Ref{a6_144}, \Ref{a6_145} ableiten
kann.
\begin{description}
\item{Zu a)}
Wir haben nach \Ref{a6_785a}
\[ \chi_L \: \tilde{s}_m^\vee \;\asymp\; \chi_L \: U_L \:\bar{s}_{0,L}^\vee
	\: U_L^{-1} \]
und damit f\"ur $\tilde{p}_m, \tilde{k}_m$
\[ \chi_L \: \tilde{C}_m \;\asymp\; \chi_L \:U_L\: \bar{C}_{0,L} \: U_L^{-1}
	\spc . \]
Nach der Ersetzung \Ref{a6_145} erh\"alt man \Ref{a6_151}.
\item{Zu b)}
Als Beitrag hat man nach \Ref{a6_61} die Entwicklungsterme
\Ref{a6_59b}-\Ref{a6_59d} von
\[ \chi_L \: \tilde{s}_m^\vee \;\asymp\; \chi_L \: U_L \: \bar{s}_{0,L}^\vee
	\: \left(\Aslsh_L + i(\Pdd U_L^{-1}) U_L \right) \:
	\bar{s}_{0,L}^\vee \: U_L^{-1} \spc . \]
F\"ur $\tilde{p}_m$, $\tilde{k}_m$ erh\"alt man folglich bei einer
asymptotischen Entwicklung die Beitr\"age \Ref{a6_141}-\Ref{a6_143} von
\begin{eqnarray*}
\chi_L \: \tilde{C}_m &\asymp& \chi_L \: U_L \left\{
	\bar{s}_{0,L} \: \left(\Aslsh_L + i(\Pdd U_L^{-1}) U_L \right) \:
		\bar{C}_{0,L} \right. \\
&& \hspace*{2.5cm}\left. + \: \bar{C}_{0,L} \: \left(\Aslsh_L +
	i(\Pdd U_L^{-1}) U_L \right) \:
	\bar{s}_{0,L} \right\} \: U_L^{-1} \spc .
\end{eqnarray*}
Da die Feldst\"arke und der Strom \Ref{a6_62}, \Ref{a6_63} des Potentials
$i(\Pdd U_L^{-1}) U_L$ verschwinden, ergeben sich nach den Ersetzungen
\Ref{a6_144}, \Ref{a6_145} die Summanden \Ref{a6_152}-\Ref{a6_154}.
\item{Zu c)}
Wir haben nach \Ref{a6_787a}
\begin{eqnarray*}
\chi_L \: \tilde{s}_m^\vee &\asymp& m \: \chi_L \:
	\left(1 - U_L \: \bar{s}_{0,L}^\vee \: (\Aslsh_L \: U_L^{-1} +
	i (\Pdd U_L^{-1})) \right) \; \sve \\
&& \hspace*{3cm} \times \; \left(1 - (U_R \: \Aslsh_R - i (\Pdd U_R)) \:
	\bar{s}_{0,R}^\vee \: U_R^{-1} \right)
\end{eqnarray*}
und somit f\"ur $\tilde{p}_m$, $\tilde{k}_m$
\begin{eqnarray*}
\lefteqn{\chi_L \: \tilde{C}_m \;\asymp\; m \: \chi_L \: C^{(1)} } \\
&&- m \: \chi_L \: U_L \: \left(\bar{s}_{0,L} \: (\Aslsh_L \: U_L^{-1} +
	i (\Pdd U_L^{-1})) \: C^{(1)} \;+\;
	\bar{C}_{0,L} \: (\Aslsh_L \: U_L^{-1} + i (\Pdd U_L^{-1})) \:
		s^{(1)} \right) \\
&&+ m \: \chi_L \: \left( C^{(1)} \: (U_R \: \Aslsh_R - i (\Pdd U_R))
	\: \bar{s}_{0,R} \;+\; s^{(1)} \: (U_R \: \Aslsh_R - i (\Pdd U_R))
	\: \bar{C}_{0,R} \right) \: U_R^{-1} \\
&&+ M \: \chi_L \: U_L \: \left( \bar{s}_{0,L} \: (\Aslsh_L \: U_L^{-1} +
	i (\Pdd U_L^{-1})) \:C^{(1)}\: (U_R \: \Aslsh_R - i (\Pdd U_R))
	\: \bar{s}_{0,L} \right. \\
&&\hspace*{2cm} +\:\bar{C}_{0,L} \: (\Aslsh_L \: U_L^{-1} +
	i (\Pdd U_L^{-1})) \:s^{(1)}\: (U_R \: \Aslsh_R - i (\Pdd U_R))
	\: \bar{s}_{0,L} \\
&&\hspace*{2cm} +\:\left. \bar{s}_{0,L} \: (\Aslsh_L \: U_L^{-1} +
	i (\Pdd U_L^{-1})) \:s^{(1)}\: (U_R \: \Aslsh_R - i (\Pdd U_R))
	\: \bar{C}_{0,L} \right) \: U_R^{-1} \;\;\; .
\end{eqnarray*}
Wir f\"uhren nun wieder die Ersetzungen \Ref{a6_144}, \Ref{a6_145}
durch. Wegen der Kommutatorbedingung \Ref{a6_150a} geht der Operator
$\bar{C}_{0,L}$ dabei auf sinnvolle Weise in
\[ \bar{C}_{0,L}(x,y) \;\longrightarrow\; \Texp\left(-i \int_x^y
	A_L^j \: \xi_j \right) \; X \; C_0(x,y) \]
\"uber. Nach asymptotischer Entwicklung mit den Formeln von Lemma
\ref{a6_lemma2} erh\"alt man \Ref{a6_155}, \Ref{a6_156}.
\end{description}
\QED

\subsection{Zus\"atzliche skalare/pseudoskalare St\"orung}
Wir betrachten nun den Fall, da{\ss} der Diracoperator \Ref{a6_150}
eine zus\"atzliche skalare/pseudoskalare St\"orung enth\"alt, also
\Equ{a6_200}
G \;=\; \chi_L \: U_R (i \Pdd + \Aslsh_R) U_R^{-1} \;+\;
	\chi_R \: U_L (i \Pdd + \Aslsh_L) U_L^{-1} \;-\;
	m \: \Xi \;-\; i \rho \:m\: \Phi
\EndEqu
mit Matrixfeldern
\[ \Xi(x) = (\Xi_{ab}(x))_{a,b=1,\ldots,n} \;\;\;,\spc
	\Phi(x) = (\Phi_{ab}(x))_{a,b=1,\ldots,n} \spc . \]
Die Massenasymmetrie und die freie Asymmetrie des Projektors $P$
ber\"ucksichtigen wir wieder durch die Matrizen $Y$ bzw. $X$.
Wir verwenden die Schreibweise
\Equ{a6_826a}
Y_L(x) \;=\; Y + \Xi(x) + i \Phi(x) \spc
Y_R(x) \;=\; Y + \Xi(x) - i \Phi(x) \spc ,
\EndEqu
haben also
\[ Y + \Xi + i \rho \Phi \;=\; \chi_R \: Y_L \:+\: \chi_L \: Y_R \spc . \]
Satz \ref{a6_thm1} l\"a{\ss}t sich leicht auf diese etwas allgemeinere
Situation \"ubertragen.
\begin{Thm}
\label{a6_thm0}
F\"ur den Diracoperator \Ref{a6_200} gilt mit der symbolischen Ersetzung
$C=p$ oder $C=k$
\begin{eqnarray}
\label{eq:a6_400}
\lefteqn{\chi_L \: \tilde{C}(x,y) \;=\; \chi_L \: U_L(x) \:
	\Texp \left(-i \int_x^y A_L^j \: \xi_j\right) \:X_L\:
	U_L^{-1}(y) \; C_0(x,y) } \\
&&-\frac{1}{2}\:\chi_L \: U_L(x) \: X_L 
	\int_x^y dz \; (2 \alpha-1) \;\T e^{-i \int_x^z
	A_L^a \: (z-x)_a} \; \xi_j \: \gamma_k \: F_L^{kj}(z) \nonumber \\
\label{eq:a6_401}
&&\hspace*{6cm} \times \; \T e^{-i \int_z^y A_L^b \: (y-z)_b} \:U_L^{-1}(y)
	\; C^{(1)}(x,y) \\
&&+\frac{1}{2}\:\chi_L \: U_L(x) \: X_L
	\int_x^y dz \;(\alpha^2-\alpha)\;\T e^{-i \int_x^z
	A_L^a \: (z-x)_a} \;\xi \slsh \: \xi_k \: j_L^k(z) \nonumber \\
\label{eq:a6_402}
&&\hspace*{6cm} \times \;\T e^{-i \int_z^y A_L^b \: (y-z)_b} \: U_L^{-1}(y)
	\; C^{(1)}(x,y) \\
&&-\frac{i}{4}\:\chi_L \: U_L(x) \: X_L
	\int_x^y dz \;\T e^{-i \int_x^z
	A_L^a \: (z-x)_a} \;\varepsilon_{ijkl} \: F_L^{ij}(z) \: \xi^k
	\; \rho \gamma^l\nonumber \\
\label{eq:a6_403}
&&\hspace*{6cm} \times \;\T e^{-i \int_z^y A_L^b \: (y-z)_b} \: U_L^{-1}(y)
	\; C^{(1)}(x,y) \\
\label{eq:a6_404}
&&+ m \: \chi_L \: U_L(x) \: \Texp \left(-i\int_x^y A_L^j \: \xi_j \right)
	\: X_L \: U_L^{-1}(y) \:Y_L(y) \; C^{(1)}(x,y) \\
&&-\frac{m}{4} \: \chi_L \: U_L(x) \inti d\lambda \nonumber \\
&&\hspace*{.5cm} \times \; \left\{ \epsilon(\lambda) \: \hat{\Pdd}_z \left(
	\T e^{-i\int_x^z A_L^j \: (z-x)_j} \: (U_L^{-1} \:Y_L\: U_R)_{|z} \:
	\T e^{-i\int_z^y A_R^k \: (y-z)_k} \right) X_R \:\xi\slsh \right.
	\nonumber \\
&&\hspace*{1cm} \left. +\:\epsilon(1- \lambda) \: X_L \: \hat{\Pdd}_z \left(
	\T e^{-i\int_x^z A_L^j \: (z-x)_j} \: (U_L^{-1}\:Y_L\: U_R)_{|z} \:
	\T e^{-i\int_z^y A_R^k \: (y-z)_k} \right) \:\xi\slsh \right\}
	\nonumber \\
\label{eq:a6_405}
&&\hspace*{0.5cm} \times \; U_R^{-1}(y) \; C^{(1)}(x,y) \;+\;
	{\cal{O}}(\ln(|\xi^2|)) \;+\; {\cal{O}}(m^2)
\end{eqnarray}
mit $F_L^{ij}$, $j_L^k$ gem\"a{\ss} \Ref{a6_62}, \Ref{a6_63}. Zur
Abk\"urzung wurde $z=\lambda y + (1-\lambda)x$ gesetzt.
F\"ur die rechtsh\"andige Komponente gilt die analoge Gleichung, wenn
man die Indizes $L$, $R$ vertauscht.
\end{Thm}
{\Beweis}
Wir verwenden f\"ur $p_0, k_0, s^\vee_0$ wieder die Schreibweise
\Ref{a6_823a}.

Zun\"achst \"uberlegen wir uns, wie sich die zus\"atzliche St\"orung
$\Xi + i \rho \Phi$ in der Formel f\"ur $s^\vee_m$ auswirkt.
Nach Lemma \ref{a6_lemma0} und Theorem \ref{theorem_sk0},
Theorem \ref{a2_theorem_sm} haben wir als einzigen zus\"atzlichen Beitrag
in Lemma \ref{a6_lemma8} den Fall zu betrachten, da{\ss} alle Faktoren
$s^\vee$ in \Ref{a6_60a} $s_0^\vee$ sind und bei
den Faktoren ${\cal{B}}$ genau einmal die St\"orung $\Xi + i \rho \Phi$
auftritt. Damit folgt
\begin{eqnarray*}
\chi_L \: \tilde{s}^\vee_m &\asymp& \chi_L \:m\: U_L \: \bar{s}^\vee_{0,L}
	\: U^{-1}_L \:(\Xi + i \rho \Phi) \: U_R \: \bar{s}^\vee_{0,R}
	\: U^{-1}_R \\
&=& \chi_L \:m\: U_L \: \bar{s}^\vee_{0,L}
	\: U^{-1}_L \:(\Xi + i \Phi) \: U_R \: \bar{s}^\vee_{0,R}
	\: U^{-1}_R \;\;\;.
\end{eqnarray*}
Als Beitrag zu $\tilde{p}_m$, $\tilde{k}_m$ erh\"alt man
\begin{eqnarray*}
\chi_L \tilde{C}_m &\asymp& \chi_L \:m\: U_L \left( \bar{s}_{0,L} \:
	U^{-1}_L \: (\Xi + i \Phi) \: U_R \: \bar{C}_{0,R} \right. \\
&& \hspace*{1cm} \left. + \bar{C}_{0,L} \: U^{-1}_L \: (\Xi + i \Phi)
	\: U_R \: \bar{s}_{0,R} \right) U^{-1}_R \;\;\;.
\end{eqnarray*}
Die Gleichung f\"ur $\tilde{C}$ folgt genau wie im Beweis von Theorem
\ref{a6_thm1} durch die Ersetzung $C_0 \:\longrightarrow\: X C_0$
\begin{eqnarray*}
\chi_L \tilde{C} &\asymp& \chi_L \:m\: U_L \left( \bar{s}_{0,L} \:
	U^{-1}_L \: (\Xi + i \Phi) \: U_R \: X_R \: \bar{C}_{0,R} \right. \\
&& \hspace*{1cm} \left. + X_L \: \bar{C}_{0,L} \: U^{-1}_L \: (\Xi + i \Phi)
	\: U_R \: \bar{s}_{0,R} \right) U^{-1}_R \;\;\;.
\end{eqnarray*}

F\"ur $C=p$ k\"onnen wir mit Hilfe von Gleichung \Ref{a3_395a} eine
asymptotische Entwicklung durchf\"uhren
\begin{eqnarray*}
\lefteqn{ \chi_L \: \tilde{p}(x,y) \;\asymp\; -\frac{1}{4} \: \chi_L \:m\:
	U_L(x) \; \inti d\lambda } \\
&&\hspace*{0.5cm} \times \; \left\{ \epsilon(\lambda) \: \Pdd_z \left(
	\T e^{-i\int_x^z A_L^j \: (z-x)_j} \right. \right. \\
&&\hspace{3cm} \left. \times \; U_L^{-1}(z) \:
		(\Xi(z) + i \Phi(z)) \: U_R(z) \:
	\T e^{-i\int_z^y A_R^k \: (y-z)_k} \right) \:X_R \:\xi\slsh  \\
&&\hspace*{1.5cm}  -\:\epsilon(1- \lambda) \: X_L\: \xi\slsh \: \Pdd_z \left(
	\T e^{-i\int_x^z A_L^j \: (z-x)_j} \: U_L^{-1}(z) \right. \\
&&\hspace{3cm} \left. \left. \times \; (\Xi(z) + i \Phi(z)) \: U_R(z) \:
	\T e^{-i\int_z^y A_R^k \: (y-z)_k} \right) \; \right\} \\
&&\hspace*{0.5cm} \times \; U_R^{-1}(y) \; p^{(1)}(x,y) \;+\;
	{\cal{O}}(\ln(|\xi^2|)) \spc .
\end{eqnarray*}
F\"ur $\tilde{k}_m$ folgt die Behauptung mit dem Analogieargument von
Lemma \ref{a6_lemma7}.
\QED

Theorem \ref{a6_thm0} gibt die Eichterme und Pseudoeichterme von $P$
bei allen f\"ur uns interessanten St\"orungen des Diracoperators an.
Wir erhalten daraus s\"amtliche in Abschnitt \ref{a6_ab1} abgeleiteten
S\"atze als Grenzfall: Lassen wir die skalare/pseudoskalare St\"orung weg,
so erhalten wir Satz \ref{a6_thm1}.
Wenn zus\"atzlich die dynamischen Eichpotentiale
$\Aslsh_{L\!/\!R}$ verschwinden, ergibt sich Satz \ref{a6_satz10} und
im Spezialfall vierkomponentiger Wellenfunktionen Satz \ref{a6_satz5}.

Besitzt $P$ auf der anderen Seite keine freie Asymmetrie (also $X=\1$), so
ist die Bedingung \Ref{a6_150a} in jedem Fall erf\"ullt. Es besteht dann
keine Notwendigkeit, die St\"orung des Diracoperators gem\"a{\ss}
\Ref{a6_150} in $\Aslsh_{L\!/\!R}$ und $U_{L\!/\!R}$ aufzuspalten. Nach
der Ersetzung
\[ \Aslsh_{L\!/\!R} \;\rightarrow\; U_{L\!/\!R} \: \Aslsh_{L\!/\!R} \: U_{L\!/\!R}^{-1}
	\;+\; i U_{L\!/\!R} \: (\Pdd U_{L\!/\!R}^{-1}) \]
k\"onnen wir annehmen, da{\ss} $U_{L\!/\!R} \equiv 1$ ist. Au{\ss}erdem lassen sich
die nichtlokalen Linienintegrale in \Ref{a6_156} nun zu dem konvexen
Integral
\[ \int_x^y d\lambda \; \hat{\Pdd}_z\left( \T e^{-i\int_x^z A_L^j \: (z-x)_j}
	\:Y\: \T e^{-i\int_z^y A_R^k \: (y-z)_k} \right) \:\xi\slsh \]
zusammenfassen.
Auf diese Weise erh\"alt man Satz \ref{a6_satz7}. Ist zus\"atzlich $Y=\1$, so
ergibt sich Satz \ref{a6_satz6}.

Man sieht an der abgeleiteten Formel f\"ur $\tilde{C}$, da{\ss} Theorem
\ref{a6_thm1} nicht auf einfache Weise weiter verallgemeinert werden
kann.
Bei einer Abschw\"achung von Bedingung \Ref{a6_150a} ist beispielsweise
die Matrix
\[ U_L(x) \: \Texp \left(-i \int_x^y A_L^j \: \xi_j\right) \:X_L \: \xi\slsh \:
	U_L^{-1}(y) \]
in \Ref{a6_400} nicht mehr selbstadjungiert. Wir m\"u{\ss}ten sie durch einen
Ausdruck
\[ {\cal{N}}(x) \: X \: {\cal{N}}^{-1}(y) \]
ersetzen, wobei ${\cal{N}}$ eine unendliche Reihe nichtlokaler
Linienintegrale der Form \Ref{a6_105} ist.

\section{Massenterme $\sim m^2$}
\label{a6_ab2}
In diesem Abschnitt wollen wir die Massenterme und Eichterme
$\sim m^2$ berechnen.
Wir ber\"ucksichtigen dabei in den Distributionen $\tilde{p}$ und
$\tilde{k}$ alle Beitr\"age bis zur Ordnung ${\cal{O}}(\xi^0)$ bzw.
${\cal{O}}(\xi^2)$.

\subsection{St\"orungsrechnung mit Massenasymmetrie}
Wir betrachten den Diracoperator \Ref{a6_108}.
Im Projektor $P$ lassen wir wieder eine Massenasymmetrie zu, die mit der
Matrix $Y$ beschrieben wird.

Zun\"achst betrachten wir die St\"orungsrechnung f\"ur $\tilde{s}^\vee_{(2)}$.
Durch Entwicklung von \Ref{2_t2} nach $m$ erh\"alt man
\begin{eqnarray*}
\tilde{s}^\vee_{(2)} &=& \sum_{p,q=0}^\infty \: (-\svn \: {\cal{B}})^p
	\:s^\vee_{(2)}\: Y^2 \: (-{\cal{B}} \: \svn)^q \\
&&+ \sum_{p,q,r=0}^\infty \: (-\svn \: {\cal{B}})^p \:\sve\:Y\:
	(-{\cal{B}} \: \svn)^q \:(-{\cal{B}} \: \sve \: Y) \:
	(-{\cal{B}} \: \svn)^r \spc ,
\end{eqnarray*}
also
\begin{eqnarray}
\label{eq:a6_254}
\lefteqn{ \chi_L \: \tilde{s}^\vee_{(2)} \;=\; \sum_{p,q=0}^\infty \: \chi_L
	\: (-\svn \: {\cal{B}}_L)^p
	\:s^\vee_{(2)}\:Y^2\: (-{\cal{B}}_L \: \svn)^q } \\
\label{eq:a6_255}
&&+ \sum_{p,q,r=0}^\infty \: \chi_L \: (-\svn \: {\cal{B}}_L)^p \:\sve\:Y\:
	(-{\cal{B}}_R \: \svn)^q \:(-{\cal{B}}_R \: \sve\:Y) \:
	(-{\cal{B}}_L \: \svn)^r
\end{eqnarray}
mit
\Equ{a6_828a}
{\cal{B}}_L = i U_L (\Pdd U_L^{-1}) \;\;\;\;\;,\spc {\cal{B}}_R = 
	i U_R (\Pdd U_L^{-1}) \spc .
\EndEqu
Wir wollen nun die Summen \Ref{a6_254} und \Ref{a6_255} vereinfachen,
dabei rechnen wir bis auf Terme der Ordnung ${\cal{O}}(\xi^2)$.
F\"ur die erste Summe hat man
\begin{eqnarray*}
\lefteqn{ \sum_{p,q=0}^\infty (-\svn \: {\cal{B}}_L)^p \: s^\vee_{(2)} \:
	Y^2 \: (-{\cal{B}}_L \: \svn)^q \: (x,y) } \\
&=& (1-U \: \svn \: U_L^{-1} \: {\cal{B}}_L) \: s^\vee_{(2)} \: Y^2 \:
	(1-{\cal{B}}_L \:U_L \: \svn \: U_L^{-1}) \\
&=& \left(1-U_L \: \svn \: [i \Pdd,\:U_L^{-1} \right) \: s^\vee_{(2)} \:
	Y^2 \: \left(1+[i \Pdd,\:U_L] \: \svn \: U_L^{-1} \right) \\
&=& U_L \: \svn \: U_L^{-1} \: (i \Pdd) s^\vee_{(2)} (i \Pdd) \: Y^2 \:
	U_L \: \svn \: U_L^{-1} \\
&=& U_L \: \svn \: U_L^{-1} \: \svn \: Y^2 \: U_L \:
	\svn \: U_L^{-1} \spc ,
\end{eqnarray*}
die zweite Summe kann man folgenderma{\ss}en umformen:
\begin{eqnarray}
\lefteqn{ \sum_{p,q,r=0}^\infty \: (-\svn \: {\cal{B}}_L)^p \:\sve\:Y\:
	(-{\cal{B}}_R \: \svn)^q \:(-{\cal{B}}_R \: \sve\:Y) \:
	(-{\cal{B}}_L \: \svn)^r } \nonumber \\
&=& \sum_{q=0}^\infty \: (1-U_L \: \svn \: U_L^{-1} \: {\cal{B}}_L)
	\: \sve \: Y \: (-{\cal{B}}_R \: \svn)^q \: (-{\cal{B}}_R \:
	\sve \: Y) \: (1-{\cal{B}}_L \: U_L \: \svn \: U_L^{-1}) \nonumber \\
&=& \sum_{q=0}^\infty \: U_L \: \svn \: U_L^{-1} \: Y \: \svn
	\: (-{\cal{B}}_R \: \svn)^q \: (-{\cal{B}}_R) \:
	\svn \: Y \:U_L\: \svn \: U_L^{-1}  \nonumber \\
&=& U_L \: \svn \: U_L^{-1} \:Y\: \svn \: (1-{\cal{B}}_R \: U_R \:
	\svn \: U_R^{-1} ) \: (-{\cal{B}}_R) \: \svn \:Y\: U_L \:
	\svn \: U_L^{-1} \nonumber \\
&=& U_L \: \svn \: (U_L^{-1} \:Y\: U_R) \: \svn \: U_R^{-1} \:
	(-{\cal{B}}_R) \: \svn \:Y\: U_L \: \svn \: U_L^{-1} \nonumber \\
&=& -U_L \: \svn \: (U_L^{-1} \:Y\: U_R) \: \svn \: [i \Pdd, \: U_R^{-1}]
	\: \svn \:Y\: U_L \: \svn \: U_L^{-1} \nonumber \\
&=& -U_L \: \svn \: U_L^{-1} \:Y^2\: \svn \: U_L \: \svn \: U_L^{-1} \nonumber \\
\label{eq:a6_835a}
&&+U_L \: \svn \: (U_L^{-1} \:Y\: U_R) \: \svn \:
	(U_R^{-1} \:Y\: U_L) \: \svn \: U_L^{-1}
\end{eqnarray}
Insgesamt erh\"alt man also
\begin{eqnarray}
\lefteqn{ \chi_L \: \tilde{s}^\vee_{(2)}(x,y) } \nonumber \\
\label{eq:a6_828b}
&=& \chi_L \: \left(
	U_L \: \svn \: (U_L^{-1} \:Y\: U_R) \: \svn \:
	(U_R^{-1} \:Y\: U_L) \: \svn \: U_L^{-1} \right)(x,y)
	\;+\; {\cal{O}}(\xi^2) \;\;\;.
\end{eqnarray}
Nun f\"uhren wir eine asymptotische Entwicklung um den Lichtkegel durch.
\begin{Lemma}
\label{a6_lemma10}
F\"ur beliebige glatte (auch matrixwertige) Funktionen $f,g$ gilt mit
der Ab\-k\"ur\-zung $z=\alpha y + (1-\alpha)x$
\begin{eqnarray*}
\lefteqn{ (\svn \:f\: \svn \:g\: \svn)(x,y)
\;=\; \int_x^y dz \; f(z) \: g(z) \; s^\vee_{(2)}(x,y) } \\
&&-\frac{i}{16 \pi} \: \int_x^y (\alpha^2-\alpha) \; \Box(f \:
	g)_{|z} \; \xi\slsh \; \Theta^\vee(\xi) \\
&&+\frac{i}{16 \pi} \: \int_x^y dz \int_x^z du \; (\Pdd f)(u) \:
	\: \gamma^j \: (z-x)_j \: (\Pdd g)(z) \: g(y) \;
	\Theta^\vee(\xi) \\
&&+\frac{i}{8 \pi} \: \int_x^y (1-\alpha) \; ((\Pdd f) \:
	g)_{|z} \; \Theta^\vee(\xi) \\
&&-\frac{i}{8 \pi} \int_x^y \alpha \; (f \: (\Pdd g))_{|z} \;
	\Theta^\vee(\xi)
	\;+\; {\cal{O}}(\xi^2) \spc .
\end{eqnarray*}
\end{Lemma}
{\Beweis}
Wir rechnen bis auf Terme der Ordnung ${\cal{O}}(\xi^2)$.
Wir berechnen das Operatorprodukt schrittweise.
Nach \Ref{a1_127a} haben wir
\begin{eqnarray}
\lefteqn{ (\svn \: f \: \svn)(x,y) \;=\; -\frac{1}{4 \pi^2} \:
	\Pdd_x \left(\lint_x^y f \right) \: \Pdd_y
\;=\; \frac{1}{4 \pi^2} \: \Pdd_x \Pdd_y \: \lint_x^y f } \nonumber \\
&=& -\frac{1}{2\pi} \: f(y) \: l^\vee(\xi) \:+\: \frac{1}{4\pi}
	\int_x^y (\Pdd f) \: \xi\slsh \: l^\vee(\xi) \nonumber \\
&&+\frac{1}{4 \pi^2} \lint_x^y dz \; \left( \Box f \:-\: \int_x^z
	\alpha \; \Box f \:-\: \frac{1}{2} \int_x^z \alpha^2
	\; (\Pdd f) \: \zeta \slsh \right) \nonumber \\
\label{eq:a6_207}
&=& -\frac{1}{2\pi} \: f(y) \: l^\vee(\xi) \:+\: \frac{1}{4\pi}
	\int_x^y (\Pdd f) \: \xi\slsh \: l^\vee(\xi) \\
\label{eq:a6_208}
&&+ \frac{1}{8 \pi} \int_x^y \alpha \; \Box f \; \Theta^\vee(\xi) \;-\;
	\frac{1}{16 \pi} \int_x^y (\alpha-\alpha^2) \: (\Pdd \Box f)
	\: \xi \slsh \; \Theta^\vee(\xi) \;\;\; .
\end{eqnarray}
Wir wollen jetzt die einzelnen Summanden in \Ref{a6_207}, \Ref{a6_208}
mit dem Operator $g\:\svn$ multiplizieren und asymptotisch
entwickeln. Mit ``$\asymp$'' bezeichnen wir die dabei jeweils auftretenden
Beitr\"age.
Die Summanden \Ref{a6_208} k\"onnen wir mit der Umformung
\begin{eqnarray*}
\lefteqn{ (\Theta^\vee \:h\: \svn)(x,y) \;=\; -\frac{1}{2\pi} \:
	(\Theta^\vee \:h\: (i \Pdd) \: l^\vee)(x,y) } \\
&=& -\frac{i}{2 \pi} \: \frac{\partial}{\partial x^j}\left( \Theta^\vee
	\:h\: l^\vee \right)(x,y) \: \gamma^j \;+\; {\cal{O}}(\xi^2)
\;=\; \frac{i}{2} \: \int_x^y \alpha \: h \; \xi\slsh
\end{eqnarray*}
behandeln, dabei ist $h$ die glatte Funktion
\Equ{a6_230a}
h(y) \;=\; \frac{1}{8 \pi} \int_x^y \alpha \; \Box f \;-\;
	\frac{1}{16 \pi} \int_x^y (\alpha-\alpha^2) \: (\Pdd \Box f)
	\: \xi \slsh \spc .
\EndEqu
Der zweite Summand in \Ref{a6_230a} f\"allt dabei weg, und wir erhalten
mit den konvexen Kombinationen
\[ z=\alpha y + (1-\alpha)x \;\;\;\;\;,\spc u=\beta y + (1-\beta) x \]
den Beitrag
\Equ{a6_210}
\Ref{a6_208} \;\asymp\; \frac{i}{16 \pi} \int_0^1 d\alpha \int_0^\alpha
	d\beta \; \beta \; (\Box f)(u) \: g(z) \; \Theta^\vee(\xi)
	\spc .
\EndEqu
Den ersten Summanden in \Ref{a6_207} k\"onnen wir mit dem Lichtkegelintegral
$\slint$ umschreiben
\begin{eqnarray}
-\frac{1}{2\pi} \: f(y) \: l^\vee(\xi) &\asymp&
\frac{1}{4 \pi^2} \: \left(\lint_x^y f \: g \right) \:
	(i \Pdd_y)
\;=\; -\frac{i}{4 \pi^2} \: \Pdd_y \lint_x^y f \: g \nonumber \\
&=& -\frac{i}{4 \pi} \: l^\vee(\xi) \: \xi\slsh \; \int_x^y f \:
	g \;-\; \frac{i}{8 \pi} \int_x^y \alpha \; \Pdd(f \: g) \;
	\Theta^\vee(\xi) \nonumber \\
&&-\frac{i}{16 \pi} \int_x^y (\alpha^2-\alpha) \: \Box(f \: g)(z)
	\; \xi\slsh \: \Theta^\vee(\xi) \nonumber \\
\label{eq:a6_211}
&=& \int_x^y f \: g \; s^\vee_{(2)}(x,y) \\
\label{eq:a6_212}
&&-\frac{i}{8 \pi} \int_x^y \alpha \; \Pdd(f \: g) \;
	\Theta^\vee(\xi) \\
\label{eq:a6_213}
&&-\frac{i}{16 \pi} \int_x^y (\alpha^2-\alpha) \: \Box(f \: g)(z)
	\; \xi\slsh \: \Theta^\vee(\xi) \spc .
\end{eqnarray}
Den zweiten Summanden in \Ref{a6_207} k\"onnen wir mit der
Abk\"urzung
\[ h^x_j(y) \;=\; \int_x^y dz \; (\partial_j f(z)) \: g(y) \]
einfach umformen
\begin{eqnarray}
\lefteqn{\frac{1}{4 \pi} \int_x^y (\Pdd f) \: \xi\slsh \; l^\vee(\xi)
	\;\asymp\; i \: \gamma^j \: (s^\vee_{(2)} \:h^x_j\: \svn)(x,y) }
	\nonumber \\
&=& i \gamma^j \int d^4z \; s^\vee_{(2)}(x,z) \: h^x_j(z) \: (i \Pdd_z)
	\: \sve(z,y) \nonumber \\
&=& -i \gamma^j \int d^4z \; \left[ (i \Pdd_z s^\vee_{(2)}(x,z)) \: h^x_j(z) \:+\:
	s^\vee_{(2)}(x,z) \: (\Pdd_z h^x_j(z)) \right] \: \sve(z,y) \nonumber \\
&=& i \: \gamma^j \: \left( \sve \:h^x_j\: \sve \;-\;
	s^\vee_{(2)} \: (i \Pdd h^x_j) \: \sve \right)(x,y) \nonumber \\
\label{eq:a6_209}
&=& \frac{i}{4 \pi^2} \lint_x^y k \;=\; \frac{i}{8 \pi} \int_x^y
	k \; \Theta^\vee(\xi) \spc ,
\end{eqnarray}
wobei $k$ f\"ur die Funktion
\[ k(y) \;=\; \gamma^j \: \left( h^x_j(y) \:+\: \frac{1}{2} \: \xi\slsh \:
	\Pdd_y h^x_j(y) \right) \]
steht.
Bei der Berechnung von $k$ erh\"alt man die Zwischenergebnisse
\begin{eqnarray*}
\Pdd_y h^x_j(y) &=& \int_x^y \alpha \; (\partial_j \Pdd f)(z)
	\: g(y) \;+\; \int_x^y (\partial_j f)(z) \:
	(\Pdd g)(y) \\
\gamma^j \: \xi\slsh \: \Pdd_y h^x_j(y) &=& 2 \int_x^y \alpha \;
	\left(\frac{d}{d\alpha} \: \Pdd f \right)\!(z) \:
	g(y) \\
&&\hspace*{-2cm} -\int_x^y \alpha \: (\Box f)(z) \: g(y) \:
	\xi\slsh \;+\; \int_x^y \alpha \: (\Pdd f)(z) \:
	\xi\slsh \: (\Pdd g)(y) \\
&&\hspace*{-3cm} =\; 2 ((\Pdd f) \: g)_{|y} \;-\; 2 \: \gamma^j h^x_j(y) \\
&&\hspace*{-2cm} -\int_x^y \alpha \: (\Box f)(z) \: g(y) \:
	\xi\slsh \;+\; \int_x^y \alpha \: (\Pdd f)(z) \: \xi\slsh \:
	(\Pdd g)(y)
\end{eqnarray*}
und somit
\begin{eqnarray*}
k(y) &=& (\Pdd f) \: g
\;-\; \frac{1}{2} \int_x^y \alpha \: (\Box f)(z) \: g(y) \: \xi\slsh
\;+\; \frac{1}{2} \int_x^y \alpha \: (\Pdd f)(z) \: \xi\slsh \: (\Pdd g)(y)
	\;\;\; .
\end{eqnarray*}
Wir setzen in \Ref{a6_209} ein und erhalten schlie{\ss}lich
\begin{eqnarray}
\label{eq:a6_214}
\lefteqn{\frac{1}{4 \pi} \int_x^y (\Pdd f) \: \xi\slsh \; l^\vee(\xi)
	\;\asymp\; \frac{i}{8 \pi} \int_x^y (\Pdd f) \: g \;
	\Theta^\vee(\xi) } \\
\label{eq:a6_215}
&&-\frac{i}{16 \pi} \int_0^1 d\alpha \int_0^\alpha d\beta \; \beta
	\; (\Box f)(u) \: g(z) \: \xi\slsh \; \Theta^\vee(\xi) \\
\label{eq:a6_216}
&&+\frac{i}{16 \pi} \int_0^1 d\alpha \int_0^\alpha d\beta \;
	(\Pdd f)(u) \: \xi\slsh \: (\Pdd g)(z) \; \Theta^\vee(\xi)
	\;\;\;.
\end{eqnarray}
Beim Aufsummieren aller Beitr\"age heben sich \Ref{a6_210}, \Ref{a6_215}
weg. Der gesuchte Ausdruck ist damit gleich
\[ \Ref{a6_211} + \Ref{a6_212} + \Ref{a6_213} + \Ref{a6_214}
	+ \Ref{a6_216} \spc . \]
\QED
Mit diesem Lemma k\"onnen wir \Ref{a6_828b} weiter berechnen.
\begin{Lemma}
\label{a6_lemma24}
F\"ur den Diracoperator \Ref{a6_108} gilt
\begin{eqnarray*}
\lefteqn{ \chi_L \: \tilde{s}^\vee_{(2)}(x,y) \;=\; \chi_L \: U_L(x) \int_x^y dz \;
	(U_L^{-1} \:Y^2\: U_L)_{|z} \; U_L^{-1}(y) \; s^\vee_{(2)}(x,y) } \\
&&- \chi_L \: \frac{i}{16 \pi} \: U_L(x) \int_x^y (\alpha^2-\alpha) \; \Box
	(U_L^{-1} \:Y^2\: U_L)_{|z} \: U_L^{-1}(y) \:
	 \xi\slsh \; \Theta^\vee(\xi) \\
&&+ \chi_L \: \frac{i}{16 \pi} \: U_L(x) \int_x^y dz \int_x^z du \;
	\Pdd (U_L^{-1} \:Y\: U_R)_{|u} \; (z-x)^j \gamma_j \;
	\Pdd (U_R^{-1} \:Y\: U_L)_{|z} \\
&& \hspace*{4cm} \times U_L^{-1}(y) \; \Theta^\vee(\xi) \\
&&+ \chi_L \: \frac{i}{8 \pi} \: U_L(x) \int_x^y (1-\alpha) \;
	\Pdd (U_L^{-1} \:Y\: U_R)_{|z} \; (U_R^{-1} \:Y\: U_L)_{|z} \;
	U_L^{-1}(y) \; \Theta^\vee(\xi) \\
&&- \chi_L \: \frac{i}{8 \pi} \: U_L(x) \int_x^y \alpha \;
	(U_L^{-1} \:Y\: U_R)_{|z} \; \Pdd (U_R^{-1} \:Y\: U_L)_{|z} \;
	U_L^{-1}(y) \; \Theta^\vee(\xi) \;+\; {\cal{O}}(\xi^2) \spc .
\end{eqnarray*}
\end{Lemma}
{\Beweis}
Folgt durch Einsetzen von \Ref{a6_828b} in die Entwicklungsformel von
Lemma \ref{a6_lemma10}.
\QED

\begin{Satz}
\label{a6_satz25}
F\"ur den Diracoperator \Ref{a6_108} gilt mit der symbolischen Ersetzung $C=p$
oder $C=k$
\begin{eqnarray*}
\lefteqn{ \chi_L \: \tilde{C}^{(2)}(x,y) \;=\; \chi_L \: U_L(x) \int_x^y dz \;
	(U_L^{-1} \:Y^2\: U_L)_{|z} \;U_L^{-1}(y) \; C^{(2)}(x,y) } \\
&&-\frac{i}{2} \chi_L \: U_L(x) \int_x^y (\alpha^2-\alpha) \; \Box
	(U_L^{-1} \:Y^2\: U_L)_{|z} \: U_L^{-1}(y) \: \xi\slsh \; C^{(3)}(x,y) \\
&&+\frac{i}{2} \chi_L \: U_L(x) \int_x^y dz \int_x^z du \;
	\Pdd (U_L^{-1} \:Y\: U_R)_{|u} \; (z-x)^j \gamma_j \;
	\Pdd (U_R^{-1} \:Y\: U_L)_{|z} \\
&& \hspace*{5.5cm} \times U_L^{-1}(y) \; C^{(3)}(x,y) \\
&&+ i\:\chi_L \:  U_L(x) \int_x^y (1-\alpha) \;
	\Pdd (U_L^{-1} \:Y\: U_R)_{|z} \; (U_R^{-1} \:Y\: U_L)_{|z} \;
	U_L^{-1}(y) \; C^{(3)}(x,y) \\
&&- i\:\chi_L \:  U_L(x) \int_x^y \alpha \;
	(U_L^{-1} \:Y\: U_R)_{|z} \; \Pdd (U_R^{-1} \:Y\: U_L)_{|z} \;
	U_L^{-1}(y) \; C^{(3)}(x,y) \\
&& \;+\; \left\{ \begin{array}{cc}
	{\cal{O}}(\xi^2) & {\mbox{f\"ur $C=k$}} \\
	{\cal{O}}(\xi^0) & {\mbox{f\"ur $C=p$}} \end{array} \right. \spc .
\end{eqnarray*}
\end{Satz}
{\Beweis}
F\"ur $C=k$ folgt die Behauptung aus \Ref{2_tm}, f\"ur $C=p$ verwendet man das
Analogieargument von Lemma \ref{a6_lemma7}.
\QED

\subsection{Zus\"atzliche skalare/pseudoskalare St\"orung}
Wir lassen nun im Diracoperator eine zus\"atzliche skalare/pseudoskalare
St\"orung zu, also
\Equ{a6_830a}
G \;=\; \chi_L \: U_R (i \Pdd U_R^{-1}) \;+\;
	\chi_R \: U_L (i \Pdd U_L^{-1}) \;-\;
	m \: \Xi \;-\; i \rho \:m\: \Phi \spc .
\EndEqu
Wir verwenden wieder die Schreibweise \Ref{a6_826a} und \Ref{a6_828a}.
Durch Entwicklung von \Ref{2_t2} erh\"alt man
\begin{eqnarray*}
\lefteqn{ \chi_L \: \tilde{s}^\vee_{(2)} \;=\; \Ref{a6_254} +
	\Ref{a6_255} } \\
&&+ \sum_{p,q,r=0}^\infty \: \chi_L \: (-\svn \: {\cal{B}}_L)^p \:\svn\:
	(\Xi + i \Phi) \: (-\svn \: {\cal{B}}_R)^q
	\: \sve\:Y \: (-{\cal{B}}_L \: \svn)^r \\
&&+ \sum_{p,q,r=0}^\infty \: \chi_L \: (-\svn \: {\cal{B}}_L)^p \:\sve \:Y\:
	(-{\cal{B}}_R \: \svn)^q \:(\Xi - i \Phi)
	\: \svn \: (-{\cal{B}}_L \: \svn)^r \\
&&+ \sum_{p,q,r=0}^\infty \: \chi_L \: (-\svn \: {\cal{B}}_L)^p \:\svn\:
	(\Xi + i \Phi) \: \svn \: (-{\cal{B}}_R \: \svn)^q \: (\Xi - i \Phi)
	\: \svn \: (-{\cal{B}}_L \: \svn)^r \spc .
\end{eqnarray*}
Wir setzen die Umformung \Ref{a6_828b} ein und behandeln die anderen
Summen analog
\begin{eqnarray*}
&=& \chi_L \: U_L \: \svn \: (U_L^{-1} \:Y\: U_R) \: \svn \:
	(U_R^{-1} \:Y\: U_L) \: \svn \: U_L^{-1} \\
&&+\chi_L \: U_L \: \svn \: U_L^{-1} \: (\Xi + i \Phi) \: U_R \:
	\svn \: (U_R^{-1} \:Y\: U_L) \: \svn \: U_L^{-1} \\
&&+\chi_L \: U_L \: \svn \: (U_L^{-1} \: Y \: U_R) \:
	\svn \: U_R^{-1} \:(\Xi-i\Phi)\: U_L \: \svn \: U_L^{-1} \\
&&+\chi_L \: U_L \: \svn \: (\Xi + i \Phi) \:
	\svn \: (\Xi - i \Phi) \: \svn \: U_L^{-1} \spc ,
\end{eqnarray*}
also
\begin{eqnarray}
\label{eq:a6_828d}
\tilde{s}^\vee_{(2)}(x,y) &=& \chi_L \: \left( U_L \: \svn \:
	(U_L^{-1} \: Y_L \: U_R) \:\svn\: (U_R^{-1} \:Y_R\: U_L) \:
	\svn \: U_L^{-1} \right)(x,y) 
\:+\: {\cal{O}}(\xi^2) \spc .
\end{eqnarray}

\begin{Thm}
\label{a6_satz26}
F\"ur den Diracoperator \Ref{a6_830a} gilt mit der symbolischen Ersetzung $C=p$
oder $C=k$
\begin{eqnarray}
\label{eq:a6_A}
\lefteqn{ \chi_L \: \tilde{C}^{(2)}(x,y) \;=\; \chi_L \: U_L(x) \int_x^y dz \;
	(U_L^{-1} \:Y_L \: Y_R\: U_L)_{|z} \;U_L^{-1}(y) \; C^{(2)}(x,y) } \\
\label{eq:a6_B}
&&-\frac{i}{2} \chi_L \: U_L(x) \int_x^y (\alpha^2-\alpha) \; \Box
	(U_L^{-1} \:Y_L \: Y_R\: U_L)_{|z} \: U_L^{-1}(y) \:
	\xi\slsh \; C^{(3)}(x,y) \\
&&+\frac{i}{2} \chi_L \: U_L(x) \int_x^y dz \int_x^z du \;
	\Pdd (U_L^{-1} \:Y_L\: U_R)_{|u} \; (z-x)^j \gamma_j \;
	\Pdd (U_R^{-1} \:Y_R\: U_L)_{|z} \nonumber \\
\label{eq:a6_C}
&& \hspace*{5.5cm} \times U_L^{-1}(y) \; C^{(3)}(x,y) \\
&&+ i \: \chi_L \:  U_L(x) \int_x^y (1-\alpha) \;
	\Pdd (U_L^{-1} \:Y_L\: U_R)_{|z} \; (U_R^{-1} \:Y_R\: U_L)_{|z}
	\nonumber \\
\label{eq:a6_D}
&& \hspace*{5.5cm} \times U_L^{-1}(y) \; C^{(3)}(x,y) \\
\label{eq:a6_E}
&&- i \: \chi_L \:  U_L(x) \int_x^y \alpha \;
	(U_L^{-1} \:Y_L\: U_R)_{|z} \; \Pdd (U_R^{-1} \:Y_R\: U_L)_{|z} \;
	U_L^{-1}(y) \; C^{(3)}(x,y) \\
&& \;+\; \left\{ \begin{array}{cc}
	{\cal{O}}(\xi^2) & {\mbox{f\"ur $C=k$}} \\
	{\cal{O}}(\xi^0) & {\mbox{f\"ur $C=p$}} \end{array} \right. \spc . \nonumber
\end{eqnarray}
\end{Thm}
{\Beweis}
Folgt aus \Ref{a6_828d} genau wie Lemma \ref{a6_lemma24} und Satz \ref{a6_satz25}
mit Hilfe der asymptotischen Entwicklungsformel von Lemma \ref{a6_lemma10}.
\QED
Die Beitr\"age \Ref{a6_C}-\Ref{a6_E} sind Massenterme;
in erster Ordnung St\"orungstheorie stimmen sie mit \Ref{a2_x6}, \Ref{a4_z6} \"uberein.

Der Summand \Ref{a6_A} beschreibt eine verallgemeinerte Phasenverschiebung von
$C^{(2)}$. Er ist der Eichterm/Pseudoeichterm $\sim m^2$.

Felix Finster\\
Harvard University, Dept. of Mathematics\\
Cambridge, MA 02138, USA \\
email: felix@math.harvard.edu

\end{document}